\documentstyle[11pt]{article}

\newtheorem {thm}{Theorem}[section]
\newtheorem {prop}[thm]{Proposition}
\newtheorem {lem}[thm]{Lemma}
\newtheorem {cor}[thm]{Corollary}
\newtheorem {defn}[thm]{Definition}
\newtheorem {example}[thm]{Example}

\def\om{\omega}
\def\g{\gamma}

\def\J{{\bf J}}
\def\h{{\bf h}}
\def\p{{\bf p}}
\def\K{{\bf K}}

\def\eps{\varepsilon}
\def\e{\eta}
\def\s{\sigma}
\def\a{\alpha}
\def\b{\beta}

\def\Om{\Omega}
\def\Z{{\bf Z}}
\def\R{{\bf R}}
\def\T{{\bf T}}

\def\de{\delta}
\def\la{\lambda}
\def\La{\Lambda}
\def\De{\Delta}
\def\L{{\cal L}}
\def\B{{\cal B}}
\def\E{{\cal E}}
\def\F{{\cal F}}

\def\text{{\rm}}
\def\leqd{\preceq_{\cal D}}
\def\geqd{\succeq_{\cal D}}
\def\Cox{\hfill{\Box}}
\def\sgn{{\rm sgn}}
\def\ti{\to\infty}
\def\rf#1{(\ref{#1})}

\begin{document}

\title{{\bf The random geometry of equilibrium phases}}

\author{{\bf Hans-Otto Georgii}, Mathematisches Institut, \\
Ludwig-Maximilians-Universit\"at, D-80333 M\"unchen,\\ Germany.
\\ Email: georgii@rz.mathematik.uni-muenchen.de
 \and {\bf Olle H\"aggstr\"om}, Department of Mathematics, \\
 Chalmers University of Technology, S-412 96 G\"oteborg, \\
 Sweden. \\ Email: olleh@math.chalmers.se \and
 {\bf Christian Maes}, Instituut voor  Theoretische Fysica, \\
K.U.Leuven, B-3001 Leuven, \\ Belgium. \\ Email:
Christian.Maes@fys.kuleuven.ac.be}

\date{ }

\maketitle

\tableofcontents
\newpage

\section{Introduction}

Equilibrium statistical mechanics intends to describe and explain
the macroscopic behavior of systems in thermal equilibrium in terms
of the microscopic interaction between their
great many constituents. As a typical example, let us take some
ferromagnetic material like iron; the constituents are then the spins of
elementary magnets at the sites of some crystal lattice. Or we may
think of a lattice approximation to a real gas, in which case the
constituents are the particle numbers in the elementary cells of
any partition of space.  The central object is the Hamiltonian
describing the interaction between these constituents. This
interaction determines the relative energies between configurations
that differ only microscopically. The equilibrium states with respect to
the given interaction are described by the associated Gibbs
measures. These are probability measures on the space of
configurations which have prescribed conditional probabilities with
respect to fixed configurations outside of finite regions. These
conditional probabilities are given by the Boltzmann factor,
the exponential of the inverse temperature times the relative
energy. This allows one to compute, at least in principle, equilibrium
expectations and spatial correlation functions following the
standard Gibbs formalism. Most important are the so called extremal
Gibbs measures since they describe the possible
 macrostates of our physical system.  In such a state, macroscopic
observables do not fluctuate while the correlation between local
observations made far apart from each other decays to zero.

Since the
early days of statistical mechanics, geometric notions have played a
role in elucidating certain aspects of the theory.  This has taken many
different forms.  Arguably, the thermodynamic formalism, as
first developed by Gibbs, already
admits some geometric interpretations primarily related to convexity.
For example, entropy is a concave function of the
specific energy, the pressure is convex as a function of the interaction
potential, the
Legendre--Fenchel
 transformation relates various fundamental thermodynamic quantities to each other, and
the set of Gibbs measures for an interaction is a simplex with vertices
corresponding to the
physically realized macrostates, the equilibrium phases.

Here, however,
 we will not be concerned
 with this kind
 of convex geometry which
 is described in detail e.g.\ in the books by Israel \cite{Isr} and
Georgii \cite{Geo}.
Rather, the geometry considered here is a way of visualizing
the structure in the typical
 realizations of the system's constituents.
%
To be more specific let us consider for a moment the case of the
standard ferromagnetic Ising model on the square lattice.  At each
site we have a spin variable taking only two possible values, $+1$
and $-1$.  The interaction is nearest-neighbor and tends to align
neighboring spins in the same direction. By the  ingenious
arguments first formulated in 1936 by Peierls \cite{Pe} (see also
\cite{D1,Si,Geo}), the phase transition in this model can be
understood from looking at the typical configurations of contours,
i.e., the broken lines separating the  domains with plus resp.\
minus spins.
 The
plus phase (the positively magnetized phase) is realized by an infinite
ocean of plus spins with finite islands of minus spins (which in
turn may contain lakes of plus spins, and so on).  On the
other hand,  above the Curie temperature (first computed by
Onsager) there is no infinite path joining nearest neighbors
with the same spin value. So, for this model the geometric
picture is rather complete (as we will show later). In general,
however, much less is known, and much less is true. Still, certain
aspects of this  geometric analysis have wide applications, at
least in certain regimes of the phase diagram.  These `certain
regimes' are, on the one hand, the high-temperature (or, in a
lattice gas setting, low-density) regime and, to the other extreme,
the low temperature behavior.

At high temperatures, all thermodynamic considerations are based on the fact
that
 entropy dominates over energy. That is,
 the interaction
between the constituents is not effective enough to enforce a
macroscopic ordering of the system.  As a result, every constituent is
more or less free to behave at random, not much influenced
by other constituents which are far apart.
 So, the system's behavior is almost like
that of a free system with independent components. This means, in
particular, that in the center of a large box we will typically encounter
more or less the same configurations no matter what boundary
conditions outside this box are imposed. That is, if we compare two
independent realizations of the system in the box with different
boundary conditions outside then, still at high temperatures, the
difference between the boundary conditions cannot be felt by the
spins in the center of the box; specifically, there should not exist any
path from the boundary to the central part of the box along which
the spins of the two realizations disagree. This picture is rather
robust and can for example
 also be applied when the interaction is random; see Sections
\ref{sect:unireg} and \ref{sect:rain}.

At low temperatures, or large densities (when the interaction is
sufficiently strong), the picture above  no longer holds.
Rather, the specific characteristics of the interaction will come
into play and determine the specific features of the low temperature phase.
In many cases, the low temperature behavior can be
described as a random perturbation of a ground state, i.e., of a
fixed configuration of minimal energy.  Then we can expect that at
low temperatures, and sometimes even up to the critical
temperature, the equilibrium phases are realized as a deterministic
ground state configuration, perturbed by finite random islands on
which the configuration disagrees with the ground
state. This means that the ground state pattern can percolate
through the space to infinity. One prominent way of confirming this
picture is provided by the so called Pirogov--Sinai theory which is
described in detail e.g.\ in \cite{Zah3}.
 In Section \ref{sect:agreep} we will discuss some other
techniques of establishing the same geometric picture.

It is evident from the above that percolation theory will
play an important role in this text. In fact, we will mainly be
concerned with dependent percolation, but one can say that
independent percolation stands as a prototype for the study of
statistical  equilibrium properties in geometric terms. In
independent percolation, the model is extremely simple: the
components are binary-valued and independent from each other. What
is hard is the type of question one asks, namely the question of
existence of infinite paths of $1$'s and their geometry. We will
introduce percolation below but refer to other publications (such
as the book by Grimmett \cite{Gr1})
for a systematic account of the theory.

Percolation will come into play here on various levels.  Its
concepts like clusters, open paths, connectedness etc.\
will be useful for describing certain geometric features of
equilibrium phases, allowing characterizations of phases
in percolation terms.  Examples will be presented
where the (thermal) phase transition goes hand in hand with a phase
transition in an associated percolation process. Next,
percolation techniques can be used to obtain specific information
about the phase diagram of the system.  For example, equilibrium
correlation functions are sometimes dominated by connectivity
functions in an associated percolation problem which is easier to
investigate. Finally, representations in terms of percolation
models will yield explicit relations between certain observables in
equilibrium models and some corresponding percolation quantities.  In
fact, the resulting percolation models, like the random-cluster
model, have some interest in their own right and will also be studied in
some detail.

This text is supposed to be self-contained.  Therefore we need to
introduce various concepts and techniques which are well-known to some
readers.  On the other hand, important related issues will not be
discussed when they are not explicitly needed. For more complete
discussions on the introductory material we will refer the interested reader
to other sources.
More seriously, we will not include here a discussion of some important
geometric concepts developed in the 1980's  for the study
of critical behavior in statistical mechanical systems, namely
random walk expansions or random current representations.
Fortunately, we can refer to an excellent book \cite{FFS} where the
interested reader will find all the relevant results and
references.  Important steps in this context include
\cite{Ai,AB,ABF,AF,BFSo,BFSp} and the
references therein.

Finally, to avoid misunderstanding, the
random geometry in the title of this work should not be confused
with stochastic geometry (or geometric probability) which, as a
branch of integral geometry,
provides a very interesting tool-box for the discussion of
morphological characteristics of random fields appearing in
statistical physics and beyond, see \cite{KR,Me,Ad}.

\medskip\noindent
{\bf Acknowledgements: } It is a pleasure to thank J.\ Lebowitz for
suggesting (and insisting) that we should write this report. We are
also grateful to L.\ Chayes, A.C.D.\ van Enter and J.\ L\H{o}rinczi
who looked at parts of the manuscript and made numerous
suggestions, and to Y.\ Higuchi for discussions on Proposition
\ref{prop:MM}.

\section{Equilibrium phases} \label{sect:eq_phases}

\subsection{The lattice} \label{sect:lattice}

Our object of study are physical systems with many constituents,
spins or particles, which
will be located at the sites of a crystal lattice $\L$. The
standard case is when ${\L}
=\Z^d$,
the $d$-dimensional hypercubic lattice.  In general, we shall
assume that $\L$ is the vertex set of a countable graph. That is,
$\L$ is an at most countable set and comes equipped with a
(symmetric) adjacency relation. Namely, we write   $x \sim y$ if
the vertices $x,y\in \L$ are adjacent, and this is visualized by
drawing an edge between $x$ and $y$. In this case, $x$ and $y$ are
also called neighbors, and the edge (or bond) between $x$ and $y$
is denoted by $\langle xy \rangle$. We write $\B$ for the set of
all edges (bonds) in $\L$. A complete description of the graph is
thus given by the pairs $(\L,\sim)$ or $(\L,\B)$.

In the case $\L=\Z^d$, the edges will usually be drawn between lattice
sites of distance one; hence $x \sim y$ whenever $|x-y|=1$. Here,
$|\cdot|$ stands for the sum-norm, i.e. $|x|=\sum_{i=1}^d |x_i|$ whenever
$x=(x_1,\ldots,x_d)\in\Z^d$. This choice is natural because then $|x-y|$
coincides with
 the graph-theoretical (or lattice) distance, viz.\ the length of
 the shortest path (consisting of consecutive edges)
connecting $x$ and $y$. For convenience, we sometimes use the same
notation in the case of a general graph.
On $\Z^d$ we will occasionally distinguish between the standard metrics
 $d_1(x,y) = \sum_i |x_i-y_i|$,
$d_2(x,y) = [\sum_{i=1}^d (x_i - y_i)^2)]^{1/2}$
and $d_{\infty}(x,y)= \max_i |x_i - y_i|$.
Given any metric $d$ on $\L$, we write $d(\La,\De)=\inf_{x\in\La,y\in\De}d(x,y)$
for the distance of two subsets $\La,\De\subset\L$.

We will always assume that the graph  $(\L,\sim)$ is {\em locally finite},
which means that each $x\in\L$ has only a finite number $N_x$ of nearest
neighbors.
In other words, $N_x$ is the number of edges emanating from $x$. $N_x$ is
also called the {\em degree} of the graph at $x$. In many cases we will even
assume that $(\L,\sim)$ is of bounded degree, which means that
$N=\sup_{x\in\L}N_x<\infty$. Common examples of such graphs, besides
$\Z^d$, are  the triangular lattice in two dimensions, and the regular
tree ${\bf T}_d$
(also known as the Cayley tree or the Bethe lattice), which is defined as
the (unique) infinite connected graph containing no circuits
in which every vertex has exactly $d+1$ nearest neighbors.

A region of the lattice, that is a subset $\La\subset{\L}$,  is called
finite if its cardinality $|\La|$ is finite. We write ${\cal E}$
for the collection of all finite regions.  The complement of a region $\La$
will be denoted by $\La^c = {\cal L} \setminus \La$.
The boundary $\partial \La$ of $\La$ is the set of all sites (vertices)
in $\La^c$ which are adjacent to some site of $\La$.

At some occasions we will need the notion of thermodynamic (or
infinite volume) limit, and we need to describe in what sense a
region $\La \in \cal E$ grows to the full lattice $\cal L$. For
our purposes, it will in general be sufficient to take an arbitrary
increasing sequence $(\La_n)$ with $\bigcup_{n\ge1}\La_n=\L$. In
the case $\L=\Z^d$, we will often make the standard choice
$\La_n = [-n,n]^d \cap \Z^d$, the lattice cubes centered around
the origin. As $\E$ is a directed set ordered by inclusion, we will
occasionally also consider the limit along $\E$. In each of these cases
we will use the notation $\La\uparrow\L$.


\subsection{Configurations}

The constituents of our systems are the spins or particles at the lattice sites. So,
at each site $x\in \L$ we have a
 variable $\sigma(x)$ taking values in a non-empty set $S$, the
{\em state space} or {\em single-spin space}. In a magnetic set-up
(to which we mostly adhere for simplicity), $\s(x)$ is interpreted
as the spin of an elementary magnet at $x$. In a lattice gas
interpretation, there is a distinguished vacuum state $0\in S$
representing the absence of any particle, and the remaining
elements correspond to the types and/or the number of the particles
at $x$. Unless stated otherwise, we will always assume that $S$ is
finite. Elements of $S$ will typically be denoted by $a,b,\ldots$.

A {\em configuration} is a
function $\sigma : {\cal L}\rightarrow S$ which assigns to each vertex
$x\in \cal L$ a spin value
$\sigma(x) \in S$. In other words, a configuration $\sigma$ is
an element of the product space $\Omega = S^{\cal L}$.
$\Omega$ is called the configuration space
and its elements are in general written as $\sigma, \eta, \xi,\ldots$.
(It is sometimes useful to visualize $a, b\ldots$ as colors.
 A configuration is then a coloring of the lattice.)
A configuration $\sigma$ is constant if for some $a\in S$,
$\sigma(x) = a$ for all $x\in \cal L$.  Two configurations $\sigma$
and $\eta$ are said to agree on a region  $\La\subset\L$,  written
as ``$\sigma \equiv\eta$ on $\La$'', if $\sigma(x) = \eta(x)$ for all
$x\in \La$. Similarly, we write ``$\sigma\equiv\eta$ off $\La$'' if
$\sigma(x) = \eta(x)$ for all $x\notin \La$.

We also consider configurations in regions $\La\subset\L$.
These are elements of $S^\La$, again denoted by letters like $\sigma, \eta,
\xi,\ldots$. Given $\sigma,\eta \in \Omega$, we write
$\sigma_\La\eta_{\La^c}$ for
the configuration $\xi\in\Omega$ with $\xi(x)=\sigma(x)$ for $x\in \La$ and
$\xi(x)=\eta(x)$ for $x\in \La^c$.  Then, obviously, $\xi\equiv\sigma$ on
$\La$. The cylinder sets
$$
{\cal N}_\La(\sigma) = \{\xi \in \Omega: \ \xi\equiv\sigma \mbox{ on }
\La\}\,, $$
$\La \in \cal E$, form a countable neighborhood basis of $\s\in\Om$; they
generate the
product topology on $\Omega$. Hence, two configurations are close to each
other if they agree on some large finite region, and a diagonal-sequence argument
shows that $\Omega$ is a compact in this topology.

We will often change a configuration $\sigma\in \Omega$ at just one site
$x\in\L$. Changing $\s(x)$ into a prescribed value $a\in S$ we obtain a
new configuration written $\s^{x,a}$.  In particular, for $S=\{-1,+1\}$
we write
\[
\sigma^x(y) =
\left\{
\begin{array}{ll}
\sigma(y) & \mbox{for }y\neq x \\
-\sigma(x) & \mbox{for } y=x
\end{array}\right.
\]
for the configuration resulting from flipping the spin at $x$.

We will also deal with automorphisms of the underlying lattice $(\L,\sim)$.
Each such automorphism defines a measurable transformation of the configuration
space $\Om$. The most interesting automorphisms are the translations
of the integer lattice ${\cal L} = \Z^d$; the associated translation group
acting on $\Om$ is given by
$\theta_x \sigma (y) = \sigma(x+y)$, $y\in \Z^d$.  In particular, any constant
configuration is translation invariant. Similarly, we can speak
about periodic configurations which are invariant under $\theta_x$ with
$x$ in some sublattice of $\Z^d$.

Later on, we will also consider configurations which refer to the lattice
bonds rather than the vertices. These are elements $\eta$ of the product
space $\{0,1\}^\B$, and a bond $b\in\B$ will be called open if
$\eta(b)=1$, and otherwise closed. The above notations apply to this
situation as well.

\subsection{Observables}

An observable is a real function on the configuration
space which may be thought of as the numerical outcome of some physical measurement.
 Mathematically, it is a measurable real function on $\Omega$.
Here, the natural underlying $\sigma$-field of measurable events in
$\Omega$ is the
product $\s$-algebra ${\cal F} = ({\cal F}_0)^{\cal L}$, where
${\cal F}_0$ is the set of all subsets of $S$. $\F$ is defined as the
smallest $\s$-algebra on $\Omega$ for which all projections $X(x):
\Omega\to S$, $X(x)(\s)=\s(x)$ with $\s\in\Omega$ and $x\in\L$, are
measurable. It coincides, in fact, with the Borel
$\s$-algebra for the product topology on $\Omega$.

We also consider events and observables depending only on some region
$\La\subset\L$. We let ${\cal F}_\La$ denote the smallest sub-$\sigma$-field
of ${\cal F}$ containing the events ${\cal N}_\Delta(\sigma)$ for $\sigma\in
S^\Delta$ and $ \Delta\in\E$ with $\Delta\subset\La$. Equivalently,
${\cal F}_\La$ is the $\s$-algebra generated by the projections
$X(x)$ with $x\in\La$. ${\cal F}_\La$ is the $\s$-algebra of events
occurring in $\La$.

An event $A$ is called {\em local\/} if it occurs in some finite region,
which means that $A\in\F_\La$ for some $\La\in\E$.
Similarly, an observable $f:\Omega\to\R$ is called {\em local\/} if it
depends on only finitely many spins, meaning that $f$ is measurable with
respect to ${\cal F}_\La$ for some $\La\in \cal E$.
%
%
 More generally, an observable
 $f$ is called {\em quasilocal\/} if it is (uniformly) continuous, i.e., if
 for all $\epsilon > 0$ there is some
$\La\in \cal E$ such that $|f(\sigma) - f(\xi)| < \epsilon$
whenever $\xi\equiv\sigma$ on $\La$. The set $C(\Omega)$ of
continuous observables is a Banach space for the supremum norm
$||f|| = \sup_\sigma |f(\sigma)|\,,$ and the local observables are
dense in it.

The local events and observables should be viewed as microscopic
quantities. On the other side we have the macroscopic quantities which
only depend on the collective behavior of all spins, but not on the
values of any finite set of spins. They are defined
in terms of the  {\em tail $\s$-algebra\/} ${\cal T}
=\bigcap_{\La\in\E}\F_{\La^c}$, which is also called the
$\s$-algebra of all events at infinity. Any tail event $A\in{\cal T}$
and any ${\cal T}$-measurable observable is called {\em macroscopic}.

As a final piece of notation we introduce the {\em indicator function\/}
 $I_A$ of an event $A$; it takes the value $1$ if
the event occurs ($I_A(\sigma) = 1$ if $\sigma \in A$) and is zero
otherwise.

\subsection{Random fields}

As the spins of the system are supposed to be random, we will consider
 suitable probability measures $\mu$ on $(\Omega, {\cal F})$. Each such
 $\mu$ is called a {\em random field}. Equivalently, the family
 $X=(X(x), x\in {\cal L})$ of random variables on the probability space
 $(\Omega, {\cal F},\mu)$ which describe the spins at all sites is called
 a random field.

Here are some standard notations concerning probability measures.
The expectation
of an observable $f$ with respect to $\mu$ is written as $\mu(f) = \int f d\mu$.
%
%
The probability of an event $A$ is
$\mu(A) = \mu(I_A) = \int_{A} d\mu$,
and we omit the set braces when $A$ is given explicitly.
For example, given any $x\in \L$ and $a\in S$ we write
$\mu(X(x)=a)$ for $\mu(A)$ with $A=\{\s\in \Om:\ \s(x)=a\}$.
Covariances are abbreviated as $\mu(f;g) = \mu(fg) - \mu(f)\mu(g)$.

Whenever we need a topology on probability measures on $\Omega$, we shall
take the {\em weak topology}. In this (metrizable) topology, a sequence of
probability measures $\mu_n$ converges to $\mu$, denoted by $\mu_n
\rightarrow \mu$, if $\mu_n(A)\rightarrow \mu(A)$ for all local events
$A\in\bigcup_{\La\in\E}\F_\La$. This holds if and only if
$\mu_n(f)\rightarrow \mu(f)$ for all local, or equivalently, all
quasilocal functions $f$. In applications, $\mu_n$ will often be an
equilibrium state in a finite box $\La_n$ tending to $\L$ as
$n\ti$, and we are interested in whether the probabilities of
events occurring in some fixed finite volume have a well-defined
thermodynamic (or bulk) limit. That is, we observe what happens
around the origin (via the local function $f$) while the boundary
of the box in which we realize the equilibrium state receeds to
infinity. As there are only countably many local events, one can
easily see by a diagonal-sequence argument that the set of all
probability measures on $\Omega$ is compact in the weak topology.

\subsection{The Hamiltonian}

We will be concerned with systems of {\em interacting\/} spins. As usual, the
interaction is described by a Hamiltonian. As the spins are located at the
sites of a graph $(\L,\sim)$, it is natural to consider the
case of homogeneous neighbor potentials. (We will deviate from homogeneity
in Section \ref{sect:rain} when considering random interactions.)
The Hamiltonian $H$ then takes the form
\begin{equation} \label{eq:Hamiltonian}
H(\sigma)= \sum_{x\sim y } U (\s(x),\s(y)) + \sum_x V
(\s(x))
\end{equation}
with a symmetric function $U: S\times S
\rightarrow {\R \cup \{\infty\}}$, the {\em neighbor-interaction},
and a {\em self-energy\/} $V: S \rightarrow\R $.
The infinite sums are formal; the summation index $x\sim y$
 means that the sum extends over all bonds $\langle x y\rangle\in\B$ of the
lattice. $U$ thus describes
the interaction between spins at neighboring sites,
while $V$ might come from the action of an external magnetic field. In a lattice gas
interpretation when $S=\{0,1\}$ (the value $1$ being assigned to
sites which are occupied by a particle), $V$ corresponds to a chemical
potential.

To make sense of the formal sums in \rf{eq:Hamiltonian} we compare the
Hamiltonian for two different configurations $\s ,\e \in \Omega$ which
{\em differ only locally} (or are ``local perturbations'' or  ``excitations''
of each other), in that
$\sigma \equiv\eta$ off some $\La\in\E$.
For such configurations we can define the
{\em relative Hamiltonian}
\begin{equation}\label{relener}
H(\s|\e) =
 \sum_{x\sim y }\ [U (\s(x),\s(y)) -
U (\e(x),\e(y)) ]
+ \sum_x\ [V(\s(x))-V(\e(x))]
\end{equation}
in which the sums now contain only finitely man non-zero terms: the first part
is over those neighbor pairs $\langle xy\rangle$ for which at least one
of the sites belongs to $\La$, and the second part is over all
$x\in\La$.

\subsection{Gibbs measures}\label{sect:Gibbsmeasures}

Gibbs measures are random fields which describe our physical spin system
when it is in {\em macroscopic equilibrium with respect to the given microscopic
interaction\/} at a fixed temperature. Here, macroscopic equilibrium means
that all parts of the system are in equilibrium with their exterior
relative to the prescribed interaction and temperature. So it is
natural to define Gibbs measures in terms of conditional probabilities.
\begin{defn}
A probability measure $\mu$ on the configuration space $\Omega$ is called
a {\bf Gibbs measure} for the Hamiltonian $H$
in (\ref{eq:Hamiltonian}) or (\ref{relener}) at inverse temperature
$\beta \sim 1/T$ if for all $\La\in \cal E$ and all $\s \in\Om$,
\begin{equation} \label{Gibbs_specification}
\mu(X\equiv\s \mbox{ \rm on } \La\,|\,X\equiv\e \mbox{ \rm off }\La) =
\mu_{\beta,\La}^\eta(\s)
\end{equation}
for $\mu$-almost all $\eta\in\Omega$. In the above,
$\mu_{\beta,\La}^\eta(\s)$ is the {\bf Boltzmann--Gibbs distribution} in
$\La$ for $\beta$ and $H$, which is  given by
\begin{equation} \label{eq:FvolGibbs}
\mu_{\beta ,\La}^{\e }(\s) =
\frac{I_{\{\s\equiv\e \mbox{ \scriptsize off } \La\}}}{Z_\La (\beta,\e)}
\exp[-\beta H(\s|\e)]\;.
\end{equation}
Here, $Z_\La (\beta,\e)$ is a normalization constant making
$\mu_{\beta ,\La}^\e$ a probability measure, and
the constraint that $\sigma$ has to coincide with $\eta$ outside $\La$ is
added because we want to realize these probability measures immediately on the
infinite lattice. Note that $\mu_{\beta ,\La}^\e$ in fact only depends on
the restriction of $\e$ to $\La^c$.
\end{defn}

\noindent
So, $\mu$ is a Gibbs measure if it has prescribed conditional distributions {\it inside}
some finite set of vertices, given that the configuration is held fixed {\it
outside}, and these conditional distributions are given by the usual
Boltzmann--Gibbs formalism. This definition goes back to
the work of Dobrushin \cite{D3} and Lanford and Ruelle \cite{LR} in the
late 1960's, whence Gibbs measures are often called {\em DLR-states}.
By this work,
equilibrium statistical physics and the study of phase transitions
made firm contact with probability theory and the study of random fields.
A thermodynamic justification of this definition can be given by the
variational principle, which states that (in the case $\L=\Z^d$) the
translation invariant Gibbs measures are precisely those translation
invariant random fields which minimize the free energy density, cf.\
\cite{LR,Geo}. For better distinction, the Gibbs distributions $\mu_{\beta ,\La}^\e$
are often called {\em finite volume} Gibbs distributions, whereas the
Gibbs measures are sometimes specified as {\em infinite volume} Gibbs measures.

 We write ${\cal G}(\beta H)$ for the set of
all Gibbs measures with given Hamiltonian $H$ and inverse
temperature $\beta$. In the special case $U \equiv 0$ of no
interaction, there is only one Gibbs measure, namely the product
measure with one-site marginals $\mu(X(x) = a)$\linebreak $=
e^{-\beta  V(a)} / \sum_{b\in S}e^{-\beta  V(b)}$. In general,
several Gibbs measures for the same interaction and temperature can
coexist. This is the fundamental phenomenon of nonuniqueness of
phases which is one of our main subjects; we return to this point
in Section \ref{sect:phases} below.

First we want to emphasize an important consequence of our assumption
that the underlying interaction $U$ involves only neighbor spins.
Due to this assumption, the Gibbs distribution $\mu_{\beta,\La}^\eta$
only depends on the restriction $\e_{\partial\La}$ of $\e$ to the
boundary $\partial\La$ of $\La$, and this implies that each Gibbs
measure $\mu\in{\cal G}(\beta H)$ is a {\em  Markov random field}.
By definition, this means that for each $\La \in \cal E$ and
$\s \in S^{\La}$
\begin{equation} \label{markov}
\mu(X\equiv \s \mbox{ on }\La\,|\,\F_{\La^c}) =
\mu(X\equiv \s \mbox{ on }\La\,|\,\F_{\partial\La}),
\end{equation}
$\mu-$almost surely. This Markov property will be an essential tool
in the geometric arguments to be discussed in this review. There is
in fact an equivalence between Markov random fields and Gibbs
measures for nearest neighbor potentials, see e.g.\ Averintsev
\cite{Av}, Grimmett \cite{Gr73} or Georgii \cite{Geo}.

As an aside, let us comment on the case when the interaction of spins is not
nearest-neighbor but only decays sufficiently fast with their
distance. The Boltzmann-Gibbs distributions in \rf{eq:FvolGibbs}, and
therefore also Gibbs measures, can then still be defined, but the
Gibbs measures fail to possess the Markov property \rf{markov}.
Rather their local conditional distributions
 $\mu_{\beta ,\La}^\e$ satisfy a weakening of the
Markov property called {\em quasilocality} or {\em almost-Markov property\/}:
for every $\La\in\E$ and $A \in {\cal F}_\La$,
$\mu_{\b,\La}^\e(A)$ is a continuous function of $\eta$. So, in this
case, Gibbs measures have prescribed continuous versions of their
local conditional probabilities. To obtain a sufficiently general
definition of Gibbs measures including this and other cases, one
introduces the concept of a specification $G = (G_\La,\ \La \in \E)$.
This is a family of probability kernels $G_\La$ from $(\Om,
\F_{\La^c})$ to $(\Om,\F)$.  $G_\La(\cdot,\eta)$
 stands for any distribution of
 spins with fixed configuration $\e_{\La^c} \in
S^{\L \setminus\La}$ outside $\La$; the standard case is the
Gibbs specification $G_\La(\cdot,\eta)
=\mu_{\beta ,\La}^{\e }$.
A Gibbs measure is then a probability measure $\mu$ on $\Omega$
satisfying $\mu(A\,|\,\F_{\La^c})= G_\La(A,\cdot)$  $\mu$-almost
surely for all $\La\in\E$ and $A\in \F$; this property can be
expressed in a condensed form by the invariance equation $\mu G_\La
=\mu$. In order for this definition to make sense the specification
$G$ needs to satisfy a natural compatibility condition for pairs of
volumes $\La \subset
\La'$ expressing the fact that if the system in $\La'$ is
in equilibrium with its exterior, then the subsystem in $\La$
is also in equilibrium with its own exterior. It is easy to see that
the Gibbs distributions in \rf{eq:FvolGibbs} are compatible in this sense.
 Details and further
discussion can be found in many books and articles dealing with
mathematical results in equilibrium statistical mechanics,
including \cite{Isr,Ru2,Geo,Pr,vEFS,Si}. In \cite{L}, the relation between Gibbs
measures and the condition of detailed balance (reversibility) in
certain stochastic dynamics is explained.

Finally, we mention an alternative and constructive approach to the
concept of Gibbs measures. Starting from the finite-volume Gibbs
distributions $\mu_{\beta ,\La}^{\e }$, one might ask what kind
of limits could be obtained if $\e$ is randomly chosen and $\La$
increases to the whole lattice $\L$. (This slightly older but still
important approach was suggested by Minlos \cite{Min}.) To make
this precise we consider the measures $\mu_{\beta,\La}^{\rho} =
\int
\mu^\eta_{\beta,\La}\, \rho(d\eta)$, where  $\rho$ is any probability measure
on $\Om$ describing a ``stochastic boundary condition''.
Any such $\mu_{\beta,\La}^\rho$ is called a (finite volume)
Gibbs distribution with respect to $H$ at inverse temperature $\beta$, and
their collection is denoted by ${\cal G}_\La(\beta H)$.
The set of all (infinite volume) Gibbs measures is then equal to
\[
{\cal G}(\beta H) = \bigcap_{\La \in {\cal E}}{\cal G}_\La(\beta H).
\]
Equivalently, a probability measure $\mu$ on $\Omega$ is a Gibbs
measure for the Hamiltonian $\beta H$ if it belongs to the closed convex hull of the
set of limit points of $\mu_{\beta,\La}^\eta$ as $\La\uparrow\L$.

One important consequence is that ${\cal G}(\b H) \neq \emptyset$. This is
because each ${\cal G}_\La(\b H)$ is obviously non-empty and compact.
Equivalently, to obtain an infinite volume Gibbs measure one can fix a particular
configuration $\e$ and take it as boundary condition. By compactness,
we obtain an infinite volume Gibbs measure $\mu_\beta^\e$
by taking the limit of (\ref{eq:FvolGibbs}) as $\La
\uparrow {\L}$, at least along suitable subsequences; for details see
e.g.\ Preston \cite{Pr} or Georgii \cite{Geo}.
We remark that in general there is no unique limiting measure $\mu^\e
_\beta$; rather there may be several such limiting measures obtained as
limits along different
subsequences.  Fortunately, however, this is not the case for a wide class of
models, either at low temperatures ($\beta$ large) when $\eta$ is a ground
state configuration (in the
realm of the Pirogov--Sinai theory),  or at high temperatures when
$\beta$ is small.

We conclude this subsection with a general remark.
As all systems in nature are finite, one may wonder why
we consider here systems with infinitely many constituents.
The answer is that sharp results for bulk
quantities can only be obtained when we make the
idealization to an infinite system.  The thermodynamic limit eliminates
finite size effects (which are always present but which are not always
relevant for certain phenomena)  and it is only in the thermodynamic
limit of infinite volume that
we can get a clean and precise picture of realistic phenomena such
as phase
transitions or phase coexistence. This is a consequence of the general
probabilistic principle of large numbers. In this sense, infinite
systems serve as an idealized
approximation to very large finite systems.

\subsection{Phase transition and phases}\label{sect:phases}

As pointed out above, in general there may exist several solutions
$\mu$ to the DLR-equation (\ref{Gibbs_specification}) for given
$U$, $V$ and $\beta$, which means that multiple Gibbs measures
exist. The system can then choose between several equilibrium
states. (In a dynamical theory this choice would depend on the
past; but here we are in a pure equilibrium setting.) {\em The
phenomenon of non-uniqueness therefore corresponds to a phase
transition}. In fact, it is then possible to construct different
Gibbs measures as infinite volume limits of Gibbs distributions
with different choices of boundary conditions \cite{Gii1,Geo}.
 Since any two Gibbs measures can be distinguished by a suitable local
 observable, a phase transition can be detected by looking at such a local observable
which is then called an {\em order para\-meter}.
Varying the external parameters such as temperature or an external magnetic
field (which can be tuned by the experimenter via some heatbath or
reservoir) one will observe different scenarios;
these are collected in the so
called {\em phase diagram\/} of the considered system.

As we have indicated in the introduction, the phase transition
phenomenon is of central interest in equilibrium statistical
mechanics.
When phase transitions occur and when they do not is also one of the primary
questions (although we will encounter many others) that we will try
to answer with the geometric methods to be developed in subsequent
sections.

If multiple Gibbs measures for a given interaction exist, the structure of the
set ${\cal G}(\beta H)$ of all Gibbs
measures becomes relevant. We only state here the most basic results;
a detailed exposition can be found in \cite{Geo}, for example.
The basic observation is that ${\cal G}(\beta H)$  is a convex
set.  Its extremal elements, the extremal Gibbs measures, have a trivial
tail $\s$-field ${\cal T}$ (which means that all events in ${\cal T}$
have probability 0 or 1). Equivalently, all macroscopic observables
are almost surely constant. In addition, the tail triviality can be
characterized by an asymptotic independence (or mixing) property.
On the other hand, any Gibbs measure $\mu$ can be decomposed into extremal
Gibbs measures; therefore
every configuration which is typical for $\mu$ is in fact
typical for some extremal Gibbs measure. This shows that the extremal Gibbs
measures correspond to what one can really see in nature as far as large systems in
equilibrium are concerned.  The extremal Gibbs measures therefore correspond to the
physical macrostates, whereas non-extremal Gibbs measures only provide a limited
description when the system's precise state is unknown. For all these reasons,
the extremal Gibbs measures are called (equilibrium) {\em phases}. The
central subject of this review is the geometric
analysis of their typical configurations, and thereby the analysis of
 the phase diagram giving the variation in the number and the nature of the
phases as one changes various control parameters (coupling, temperature,
external fields, etc.).

Often it is natural to consider automorphisms
of the graph $(\L,\sim)$.  For example,
if $\L=\Z^d$ we consider the translation group $(\theta_x)_{x\in \Z^d}$.
A homogeneous phase is then an
extremal Gibbs measure which is also translation invariant.
On the other hand, we can regard the extremal points of
the convex set of all translation invariant Gibbs measures.
 These are ergodic, which means that they cannot
be decomposed into distinct translation
invariant probability measures, and are trivial on the $\s$-algebra
of all translation invariant events. However, these extremal translation
invariant Gibbs measures need not
be homogeneous phases; they are {\it only\/} ergodic.  Yet, ergodic
measures $\mu$ satisfy a law of large numbers:  for
any observable $f$ and any sequence of increasing cubes
$\La$,
$$\lim_{\La\uparrow \Z^d}
\frac 1{|\La|}
\sum_{x\in \La} f\circ \theta_x = \mu(f)\quad \mu\mbox{-almost surely.}$$
Hence, ergodic Gibbs measures
are suitable for modelling macrostates  in equilibrium if
one limits oneself to measuring certain bulk observables or
macroscopic quantities with additivity properties.
Notice, however, that there exists a certain non-uniformity in the literature
concerning the nomenclature.  Sometimes these ergodic Gibbs measures are called
(pure) phases.  It is then argued that
it might happen that two phases (as defined above) for a system can by no
means be macroscopically distinguished (for example if one is a translation
of the other).  We do not wish to enter into a detailed discussion of these
points.

\section{Some models}
In this section we discuss briefly the phase transition behaviour of
some prototypical examples of Gibbs systems. Although these examples
are fairly standard and well-known to most of our readers, we need to
include them here to set up the stage. They will be studied in
detail in the later sections. An account of phase transition phenomena
in more general lattice models can be found in many other sources,
including \cite{Kot,Geo,Sim,Si,SM}.

\subsection{The ferromagnetic Ising model} \label{sect:Ising}

The Ising model was introduced in the 1920's by Wilhelm Lenz \cite{Lenz}
and his student Ernst Ising \cite{I} as a simple model for magnetism
and, in particular, ferromagnetic phase transitions.
Each site $x\in \L$ can take either of two spin values, $+1$
(``spin up'') and
$-1$ (``spin down''), so that the state space is equal to
$S=\{-1,+1\}$. The Hamiltonian is given by (\ref{eq:Hamiltonian}) with
$U(\s(x),\s(y))=-\s(x)\,\s(y)$ and $V(\s(x))=-h\,\s(x)$. The
parameter $h\in \R$ describes an external field. The finite volume
Gibbs distribution in a box $\La$ with external field $h$ at inverse
temperature $\b>0$ with
boundary condition $\eta$ is thus the probability measure
$\mu_{h,\b,\La}^\eta$ on $\Omega=\{-1,+1\}^\L$ which to each
$\sigma\in\Omega$ assigns probability proportional to
\[
I_{\{\s\equiv\e \mbox{ \scriptsize on } \La^c\}} \
\exp\Bigg[\;\b\,
\Big( \sum_{x\sim y \atop x \in \La\,
 \mbox{\tiny or}\,y\in\La}\sigma(x)\,\sigma(y)
+h \sum_{x\in \La}\sigma(x)\Big)\Bigg]\,.
\]
For $\b=0$ (``infinite temperature'') the spin
variables are independent under $\mu_{h,\b,\La}^\eta$,
but as soon as $\b>0$
the probability distribution starts to favour configurations with
many neighbor pairs of aligned spins. This
tendency becomes stronger and stronger as $\b$ increases.

In the case $h=0$ of no external field, the model is symmetric
under interchange of the spin values $-1$ and $+1$, so that there
is an equal chance of having many pairs of plus spins or having
many pairs of minus spins. This dichotomy gives rise to the
following interesting behavior. Suppose that $\L=\Z^d$, $d\geq 2$.
If $\b$ is sufficiently small (i.e., in the high temperature
regime), the interaction is not strong enough to produce any long
range order, so that the boundary conditions become irrelevant in
the infinite volume limit and the Gibbs measure is uniquely
determined. By ergodicity and the $\pm$ symmetry, the limiting
fraction of plus spins will almost surely\ be ${1}/{2}$ under
this unique Gibbs measure. In contrast, when $\b$ is
sufficiently large (in the low temperature regime), the interaction
becomes so strong that a long range order appears: the bias towards
neighbor pairs of equal spin then implies that Gibbs measures
prefer configurations with either a vast majority of plus spins or
a vast majority of minus spins, and this preference even survives
in the infinite volume limit. The system thus undergoes a phase
transition which manifests itself in a non-uniqueness of Gibbs
measures. Specifically, there exist two particular Gibbs measures
$\mu^+$ and $\mu^-$, obtained as infinite volume limits with
respective boundary conditions $\eta\equiv +1$ and $\eta\equiv
-1$, which can be distinguished by their overall density of $+1$'s:
the density is greater than ${1}/{2}$ under $\mu^+$ and (by
symmetry) less than ${1}/{2}$ under $\mu^-$. This is the
spontaneous magnetization phenomenon that Lenz and Ising were
looking for but were discouraged by not finding it in one
dimension.
In higher dimensions, the uniqueness regime and the phase transition regime
are separated by a sharp critical value $\b_c$, as is
summarized in the following classical theorem \cite{Pe,D1,D2}:
\begin{thm} \label{thm:Ising}
For the ferromagnetic Ising model on the integer lattice ${\Z}^d$
of dimension $d\geq 2$ at zero external field,
there exists a critical inverse temperature $\b_c\in(0,\infty)$
(depending on $d$)
such that for $\b<\b_c$ the model has a unique Gibbs measure
while for
$\b>\b_c$ there are multiple Gibbs measures.
\end{thm}
A stochastic-geometric proof of this result will be given in
Section \ref{sect:random-cluster}. In fact, the result (as well as
its proof) holds for any graph $(\L,\sim)$ in place of $\Z^d$,
except that $\b_c$ may then take the values $0$ or $\infty$. For
instance, on the one-dimensional lattice ${\bf Z}^1$ we have
$\b_c=\infty$, which means that there is a unique Gibbs measure for
all $\b$. For $\Z^2$, the critical value has been found to be
$\b_c=\frac{1}{2}\log(1+\sqrt{2})$. This calculation is a
remarkable achievement which began with Onsager \cite{Ons} in 1944.
An account of various (algebraic and/or combinatorial) methods can
be found e.g. in \cite{Tho,CW}. Let us also mention the work done
in 1973 by Abraham and Martin-L\"of \cite{AML} relating these exact
computations to the real magnetization in the appropriate Gibbs
measures; it also gives the result that there is a unique Gibbs
measure {\em at} the critical value $\b=\b_c$. A rigorous
calculation of the critical value in higher dimensions is beyond
current knowledge. It is believed that uniqueness holds at
criticality in all dimensions $d\geq 2$, but so far this is only
known for $d=2$ and $d\geq 4$ \cite{AF}.

The case of a nonzero external field $h\neq 0$ is less interesting,
in that one finds a unique Gibbs measure for all $\b$ and $d$.
The intuitive explanation is that for $h\neq 0$ there is no $\pm$
symmetry  which could be broken;  depending on the sign of $h$,
the system is forced to prefer either $+1$'s or $-1$'s.  This comes
from the fact that the magnetic field acts on the whole volume,
whereas the influence of a boundary condition is of smaller order
as the volume increases. In contrast, a phase transition for $h\neq 0$
does occur when ${\bf Z}^d$ is replaced by
certain nonamenable graph structures for which the boundary of a
volume is of the same order of magnitude as the volume itself
(which makes them physically perhaps less realistic) --
an example is the regular tree ${\bf T}_d$ with $d\geq 2$;
we refer to \cite{Sp,Geo,JS}.  A phase
transition can also occur for a non-zero external field for the
Ising model on a half-space where it is due to the so-called Basuev
phenomenon \cite{Bas, Ba2}.

Because of the simplicity of its model assumptions, the standard
Ising model has inspired a variety of techniques for analyzing
interacting random fields. Its ferromagnetic structure suggests
various monotonicity properties which can be checked by the
coupling methods to be described in Section \ref{sect:coupling_SD},
and the assumption of neigbor interaction implies the spatial
Markov property  \rf{markov} which plays a fundamental role in
the geometric analysis of typical configurations.
Many techniques which were developed on this testing ground turned
out to be fruitful also in more general cases.

\subsection{The antiferromagnetic Ising model} \label{sect:anti}

The Ising antiferromagnet is defined quite similarly to the
ferromagnetic case, except that $U(\sigma(x),\sigma(y))$ is
taken to be $+\sigma(x)\,\sigma(y)$ rather than $-\sigma(x)\,\sigma(y)$.
This means that neighboring sites now prefer to take {\em opposite} spins.

Suppose that $h=0$ and that the underlying graph is bipartite; this
means that $\L$ can be partitioned into two sets $\L_{even}$ and
$\L_{odd}$ such that sites in $\L_{even}$ only have edges to sites
in $\L_{odd}$, and vice versa. Clearly,
$\Z^d$ is an example of a bipartite graph. In this situation, we can
reduce the antiferromagnetic Ising model to the ferromagnetic case by
a simple spin-flipping trick: The bijection $\s\leftrightarrow\tilde{\s}$
of $\Omega$ defined by
\begin{equation} \label{flipping_the_odd_lattice}
\tilde{\s}(x)=\left\{
\begin{array}{ll}
\s(x) & \mbox{if } x\in \L_{even}\,, \\
-\s(x) & \mbox{if } x\in \L_{odd}
\end{array} \right.
\end{equation}
maps any Gibbs measure  for the antiferromagnetic
Ising model to a Gibbs measure for the ferromagnetic
Ising model with the same parameters, and vice versa. As a
consequence, a phase transition in the
antiferromagnetic model is equivalent to a phase transition in the
ferromagnetic model with the same parameters. Hence, Theorem
\ref{thm:Ising} immediately carries over to the antiferromagnetic
case.

The model becomes more interesting (or, at least, more genuinely
antiferromagnetic) if either $h\neq 0$ or the graph is taken to be
non-bipartite. Suppose first that $h\neq 0$ but still $\L=\Z^d$. If
$|h|$ is small and $\b$ sufficiently large, we have the same
picture as in the case $h=0$: there exist two distinct phases, one
having a majority of plus spins on the even sublattice and a
majority of minus spins on the odd sublattice, the other one having
a majority of plus spins on the odd sublattice and a majority of
minus spins on the even sublattice. We will show this in Section
\ref{ssect:ground_energy_perc}, Example \ref{ex:IsingAF}; see also
\cite{D2,Geo}. Note that this phase transition is somewhat
different in flavor compared to that in the ferromagnetic Ising
model: whereas in the Ising ferromagnet the phase transition
produces a breaking of a state-space symmetry, the phase transition
in the Ising antiferromagnet instead breaks the translation
symmetry between the sublattices $\L_{even}$ and $\L_{odd}$.

To see what happens in the case of a non-bipartite graph we consider
the triangular lattice which can be obtained
by taking the usual square lattice $\Z^2$ and adding
an edge between each vertex $x$ and its north-east neighbor $x+(1,1)$.
In this case,  one expects uniqueness when $h=0$, and
existence of three distinct phases when $|h|\neq 0$ is small and
$\b$ is large. Phase transitions in these models were studied in
\cite{D2,Hei}, for example.

\subsection{The Potts model}  \label{sect:Potts}

A natural generalization of the ferromagnetic Ising model is the
(ferromagnetic) Potts model \cite{Po}, in which each spin may take
$q\ge2$ (rather than only two) different values. The state space is
then $S=\{1,2,\ldots, q\}$, and the pair interaction is given by
\[
U(\s(x),\s(y))= 1-2I_{\{\s(x)=\s(y)\}}.
\]
We confine ourselves to the case of zero external field, so that
$V(\sigma(x))\equiv 0$. Taking $q=2$ and identifying the state
space $\{1,2\}$ with $\{-1,+1\}$ we reobtain the ferromagnetic
Ising model with zero external field. Just as in the latter case,
the Potts interaction favours configurations where many neighbor
pairs agree, and Theorem \ref{thm:Ising} can be extended to the
Potts model as follows, as we will show in Section
\ref{sect:PT_Potts}.
\begin{thm} \label{thm:Potts}
For the $q$-state Potts model on ${\bf Z}^d$, $d\geq 2$,
there exists a critical inverse temperature $\b_c\in(0,\infty)$
(depending on $d$ and $q$) such that for $\b<\b_c$ the model
has a unique Gibbs measure while for
$\b>\b_c$ there exist $q$ mutually singular Gibbs measures.
\end{thm}
In the same way as Theorem \ref{thm:Ising}, this theorem also holds
on general graphs provided we allow $\b_c$ to be 0 or $\infty$.
The  Potts model differs from the Ising model in that, for $q$
large enough, there are multiple Gibbs measures also at the
critical value $\b=\b_c$, as demonstrated by Koteck\'y and
Shlosman \cite{KS}; an outline of a proof will be given in Example
\ref{ex:Potts_first_order}. The Onsager critical value for the
two-dimensional Ising model is believed to extend to the
 Potts model on $\Z^2$ through the formula
$\b_c(q)=\frac{1}{2}\log(1+\sqrt{q})$; see Welsh
\cite{W1}, for example. This has so far only been established when $q$ is
sufficiently large \cite{LMMRS}.

\subsection{The hard-core lattice gas model} \label{sect:hard-core}

The hard-core lattice gas model (or hard-core model for short)
describes a gas of particles which can only sit on the lattice sites
but are so large that adjacent sites cannot be occupied
simultaneously. The state space is $S=\{0,1\}$, the pair interaction
\[
U(\s(x),\s(y))=\left\{
\begin{array}{ll}
\infty & \mbox{if }\s(x)=\s(y)=1 \\
0 & \mbox{otherwise,}
\end{array} \right.
\]
describes the hard core of the particles,
and the chemical potential is $V(\s(x))=-(\log \la)\;\s(x)$. Here
$\la>0$ is the so-called activity parameter.
The hard-core model shows some similarities to the Ising
antiferromagnet in an external field and can, in fact, obtained from
it by a limiting procedure ($\b\rightarrow
\infty$, $h\rightarrow 2d$, $\b (h-2d)=$ const, \cite{DKS}). Since
$U$ is either 0 or $\infty$, the inverse temperature $\b$ is
irrelevant and will thus be fixed as $1$, and we can vary only the
parameter $\la$. Finite volume Gibbs distribitions can
then be thought of as first letting all spins be independent, taking
values $0$ and $1$ with respective probabilities $\frac{1}{1+\la}$
and $\frac{\la}{1+\la}$, and then conditioning on the event that
no two $1$'s sit next to each other anywhere on the lattice.

The phase transition behavior of the hard-core model on $\Z^d$,
$d\geq 2$, is as follows. For $\la$ sufficiently close to 0,
the particles are spread out rather sparsely
on the lattice, and we get a unique Gibbs measure, just as in the
Ising antiferromagnet at high temperatures. When $\la$ increases,
the particle density also increases, and  the
system  finally starts looking for optimal packings of particles.
There are two such
optimal packings, one where all sites in $\L_{even}$ are occupied and
those in $\L_{odd}$ are empty, and one vice versa; we denote these
configurations by $\eta_{even}$ and $\eta_{odd}$, respectively.
(These chessboard configurations look similar to those favoured
in the Ising  antiferromagnet.) For sufficiently large $\la$,
the infinite volume construction of Gibbs measures with these
two choices of boundary condition produces different Gibbs measures,
so we get a phase transition \cite{D2}.
\begin{thm} \label{thm:hard-core}
For the hard-core model on ${\bf Z}^d$, $d\geq 2$,
there exist two constants $0<\la_c\leq \la_c'<\infty$
(depending on $d$) such that for $\la<\la_c$
the model has a unique Gibbs measure while for
$\la>\la_c'$ there are multiple Gibbs measures.
\end{thm}
This result will be proved in Section \ref{sect:WR_RC}. From a
computer-assisted proof \cite{Ra} we know that $\la_c
\geq 1.50762$.  It is widely believed that one should be able to take
$\la_c=\la_c'$ in this result, which would mean that the occurrence of
phase transition is increasing in $\la$. Such a result,
however, would (unlike Theorems \ref{thm:Ising} and
\ref{thm:Potts}) {\em not} extend to arbitrary graph structures;
some counterexamples were recently provided by Brightwell,
H\"aggstr\"om and Winkler \cite{BHW}.

The hard-core model analogue of introducing an external field in
the Ising model on $\Z^d$ is obtained by replacing the single activity
parameter $\la$ by two different activities
$\la_{even}$ and $\la_{odd}$, one for sites in $\L_{even}$
and the other for sites in $\L_{odd}$. By analogy with the Ising
model, one would expect to have a unique Gibbs measure as soon as
$\la_{even}\neq\la_{odd}$; this was conjectured by Van den
Berg and Steif \cite{vdBM} and proved for the case $d=2$ by
H\"aggstr\"om \cite{H3}.

\subsection{The Widom--Rowlinson lattice model} \label{sect:WR}

The Widom--Rowlinson model is another lattice gas model, where this
time there are two types of particles, and two particles are
allowed to sit on neighboring sites only if they are of the same
type. Actually, Widom and Rowlinson \cite{WR} originally introduced
it as a continuum model of particles living in $\R^d$; see Section
\ref{sect:cont_WR} below. The lattice variant described here was
first studied by Lebowitz and Gallavotti \cite{LG}. The state space
is $S=\{-1,0,+1\}$, where $-1$ and $+1$ are the two particle types,
and $0$'s correspond to empty sites. The pair interaction is given
by
\[
U(\s(x),\s(y))=\left\{
\begin{array}{cl}
\infty & \mbox{if } \s(x)\,\s(y)=-1\,, \\
0 & \mbox{otherwise,}
\end{array} \right.
\]
and the chemical potential by
\[
V(\s(x))=\left\{
\begin{array}{cll}
 -\log\la_- &\mbox{if} &\s(x)=-1 \,,\\
0 & \mbox{if} &\s(x)=0\,,\\
-\log\la_+ &\mbox{if} &\s(x)=+1\,.
\end{array} \right.
\]
Here $\la_-,\, \la_+>0$ are the activity parameters for the two
particle types $-1$ and $+1$. As in the hard-core model, we fix the
inverse tempreature $\b=1$ and only vary the activity parameters.
Gibbs measures can then be thought of as first picking all spins
independently, taking values $-1$, $0$ or $+1$ with probabilities
proportional to $\la_-$, $1$, and $\la_+$, and then conditioning
on the event that no two particles of different type sit next to each
other in the lattice.
We are mainly interested in the symmetric case $\la_-=\la_+=\la$,
where the phase transition behavior on ${\bf Z}^d$, $d\geq 2$
is similar to the Ising model: For $\la$ small, there is a unique
Gibbs measure in which the overall density of plus-particles is almost
surely equal to that of the minus-particles. For $\la$ sufficiently
large, the system wants to pack the particles so densely that the $\pm 1$
symmetry is broken. As for the Ising model, one can construct two
particular Gibbs measures
$\mu^+$ and $\mu^-$ using boundary conditions $\eta\equiv +1$ and
$\eta\equiv -1$; for small $\la$ we get $\mu_+=\mu_-$ whereas for
large $\la$ the two measures are different (and distinguishable through
the densities of the two particle types), producing a phase transition.
\begin{thm} \label{thm:WR}
For the Widom--Rowlinson model on ${\bf Z}^d$, $d\geq 2$,
with activities $\la_-=\la_+=\la$,
there exist $0<\la_c\leq \la_c'<\infty$
(depending on $d$) such that for $\la<\la_c$
the model has a unique Gibbs measure while for
$\la>\la_c'$ there are multiple Gibbs measures.
\end{thm}
As in the hard-core model, we expect that one should be able to take
$\la_c=\la_c'$, but such a monotonicity is not known.
Examples of graph structures where the desired monotonicity fails
can be found in \cite{BHW}.

We furthermore expect that the asymmetric Widom--Rowlinson model on
${\bf Z}^d$ with $\la_-\neq \la_+$ always has a unique Gibbs measure
(similarly to the Ising model with a nonzero external field), but this also is
not rigorously known.

\section{Coupling and stochastic domination} \label{sect:coupling_SD}

Geometry alone will not be sufficient for our analysis of equilibrium
phases.
We also need some probabilistic tools which allow us to compare different
configurations and different probability measures.
So we need to include another preparatory section
describing these tools and their basic applications to our setting.

{\em Coupling} is a probabilistic technique which has turned out to
be immensely useful in virtually all areas of probability theory,
and especially in its applications to statistical mechanics. The
basic idea is to define two (or more) stochastic processes jointly
on the same probability space so that they can be compared
realizationwise. This direct comparison often leads to conclusions
which would not be easily available by considering the processes
separately. Although an independent coupling is sometimes already
quite useful (as we will see in Section \ref{sect:appli}, for
example), it is usually more efficient to introduce a dependence
which relates the two processes in an efficient way. One such
particularly nice relationship is that one process is pointwise
smaller than the other in some partial order. This case is related
to the central concept of {\em stochastic domination}, via
Strassen's Theorem (Theorem \ref{thm:strassen}) below. We will
confine ourselves to those parts of coupling theory that are needed
for our applications; a more general account can for example be
found in the monograph by Lindvall \cite{Lin}.

\subsection{The coupling inequality}\label{sect:ineq}

In this and the next subsection of general character,
$\L$ will be an arbitrary finite or countably infinite set.
As the notation indicates, we think of the standard case that $\L$ is
the lattice introduced in Section \ref{sect:lattice}, but the
following results will  also be applied to the case when $\L$ is replaced
by its set $\B$ of bonds. We consider again the product space $\Om=S^\L$,
where for the moment $S$ is an arbitrary measurable space.

Suppose $X$ and $X'$ are random elements of $\Om$, and let
$\mu$ and $\mu'$ be their respective
distributions. We define the (half)
{\em total variation distance} $\|\mu-\mu'\|$
between $\mu$ and $\mu'$ by
\begin{equation}\label{vardi}
\|\mu- \mu'\| = \sup_{A\subset \Om} |\mu(A)-\mu'(A)|
\end{equation}
where $A$ ranges over all measurable subsets of $\Om$. The
coupling inequality below provides us with a convenient upper bound on
this distance. To state it we first need to
define what we mean by a coupling of $X$ and $X'$.
\begin{defn}
A {\bf coupling} $P$ of two $\Om$-valued random variables $X$ and $X'$, or of
their distributions $\mu$ and $\mu'$, is a probability measure
on $\Om\times \Om$ having marginals $\mu$ and $\mu'$,
in that for every event $A\subset \Om$
\begin{equation} \label{first_marginal}
P((\xi, \xi'):\, \xi\in A) = \mu(A)
\end{equation}
and
\begin{equation} \label{second_marginal}
P((\xi, \xi'):\, \xi'\in A) = \mu'(A)\;.
\end{equation}
\end{defn}
We think of a coupling as a redefinition of the random variables $X$ and
$X'$ on
a new common probability space such that their distributions are preserved.
Sometimes it will be convenient to keep the underlying probability space
implicit, but in general, as in  (\ref{first_marginal})
and (\ref{second_marginal}), we make the canonical choice, which is the
product space $\Om\times \Om$; $X$ and $X'$ are then simply the
projections on the two coordinate spaces. With this in mind, we write
$P(X\in A)$ and $P(X'\in A)$ for
 the left hand sides of (\ref{first_marginal})
and (\ref{second_marginal}), respectively.
In the same spirit,  $P(X = X')$ is a short-hand for $P((\xi, \xi'):\, \xi
= \xi')$.
\begin{prop}[The coupling inequality]
Let $P$ be a coupling of two $\Om$-valued random variables $X$ and $X'$,
with distributions $\mu$ and $\mu'$. Then
\begin{equation} \label{eq:coupling_inequality}
\|\mu- \mu'\| \leq P(X\neq X').
\end{equation}
\end{prop}
{\bf Proof:} For any $A\in \Om$, we have
\begin{eqnarray*}
\mu(A) - \mu'(A) & = & P(X\in A) - P(X' \in A) \\
& = & P(X\in A, X'\not\in A) - P(X\not\in A, X' \in A) \\
& \leq & P(X\in A, X'\not\in A) \\
& \leq & P(X\neq X')\;,
\end{eqnarray*}
whence (\ref{eq:coupling_inequality}) follows by symmetry. $\Cox$

\medskip\noindent
The next result states that there always exists some coupling which
achieves equality in (\ref{eq:coupling_inequality}). We call such a
coupling {\em optimal}.
\begin{defn}   \label{def:opt}
A coupling $P$ of two $\Om$-valued random variables $X$ and $X'$,
with distributions $\mu$ and $\mu'$, is said to be an {\bf optimal
coupling} if
\[
 \|\mu-\mu'\|=P(X\neq X')\,.
\]
\end{defn}
\begin{prop}
For any two $\Om$-valued random variables $X$ and $X'$, there exists
an optimal coupling of $X$ and $X'$.
\end{prop}
To construct an optimal coupling, one simply puts the common mass
$\mu\wedge\mu'$ of $\mu$ and $\mu'$ on the diagonal of $\Om\times\Om$ and
adds any measure with marginals $\mu-\mu\wedge\mu'$ and
$\mu'-\mu\wedge\mu'$; the simplest choice of such a
measure is the appropriately scaled product measure.
For details we refer to \cite{Lin}, where such a coupling is called the
$\gamma$-coupling of $\mu$ and $\mu'$. This construction shows, in
particular, that
 the optimal coupling is in general not unique. Applications of
optimal couplings to
interacting particle systems can be found in \cite{L,Ma}, for example.

\subsection{Stochastic domination} \label{sect:stoch_dom}

Suppose now that $S$ is a closed subset of $\R$, so that $S$ is linearly
ordered.
The product space $\Om$ is then equipped with a natural partial order
$\preceq$
which is defined coordinatewise: For $\xi, \xi'\in\Om$, we write
$\xi\preceq \xi'$ (or $\xi' \succeq \xi$) if
$\xi(x) \leq \xi'(x)$ for every $x\in \L$.  A function
$f:\, \Om\rightarrow \R$ is said to be {\em increasing}
(or, non-decreasing)
if
$f(\xi)\leq f(\xi')$
whenever $\xi\preceq \xi'$. An event $A$ is said to be increasing if its
indicator function $I_A$ is increasing. The following standard definition
of stochastic domination expresses the fact that $\mu'$ prefers larger
elements of $\Om$ than $\mu$.
\begin{defn}
Let   $\mu$ and $\mu'$ be two probability measures on $\Om$. We say that
$\mu$ is {\bf  stochastically dominated} by $\mu'$, or
$\mu'$ is stochastically larger than $\mu$, writing  $\mu\leqd \mu'$, if
for every bounded increasing observable
$f:\, \Om\rightarrow \R$ we have $\mu(f)\leq \mu'(f)$.
\end{defn}
In the one-dimensional case when  $|\L|=1$ and $\Om=S\subset\R$, this
definition is equivalent to the classical requirement that
$\mu([r,\infty)) \leq  \mu'([r,\infty))$ for all $r\in \R$.
The following fundamental result of Strassen \cite{Str} characterizes
stochastic domination in coupling terms.
\begin{thm}[Strassen] \label{thm:strassen}
For any two probability measures $\mu$ and $\mu'$ on $\Om$, the following
statements are equivalent.
\begin{description}
\item{\rm (i)} $\mu\leqd \mu'$
\item{\rm (ii)} For all {\em continuous} bounded increasing functions
$f:\Om\to\R$,  $\mu(f)\leq \mu'(f)$.
\item{\rm (iii)} There exists a coupling $P$ of $\mu$ and $\mu'$ such that
$P(X\preceq X')=1$.
\end{description}
\end{thm}
{\bf Sketch of proof: } While the implications (i) $\Rightarrow$
(ii) and (iii) $\Rightarrow$ (i) are trivial, the assertion (ii)
$\Rightarrow$ (iii) is too deep to be explained here in detail. To
start one should note that $\mu$ and $\mu'$ may be considered as
measures on the compact space $[-\infty,\infty]^\L$, and one can
then follow the arguments outlined in \cite{L}, pp.\ 72 ff., or
\cite{Lin}. Further discussion of Strassen's theorem can be found
in \cite{Lin,KKO}. $\Cox$ \medskip

\noindent The  equivalence (i) $\Leftrightarrow$ (ii) in Theorem
\ref{thm:strassen} readily implies the following corollary.
\begin{cor}\label{cor:stoch_dom_preserved_under_limits}
The relation $\leqd$ of stochastic domination is preserved under weak
limits.
\end{cor}
Next we recall a famous sufficient condition for stochastic domination.
This condition is (essentially) due
 to Holley \cite{Ho} and refers to the finite-dimensional case when
$|\L|<\infty$. We also assume for simplicity that $S\subset\R$ is finite.
Hence $\Om$ is finite. In this case, a probability measure $\mu$ on $\Om$
is called
{\em irreducible} if the set $\{\eta\in\Om:\, \mu(\eta)>0\}$ is connected
in the sense that any element of $\Om$ with positive $\mu$-probability
can be reached from any other via successive coordinate changes without
passing through elements with zero $\mu$-probability.
\begin{thm}[Holley] \label{thm:holley}
Let $\L$ be finite, and let $S$ be a finite subset of ${\bf R}$. Let
$\mu$ and $\mu'$ be probability measures on $\Om$. Assume that $\mu'$ is
irreducible and assigns
positive probability to the maximal
element of $\Om$ (with respect to $\preceq$). Suppose further
that
\begin{equation} \label{eq:conditional_monotone}
\mu(X(x)\geq a\, | \, X=\xi\mbox{ \rm off } x) \leq
\mu'(X(x)\geq a\, | \, X=\eta\mbox{ \rm off } x)
\end{equation}
whenever  $x\in \L$, $a \in S$, and $\xi, \eta\in S^{\L\setminus\{x\}}$
are such that $\xi\preceq\eta$,
$\mu(X=\xi\mbox{ \rm off } x)>0$ and $\mu'(X=\eta \mbox{ \rm off } x)>0$.
Then $\mu\leqd \mu'$.
\end{thm}
{\bf Proof:} Consider a Markov chain $(X_k)_{k=0}^\infty$ with
state space $\Om$ and transition probabilities defined as follows.
At each integer time $k\geq 1$, pick a random site $x\in \L$
according to the uniform distribution. Let $X_k=X_{k-1}$ on
$\L\setminus\{x\}$, and select $X_k(x)$ according to the
single-site conditional distribution prescribed by $\mu$. This is a
so-called Gibbs sampler for $\mu$, and it is immediate that if the
initial configuration $X_0$ is chosen according to $\mu$, then
$X_k$ has distribution $\mu$ for each $k$. Define a similar Markov
chain $(X'_k)_{k=0}^\infty$ with $\mu$ replaced by $\mu'$.

Next, define a coupling of $(X_k)_{k=0}^\infty$ and
$(X'_k)_{k=0}^\infty$ as follows. First pick the initial values
$(X_0,X^\prime_0)$ according to the product measure
$\mu\times \mu'$. Then, for each
$k$, pick a site $x\in \L$ at random and let $U_k$ be an independent
random variable, uniformly distributed on the interval $[0,1]$. Let
$X_k(y)=X_{k-1}(y)$ and $X'_k(y)=X'_{k-1}(y)$ for each
site $y \neq x$, and update the values at site $x$ by letting
\[
X_k(x)=\max\{a\in S: \mu(X(x)\geq a\, | \, X=\xi\mbox{ off } x)\geq
U_k\} \]
and
\[
X'_k(x)=\max\{a\in S: \mu'(X'(x)\geq a\, | \,
X'=\eta\mbox{ off } x)\geq U_k\}
\]
where $\xi=X_{k-1}(\L\setminus\{x\})$ and
$\eta=X'_{k-1}(\L\setminus\{x\})$. It is clear that this construction
gives the correct marginal behaviors of $(X_k)_{k=0}^\infty$ and
$(X'_k)_{k=0}^\infty$.
The assumption (\ref{eq:conditional_monotone}) implies that
$X_k\preceq X'_k$ whenever $X_{k-1}\preceq X'_{k-1}$.
By the irreducibility of $\mu'$,
the  chain $(X'_k)_{k=0}^\infty$ will almost surely hit the maximal
state of $\Om$ at some finite (random) time, and from
this time on we will thus have $X_k\preceq X'_k$. Since the coupled chain
 $(X_k,X'_k)_{k=0}^\infty$ is a finite state aperiodic
Markov chain,  $(X_k,X'_k)$ has a limiting distribution as
$k\rightarrow\infty$. Picking $(X,X')$ according to this limiting
distribution gives a coupling of $X$ and $X'$ such that
$X\preceq X'$ almost surely, whence $\mu \leqd \mu'$ by Theorem
\ref{thm:strassen}. $\Cox$

\medskip\noindent
A non-dynamical proof of Holley's inequality by induction on $|\L|$,
together with an extension to non-finite $S$, was given by Preston
\cite{Pr2}; the simplest induction proof of an even more general result
can be found in \cite{BaBo}.

As a consequence of Holley's inequality we obtain the celebrated FKG
inequality (Theorem \ref{thm:FKG} below)
of  Fortuin, Kasteleyn and Ginibre \cite{FKG}, who
stated it under slightly different conditions. It concerns the correlation
structure in a single probability measure
rather than a comparison between two probability measures.
\begin{defn}  \label{def:mon}
A probability measure $\mu$ on $\Om$ is called  {\bf monotone} if
\begin{equation} \label{eq:FKG_property}
\mu(X(x)\geq a\, | \, X=\xi\mbox{ \rm off } x) \leq
\mu(X(x)\geq a\, | \, X=\eta\mbox{ \rm off } x)
\end{equation}
whenever $x\in \L$, $a\in S$, and $\xi, \eta\in S^{\L\setminus\{x\}}$ are
such that $\xi\preceq\eta$,
$\mu(X=\xi\mbox{ \rm off } x)>0$ and $\mu(X=\eta\mbox{ \rm off } x)>0$.
\end{defn}
Intuitively, $\mu$ is monotone if the spin at a site $x$ prefers to take
large values
whenever its surrounding sites do.
\begin{defn} \label{defn:positive_correlations}
A probability measure $\mu$ on $\Om$ is said to have {\bf positive
correlations} if for all bounded increasing functions
$f,g:\, \Om \rightarrow \R$ we have
\[
\mu(fg) \geq \mu(f)\,\mu(g)\,.
\]
\end{defn}
Since the preceding inequality is preserved under rescaling and
addition of constants to $f$ and $g$, $\mu$ has positive
correlations whenever $\mu \leqd \mu'$ for any probability measure
$\mu'$ with bounded increasing Radon--Nikodym density  relative to
$\mu$. Theorem \ref{thm:strassen} thus
shows that $\mu$ has positive correlations whenever $\mu(fg) \geq
\mu(f)\,\mu(g)$ for all {\em continuous} bounded increasing
functions $f$ and $g$. Hence, the property of positive correlations
is also preserved under weak limits.
\begin{thm}[The FKG inequality] \label{thm:FKG}
Let $\L$ be finite,  $S$ a finite subset of ${\bf R}$, and $\mu$ a
probability measure on $\Om$ which is irreducible and assigns
positive probability to the maximal element of $\Om$ (relative to
$\preceq$). If $\mu$ is monotone, it also has positive correlations.
\end{thm}
{\bf Proof:}
Suppose $\mu'$ is a second probability measure on $\Om$ such that
$\mu'(\eta)={\mu(\eta)g(\eta)}$
for all $\eta\in \Om$ and some positive increasing function $g$. For $x\in
\L$, $a\in S$ and $\xi\in S^{\L\setminus\{x\}}$ such
that $\mu(X=\xi\mbox{ off } x)>0$ we write
$q_x(a,\xi)= \mu(X(x)\geq a \, | \, X=\xi\mbox{ off } x)$ and define
$q_x'(a,\xi)$ similarly in terms of $\mu'$. Then
\begin{eqnarray*}
{ q_x'(a,\xi)}\,\big/\,\big({1-q_x'(a,\xi)\big)} & = &
 \sum_{s\geq a}\mu(\xi^{x,s})\,g(\xi^{x,s})\bigg/
{\sum_{s< a}\mu(\xi^{x,s})\,g(\xi^{x,s})} \\
& \geq & {\sum_{s\geq a}\mu(\xi^{x,s})}\bigg/{\sum_{s<a}
\mu(\xi^{x,s})} \\
& = & {q_x(a,\xi)}\,\big/\,\big({1-q_x(a,\xi)}\big).
\end{eqnarray*}
Together with  assumption (\ref{eq:FKG_property}), this implies that $\mu$
and $\mu'$ satisfy (\ref{eq:conditional_monotone}).
Theorem \ref{thm:holley} thus implies that
$\mu \leqd \mu'$, and the corollary follows. $\Cox$

\bigskip\noindent
Finally we state a simple observation showing that, under the condition of
stochastic domination, the equality of the
single-site marginal distributions already implies the equality of the
whole probability measures.
\begin{prop} \label{prop:single_site+dom}
Let $\cal L$ be finite or countable, and let $\mu$ and $\mu'$ be
two probability measures on $\Om=\R^\L$ satisfying $\mu\leqd\mu'$. If,
in addition,
$\mu(X(x)\leq r)=\mu'(X(x)\leq r)$ for all $x\in \L$ and $r\in \R$ then
$\mu=\mu'$.
\end{prop}
{\bf Proof:}
Let $P$ be a coupling of $\mu$ and $\mu'$
such that $P(X\preceq X')=1$ which exists by Theorem \ref{thm:strassen}.
Writing ${\bf Q}$ for the set of
rational numbers, we have for each $x\in \L$
\begin{eqnarray*}
P(X(x)\neq X'(x)) & = & P(X(x) < X'(x)) \\[1ex]
& \leq & \sum_{r\in{\bf Q}} P(X(x)\leq r ,\ X'(x)> r) \\
& = & \sum_{r\in{\bf Q}} \bigg(P(X(x)\leq r)-P(X'(x)\leq r)\bigg) \\
& = & 0.
\end{eqnarray*}
Summing over all $x\in\L$ we get $P(X\neq X')=0$,
whence $\mu=\mu'$ by \rf{eq:coupling_inequality}. $\Cox$

\subsection{Applications to the Ising model} \label{sect:Ising_monotone}

We will now apply the results of the previous subsection to the
ferromagnetic Ising model. Let $(\L,\sim)$ be any infinite
locally finite graph.
For definiteness, one may think of the case $\L=\Z^d$; the
arguments are, however, independent of the particular graph structure.

As in Section \ref{sect:Ising}, we write $\mu^\eta_{h,\beta,\La}$ for the
Gibbs distribution in  a finite region $\La$ with boundary condition
$\eta\in\Om$ and external field $h\in\R$ at inverse temperature $\beta>0$.
Our first application of Holley's
theorem
asserts that if one boundary condition dominates another, then we also have
stochastic domination between the corresponding finite volume Gibbs
distributions.
\begin{lem} \label{lem:boundary_domination}
If the boundary conditions $\xi, \eta\in \Omega$ satisfy $\xi \preceq
\eta$,
then
\[
\mu^\xi_{h,\beta,\Lambda} \leqd \mu^\eta_{h,\beta,\La} \, .
\]
Also, each $\mu^\xi_{h,\beta,\Lambda}$ has positive correlations.
\end{lem}
{\bf Proof:} The conditional probability of having a plus spin at a given
site $x$ given the configuration $\xi$ everywhere else is equal to
$(1+\exp[-2\b(h+\sum_{y:y\sim x}\xi(y))])^{-1}$, which is an increasing
function of $\xi$.  Theorem \ref{thm:holley} and Theorem \ref{thm:FKG}
thus imply stochastic domination between the projections of the Gibbs
distributions to $S^\La$ and the positive correlations property.
As their behavior outside $\La$ is
deterministic, the lemma follows. $\Cox$

\medskip\noindent
We  write $\mu^+_{h,\beta,\La}$ and $\mu^-_{h,\beta,\La}$ for the
finite volume Gibbs distributions obtained with respective boundary
conditions
$\eta\equiv +1$ and $\eta \equiv -1$.  Lemma
\ref{lem:boundary_domination} then shows that
\begin{equation}\label{eq:sandwich}
\mu^-_{h,\beta,\La} \leqd
\mu^\eta_{h,\beta,\La} \leqd \mu^+_{h,\beta,\La}
\end{equation}
for any $\eta\in\Omega$. This sandwich inequality reveals the special role
played by the ``all plus'' and ``all minus''
boundary conditions.
Next we establish the existence of the limiting ``plus measure''
discussed in Section \ref{sect:Ising}. We will say that a measure $\mu$ on
$\Om$ is {\em homogeneous} if it is invariant under all graph
automorphisms of $(\L,\sim)$. In the case $\L=\Z^d$, a homogeneous measure
is thus invariant under translations,
lattice rotations, and reflections in the axes. We write $\La\uparrow\L$
for the limit along an arbitrary increasing sequence of finite regions
which exhaust the full graph $\L$.

\begin{prop} \label{prop:plus_measure}
The limiting probability measure
\[
\mu^+_{h,\beta}=\lim_{\La\uparrow\L}\mu^+_{h,\beta,\La}
\]
exists. $\mu^+_{h,\beta}$ is a homogeneous Gibbs measure for the Ising
model on $(\L,\sim)$ with external field $h$ and inverse temperature
$\beta$ and has positive correlations.
\end{prop}
{\bf Proof:} By the general theory in Section \ref{sect:Gibbsmeasures},
the limit
is a Gibbs measure whenever it exists. Also, by Lemma
\ref{lem:boundary_domination} the limit must have positive correlations.
To show the existence of the limit we note that
\begin{equation} \label{eq:decreasing}
\mu^+_{h,\beta,\La} \geqd \mu^+_{h,\beta,\Delta}\quad\mbox{ whenever
}\La\subset\Delta\;.
\end{equation}
This follows from Lemma \ref{lem:boundary_domination} because
$\mu^+_{h,\beta,\La} $ is obtained from $\mu^+_{h,\beta,\Delta}$ by
conditioning on the increasing event that $X\equiv +1$ on
$\Delta\setminus\La$.

Now, for any finite $A\subset\La$, if $\La$ increases, then, by
(\ref{eq:decreasing}), $\mu^+_{h,\beta,\La}(X\equiv +1\mbox{ on
}A)$ decreases and therefore converges to $\
\inf_{\Delta}\mu^+_{h,\beta,\Delta} (X\equiv +1\mbox{ on }A)\,$.
Note that this limit is obviously invariant under any automorphism
of $(\L,\sim)$. By inclusion-exclusion it follows that, for any
local observable $f$,  $\mu^+_{h,\beta,\La}(f)$ converges to an
automorphism invariant limit as $\La\uparrow\L$. These limits
determine a unique homogeneous probability measure
$\mu^+_{h,\beta}$ which, as a weak limit of finite volume Gibbs
distributions, is a Gibbs measure. The lemma is thus proved.
$\Cox$

\medskip\noindent
Obviously, replacing the ``all plus'' boundary condition by the
``all minus'' boundary condition, we obtain in the same way an
automorphism invariant infinite volume Gibbs measure
$\mu^-_{h,\beta}$ with positive correlations.   In the same way,
Lemma \ref{lem:boundary_domination} shows that {\em any extremal Gibbs
measure has positive correlations\/} (since it is a weak limit of $
\mu^\eta_{h,\beta,\La}$ for suitable $\eta$). However, positive
correlations may fail for suitable $(\L,\sim)$ and some particular
non-extremal Gibbs measures. For instance, when $\L={\bf Z}^3$ and
$\beta$ is sufficiently large, one can take a convex combination of
two different so-called Dobrushin states; see \cite{D72,vBe}.

We now take the limit in the sandwich inequality \rf{eq:sandwich}. Let
$\mu$ be any Gibbs measure for the Ising model with parameters $h$ and
$\beta$.
Taking the mean $\int\mu(d\eta)$ in \rf{eq:sandwich}, we obtain that
$\mu^-_{h,\beta,\La} \leqd\mu \leqd \mu^+_{h,\beta,\La}$, and since
stochastic domination is preserved under weak limits, we end up with
\begin{equation} \label{eq:Ising_sandwich}
\mu^-_{h,\beta} \leqd
\mu \leqd \mu^+_{h,\beta}
\end{equation}
when $\mu$ is any Ising-model Gibbs measure for $\beta,\,h$.
On the one hand, this shows that $\mu^-_{h,\beta}$ and $\mu^+_{h,\beta}$
are extremal, and thus equilibrium phases in the sense of Section
\ref{sect:phases}.
On the other hand, we obtain
 an efficient criterion for the existence of a phase transition which was
first observed by Lebowitz and Martin-L\"of \cite{LML} and Ruelle
\cite{Ru3}.
\begin{thm} \label{thm:Ising_sandwich}
For the Ising model on an infinite locally finite graph $(\L,\sim)$ with
external field $h\in\R$ and inverse temperature $\beta$, the following
statements are equivalent.
\begin{description}
\item{\rm (i)} There is a unique infinite volume Gibbs measure.
\item{\rm (ii)} $\mu^-_{h,\beta}=\mu^+_{h,\beta}$
\item{\rm (iii)} $\mu^-_{h,\beta}(X(x)=+1)=\mu^+_{h,\beta}(X(x)=+1)$ for
all $x\in\L$.
\end{description}
\end{thm}
{\bf Proof:} The implications (i) $\Rightarrow$ (ii) $\Rightarrow$ (iii)
are immediate. (iii) $\Rightarrow$ (ii) follows directly from
\rf{eq:Ising_sandwich} and Proposition
\ref{prop:single_site+dom}, and (ii) $\Rightarrow$ (i) from
\rf{eq:Ising_sandwich}. $\Cox$

\medskip\noindent
{\bf Remarks:} (a) In the case $h=0$ of no external field,
assertion (iii) is equivalent to $\mu^+_{h,\beta}(X(x)=+1)=1/2$ for
all $x\in\L$, by the $\pm$ symmetry of the model. An extension of
Theorem \ref{thm:Ising_sandwich} in this case to the $q$-state
Potts model will be given in Theorem \ref{thm:Potts_sandwich}.

(b) If the graph automorphisms act transitively on $(\L,\sim)$
then, by homogeneity, assertion (iii) is equivalent to having the
equation only for {\em some} $x\in\L$. Using for example the
random-cluster methods of Section \ref{sect:random-cluster} one can
obtain the same equivalence also in the general case, assuming only
that $\L$ is connected.

(c) If $\L=\Z^d$, (iii) is equivalent to the condition that the free energy
 density is differentiable with respect to $h$ at the given values of $h$
and $\beta$ \cite{LML}. By the celebrated Lee-Yang circle theorem
\cite{Ru4}, this is the case whenever $h\ne0$. Alternatively, one can use
the so-called GHS inequality to
establish (iii) for $h\ne0$ \cite{Pr3}. Hence, for non-zero external
field the Ising model on $\Z^d$ does not exhibit a phase transition.

\medskip\noindent
We conclude this subsection comparing the Ising-model ``plus'' measures
for different values of the parameters. (For similar results in a lattice
gas setting see \cite{LePe}.)
\begin{prop} \label{prop:Ising_external_field}
Consider the Ising model on an arbitrary graph $(\L,\sim)$ at two inverse
temperatures
$\beta_1,\,\beta_2$ and two external fields $h_1,\, h_2$. Suppose that
either $\beta_1=\beta_2$ and $h_1\le h_2$, or $(\L,\sim)$ is of bounded
degree $N=\sup_{x\in\L}N_x$ and
$\beta_2h_2\ge\beta_1h_1+N|\beta_1-\beta_2|$. Then
\[
\mu^+_{h_1,\beta_1} \leqd \mu^+_{h_2,\beta_2} \, .
\]
\end{prop}
{\bf Proof:} The stated conditions imply that
$\mu^\xi_{h_1,\beta_1,\{x\}}(+1) \le \mu^\eta_{h_2,\beta_2,\{x\}}(+1)$
whenever $\xi\preceq\eta$. Hence, by Theorem \ref{thm:holley},
$
\mu^+_{h_1,\beta,\Lambda} \leqd \mu^+_{h_2,\beta,\Lambda}
$
for all $\La$, and the proposition follows by letting $\La\uparrow\L$.
$\Cox$

\medskip\noindent
If $\L$ has bounded degree $N$, we obtain in particular a comparison of
Ising and Bernoulli measures. Let $\psi_p$
denote the Bernoulli measure on $\{-1,1\}^\L$ with density $p$. Then
$\mu^+_{h,\beta}\to\psi_p$ as $\beta\to 0$ and $\beta
h\to\frac{1}{2}\log\frac{p}{1-p}$. The preceding proposition and Corollary
\ref{cor:stoch_dom_preserved_under_limits} thus show that
\[
\mu^+_{h,\beta}\geqd\psi_p\ \mbox{ for }\
h\ge\frac{1}{2\beta}\log\frac{p}{1-p}+N
\]
and
\[
\mu^+_{h,\beta}\leqd\psi_p\ \mbox{ for }\
h\le\frac{1}{2\beta}\log\frac{p}{1-p}-N\;.
\]

\subsection{Application to other models} \label{sect:other_appl}

Do the arguments of the previous subsection extend to the other models of
Section 3? The answer to this question
is different for the different models. For the {\em Potts model\/} with
$q\geq 3$, Lemma \ref{lem:boundary_domination} fails because the
conditional
probability that the spin at a site $x$ takes a large value is not
increasing in the surrounding spin configuration. Nevertheless,
the Potts model admits some analogues of Proposition
\ref{prop:plus_measure} and
Theorem \ref{thm:Ising_sandwich}. These results are deeper than their
Ising counterparts, and will be demonstrated by random-cluster
arguments in Section \ref{sect:random-cluster}.

The {\em Widom--Rowlinson model \/} exhibits the same monotonicity properties
as the Ising model. We thus obtain
Widom--Rowlinson analogues of Lemma \ref{lem:boundary_domination} and
Proposition
\ref{prop:plus_measure}. Fixing any two activity parameters
$\lambda_+,\,\lambda_->0$ and
 writing $\mu^+$ (resp.\ $\mu^-$)
for the associated limiting Gibbs measures with ``all plus'' (resp.\
``all minus'') boundary conditions, we find that any other Gibbs measure
for the same parameters is sandwiched
(in the sense of (\ref{eq:Ising_sandwich})) between
$\mu^+$ and $\mu^-$. The analogue of Theorem \ref{thm:Ising_sandwich}
reads as
follows.
\begin{thm} \label{WR_sandwich}
For the Widom--Rowlinson model on an infinite locally finite graph
$(\L,\sim)$
with activity parameters $\lambda_+,\,\lambda_->0$, the following
statements are
equivalent.
\begin{description}
\item{\rm (i)} There is a unique infinite volume Gibbs measure.
\item{\rm (ii)} $\mu^+=\mu^-$
\item{\rm (iii)}
$\mu^+(X(x)=+1)=\mu^-(X(x)=+1)$
for all $x\in\L$.
\end{description}
\end{thm}

\noindent
The {\em Ising antiferromagnet \/} and the {\em  hard-core lattice
gas model \/} are far from
satisfying the monotonicity properties needed for the
arguments in Section \ref{sect:Ising_monotone};  the conditional
probability that a site $x$ takes the value $+1$
is {\em decreasing}, rather than
increasing, in the surrounding configuration. However, when $(\L,\sim)$ is
bipartite (as in the case $\L=\Z^d$) we can use again the trick of
(\ref{flipping_the_odd_lattice}) to flip all spins on the odd sublattice.
The Ising antiferromagnet is then mapped onto the Ising ferromagnet with a
staggered external field (having
alternating signs on the even and the odd sublattices).
Similarly, the hard-core model is mapped into a model which also exhibits
the  necessary monotonicity properties.
We thus obtain analogous results which we spell out only for the hard-core
model  with activity $\lambda>0\,$: There exist two particular Gibbs
measures $\mu_\lambda^{even}$ and $\mu_\lambda^{odd}$, obtained as infinite
volume limits of finite volume Gibbs distributions with respective boundary
conditions $\eta_{even}$ and $\eta_{odd}$, defined in
Section \ref{sect:hard-core}. In terms of these two Gibbs measures,
the existence of a phase transition can be characterized as follows.
\pagebreak
\begin{thm} \label{hard-core_sandwich}
For the hard-core model on an infinite locally finite bipartite graph
$(\L,\sim)$
with activity parameter $\lambda$, the following are
equivalent.
\begin{description}
\item{\rm (i)} There is a unique infinite volume Gibbs
measure.
\item{\rm (ii)} $\mu^{even}_\lambda=\mu^{odd}_\lambda$
\item{\rm (iii)}
$\mu^{even}_\lambda(X(x)=1)=\mu^{odd}_\lambda(X(x)=1)$
for all $x\in\L$.
\end{description}
\end{thm}

\section{Percolation} \label{sect:percolation}

We will now introduce the ideas of random geometry we referred to
in the title of this paper. As these were developed first in the framework
of percolation theory, we devote this section to a description
of this subject. We will start with the classical case of
independent, or Bernoulli, percolation, and will then proceed to the case
of dependent percolation. In the subsequent sections we will see how
these results and ideas can be used for the geometric analysis of
equilibrium phases.

\subsection{Bernoulli percolation} \label{sect:bond_percolation}

Bernoulli percolation was introduced in the 1950's
in papers by Broadbent and Hammersley \cite{BH,Ham1,Ham2} as a model for
 the passage of a fluid through a porous medium. In fact, the model has
appeared first in \cite{Flory} in the context of polymerization.
We give here only a brief introduction; a thorough treatment
can be found in the
the books and lectures by Grimmett \cite{Gr1,Gr2}, Chayes and Chayes \cite{CCh},
and Kesten \cite{Kes}.

The porous medium is modelled by a graph
$(\L,\sim)$, and either the sites or the bonds of this graph are considered
to be randomly open or closed (blocked).
We begin with the case of site percolation;
the alternative case of random bonds will be discussed at the end of this
subsection. 

The basic question of percolation theory is how a fluid can spread through
the medium. This involves the connectivity properties of the set of
open vertices. To describe this we introduce some terminology.
A finite {\em path} is a sequence  $(v_0, e_1, v_1, e_2, $ $\ldots,e_k,
v_k)$,
where $v_0, \ldots, v_k\in\L$ are pairwise distinct vertices and $e_1,
\ldots, e_k\in\B$
are pairwise distinct edges such that, for each $i\in \{1,\ldots, k\}$, the
edge $e_i$ connects the vertices $v_{i-1}$ and $v_i$. Obviously, a path is
equivalently described by its sequence $(v_0, v_1, \ldots, v_k)$ of
vertices or
its sequence $(e_1, e_2, \ldots, e_k)$ of edges. The number $k$ is called
the {\em length} of the path.
In the same way, we can also speak of infinite paths
$(v_0, e_1, v_1, e_2, \ldots)$ and doubly infinite paths
$(\ldots, v_{-1}, e_0, v_0, e_1, v_1, \ldots)$. A region $C\subset\L$
is called {\em connected} if for any $x,y\in C$ there exists a path which
starts at $x$, ends at $y$, and which only contains vertices in $C$.

An {\em open path} is a path  on which all vertices are open. An
{\em open cluster} is a maximal connected set $C$ in which all
vertices are open; here maximal means that there is no larger
region $C'\supset C$ which is connected and only contains open
vertices. An infinite open cluster (or infinite cluster, for short)
is an open cluster containing infinitely many vertices. Using these
terms, we may say that the existence of an infinite open cluster is
equivalent to the fact that a fluid can wet a macroscopic part of
the medium.

We now turn to the classical case of Bernoulli site percolation with
retention parameter $p\in [0,1]$. In this case, each
vertex of $\L$,  {\em independently of all others}, is declared to be open
(and represented by the value $1$) with probability $p$ and closed
(with value $0$) with probability $1-p$. We write $\psi_p$ for the
associated (Bernoulli) probability measure on the configuration space
$\{0,1\}^\L$.

The first question to be asked is whether or not
infinite clusters can exist. This depends, of course,  on both the graph
$(\L,\sim)$ and the parameter $p$. The basic observation is the following.
\begin{prop} \label{prop:perc_critical_value}
For Bernoulli site percolation on an infinite locally finite graph
\linebreak$(\L,\sim)$, there exists a critical value $p_c\in [0,1]$ such
that
\[
\psi_p(\exists \mbox{ \rm an infinite open cluster})= \left\{
\begin{array}{ll}
0 & \mbox{if } p<p_c \\
1 & \mbox{if } p>p_c.
\end{array} \right.
\]
At the critical value $p=p_c$, the $\psi_p$-probability of having an
infinite
open cluster is either $0$ or $1$.
\end{prop}
{\bf Proof:} A moment's thought reveals that the existence of infinite
clusters is invariant under a change of the status of finitely many
vertices. By Kolmogorov's zero-one law, the $\psi_p$-probability of having
infinite clusters is therefore either
$0$ or $1$. It remains to show that this probability is increasing in $p$.
The existence of an infinite open cluster is obviously an increasing event,
so we are done if we can show that
\begin{equation} \label{eq:perc_monotone}
\psi_{p_1} \leqd \psi_{p_2} \ \mbox{ whenever } p_1\leq p_2.
\end{equation}
 This is intuitively obvious and can be proved
by the following elementary coupling argument.
Let $Y=(Y(x))_{x\in\L}$ be a family of i.i.d.\ random variables with
uniform distribution on $[0,1]$, and for $p\in[0,1]$ let
$X_p=(X_p(x))_{x\in\L}$
be defined by $X_p(x)=I_{\{Y(x)\le p\}}$. It is then clear that $X_p$ has
distribution $\psi_p$ and
$X_{p_1} \preceq X_{p_2}$ whenever $p_1\leq p_2$. This implies
(\ref{eq:perc_monotone}) by (the trivial part of)
Theorem \ref{thm:strassen}. (Note that we have in fact constructed a
simultaneous coupling of all $\psi_p$'s, and that this construction would
be used in Monte Carlo simulations of Bernoulli percolation.) $\Cox$

\medskip\noindent
We write $\{x \leftrightarrow \infty\}$ for the event that $x\in \L$
belongs to an
infinite cluster, and set $\theta_x(p)=\psi_p(x \leftrightarrow \infty)$.
For homogenous graphs such as $\L=\Z^d$,
 $\theta_x(p)$ does not depend on $x$, and then  we write simply
$\theta(p)$.
Equation \rf{eq:perc_monotone}
shows that $\theta_x(p)$ is
increasing in $p$. We also make the following observation.
\begin{prop} \label{prop:perc_single_site}
For any infinite locally finite connected graph $(\L,\sim)$, any $x\in \L$
and any $p\in (0,1)$,
we have $\theta_x(p)>0$ if and only if
$\psi_p(\exists \mbox{ \rm an infinite open cluster})=1$.
\end{prop}
{\bf Proof:} From Proposition \ref{prop:perc_critical_value} we
know that  an infinite cluster exists with probability $0$ or $1$.
The implication ``only if''  is therefore immediate. For the ``if''
part, we note that if  an infinite cluster exists with positive
probability, then there is some $N$ such that $\psi_p(A_N)>0$,
where $A_N$ is the event that some vertex within distance $N$ from
$x$ belongs to an infinite cluster. On the other hand we have
$\psi_p(B_N)>0$, where $B_N$ is the event that all vertices within
distance $N-1$ from $x$ are open. The event $A_N\cap B_N$ implies that $x$
belongs to an infinite cluster.
 But $A_N$ and $B_N$ are increasing events (see Section
\ref{sect:stoch_dom}),
so that we can apply Theorem \ref{thm:FKG} to obtain $\theta_x(p)\geq
\psi_p(A_N\cap B_N) \geq\psi_p(A_N)
\psi_p(B_N)>0$.
 $\Cox$

\medskip\noindent
Note that the proof above applies to the
much broader class of all measures  with positive
correlations (recall Definition \ref{defn:positive_correlations}), rather
than only the Bernoulli measures.

Next we ask whether both possibilities in Proposition
\ref{prop:perc_single_site}
really occur, that is, if $0<p_c<1$. For, only in this case we really have
a {\em nontrivial} critical
phenomenon at $p_c$.  The answer depends on the graph.
For $\L=\Z^d$ with dimension $d\geq 2$, the threshold $p_c$ is indeed
nontrivial, as is
stated in the theorem below. This nontriviality
of $p_c$ is a fundamental ingredient of many of the stochastic-geometric
arguments employed later on. On the other hand, it is easy
to see that $p_c=1$ for $\L=\Z^1$.
\begin{thm} \label{thm:perc_nontrivial}
The critical value $p_c=p_c(d)$ for site percolation on $\L=\Z^d$, $d\geq
2$,
satisfies the inequalities
\begin{equation} \label{eq:perc_nontrivial}
\frac{1}{2d-1} \leq p_c \leq \frac{6}{7}.
\end{equation}
\end{thm}
{\bf Proof:}
We begin with the lower bound on $p_c$.  For $k=1,2,\ldots$, we write
$N_k$ for the random number of open paths of length $k$ starting at $0$.
On the event $\{0\leftrightarrow\infty\}$ we have $N_k\geq 1$ for each $k$,
whence
\begin{equation} \label{eq:percolation_needs_open_paths}
\theta(p)\leq \psi_p(N_k) \, .
\end{equation}
The number of all paths of length $k$ starting at $0$ is at most
$2d(2d-1)^{k-1}$, and each path is open with probability $p^k$. Hence
\[
\psi_p(N_k) \leq 2d(2d-1)^{k-1}p^k
\]
which tends to $0$ as $k\rightarrow\infty$ whenever $p<{1}/({2d-1})$.
In combination with (\ref{eq:percolation_needs_open_paths}) this implies
that $\theta(p)=0$ for $p< {1}/({2d-1})$, and the first half of
(\ref{eq:perc_nontrivial}) is established.

The second half of (\ref{eq:perc_nontrivial}) only needs to be
proved for $d=2$; this is because $\Z^2$ can be embedded into
$\Z^d$ for any $d \geq 2$, so that $p_c(d)\leq p_c(2)$. So let
$d=2$. We first need some additional terminology. A $*$-path in
$\Z^2$ is a sequence $(v_0, v_1, \ldots, v_k)$ of distinct vertices
such that $d_\infty(v_{j-1},v_{j})=1$ for $j=1,\ldots k$. Note that
two consecutive vertices in a $*$-path need not be
nearest-neighbors; they may also be ``diagonal neighbors''. A {\em
$*$-circuit} is a sequence $(v_0, v_1, \ldots,$ $v_k, v_0)$ such
that $(v_0, v_1, \ldots, v_k)$ is a $*$-path and $d_\infty(v_k,
v_0)=1$. Informally, a $*$-circuit is a $*$-path which ends where
it starts.  $*$-circuits  with the same set of sites  are
identified. A {\em closed $*$-circuit} is a $*$-circuit in which
all vertices are closed.

Now let $M$ be the number of closed $*$-circuits that surround the origin
$0$. As $\Z^2$ is a planar graph, the ``outer boundary'' of a finite open
cluster containing $0$ defines a closed $*$-circuit around $0$. Hence, the
event
$\{0 \leftrightarrow\infty\}$ occurs if and only if $M=0$.

The number of {\em all\/} (not necessarily closed) $*$-circuits of a given
length $k$ surrounding $0$   allows the following crude estimate. Consider
the
leftmost crossing of the $x$-axis of such a circuit; the location of such
a crossing is at distance at most $k$ from the origin, so
there are at most $k$ such locations to choose from. Starting at this
location, we may trace the $*$-circuit clockwise (say), and at each step we
have at most $7$ $d_\infty$-neighbors to choose from. Hence, the number
of $*$-circuits of length $k$ around $0$ is at most $k\,7^{k-1}$.
Each one is closed with probability $(1-p)^k$, so
\[
\psi_p(M) \leq \sum_{k=1}^\infty k\,7^{k-1}(1-p)^k.
\]
The last sum is finite for $p>{6}/{7}$. Hence, for such a $p$ and $n$
large enough, there is a positive probability for having no closed
$*$-circuit around $0$ which contains a site of distance at least $n$ from
$0$. By the argument in the proof of Proposition
\ref{prop:perc_single_site},
 it follows  that
 $\psi_p(M=0)>0$ for such $p$.
Hence $\theta(p)>0$ for $p>{6}/{7}$, and the
second half of (\ref{eq:perc_nontrivial}) follows.
$\Cox$

\medskip\noindent
The exact value of the percolation threshold $p_c(d)$ of $\Z^d$ is not
known for any $d\geq 2$. The best rigorous
bounds for $d=2$ are presently
\begin{equation}\label{erma}
0.556 < p_c(2) < 0.680
\end{equation}
where the first inequality is due to van den Berg and Ermakov
\cite{vdBE} (inspired by \cite{MeP}) and the second to Wierman
\cite{Wie}. For high dimensions, it is known that
\begin{equation} \label{eq:perc_high_dimensions}
\lim_{d\rightarrow\infty} {2d} \,{p_c(d)}=1\, ,
\end{equation}
see Kesten \cite{K90}, Gordon \cite{Gor} and Hara and Slade \cite{HS2} for
this and finer asymptotics. On some particular graphs, $p_c$ can be
determined exactly. For example, for $\L=\T_d$, the regular (Cayley or
Bethe) tree,   branching-process arguments immediately show that
$p_c(\T_d)=1/d$, and for the triangular lattice it follows from planar
duality that $p_c=1/2$ \cite{Kes}.

The preceding considerations do not tell us what happens {\em at}
the critical value $p_c$. It is believed that, for the integer lattice
$\Z^d$ of any dimension $d\geq 2$,  there is
$\psi_{p_c}$-a.s.\ no infinite cluster, which means that $\theta(p_c)=0$;
 so far this is only known
for $d=2$ and $d\geq 19$, see Russo \cite{Russo} and Hara and Slade
\cite{HS}. One can show that the relation $\theta(p_c)=0$ implies
continuity of
$\theta(p)$ at $p=p_c$, so in combination with the trivial continuity
of $\theta(p)$ in the subcritical regime and the following more interesting
result (which can be found e.g.\ in \cite{Gr1}) we get continuity of
$\theta(p)$ throughout $[0,1]$ as soon as absence of infinite clusters
at criticality is established.
\begin{thm} \label{thm:perc_continuity}
For Bernoulli site percolation on $\Z^d$, $d\geq 2$, the function
$\theta(p)$ is
continuous throughout the supercritical regime $(p_c,1]$.
\end{thm}

So far we were interested in the existence of infinite clusters. In the
subcritical regime $p<p_c$ when no infinite cluster exists, one may ask
for the size of a typical cluster.
 Let $|C_0|$ be the random number of vertices in the open
cluster containing the origin; we set $|C_0|=0$ if the origin is closed. By
the definition of $p_c$,   $\psi_p(|C_0|=\infty)=0$. For $\L=\Z^d$, we
even have the stronger statement that the expected value of $|C_0|$ is
finite.
\begin{thm}  \label{thm:expected_cluster_size}
For Bernoulli site percolation on $\Z^d$ with retention parameter $p<p_c$,
we have $\psi_p(|C_0|)<\infty$.
\end{thm}
This was proved independently by Menshikov \cite{Men} and by Aizenman and
Barsky \cite{AB}. The proofs  are rather involved, so we refer the reader
to the original articles and \cite{Gr1}.
It is worth noting that Theorem \ref{thm:expected_cluster_size}
fails in the setting of general graphs; a striking counterexample is the
``three-one-tree'' discussed on p.\ 936 of Lyons \cite{Ly}.

Menshikov even showed that the distribution of the radius of the open
cluster containing the origin decays exponentially. The following even
stronger result
states that the same is true for the distribution of $|C_0|$; see Grimmett
\cite{Gr1} for a proof (in the case of bond percolation).
\begin{thm}    \label{thm:exponential_tail}
For Bernoulli site percolation on $\Z^d$ with retention parameter $p<p_c$,
there exists
a constant $c$ (depending on $p$) such that
\[
\psi_p(|C_0| \geq n) \leq e^{-c\,n}
\]
for all $n$.
\end{thm}
Looking at a fixed path of length $n$ starting at $0$, we immediately
obtain the lower bound $\psi_p(|C_0| \geq n) \geq p^n$. So the preceding
upper bound is best possible, except that the optimal
constant $c$ is unknown.

We conclude this section with some remarks on {\em Bernoulli bond
percolation}.
 The model is similar, except that now the edges
rather than the vertices in $(\L,\sim)$ are independently open
(described by the value $1$)   or closed (with value $0$)  with
respective probabilities $p$ and $1-p$. The associated
configuration space is thus $\{0,1\}^\B$. We write $\phi_p$ for the
associated (Bernoulli) probability measure on $\{0,1\}^\B$. In the
present context of bond percolation, an {\em open path} is a path
in which all edges are open, and an {\em open cluster} is a maximal
region $C\subset\L$ which is connected, in that  for any $x,y\in C$
there is an open path in $C$ from $x$ to $y$.

All results for site percolation discussed so far extend to the
bond percolation set-up. This is no surprise because bond
percolation is equivalent to site percolation on the so-called
covering graph for which $\B$  is taken as set of vertices, and
edges are drawn between any two coincident elements of $\B$. In
particular, there exists again a critical value $p_c$ for the
occurrence of infinite open clusters, and Propositions
\ref{prop:perc_critical_value} and \ref{prop:perc_single_site},
Theorems \ref{thm:perc_nontrivial},
\ref{thm:expected_cluster_size}, \ref{thm:exponential_tail} and
\ref{thm:perc_continuity}, and the asymptotic formula
(\ref{eq:perc_high_dimensions}) are still true in the case of bond
percolation. What is generally different, are the critical values
for site and bond Bernoulli percolation on a given graph. One
remarkable case is that of bond percolation on $\Z^2$, where
(again by planar duality) $p_c= {1}/{2}$; this is a famous result
of Kesten \cite{K80}. In the specific case of trees, however, site
and bond percolation are equivalent. In particular, for $\L={\bf
T}_d$ we have $p_c={1}/{d}$ for both site and bond percolation; for
more general trees a formula for $p_c$ was given by Lyons
\cite{Ly}.

\subsection{Dependent percolation: the role of the density}
\label{sect:density}

Our main subject is the analysis of equilibrium phases by means of
percolation methods. In this case, a site will be considered as open if,
for example, the configuration in a neighborhood of this site shows a
specified pattern, and the
events ``site $x$ is open'' with $x\in\L$ are then far from being
independent. This leads us to considering the case of dependent
percolation.
In this subsection we do some first steps in this direction.

Our starting point is the following question.
In the case of Bernoulli percolation,  there is a unique
parameter, the occupation probability or density $p$, which governs the
phase
diagram and allows to distinguish between subcritical (``no infinite
cluster'') and supercritical (``at least one infinite cluster'') behavior.
Does this also
hold in general? Specifically, is it true that for any translation
invariant probability measure $\mu$ on $\{0,1\}^{\Z^d}$, the occurrence of
an infinite cluster only depends
on the density $p(\mu)=\mu(X(x)=1)$ of open sites $x\in\Z^d$? In general,
the answer is obviously ``no'', as we will now show by two simple examples:
  there exist translation invariant measures $\mu$ on $\{0,1\}^{\Z^d}$ with
arbitrarily small densities such that infinite clusters exist almost
surely,
and also translation invariant measures with densities arbitrarily close to
$1$ for which no infinite clusters exist with probability $1$.
\begin{example}\label{ex:percolation_at_low_density}
{\rm For $q\in (0,1)$, let $(Y(x),\, x\in \Z)$ be i.i.d.\ random variables
taking values $0$ and $1$ with probability $1-q$ resp.\ $q$.
We define a translation invariant random field $(X(x), \, x\in \Z^d)$,
$d\geq 2$, by
setting $X(x)=X(x_1,\ldots, x_d)=Y(x_1)$ for each $x\in \Z^d$. Writing
$\mu$ for the distribution of $(X(x), \, x\in \Z^d)$, we have $p(\mu)=q$,
but with $\mu$-probability $1$ there exist infinitely many infinite
clusters, even if $q$ is arbitrarily small. }
\end{example}
\begin{example} \label{ex:nonpercolation_at_high_density}
{\rm Again let $q\in (0,1)$, $d\geq 2$, and
$(Y(x,i),\, x\in \Z,\,i\in\{1,\ldots,d\})$ be i.i.d.\ random variables
taking values $0$ and $1$ with probabilities $1-q$ and $q$.
We define a translation invariant random field $(X(x), \, x\in \Z^d)$ by
setting
\[
X(x)= X(x_1,\ldots,x_d)=\prod_{i=1}^d Y(x_i,i).
\]
A moment's thought reveals that $\mu$-a.s.\ there exist no infinite
clusters, despite the fact that $p(\mu)=q^d$ may be arbitrarily close to
$1$. }
\end{example}
These examples suggest to look for additional assumptions under which
high (resp.\ low) density guarantees existence (resp.\ nonexistence) of
infinite clusters. Positive correlations (in the sense of Definition
\ref{defn:positive_correlations}) does not suffice, because both
 examples above obviously have positive correlations.

An alternative might be to assume $R$-independence in the sense that
$X(\Lambda)$ and $X(\Delta)$ are independent
for any two finite regions $\Lambda, \Delta \subset \L$ such that
\[
\min_{x\in\Lambda, y\in \Delta}|x-y| >R
\]
for some given $R$.
For $\Z^d, d\geq 2$, this gives nontrivial thresholds $0<p_1<p_2<1$
(depending on $R$) such that existence
(resp.\ non-existence) of infinite clusters is guaranteed as long as
$p(\mu)>p_2$ (resp.\ $p(\mu)<p_1$); see e.g.\ Liggett, Schonmann and
Stacey \cite{LSS}.  However, $R$-independence rarely holds
in Gibbs models. For instance, for plus phase of
the Ising model with vanishing external
field and inverse temperature $\beta>0$, the spins at any two vertices are
always strictly positively
correlated no matter how far apart they are (although the correlation does
tend to $0$ in the distance).

However, in contrast to what we just saw in the case of the cubic
lattices, the density does play a significant role for the regular
trees $\T_d$. To show this we consider a  probability measure $\mu$
on $\{0,1\}^\L$, where now
  $\L=\T_d$ with $d\geq 2$.
The natural analogue of translation invariance in this setting is {\em
automorphism invariance} of $\mu$, which means that $\mu$ inherits all the
symmetries of
$\T_d$. In particular, this implies that $\mu(X(x)=1)$ is independent of
$x$,
so that the density $p(\mu)$ is well-defined. As opposed to the $\Z^d$
case,
having $p(\mu)$ sufficiently close to $1$ now does guarantee that
 an infinite cluster exists with positive probability. This is also true in
the bond percolation case, where $p(\mu)$ is defined as the
probability that a given edge is open. The following result is due to
H\"aggstr\"om \cite{H6}.
\begin{thm} \label{thm:dep_perc_on_tree}
For any automorphism invariant site percolation model $\mu$ on $\T_d$ with
density
$p(\mu)\geq \frac{d+1}{2d}$, we have
$\mu(\exists \mbox{ \rm an infinite open cluster})>0$. The same is true for
bond percolation on $\T_d$ with density $p(\mu)\geq \frac{2}{d+1}$.
\end{thm}
These bounds are in fact sharp, in that for any $p<\frac{d+1}{2d}$
there exists some automorphism invariant probability measure on
$\{0,1\}^{\T_d}$ with density $p$, which does not allow an infinite
cluster with probability $1$, and similarly for the case of bond
percolation; see \cite{H6}. It follows from Example
\ref{ex:nonpercolation_at_high_density} that the corresponding
threshold for $\Z^d$ is trivial: only density $1$ is enough to rule
out the nonexistence of infinite clusters.  The intuitive reason is
the following. On $\Z^d$, one can find finite regions
$\La\subset\L$ with arbitrarily small surface-to-volume ratio,
which means that a vast majority of sites is not adjacent to a
vertex outside $\La$; we can simply take
$\La=\Lambda_n=[-n,n]^d\cap \Z^d$ with large $n$; this property of
$\Z^d$ is known as {\em amenability}. Hence, a relatively small
number of closed vertices may easily ``surround'' a large number of
open sites. In contrast, every region in $\T_d$ has a surface of
the same order of magnitude as its volume; this makes it impossible
for a small minority of closed vertices to surround a large number
of open vertices. This intuition can be turned into a proof using
the so called mass-transport method sketched below. Benjamini,
Lyons, Peres and Schramm \cite{BLPS} have recently extended this
method to derive a similar dichotomy for a large class of graphs,
including Cayley graphs of finitely generated groups.

\medskip\noindent
{\bf Sketch proof of Theorem \ref{thm:dep_perc_on_tree}:}
For simplicity we confine ourselves to the case of bond percolation on
$\T_2$.
We want to show that if $p(\mu)\geq 2/3$, then an infinite cluster exists
with positive $\mu$-probability.
Imagine the following allocation of mass to the edges of $\T_2$.
Originally every edge receives mass $1$. Then the mass is redistributed,
 or transported, as follows. If an edge $e$ is
open and is contained in a finite open cluster, then it distributes
all its mass equally among those closed edges that are adjacent to
the open cluster containing $e$. If $e$ is open and contained in an
infinite open cluster, then it keeps its mass. Closed edges,
finally, keep their own mass and happily accept any mass that open
edges decide to send them. The expected mass at each edge before
transport is obviously $1$, and one can show --- this is an
instance of the mass-transport principle \cite{BLPS} --- that the
expected mass at a given edge is $1$ also after the transport.
Suppose now, for contradiction, that $p(\mu)\geq {2}/{3}$ and that all
open clusters are finite $\mu$-a.s. Then all open edges have mass
$0$ after transport. Furthermore, since each open cluster
containing exactly $n$ edges has exactly $n+3$ adjacent closed
edges (as is easily shown by induction --- it is here that the tree
structure
is used), the mass after transport at a closed edge adjacent to two
open clusters of sizes $n_1$ and $n_2$ has mass
\[
1+\frac{n_1}{n_1+3}+\frac{n_2}{n_2+3}<3
\]
Hence the expected mass after transport at a
given edge $e$ is strictly less than
\[
3\,\mu(X(e)=0)=3(1-p(\mu))\leq 1,
\]
contradicting the mass-transport principle.
$\Cox$

\subsection{Examples of dependent percolation}

>From the previous subsection the reader might get a rather
pessimistic view of the possibilities of establishing  existence
(or non-existence) of infinite clusters for dependent percolation
models on $\Z^d$. This is certainly {\em not}\/ the case, and a lot
can be done. One standard way of determining the percolation
behavior of a dependent model is by {\em stochastic comparison with
a suitable Bernoulli percolation model\/}: For the existence of
infinite clusters, it is sufficient to show that the given
dependent model is stochastically larger than the Bernoulli model
for some parameter  $p>p_c$, and the absence of infinite clusters
will follow if the model at hand is stochastically dominated by the
Bernoulli model for some  $p<p_c$. Let us demonstrate this
technique for the Ising model on $\Z^d$.

Consider percolation of plus spins in the
plus measure $\mu^+_{h,\beta}$, defined in Section
\ref{sect:Ising_monotone}.
If we keep $\beta$ fixed then Proposition \ref{prop:Ising_external_field}
tells us that $\mu^+_{h,\beta}$ is stochastically increasing in
$h$. Consequently, both the probability of having
an infinite cluster of plus spins, as well as the probability that a given
vertex is in such an infinite cluster, are increasing in $h$.
Furthermore, as $\Z^d$ is of bounded degree $N=2d$, the remarks after the
same proposition imply that,
for any given $p\in(0,1)$ and $\beta$, the Ising measure $\mu^+_{h,\beta}$
stochastically dominates the Bernoulli measure $\psi_p$ when $h$ is large
enough, and is dominated by $\psi_p$ for $h$ below some bound. We may combine
this observation with Proposition \ref{prop:perc_critical_value} to deduce
the following critical phenomenon:
\begin{thm}\label{thm:Ising_percolation_monotone_in_h}
For the Ising model on $\Z^d$, $d\geq 2$, at a fixed temperature $\beta$,
there exists a critical value $h_c\in \R$ (depending on $d$ and $\beta$)
for
the external field, such that
\[
\mu^+_{h,\beta}(\exists \mbox{ \em an infinite cluster of plus
spins})=\left\{
\begin{array}{ll}
0 & \mbox{\em if } h<h_c \\
1 & \mbox{\em if } h>h_c \, .
\end{array} \right.
\]
\end{thm}
As we shall see later in Theorem \ref{thm:Ising_plus_agreep}, we have
$h_c=0$ when $d\ge2$ and $\beta>\beta_c$. Higuchi \cite{Hi2}
has shown that the percolation transition at $h_c$ is sharp, in that the
connectivity function decays exponentially when $h<h_c$, and that the
percolation probability is continuous in $(\b,h)$, except on the critical
half-line $h=0,\,\b>\b_c$. In Section \ref{sect:random-cluster} below,
we will make a similar use of stochastic comparison arguments for
random-cluster measures, cf.\ Proposition \ref{prop:Potts_magnetization}.
The stochastic domination approach works also in the framework of
lattice gases with attractive potential; see Lebowitz and Penrose
\cite{LePe}.\medskip

In the rest of this subsection we shall give some examples of strongly
dependent systems where other approaches to the question of percolation
are needed. Typically in these examples, the probability that all vertices
in a finite region $\La$ are open (or closed) fails to decay exponentially
in the volume of $\La$, and as a consequence, the random field neither
dominates nor is stochastically dominated by any nontrivial Bernoulli
model.

The geometry of level heights of a random field  forms an important
object of study both from  the theoretical and the applied side.
For example, it relates to the presence of hills and valleys on a
rough surface, or to the random location of potential barriers in a
doped semi-conductor.  To fix the ideas we consider a random field
$X=(X(x), x\in \Z^d)$ with values $X(x) \in S\subset\R$ which are
not necessarily discrete.  It is often interesting to divide $S$
into two parts $S_1$ and $S_0$ and to define a new discrete random
field $Y$ via $Y(x) = I_{\{X(x)\in S_1\}}$. For $S=\R$ one typically
considers $S_1=[\ell,\infty)$ for some level $\ell\in\R$. In
this way we obtain a coarse-grained description of a system of
continuous spins. One question is to which extent one
can reconstruct the complete image from this information.  We
consider here a different question: what is the
geometry of the random set $\{x\in \Z^d: Y(x) =1\}$? This set is
called the excursion or exceedance set when it corresponds, as in
the example above, to the set on which the original
random field exceeds a given level.  For a recent
review of this subject we refer to \cite{Ad}.

We now give four
examples of equilibrium systems with continuous spins
where one can show (the absence of) percolation of an excursion
(exceedance) set. Here we only state the results.
 Some hints on the proofs will be given later in
  Section \ref{sect:agreep} via Theorem \ref{thmrusso}.
Details can be found in the paper by Bricmont, Lebowitz and Maes
\cite{BLM}.
\begin{example}\label{ex:BLM1}
{\rm Consider a general model of real-valued  spins $(\s(x),\,x\in\Z^d)$
with ferromagnetic nearest-neighbor
 interaction. The formal Hamiltonian is given by
\begin{equation}  \label{eq:nn_ferromagnet}
 H(\s) = -\sum_{x\sim y} \s(x)\s(y)\,.
\end{equation}
The reference (or  single-spin) measure $\lambda \neq
\delta_0$ on $\R$ is assumed to be even and to decay fast enough at
$\pm\infty$
 so that the model is well defined.  Then, for any Gibbs measure $\mu$
relative to (\ref{eq:nn_ferromagnet})  with $\mu(\sgn(\s(0)))>0$,
 there will be percolation of all sites $x$ with $\s(x) \geq 0$.
Such Gibbs measures always exist at sufficiently low
temperatures when $d\geq 2$.}
\end{example}
\begin{example} \label{ex:mas}
{\rm Consider  again a  spin system $(\s(x),\,x\in\Z^d)$, where now the
Hamiltonian has the `massless' form
$$H(\s)=-\sum_{x\sim y} \psi(\s(x) - \s(y))$$
with $\s(x)\in \R$ or $\Z$ and $\psi$ an even convex function.
The single-spin measure $\lambda$ is either Lebesgue measure on $\R$ or
counting measure on $\Z$. The case   $\s(x)\in\Z$ and $\psi(t) = |t|$
corresponds to
the so-called solid-on-solid (SOS) model of a $d$-dimensional surface in
$\Z^{d+1}$; the choice $\s(x)\in \R$ and $\psi(t) = t^2$  gives the
harmonic crystal.  Let $\mu$ be a Gibbs measure which is obtained as
infinite volume limit of finite volume Gibbs distributions with zero
boundary condition.
(In the continuous-spin case,
such Gibbs measures exist for any temperature when $d\geq 3$ and $\psi(t)
= \alpha t^2 + \phi(t)$, where $\alpha > 0$ and $\phi$ is convex
 \cite{BL}.) Then, for any $\ell <0$,
there is percolation of the sites $x\in\Z^d$ with $\s(x) \geq
\ell$ .}
\end{example}
\begin{example}
{\rm Consider next a model of two-component spins $\s(x) \in \R^2$,
$x\in\Z^d$, $ \s(x) = (r_x \cos \phi_x,r_x \sin \phi_x)$, with formal
Hamiltonian
$$H(\s) = -\sum_{x\sim y} \s(x)\cdot \s(y)$$
and some rotation-invariant and suitably decaying reference measure
$\lambda$ on $\R^2$. Then, for any Gibbs
measure $\mu$ with  $\mu(\cos \phi_0) > 0$, there is percolation of  the
sites $x\in\Z^d$ with $\cos \phi_x \geq 0$.  Such Gibbs measures exist at
low temperatures if $d\geq 3$.}
\end{example}
\begin{example}\label{ex:harmonic}
{\rm Consider again the massless harmonic crystal of Example \ref{ex:mas}
above
(with $\psi(t) = t^2$) in $d=3$ dimensions.  There exists a value
$\ell_c<\infty$ so that for all $\ell\geq \ell_c$ there is {\em no} percolation of
sites $x\in\Z^d$ with
$\s(x) \geq \ell$.}
\end{example}

Finally, we give an example of a strongly correlated system,
sharing some properties with the harmonic crystal of Example
\ref{ex:harmonic}, where at present there is no proof of a
percolation transition. The model is one of the simplest examples
of an interacting particle system.  What makes the problem
difficult is that the random field is not Markov (not even
Gibbsian) and not explicitly described in terms of a family of
local conditional distributions.

\begin{example}\label{ex:BLM4}
{\rm The voter model is a stochastic dynamics in which individuals
(voters) sitting at the vertices of a graph update their position
(yes/no) by randomly selecting a neighboring vertex and adopting
its position, see Liggett \cite{L} for an introduction. Using
spin language and putting ourselves on $\Z^3$, the time evolution
of this voter model is specified by giving the rate $c(x,\s)$ for a
spin flip at the site $x$ when the spin configuration is $\s\in
\{+1,-1\}^{\Z^3},$
\[
c(x,\s) =  \frac{1}{6}  \sum_{y\sim x} \Big( 1 -\s(x)\s(y)\Big)\;.
\]
There is a one-parameter family of extremal invariant measures
$\mu_p$ each obtained asymptotically (in time) from taking  the
Bernoulli measure $\psi_p$ with density $p$ as initial condition.
These stationary states $\mu_p$ are strongly correlating.  The
spin-spin correlations decay as  the inverse $1/r$ of the
spin-distance $r$ on $\Z^3$. It is an open question whether for $p$
sufficiently close to 1 the plus spins percolate, and whether for
sufficiently small $p$ there is no percolation. Simulations by
Lebowitz and Saleur \cite{LS} indicate that there is indeed a
non-trival percolation transition with critical value $p_c \approx
0.16$.

The same problem may be considered for $d\geq 4$. For $d=1,\,2$, however,
the problem is not interesting because in these cases $\mu_p$ is known
to put mass $p$ on the ``all $+1$'' configuration and mass $1-p$ on
the ``all $-1$'' configuration.
Alternatively,
one may consider the same model with $\Z^3$ replaced by $\T_d$;
Theorem \ref{thm:dep_perc_on_tree} can then be applied to show that
the plus spins do percolate for $p\geq \frac{d+1}{2d}$. }
\end{example}

\subsection{The number of infinite clusters} \label{sect:uniqueness}

Once infinite clusters have been shown to exist with positive
probability in some percolation model, the next natural question
is: {\em How many infinite clusters can exist simultaneously?}  For
Bernoulli site or bond percolation on $\Z^d$, Aizenman, Kesten and
Newman \cite{AKN} obtained the following, now classical uniqueness
result: with probability 1, there exists at most one infinite
cluster. Simpler proofs were  found later by Gandolfi, Grimmett and
Russo \cite{GGR} and by Burton and Keane \cite{BuK}. The argument
of Burton and Keane is not only the shortest (and, arguably, the
most elegant) so far. Also, it requires much weaker assumptions on
the percolation model, namely: translation invariance and the
finite-energy condition below, which is a strong way of stating
that all local configurations are really possible; its significance
for percolation theory had been discovered before by Newman and
Schulman \cite{NS}.
\begin{defn}
A probability measure $\mu$ on $\{0,1\}^\L$, with $\L$ a countable set,
is said to have {\bf finite energy} if, for every finite region
$\La\subset \L$,
\[
\mu(X\equiv \eta \mbox{ \rm on } \La \, | \, X\equiv \xi \mbox{ \rm off }
\La)>0
\]
for all $\eta\in \{0,1\}^\La$ and $\mu$-a.e.\ $\xi\in \{0,1\}^{\La^c}$.
\end{defn}
\begin{thm}[The Burton--Keane uniqueness theorem] \label{thm:Burton_Keane}
Let $\mu$ be a probability measure on $\{0,1\}^{\Z^d}$ which is translation
invariant and has finite energy. Then, $\mu$-a.s., there exists at most one
infinite open cluster.
\end{thm}
{\bf Sketch proof:} Without loss of generality we can assume that
$\mu$ is ergodic with respect to translations. For, one can easily show
that the measures in the ergodic decomposition of $\mu$ admit the same
conditional probabilities, and thus inherit the finite-energy property.
Since the number $N$ of infinite clusters is obviously invariant under
translations, it then follows that $N$ is almost surely equal to some
constant $k\in \{0,1,\ldots,\infty\}$. In fact, $k\in \{0,1,\infty\}$.
Otherwise, with positive probability each of the $k$ clusters would
meet a sufficiently large cube $\La$; by the finite-energy property,
this would imply
that with positive probability all these clusters
are connected within $\La$, so that in fact $k=1$, in contradiction to the
hypothesis. (This part of the argument goes back to \cite{NS}.)

We thus only need to exclude the case $k=\infty$. In this case,
$\mu(N\ge3)=1$,
and the finite-energy property implies again that $\mu(A_x)=\delta>0$,
where $A_x$ is the event that $x$ is a triple point, in that there exist
three disjoint infinite open paths with starting point $x$. By the (norm-)
ergodic theorem, for any sufficiently large cubic box
$\La$ we have
\begin{equation}\label{eq:Burton_Keane}
\mu\Bigl(|\La|^{-1}\sum_{x\in\La}I_{A_x}\ge\delta/2\Bigr)\ge 1/2\,.
\end{equation}
On the other hand, for geometrical reasons (which are intuitively obvious
but need some work when made precise), there cannot be more triple points
in $\La$
than  points in the boundary $\partial\La$ of $\La$. Indeed,
each of the three paths leaving a triple point meets $\partial\La$, which
gives three boundary points associated to each triple point in $\La$. If
one identifies these boundary points successively for one triple point
after the other one sees that, at each step, at least one of the boundary
points must be different from those obtained before. Hence,
\[
|\La|^{-1}\sum_{x\in\La}I_{A_x}\le |\La|^{-1}|\partial\La|<\delta/2
\]
when $\La $ is large enough. Inserting this into \rf{eq:Burton_Keane} we
arrive at the contradiction $\mu(\emptyset)\ge1/2$, and the theorem is
proved. For more details we refer
to the original paper \cite{BuK}.
$\Cox$

\medskip\noindent
We stress that the last argument relies essentially on the amenability
property of $\Z^d$
discussed in Section \ref{sect:density}. The finite-energy condition is also
indispensable: In another paper \cite{BuK2}, Burton and Keane construct,
for any $k\in \{2,3,\ldots,\infty\}$,
translation invariant percolation models on $\Z^2$ for
which finite energy fails and which have exactly $k$
infinite open clusters.
 For example, we have $k=\infty$ in Example 6.1.
Fortunately, the finite-energy condition  holds in most of the dependent
percolation models which show up in stochastic-geometric studies of
Gibbs measures.

The situation becomes radically
different when $\Z^d$ is replaced by  the non-amenable tree $\T_d$.
Instead of having a unique infinite cluster, supercritical
percolation models on $\T_d$ tend to produce infinitely many infinite
clusters. It is not hard to verify that this is indeed the case for
supercritical Bernoulli site
or bond percolation
(except in the trivial case when the retention probability
$p$ is $1$), and a corresponding result for automorphism invariant
percolation
on $\T_d$ can be found in \cite{H6}. On more general nonamenable graph
structures, the uniqueness of the infinite cluster property can fail in
more interesting ways than on trees;
see e.g.\ Grimmett and Newman \cite{GN} and
H\"aggstr\"om and Peres \cite{HP}.

Let us next consider the particular case of (possibly dependent)
site percolation on $\Z^2$. We know from Theorem \ref{thm:Burton_Keane}
that
under fairly general assumptions
there is almost surely at most one infinite open cluster. Under the same
asumptions
there is almost surely at most one infinite {\em closed} cluster (i.e., at
most
one infinite connected component of closed vertices). In fact, the  proof
of Theorem \ref{thm:Burton_Keane} even shows that almost surely there is
at most one infinite closed $*$-cluster. (Here, a closed $*$-cluster is a
maximal set $C$ of closed sites which is $*$-connected, in that any two
$x,y\in C$ are connected by a $*$-path in $C$; $*$-paths were introduced
in the proof of Theorem \ref{thm:perc_nontrivial}. Any closed cluster is
part of some closed $*$-cluster.)
 But perhaps an infinite open
cluster and an infinite closed $*$-cluster can coexist? Theorem
\ref{thm:GKR} below asserts
that under reasonably general circumstances this cannot happen.
Under slightly different conditions (replacing the finite-energy
assumption by
separate ergodicity under translations in the two coordinate directions),
it was proved by Gandolfi, Keane and Russo \cite{GKR}.
\begin{thm} \label{thm:GKR}
Let $\mu$ be an automorphism
invariant and ergodic probability measure on $\{0,1\}^{\Z^2}$ with finite
energy and positive correlations. Then
\[
\mu(\exists \mbox{ \rm infinite open cluster, }\exists \mbox{ \rm infinite
closed  $*$-cluster})= 0\,.
\]
\end{thm}
Note that automorphism invariance in the $\Z^2$-case means that, in
addition to
translation invariance, $\mu$ is also invariant under reflection in and
exchange of coordinate axes. Under the conditions of the theorem, we have
in fact some information on the geometric shape of infinite clusters: If an
infinite open cluster exists and thus
all closed $*$-clusters are finite, each finite box of $\Z^2$ is surrounded
by an open circuit, and all these circuits are part of the
(necessarily unique) infinite open cluster. Hence the infinite open cluster
is a sea, in the sense that all ``islands'' (i.e., the $*$-clusters of its
complement) are finite. Similarly, if a closed $*$-cluster in $\Z^2$ exists,
it is necessarily a sea (and in particular unique).

The corresponding result is
false in higher dimensions. To see this, consider Bernoulli site percolation
on $\Z^3$. The critical value $p_c$ for this model is strictly less than
${1}/{2}$ (see Campanino and Russo \cite{CR}), whence for $p={1}/{2}$
there exist almost surely both an infinite open cluster and an infinite
closed cluster.

  The proof of
Theorem \ref{thm:GKR} below is based on a geometric argument of Yu
Zhang who gave a new proof of Harris' \cite{Har} classical result
that the critical value $p_c$ for bond percolation on $\Z^2$ is at
least ${1}/{2}$. (Recall that this bound is actually sharp.)
Zhang's proof appeared first in \cite{Gr1} and was exploited later
in other contexts in \cite{H3,HJ}.

\medskip\noindent
{\bf Proof of Theorem \ref{thm:GKR}:}
Let $A$ be the event that there exists an infinite open cluster, let
$B$ be the event that there exists an infinite closed $*$-cluster,
and assume by contradiction that $\mu(A\cap B)>0$.
Then, by ergodicity,  $\mu(A\cap B)=1$. (This is the only use of
ergodicity we make, and ergodicity could clearly be replaced by tail
triviality or some other mixing condition.)
Next we pick $n$ so large that
\[
\mu(A_n)>1-10^{-3}\ \mbox{ and }\ \mu(B_n)>1-10^{-3}\;,
\]
where $A_n$ (resp.\ $B_n$) is the event that some
infinite open cluster
(resp.\ some infinite closed $*$-cluster)
 intersects $\Lambda_n =[-n,n]^2\cap \Z^2$.
Let $A_n^L$ (resp.\ $A_n^R$, $A_n^T$ and $A_n^B$)
be the event that some vertex in the left (resp.\ right,
top and bottom) side of the square-shaped vertex set
$\Lambda_n\setminus\Lambda_{n-1}$ belongs to some infinite open path
which contains no other vertex of $\Lambda_n$, and
define $B_n^L$, $B_n^R$, $B_n^T$ and $B_n^B$ analogously.
Then
\[
A_n=A_n^L\cup A_n^R \cup A_n^T \cup A_n^B\,.
\]
Since all four events in the right hand side are increasing and
$\mu$ has positive correlations,
\begin{eqnarray*}
\mu(A_n) & = & \mu(A_n^L\cup A_n^R \cup A_n^T \cup A_n^B) \\
& = & 1-\mu(\neg A_n^L\cap \neg A_n^R \cap \neg A_n^T \cap \neg A_n^B) \\
& \leq & 1-\mu(\neg A_n^L)\mu(\neg A_n^R)\mu(\neg A_n^T)
\mu(\neg A_n^B)\;,
\end{eqnarray*}
where $\neg$ indicates the complement of a set (for typographical
reasons).
By the automorphism invariance of $\mu$,
$A_n^L$, $A_n^R$, $A_n^T$ and $A_n^B$ all have the
same $\mu$-probability, so that
\[
\mu(\neg A_n^L)\leq(1-\mu(A_n))^{{1}/{4}}
\]
and therefore
\begin{equation} \label{eq:A-events}
\mu(A_n^L)=\mu(A_n^R)\geq 1-(1-\mu(A_n))^{{1}/{4}}=
1-10^{-{3}/{4}}> 0.82\;.
\end{equation}
In the same way, we get
\begin{equation} \label{eq:B-events}
\mu(B_n^T)=\mu(B_n^B)> 0.82\;.
\end{equation}
Now define the event $D=A_n^L\cap A_n^R\cap B_n^T\cap B_n^B$. From
(\ref{eq:A-events}) and (\ref{eq:B-events}) we obtain
\[
\mu(D)\geq 1-4(1-0.82)=0.28>0\;.
\]
When $D$ occurs, both the left side and
the right side of $\Lambda_n$ are intersected by some infinite open
cluster.
By Theorem \ref{thm:Burton_Keane}, these
infinite open clusters are identical and separate their (common)
complement into (at least) two pieces, preventing
the infinite closed $*$-clusters intersecting the top and bottom
sides of $\Lambda_n$ from reaching each other (see the picture on p.\
196 of \cite{Gr1}). Consequently, there exist two infinite closed
$*$-clusters, in contradiction to Theorem \ref{thm:Burton_Keane}.
$\Cox$

\medskip\noindent
Theorem \ref{thm:GKR} admits some variants. First,
the assumption of ergodicity can be avoided if the assumption of
positive correlations is strenghtened to the condition that $\mu$
is monotone in the sense of Definition \ref{def:mon}. (This
is because monotonicity is preserved under ergodic
decomposition, so that Theorem \ref{thm:FKG} implies
positive correlations for each ergodic component.) Moreover, as
the preceding proof shows, ergodicity is only needed to show that
infinite clusters exist with probability either 0 or 1, and
translation invariance and finite energy are only used for the
uniqueness of infinite clusters. We also need only the invariance
under lattice rotations rather than all reflections, and closed
$*$-clusters can be replaced by closed clusters. We may thus state
the following result.
\begin{prop}\label{prop:GKR2}
There exists no probability measure $\mu$ on $\{0,1\}^{\Z^2}$ which
has positive correlations, is invariant under lattice rotations and
the interchange of the states 1 (``open'') and 0 (``closed''), and
satisfies
\[
\mu(\,\exists \mbox{ \em a unique infinite open cluster})=1\;.
\]
\end{prop}
This proposition will be applied to the ferromagnetic Ising model in
Section \ref{subs:ferro}.

\section{Random-cluster representations} \label{sect:random-cluster}

In the previous section we saw a number of dependent percolation models.
Here we shall focus on a particular class of such models, namely the
(Fortuin--Kasteleyn)
random-cluster model (and some of its relatives), which has turned out
to be of great value in analyzing the phase transition behavior of
Ising and Potts models. An alternative source for much of the material
in the present section is H\"aggstr\"om \cite{H7}. In Sections
\ref{sect:FK} and \ref{sect:Inf_vol_limits},
we introduce the random-cluster model and discuss its relation
to Ising and Potts models. This relation is then applied in Section
\ref{sect:PT_Potts} to prove Theorems \ref{thm:Ising} and \ref{thm:Potts}.
Despite the fact that Theorems \ref{thm:Ising} and \ref{thm:Potts} concern
infinite systems, these applications only require defining
finite volume random-cluster measures. However, it may be interesting
in its own right to study infinite volume random-cluster measures on graphs
such as $\Z^d$; this is done in Section \ref{sect:infinite_FK}.
In Section \ref{sect:Swendsen_Wang}, we describe how the random-cluster
representation of Ising and Potts models can be used to construct
highly efficient Monte Carlo simulation algorithms.
Finally,
in Section \ref{sect:WR_RC}, we discuss a variant of the random-cluster
model
which is applicable to the Widom--Rowlinson model rather than to
Ising and Potts models.

\subsection{Random-cluster and Potts models}
\label{sect:FK}

The random-cluster model, also known as the Fortuin--Kasteleyn
(FK) model after
its inventors \cite{FK1,FK2,FK3}, is a two-parameter family
of dependent bond percolation models living on a finite graph. Let
$G=(\L,\sim)$ be
a finite graph with vertex set $\L$ and edge set $\B$. For a bond
configuration $\eta\in\{0,1\}^\B$, we write $k(\eta)$ for the number of
connected components (including isolated vertices) in the subgraph of $G$
containing all vertices but only the open edges (i.e.\ those $e\in \B$
for which $\eta(e)=1$).
\begin{defn} \label{defn:random-cluster}
The {\bf random-cluster measure} $\phi^G_{p,q}$ for $G$ with parameters
$p\in[0,1]$ and $q>0$ is the probability measure on $\{0,1\}^\B$ which
to each $\eta\in\{0,1\}^\B$ assigns probability
\[
\phi^G_{p,q}(\eta)=\frac{1}{Z^G_{p,q}}\left\{\prod_{e\in \B}
p^{\eta(e)}(1-p)^{1-\eta(e)}\right\}q^{k(\eta)},
\]
where $Z^G_{p,q}$ is a normalizing constant.
\end{defn}
Note that taking $q=1$ yields the Bernoulli bond percolation measure
$\phi_p$ defined at the end of Section \ref{sect:bond_percolation}. All
other
choices of $q$ give rise to dependencies between edges (as long
as $p$ is not $0$ or $1$, and $G$ is not a tree).

Taking $q\in\{2,3,\ldots\}$ yields a model which is intimately related to
the $q$-state Potts model, in a way which we will explain now.
Let $\mu^G_{\beta,q}$ be the Gibbs measure
for the $q$-state Potts model on $G$ at inverse temperature $\beta$,
i.e.\ $\mu^G_{\beta,q}$ is the measure on $\{1,\ldots, q\}^\L$ which
to each $\sigma \in \{1,\ldots, q\}^\L$ assigns probability
\[
\mu^G_{q,\beta}(\sigma) = \frac{1}{Z^G_{\beta,q}} \exp\left(
-2\beta\sum_{x\sim y}I_{\{\sigma(x)\neq \sigma(y)\}}\right)\;,
\]
where again $Z^G_{\beta,q}$ is a normalizing constant.

For $q\in\{2,3,\ldots\}$ and $p\in[0,1]$,
let $P^G_{p,q}$ be the probability measure on
$\{1,\ldots,q\}^\L \times \{0,1\}^\B$ corresponding to picking a random
element of $\{1,\ldots,q\}^\L \times \{0,1\}^\B$ according
to the following two-step procedure.
\begin{enumerate}
\setlength{\itemsep}{0ex}
\item Assign each vertex a spin value chosen from $\{1,\ldots,q\}$
according to uniform distribution, assign each edge value $1$ or $0$ with
respective probabilities $p$ and $1-p$, and do this independently
for all vertices and edges.
\item Condition on the event that no two vertices with different spins
have an edge with value $1$ connecting them.
\end{enumerate}
In other words, $P^G_{p,q}$ is the measure which to each element
$(\sigma,\eta)$ of $\{1,\ldots,q\}^\L \times \{0,1\}^\B$ assigns
probability proportional to
\[
\prod_{e=\langle x y \rangle\in \B}\left( p^{\eta(e)}(1-p)^{1-\eta(e)}\;
 I_{\{(\sigma(x)-\sigma(y))
\eta(e)=0\}} \right)\;.
\]
Here, $\langle x y \rangle$ denotes the edge linking $x$ and $y$.
The measure $P^G_{p,q}$ was introduced by Swendsen and Wang
\cite{SW} and made more explicit by Edwards and Sokal \cite{ES},
and is therefore called the {\em Edwards--Sokal measure}. The
following theorem states that the edge marginal of $P^G_{p,q}$ is a
random-cluster measure, and the vertex marginal is a Gibbs measure
for the Potts model, meaning that $P^G_{p,q}$ is a coupling of
$\mu^G_{\beta,q}$ and $\phi^G_{p,q}$.
\begin{thm} \label{thm:RC_coupling}
Let $P^{G,vertex}_{p,q}$ and $P^{G,edge}_{p,q}$ be the probability measures
obtained by projecting $P^G_{p,q}$ on $\{1,\ldots,q\}^\L$ and
$\{0,1\}^\B$, respectively. Then
\begin{equation} \label{eq:spinmarginal}
P^{G,vertex}_{p,q}=\mu^G_{\beta,q} \,
\end{equation}
with $\beta=\frac{1}{2}\log(1-p)$, and
\begin{equation} \label{eq:edgemarginal}
P^{G,edge}_{p,q} =\phi^G_{p,q} \, .
\end{equation}
\end{thm}
{\bf Proof:} The proof is just a matter of summing out the marginals.
Letting $Z$ be the normalizing constant in $P^G_{p,q}$,
fixing $\sigma\in\{1,\ldots,q\}^\L$, and summing over all
$\eta\in\{0,1\}^\B$ we find
\begin{eqnarray*}
P^{G,vertex}_{p,q}(\sigma)
& = & \sum_{\eta\in\{0,1\}^\B}P^G_{p,q}(\sigma,\eta) \\
& = & \frac{1}{Z}\sum_{\eta\in\{0,1\}^\B}\,
\prod_{e=\langle x y \rangle\in \B} p^{\eta(e)}
(1-p)^{1-\eta(e )}
I_{\{(\sigma(x)-\sigma(y))\eta(e)=0\}} \\
& = & \frac{1}{Z}\prod_{e=\langle x y \rangle\in \B}
(1-p)^{I_{\{\sigma(x)\neq \sigma(y)\}}} \\
& = & \frac{1}{Z}\exp\bigg(-2\beta\sum_{x\sim y}
I_{\{\sigma(x)\neq \sigma(y)\}}\bigg) \\
& = & \mu_{\beta,q}^G(\sigma)\;,
\end{eqnarray*}
since $Z$ must be equal to $Z^{\beta,q}_G$ by normalization. This proves
(\ref{eq:spinmarginal}). To verify (\ref{eq:edgemarginal}) we
proceed similarly, fixing $\eta\in\{0,1\}^\B$ and summing over
$\sigma\in\{1,\ldots,q\}^\L$. Note that,
given $\eta$, there are exactly $q^{k(\eta)}$
spin configurations $\sigma$ that are allowed, in that any two
neighboring vertices $x\sim y$ with $\eta(\langle x y\rangle )=1$
have the same spin. We get
\begin{eqnarray*}
P^{G,edge}_{p,q}(\eta)
& = & \sum_{\sigma\in\{1,\ldots,q\}^\L}P^G_{p,q}(\sigma,\eta) \\
& = & q^{k(\eta)}\frac{1}{Z}
\prod_{e\in \B}p^{\eta(e)}(1-p)^{1-\eta(e)} \\
& = & \mu^G_{p,q}(\eta)\;,
\end{eqnarray*}
again by normalization. $\Cox$

\medskip\noindent
The Edwards--Sokal coupling $P^G_{p,q}$ of $\mu^G_{\beta,q}$ and
$\phi^G_{p,q}$ is the
key to using the random-cluster model in analyzing the Potts model. The
following two results are each other's dual, and are
immediate consequences of Theorem
\ref{thm:RC_coupling} and the definition of $P^G_{p,q}$.
\begin{cor} \label{cor:RC_to_Potts}
Let $p=1-e^{-2\beta}$, and suppose we pick a
random spin configuration $X\in\{1,\ldots,q\}^\L$ as follows:
\begin{enumerate}
\item Pick a random edge configuration $Y\in\{0,1\}^\B$ according to
the random-cluster measure $\phi^{G}_{p,q}$.
\item For each connected component $C$ of $Y$, pick a spin at
random (uniformly) from $\{1,\ldots,q\}$, assign this spin to every vertex
of $C$, and do this independently for different connected components.
\end{enumerate}
Then $X$ is distributed according to the Gibbs measure $\mu_{\beta,q}^G$.
\end{cor}
\begin{cor} \label{cor:Potts_to_RC}
Let $p=1-e^{-2\beta}$, and suppose we pick a
random edge configuration $Y\in\{0,1\}^\B$ as follows:
\begin{enumerate}
\item Pick a random spin configuration $X\in\{1,\ldots,q\}^\L$ according to
the Gibbs measure $\mu_{\beta,q}^G$.
\item Given $X$, assign each edge $e=\langle x y \rangle$
independently value $1$ with probability
\[
\left\{
\begin{array}{ll}
p & \mbox{if $X(x)=X(y)$} \\
0 & \mbox{if $X(x)\neq X(y)$\, ,}
\end{array} \right.
\]
and value 0 otherwise.
\end{enumerate}
Then $Y$ is distributed according to the random-cluster measure
$\phi^{G}_{p,q}$.
\end{cor}
As a warm-up for the phase transition considerations in Section
\ref{sect:PT_Potts}, we give the following result as
a typical application of the random-cluster representation.
\begin{cor}  \label{cor:Potts_pos_cor}
If we pick a random spin configuration
$X\in\{1,\ldots,q\}^\L$ according to the Gibbs measure $\mu_{\beta,q}^G$,
then
for $i\in\{1,\ldots,q\}$ and two vertices $x, y\in \L$, the two events
$\{X(x)=i\}$ and $\{X(y)=i\}$ are positively correlated, i.e.\
\[
\mu_{\beta,q}^G(X(x)=i,X(y)=i)\geq
\mu_{\beta,q}^G(X(x)=i)\,\mu_{\beta,q}^G(X(y)=i)\,.
\]
\end{cor}
{\bf Proof:}
The measure $\mu_{\beta,q}^G$ is invariant under permutation of the
spin set $\{1,\ldots,q\}$, so that
\[
\mu_{\beta,q}^G(X(x)=i)=\mu_{\beta,q}^G(X(y)=i)=\frac{1}{q}\, .
\]
We therefore need to show that
\[
\mu_{\beta,q}^G(X(x)=i,X(y)=i)\geq\frac{1}{q^2} \, .
\]
We may now think of $X$ as being obtained as in Corollary
\ref{cor:RC_to_Potts}
by first picking an edge configuration $Y\in\{0,1\}^\B$ according to
the random-cluster measure $\phi^{G}_{p,q}$ and then assigning i.i.d.\
uniform spins to the connected components. Given $Y$, the conditional
probability that $X(x)=X(y)=i$ is ${1}/{q}$ if $u$ and $v$ are in
the same connected component of $Y$, and ${1}/{q^2}$ if they are
in different connected components. Hence, for some $\alpha\in[0,1]$,
\[
\mu_{\beta,q}^G(X(x)=i,X(y)=i)=
\alpha\frac{1}{q}+(1-\alpha)\frac{1}{q^2}\geq \frac{1}{q^2} \, .
\]
$\Cox$

\medskip\noindent
An easy modification of the above proof shows that
if $G$ is connected and $\beta>0$, then
the correlation between $I_{\{X(x)=i\}}$ and $I_{\{X(y)=i\}}$
is in fact {\em strictly} positive.

Note that the relation between the random-cluster model and the Potts model
depends crucially on the fact that all spins in $\{1, \ldots, q\}$ are
{\em a priori} equivalent. This is no longer the case when a nonzero
external field is present in the Ising model. Several attempts to find useful
random-cluster representations of the Ising model with external field
have been made, but progress has been limited. Perhaps the recent duplication
idea of Chayes, Machta and Redner \cite{CMR} represents a breakthrough on
this problem.

\subsection{Infinite-volume limits}
\label{sect:Inf_vol_limits}

In this subsection we will exploit some stochastic monotonicity properties
of random-cluster distributions on finite subgraphs of $\Z^d$. This will
give us the existence of certain limiting random-cluster distributions,
and also the existence of certain limiting Gibbs measures for the Potts
model.

The basic observation is stated in the lemma below which follows directly
from definitions.
\begin{lem} \label{lem:single_edge}
Consider the random-cluster model with parameters $p$ and $q$
on a finite graph $G$ with edge set $\B$.
For any edge $e=\langle x y \rangle \in \B$, and any configuration
$\eta\in\{0,1\}^{\B\setminus\{e\}}$, we have that
\begin{equation} \label{eq:single_edge}
\phi^G_{p,q}(e \mbox{ \em is open}\, | \, \eta)= \left\{
\begin{array}{ll}
p & \mbox{\em if $x$ and $y$ are connected via open edges in $\eta$} \\
\frac{p}{p+(1-p)q} & \mbox{\em otherwise.}
\end{array} \right.
\end{equation}
\end{lem}
\medskip\noindent
For $q\geq 1$, Lemma \ref{lem:single_edge} means in particular that
the conditional probability in (\ref{eq:single_edge})
is increasing in $\eta$ (and also in $p$). This allows us to use Holley's
Theorem and the FKG inequality to prove the following very useful result.
We write $\phi^G_p$ for Bernoulli bond percolation on $G$ with parameter $p$.
\begin{cor}  \label{cor:FKG_for_RC}
For a finite graph $G$ and the random-cluster measure $\phi^G_{p,q}$ with
$p\in[0,1]$ and $q\geq 1$, we have
\begin{description}
\item{\rm (a) }
$\phi^G_{p,q}$ is monotone, and therefore it has positive correlations,
\item{\rm (b) }
$\phi^G_{p,q}\leqd \phi^G_p \,$,
\item{\rm (c) }
$\phi^G_{p,q} \geqd \phi^G_{\frac{p}{p+(1-p)q}}\,$ .
\end{description}
Furthermore, for $0\leq p_1\leq p_2\leq 1$ and $q\geq 1$, we have
\begin{description}
\item{\rm (d) }
$\phi^G_{p_1,q} \leqd \phi^G_{p_2, q}\,$.
\end{description}
\end{cor}
{\bf Proof:} The monotonicity in (a) is just the observation that the
conditional probability in (\ref{eq:single_edge}) is increasing in $p$
and in $\eta$. Positive correlations then follows from Theorem
\ref{thm:FKG}. Next, note that (\ref{eq:single_edge}) implies that
\begin{equation}  \label{eq:single_edge_inequality}
\frac{p}{p+(1-p)q} \leq \phi^G_{p,q}(e\mbox{ is open}\, | \, \eta) \leq p
\end{equation}
for all $\eta$ as in Lemma \ref{lem:single_edge}. Theorem \ref{thm:holley}
in conjunction with the second (resp.\ first) inequality in
(\ref{eq:single_edge_inequality}) implies (b) (resp.\ (a)). Finally,
(d) is just another application of (d) and Theorem \ref{thm:holley}.
$\Cox$

\medskip\noindent
Consider now the integer lattice $\Z^d$ (for definiteness and simplicity) with
its usual graph structure. We associate with any finite region
$\La\subset\Z^d$ two specific random-cluster distributions which
correspond to two different choices of boundary condition. The latter will
we distinguished by a parameter $b\in\{0,1\}$.
Let $\B$ be the set of all nearest-neighbor bonds in $\Z^d$, $\B_\La^0$
the set of all edges of $\B$ that are contained in $\La$, and $\B_\La^1$
the set of edges with at least one endpoint in $\La$. (The difference
$\B_\La^1\setminus\B_\La^0$ thus consists of all edges leading from a
point of
$\La$ to a point of $\La^c$.) We then let $\phi^b_{p,q,\La}$ be the
probability
measure on $\{0,1\}^\B$ in which each $\eta\in \{0,1\}^\B$ is assigned
probability proportional to
\[
I_{\{\eta\equiv b \mbox{ \scriptsize off } \B_\La^b\}}
\bigg\{\prod_{e\in \B_\La^b}
p^{\eta(e)}(1-p)^{1-\eta(e)}\bigg\}\; q^{k(\eta,\La)}\;,
\]
where $k(\eta,\La)$ is the number of all $\eta$-open clusters meeting
$\La$.
We call $\phi^b_{p,q,\La}$ the random-cluster distribution in
$\Lambda$ with parameters $p$ and $q$ and boundary condition $b$. In the
case
$b=0$,  $k(\cdot,\La)$ is simply the number of all clusters that are
contained in $\La$; this corresponds to forgetting all sites in $\La^c$ and
is therefore referred to as the {\em free boundary condition}. On the
other hand, suppose that $\Lambda$ has no holes, in the sense that
$\Lambda^c$ has no finite connected components; since we can always
assume without loss of generality that $\La$ is connected, we call such a $\La$
{\em simply connected}.
Then, in the case $b=1$, all sites of $\La^c$ may be thought of as
being firmly wired together, whence this is called the {\em wired boundary
condition}.

Suppose now that $\La\subset\De$ are two finite regions in $\Z^d$. Then
$\phi^b_{p,q,\La}$ is obtained from $\phi^b_{p,q,\De}$ by conditioning on
the event $\{\eta\equiv b \mbox{ on } \B^b_\De\setminus\B^b_\La\}$ which
is increasing for $b=1$ and decreasing for $b=0$. Hence, if $q\ge1$ then
Corollary \ref{cor:FKG_for_RC} (a) implies that
\begin{equation}\label{eq:RC_monotone}
\phi^0_{p,q,\La} \leqd \phi^0_{p,q,\De} \mbox{ and } \phi^1_{p,q,\La}
\geqd\phi^1_{p,q,\De} \mbox{ when } \La\subset\De\;,
\end{equation}
in complete analogy to (\ref{eq:decreasing}). Moreover, we obtain the
following counterpart of Proposition \ref{prop:plus_measure}
on the existence of infinite-volume limits. We write
$\La\uparrow\Z^d$ for the limit along some (any) increasing
sequence of finite simply
connected subsets of ${\bf Z}^d$, converging to ${\bf Z}^d$ in the usual way.
\begin{lem} \label{lem:free_and_wired_measure}
For $p\in[0,1]$ and $q\ge1$, the limiting measures
\[
\phi^b_{p,q} =\lim_{\La\uparrow\Z^d} \phi^b_{p,q,\La}\;,\quad
b\in\{0,1\}\,,
\]
exist and are translation invariant.
\end{lem}

\medskip\noindent
This convergence result has consequences for the convergence of
Gibbs distributions for the Potts model, as we will show next. Let
$q\in\{2,3,\ldots\}$, and for $i\in\{1,\ldots,q\}$ and any finite
region $\Lambda$ in $\Z^d$ let $\mu_{\beta,q,\Lambda}^i$ denote the
Gibbs distribution in $\Lambda$ for the Potts model at inverse
temperature $\beta$ with boundary condition $\eta\equiv i$ on
$\Lambda^c$. For $i=0$, let $\mu_{\beta,q,\Lambda}^0$ be the
corresponding Gibbs distribution with free boundary condition,
which is defined by letting $\L=\La$ in
\rf{eq:FvolGibbs}, i.e., by ignoring all sites outside $\La$; we think of
$\mu_{\beta,q,\Lambda}^0$
as a probability measure on the full configuration space
$\{1,\ldots,q\}^{\Z^d}$ by using an arbitrary extension.

Still for $i=0$, Theorem \ref{thm:RC_coupling} shows that
$\mu_{\beta,q,\Lambda}^0$ and $\phi^0_{p,q,\La}$ admit an Edwards--Sokal
coupling (on $\La$) when $p=1-e^{-2\b}$. A similar Edwards--Sokal
coupling is possible for $i\in\{1,\ldots, q\}$ when $\La$ is simply
connected. Indeed, let $P^{i}_{p,q,\La}$ be the
probability measure on $\{1,\ldots,q\}^\L \times \{0,1\}^\B$
corresponding to picking a random site-and-bond configuration according
to the following procedure.
\begin{enumerate}
\item
Assign to each vertex of $\Lambda^c$ value $i$, and to all edges of
$\B\setminus\B_\La^1$ value $1$.
\item
Assign to each vertex in $\Lambda$ a spin value chosen from
$\{1,\ldots,q\}$
according to uniform distribution, assign to each edge
in $\B_\La^1$ value $1$ or $0$ with
respective probabilities $p$ and $1-p$, and do this independently
for all vertices and edges.
\item
Condition on the event that no two vertices with different spins
have an edge with value $1$ connecting them.
\end{enumerate}
It is now a simple modification of the proof of Theorem
\ref{thm:RC_coupling}
to check that the vertex and edge marginals of $P^{i}_{p,q,\La}$ are
$\mu^i_{\beta,q,\La}$ and $\phi^1_{p,q,\La}$, respectively.
(Note that by the simple connectedness of $\La$ there is always a
unique component containing $\La^c$.)
Analogues of Corollaries \ref{cor:RC_to_Potts} and \ref{cor:Potts_to_RC}
follow easily. This leads us to the following result extending
Proposition \ref{prop:plus_measure} to the Potts model.
\begin{prop} \label{prop:free_and_ordered_Potts_state}
For any $i\in\{0,1,\ldots,q\}$, the limiting probability measure
\[
\mu_{\beta,q}^i =\lim_{\La\uparrow\Z^d} \mu_{\beta,q,\Lambda}^i
\]
on
$\{1,\ldots,q\}^{\Z^d}$ exists and is a translation invariant Gibbs
measure for the
$q$-state Potts model on $\Z^d$ at inverse temperature $\beta$.
\end{prop}
{\bf Proof:} In view of the general facts reported in Section 2.6, the
limits are Gibbs measures whenever they exist.
We thus need to show that $\mu_{\beta,q,\Lambda}^i(f)$ converges as
$\La\uparrow\Z^d$, for any
local observable $f$. For definiteness, we do this for
$i\in\{1,\ldots,q\}$; the case $i=0$ is completely similar.

Fix an $f$ as above, and let $\De\subset \Z^d$ be the finite
region on which $f$ depends. As shown above, for a simply connected
$\La$ we may think of a $\{1,\ldots,q\}^{\Z^d}$-valued random
element $X$ with distribution $\mu_{\beta,q,\Lambda}^i$ as arising
by first picking an edge configuration $Y\in \{0,1\}^\La$ according to
$\phi^1_{p,q,\La}$ (with $p=1-e^{-2\beta}$) and then assigning random spins
to the connected components, forcing spin $i$ to the (unique) infinite
cluster. For $x, y\in \De$, we write $\{x\leftrightarrow y\}$ for the
event that $x$ and $y$ are in the same connected component in $Y$, and
$\{x\leftrightarrow \infty\}$ for the event that $x$ is in an infinite
cluster. Clearly, the conditional distribution of $f$ given $Y$ depends
only on the indicator functions $(I_{\{x\leftrightarrow y\}})_{x,y\in \De}$
and $(I_{\{x\leftrightarrow \infty\}})_{x \in \De}$, since the conditional
distribution of $X$ on $\De$ is uniform over all elements of
$\{1,\ldots,q\}^\De$
such that firstly $X(x)=X(y)$ whenever $x\leftrightarrow y$, and secondly
$X(x)=i$ whenever $x\leftrightarrow\infty$. Hence, the desired convergence
of
$\mu_{\beta,q,\Lambda}^i(f)$ follows if we can show that the joint
distribution of $(I_{\{x\leftrightarrow y\}})_{x,y\in \De}$
and $(I_{\{x\leftrightarrow \infty\}})_{x \in \De}$ converges as
$n\rightarrow\infty$. This, however, follows from Lemma
\ref{lem:free_and_wired_measure} upon noting that
$(I_{\{x\leftrightarrow y\}})_{x,y\in \De}$
and $(I_{\{x\leftrightarrow \infty\}})_{x \in \De}$
are increasing functions. $\Cox$

\subsection{Phase transition in the Potts model}
\label{sect:PT_Potts}

As promised, this subsection is devoted to proving Theorems
\ref{thm:Ising} and \ref{thm:Potts}, using random-cluster arguments.
The original source for the material in this subsection is
Aizenman, Chayes, Chayes and Newman \cite{ACCN1}; see also \cite{H7}
for a slightly different presentation.

We consider the Potts model on $\Z^d$, $d\geq 2$. All the arguments to be
used
here,
except those showing that the critical inverse temperature $\beta_c$ is
strictly between $0$ and $\infty$, go through on arbitrary infinite
graphs; we stick to the $\Z^d$ case for definiteness and simplicity of
notation.
We consider the limiting Gibbs measures $\mu_{\beta,q}^i$ obtained in
Proposition \ref{prop:free_and_ordered_Potts_state}.
For $i\in\{1,\ldots,q\}$, these play a role similar to that
of the ``plus'' and ``minus'' measures $\mu^+_\beta$ and $\mu^-_\beta$ for
the Ising model. In fact, we have the following result which extends
Theorem \ref{thm:Ising_sandwich} to the
Potts model and also gives a characterization of phase transition in terms
of percolation in the random-cluster model.
\begin{thm} \label{thm:Potts_sandwich} Let $\b>0$ and $p=1-e^{-2\b}$.
For any $x\in \Z^d$ and any $i\in\{1,\ldots,q\}$, the following statements
are equivalent.
\begin{description}
\item{\rm (i)} There is a unique  Gibbs measure for the
$q$-state Potts model on $\Z^d$ at inverse temperature $\beta$.
\item{\rm (ii)} $\mu_{\beta,q}^i(X(x)=i)={1}/{q}$.
\item{\rm (iii)} $\phi^1_{p,q}(x\leftrightarrow \infty)=0$.
\end{description}
\end{thm}

\medskip\noindent
As we will see in a moment, it is the percolation criterion (iii) which is
most convenient to apply. In this context  we note that
\begin{equation}\label{limit_interchange}
\phi^1_{p,q}(x\leftrightarrow \infty)=
\inf_{\La,\De}\phi^1_{p,q,\La}(x\leftrightarrow \De^c) =
\lim_{\La\uparrow\Z^d}\phi^1_{p,q,\La}(x\leftrightarrow \La^c) \;,
\end{equation}
where $\{x\leftrightarrow \De^c\}$ stands for the event that there exists
an open path from $x$ to some site in $\De^c$. This follows from
\rf{eq:RC_monotone} and the fact that $\{x\leftrightarrow \De^c\}$
decreases to $\{x\leftrightarrow \infty\}$ as $\De\uparrow\Z^d$.

The usefulness of the percolation criterion is demonstrated by the next
result which extends the scenario for Bernoulli percolation to the
random-cluster model. Together with Theorem \ref{thm:Potts_sandwich}, this
gives
Theorem \ref{thm:Potts} with $\b_c=\frac{1}{2}\log(1-p_c)$.
\begin{prop} \label{prop:Potts_magnetization} For the random-cluster model
on $\Z^d$, $d\ge2$, and any fixed $q\ge1$,
there exists a percolation threshold $p_c\in(0,1)$
(depending on $d$ and $q$)
such that
\[
\phi^1_{p,q}(x\leftrightarrow \infty)\left\{\begin{array}{cl}=0 &\mbox{for
} p<p_c\;,\\ >0 &\mbox{for } p>p_c\;.\end{array}\right.
\]
\end{prop}
{\bf Proof:}
The statement of the proposition consists of the following three
parts:
\begin{description}
\item{(i)} $\phi^1_{p,q}(x\leftrightarrow \infty)=0$ for $p$
sufficiently small,
\item{(ii)} $\phi^1_{p,q}(x\leftrightarrow \infty)>0$ for $p$
sufficiently close to 1, and
\item{(iii)} $\phi^1_{p,q}(x\leftrightarrow \infty)$ is increasing in $p$.
\end{description}
We first prove (i).
Suppose $p<p_c(\Z^d, \mbox{bond})$, the critical value for Bernoulli
bond percolation on $\Z^d$. For $\eps>0$, we can then pick $\De$ large
enough
so that
\[
\phi_p(0\leftrightarrow\De^c)\leq \eps \, .
\]
By Corollary \ref{cor:FKG_for_RC} (b), we have that
the projection of $\phi_p$ on $\{0,1\}^{\B_\De^1}$ stochastically
dominates the
projection of $\phi_{p,q,\La}^1$ on $\{0,1\}^{\B_\De^1}$ for any
$\La\supset\De$, so that
\[
\phi^1_{p,q,\La}(0\leftrightarrow\infty) \leq
\phi^1_{p,q,\La}(0\leftrightarrow\De^c)  \leq \eps
\]
for any $\La\supset\De$. Since $\eps$ was arbitrary, we find
\[
\lim_{\La\uparrow\Z^d}\phi^1_{p,q,\La}(0\leftrightarrow\infty)=0
\]
which in conjunction with \rf{limit_interchange}
implies (i).

Next, (ii) can be established by a similar argument: Let $p$ be such that
$p^\ast={p}/[{p+(1-p)q}]>p_c(\Z^d,\mbox{bond})$.
Corollary \ref{cor:FKG_for_RC} (c) then shows that
$
\phi^1_{p,q,\La} \geqd \phi_{p^\ast}
$
for every $\La$, so that
\[
\lim_{\La\uparrow\Z^d}\phi^1_{p,q,\La}(0\leftrightarrow\infty)>0,
\]
proving (ii).

To check (iii) we note that
Corollary \ref{cor:FKG_for_RC} (d) implies that, for any $\La$,
\begin{equation} \label{eq:RC_domination}
\phi^1_{p_1,q,\La} \leqd \phi^1_{p_2,q,\La} \mbox{ whenever $p_1\leq p_2$.}
\end{equation}
This proves (iii) and thereby the proposition.  $\Cox$

\medskip\noindent
Before the proof of Theorem \ref{thm:Potts_sandwich} we need another
definition and a couple of lemmas.
For a finite box $\Lambda$ in ${\bf Z}^d$ and a spin configuration
$\xi\in \{1, \ldots, q\}^{\partial\Lambda}$, let
\[
A^i_\xi= \{x\in \partial \Lambda : \, \xi(x)=i \}
\]
for $i=1, \ldots, q$. We now define the {\em random-cluster
distribution
$\phi^\xi_{p,q, \Lambda}$ for $\Lambda$ with boundary condition $\xi$},
as the probability measure on $\{0,1\}^{\B^1_\Lambda}$ which to each
$\eta\in \{0,1\}^{\B^1_\Lambda}$ assigns probability proportional to
\[
I_{D(\xi, \eta)}\, \bigg\{\prod_{e\in \B_\La^1}
p^{\eta(e)}(1-p)^{1-\eta(e)}\bigg\} \, q^{k_\xi(\eta)}
\]
where $k_\xi(\eta)$ is the number of connected components in $\eta$
that do not intersect $\partial\Lambda$, and $D(\xi, \eta)$ is the event
that there is no open path in $\eta$ connecting any two vertices
in $A^i_\xi$ and $A^j_\xi$ for any $i\neq j$.
\begin{lem} \label{lem:mixed_RC_to_Potts}
Let $p=1-e^{-2\beta}$, let $\Lambda$ be a finite region in $\Z^d$,
and fix some boundary condition  $\xi\in \{1, \ldots,
q\}^{\partial\Lambda}$. Suppose that we pick a random spin
configuration $X\in\{1, \ldots, q\}^\Lambda$ as follows.
\begin{enumerate}
\item
Pick $Y\in \{0,1\}^{\B^1_\Lambda}$ according to $\phi^\xi_{p,q, \Lambda}$.
\item
For each $i\in\{1, \ldots, q\}$ and each connected component $C$ in $\eta$
intersecting $A^i_\xi$, assign spin $i$ to every vertex in $C$.
\item For all other connected components $C$ in $Y$, pick a spin at
random (uniformly) from $\{1,\ldots,q\}$, assign this spin to every vertex
of $C$, and do this independently for different connected components.
\end{enumerate}
Then $X$ is distributed according to the Gibbs distribution $\mu^\xi_{\b,q,\La}$
for the Potts model on $\La$ with boundary condition $\xi$.
\end{lem}
{\bf Proof:} This is a straightforward generalization of the proofs of
Theorem \ref{thm:RC_coupling} and Corollary \ref{cor:RC_to_Potts}.
$\Cox$
\begin{lem} \label{lem:mixed_RC_domination}
With notation as above, we have, for any
$\xi\in\{1, \ldots, q\}^{\partial\Lambda}$, that the projection of
$\phi^1_{p,q,\Lambda}$ on $\{0,1\}^{\B^1_{\Lambda}}$ stochastically
dominates $\phi^\xi_{p,q,\Lambda}$.
\end{lem}
{\bf Proof:}
Just write down single-edge conditional distributions for
$\phi^\xi_{p,q, \Lambda}$ and $\phi^1_{p,q, \Lambda}$ (as in
Lemma \ref{lem:single_edge}) and invoke Theorem \ref{thm:holley}.
$\Cox$

\medskip\noindent
We are finally ready for the proof of Theorem \ref{thm:Potts_sandwich}.

\medskip\noindent
{\bf Proof of Theorem \ref{thm:Potts_sandwich}:} We begin with the
implication (i) $\Rightarrow$ (ii). If there is a unique Gibbs measure for
the $q$-state
Potts model on $\Z^d$ at inverse temperature $\beta$, then we have in
particular that
\[
\mu^1_{\beta,q} = \ldots = \mu^q_{\beta,q} \, .
\]
But since by symmetry $\mu^i_{\beta,q}(X(x)=i)=\mu^j_{\beta,q}(X(x)=j)$
for any $i,j\in\{1,\ldots, q\}$, we must have
$\mu^i_{\beta,q}(X(x)=i)= {1}/{q}$, and (ii) is established.

Next we turn to the implication (ii) $\Rightarrow$ (iii). By the
Edwards--Sokal coupling
of edge and site processes introduced before
Proposition \ref{prop:free_and_ordered_Potts_state}, we have
\begin{eqnarray}
\mu^i_{\beta,q}(X(0)=i) & = &
\lim_{\La\uparrow\Z^d}\mu^i_{\beta,q, \Lambda}(X(0)=i) \nonumber
\\
& = & \frac{1}{q} + \frac{q-1}{q}
\lim_{\La\uparrow\Z^d}\phi^1_{p,q,\La}(0\leftrightarrow\La^c)\, .
\label{eq:Potts_magnetization_connectivity}
\end{eqnarray}
Together with \rf{limit_interchange}, the result follows.

Most of the work is needed for the implication (iii) $\Rightarrow$
(i). Roughly speaking, the absence of percolation in the
random-cluster model implies that every finite region is cut off
from infinity by a set of closed edges. Thus, independently of what
happens macroscopically, the local spins feel as if they are in a
system with free boundary condition. This makes a phase transition
impossible.

To make this intuition precise we let $\mu$ be an arbitrary Gibbs measure
for the Potts model at inverse temperature $\b$. We will show that
$\mu=\mu^0_{\b,q}$, the limiting measure with free boundary condition
obtained in Proposition \ref{prop:free_and_ordered_Potts_state}. We fix
any local observable $f$ and some $\eps>0$. We then can
find a finite box $\De\subset\Z^d$ such that
\begin{equation}\label{free_bc_approximation}
\left|\;\mu^0_{\b,q,\Gamma}(f)-\mu^0_{\b,q}(f)\,\right|<\eps \ \mbox{ for
all finite  } \Gamma\supset\De\;.
\end{equation}
In view of \rf{limit_interchange} we can also choose a finite box
$\La\supset\De$ satisfying
$\phi^1_{p,q,\La}(x\leftrightarrow\La^c)<\eps/|\De|$ for all $x\in\De$,
and thus
\[
\phi^1_{p,q,\La}(\De\leftrightarrow\La^c)<\eps\;.
\]
Here, $\{\De\leftrightarrow\La^c\}$ is the event that there exists an open
path from $\De$ to $\La^c$.

Consider the complementary event $C=\{\De\leftrightarrow\La^c\}^c$. For
any edge configuration $\eta\in C$, $\De$ is cut off from $\La^c$ by a set
of closed edges. Indeed, let $\Gamma$ be the union of all $\eta$-open
clusters meeting $\De$. Then
\begin{description}
\item{(a)} $\De\subset\Gamma\subset\La$, and
\item{(b)} $\eta(e)=0$ for all edges from $\Gamma$ to $\Gamma^c$, i.e.,
$\eta\equiv 0$ on $\B^1_\Gamma\setminus\B^0_\Gamma$.
\end{description}
Since these properties are stable under finite unions, there exists a
largest set $\Gamma(\eta)$ satisfying (a) and (b). The maximality implies
that, for each fixed region $\Gamma$, the event
$\{\eta:\Gamma(\eta)=\Gamma\}$ only depends on the status of the edges off
$\B^0_\Gamma$.

Coming to the core of the argument, we fix any boundary spin configuration
$\xi\in \{1, \ldots, q\}^{\partial\La}$ and consider the Edwards--Sokal
coupling $P=P^\xi_{p,q,\La}$ of
$\mu^\xi_{\b,q,\La}$ and $\phi^\xi_{p,q, \La}$ introduced in
Lemma \ref{lem:mixed_RC_to_Potts}.
We write $X$ for the random spin configuration in $\{1,\ldots,q\}^{\Z^d}$
and $Y$ for the random edge configuration in $\{0,1\}^\B$. Then
$P(f\circ X)=\mu^\xi_{\b,q,\La}(f) $, and
\[
P(Y\in C)=\phi^\xi_{p,q, \La}(C)\ge \phi^1_{p,q, \La}(C)>1-\eps
\]
by Lemma \ref{lem:mixed_RC_domination} and the choice of $\La$.
Assuming without loss of generality that $\|f\|\le 1$, we can
therefore conclude that
\begin{equation}\label{eq:cutset_approx}
\Big|P\big(f\circ X \,\big|\,Y\in C\big) -\mu^\xi_{\b,q,\La}(f)\Big|<2\eps\;.
\end{equation}
However, the conditional expectation on the left is an average of the
conditional expectations $P(f\circ X \,|\,\Gamma(Y)=\Gamma)$ with
$\De\subset\Gamma\subset\La$, and these in
turn are averages of conditional expectations of the form
\[
P(f\circ X \;|\;Y=\eta \mbox{ off }\B^1_\Gamma,\; Y\equiv 0 \mbox{ on
}\B^1_\Gamma\setminus\B^0_\Gamma )
\]
which, by construction of $P$, are equal to $\mu^0_{\b,q,\Gamma}(f)$.
Together with \rf{free_bc_approximation} and \rf{eq:cutset_approx}, we
conclude that
\[
\left|\;\mu^\xi_{\b,q,\La}(f)-\mu^0_{\b,q}(f)\,\right|<3\eps\;.
\]
Taking the $\mu$-average over $\xi$ and letting $\eps\to0$ we finally get
$\mu(f)=\mu^0_{\b,q}(f)$, and the proof is complete.  $\Cox$

\subsection{Infinite volume random-cluster measures}
\label{sect:infinite_FK}

The random-cluster arguments used
in the previous subsections for studying infinite volume Ising and Potts
models only required defining finite volume random-cluster distributions,
although we have seen already the limiting random-cluster measures
$\phi^b_{p,q}$. (Their existence was convenient in the formulation of
Theorem \ref{thm:Potts_sandwich}, but not really needed
for the arguments.) Recent
years have nevertheless witnessed a rapid development of a theory for
infinite volume random-cluster measures, defined in the DLR spirit. Here
we shall discuss the basics of such a theory. A similar theory of
infinite volume Edwards--Sokal measures for joint spin and edge
distributions was recently developed in \cite{BCCK}.

Let $G$ be an infinite (locally finite) graph with vertex set $\L$ and
edge set $\B$. Fix
$p\in [0,1]$ and $q>0$, and let $B\subset \B$ be a finite simply
connected region; since
the random-cluster model lives on edges rather than on vertices, we let
``region'' refer to edge sets rather than vertex sets in this subsection.
Let $V(B)=\{x\in \L : \, \exists\, e \in B \mbox{ incident to } x\}$. For
an edge configuration $\xi \in \{0,1\}^{B^c}$, define the {\em
random-cluster distribution} $\phi^{B, \xi}_{p,q}$ in $B$ as the
probability
measure on $\{0,1\}^\B$ in which each $\eta\in \{0,1\}^\B$ is assigned
probability proportional to
\begin{equation} \label{eq:finite_volume_RC}
I_{\{\eta=\xi \mbox{ \scriptsize off } B\}}
\left\{\prod_{e\in B}
p^{\eta(e)}(1-p)^{1-\eta(e)}\right\}q^{k(\eta,B)}\;,
\end{equation}
where $k(\eta,B)$ is the number of connected components of $\eta$
which intersect $V(B)$. This is a generalization of the  random-cluster
distributions $\phi^b_{p,q,\La}$
defined in Section \ref{sect:Inf_vol_limits},
which are recovered by taking $B=\B^b_\La$ and $\xi\equiv b$,
$b\in\{0,1\}$. It is easy to see that the random-cluster distributions
are consistent in the sense that conditioning on a configuration in some
$B'\subset B$ yields the corresponding random-cluster distribution in
$B\setminus B'$.
\begin{defn} \label{defn:RC_type_I}
A probability measure $\phi$ on $\{0,1\}^\B$ is said to be a
random-cluster measure with parameters $p$ and $q$ if its conditional
probabilities
satisfy
\[
\phi(\eta \, | \, \xi) \equiv \phi(\eta \mbox{ \rm in }B\, | \, \xi\mbox{ \rm off }B) = \phi^{B, \xi}_{p,q}(\eta)
\]
for all finite simply connected
$B\subset \B$, $\phi$-almost all $\xi$ and all $\eta$ such
that $\eta=\xi$ off $B$.
\end{defn}
This is the direct analogue of the definition of a Gibbs measure in the
random-cluster setting. There is, however, also another possibility which
differs from the preceding one for graphs in which the complements of
finite regions are not connected. The idea is to connect all infinite
clusters  at infinity. This corresponds to a one-point compactification of
$G$. Accordingly, we shall use a prefix `C' which stands for
``compactified''.
Thus, we define a {\em C-random-cluster distribution} $\hat{\phi}^{B,
\xi}_{p,q}$ as in
(\ref{eq:finite_volume_RC}), except that $k(\eta,B)$ is replaced by
$\hat{k}(\eta,B)$, defined as the number of all {\em finite} connected
components
of $\eta$ intersecting $V(B)$.
\begin{defn} \label{defn:RC_type_II}
A probability measure $\phi$ on $\{0,1\}^\B$ is said to be a
C-random-cluster measure for $p$ and $q$ if its conditional probabilities
satisfy
\[
\phi(\eta \, | \, \xi) = \hat{\phi}^{B, \xi}_{p,q}(\eta)
\]
for all finite simply connected
$B\subset \B$, $\phi$-almost all $\xi$ and all $\eta$ such
that $\eta=\xi$ off $B$.
\end{defn}
The study of random-cluster measures in the case $G=\Z^d$ was
initiated by Grimmett \cite{Gr5}, and about
simultaneously by Pfister and Vande Velde \cite{PV} and Borgs and
Chayes \cite{BC}; see also Sepp\"al\"ainen \cite{Sep}
for some even more recent developments.
The C-variant (Definition \ref{defn:RC_type_II}) was introduced in the
regular tree case by H\"aggstr\"om \cite{H2,H4}, and further studied
in a general graph context
by Jonasson \cite{J}. In the following, we will try to convince the reader
that random-cluster measures of both types are of interest, and also
discuss their relation to each other.

We shall concentrate mainly on the case $q\geq 1$. The reason for this is
that it is only for $q\geq 1$ that the conditional probability in
(\ref{eq:single_edge}) is increasing in $\eta$, which allows the use
of the stochastic domination and correlation results in Section
\ref{sect:coupling_SD} (Theorems \ref{thm:holley} and \ref{thm:FKG}).
For $q<1$, these tools are not available, and for this reason the
random-cluster model with $q<1$ is much less understood than
the $q\geq 1$ case, although in Grimmett \cite{Gr5}, H\"aggstr\"om
\cite{H1}
and Sepp\"al\"ainen \cite{Sep} one can find at least
some results in the $q<1$ regime of the parameter space.

We write $\phi^{B, 0}_{p,q}$ for $\phi^{B, \xi}_{p,q}$ with $\xi \equiv 0$,
and  $\hat{\phi}^{B,1}_{p,q}$ for $\hat{\phi}^{B,\xi}_{p,q}$ with
$\xi \equiv 1$.
We can omit the hat when  $G=\Z^d$, because there is always exactly
one infinite cluster regardless of the configuration on $B$ (this is
related to Proposition \ref{prop:types_I_II_identical} below).
On the other hand, the two different ways
of counting clusters with wired boundary condition are not equivalent for
all graph structures; a simple counterexample is $G=\T_d$ in which the
wired
boundary condition gives rise to several infinite clusters.

Still in the context of general infinite graphs, we write $B\uparrow\B$
for the limit along some (any) sequence of finite simply connected regions
increasing to $\B$ in the usual way.
In complete analogy to  \rf{eq:RC_monotone} and Lemma
\ref{lem:free_and_wired_measure} we then obtain the following monotonicity
and convergence result.
\begin{lem} \label{lem:wired/free_monotone}
For $p\in[0,1]$, $q\geq 1$, and any two finite bond sets $B_1\subset B_2$,
we have
\begin{description}
\item{\rm(i) } $\phi_{p,q}^{B_1, 0}\leqd\phi_{p,q}^{B_2, 0} $, so that the
limit
$\phi_{p,q}^{0}=\lim\limits_{B\uparrow\B}
\phi_{p,q}^{B, 0}$ exists; and
\item{\rm(ii)} $\hat{\phi}_{p,q}^{B_1, 1}\geqd\hat{\phi}_{p,q}^{B_2,1} $,
so that the limit
$\hat{\phi}_{p,q}^{1}=\lim\limits_{B\uparrow\B}\hat{\phi}_{p,q}^{B, 1}$
exists.
\end{description}
\end{lem}

\medskip\noindent
Now let $\phi$
be any random-cluster measure of either type, compactified or not,
with the given parameters $p$ and $q$. Further application of
Theorem \ref{thm:holley} (or Corollary \ref{cor:FKG_for_RC}) implies that
\[
\phi_{p,q}^{0} \leqd \phi  \leqd \hat{\phi}_{p,q}^{1} \, .
\]
This is analogous to the sandwiching relation (\ref{eq:Ising_sandwich})
for the Ising model.

Furthermore, the arguments of Section \ref{sect:PT_Potts} go through to
show that the $q$-state Potts model on $G$ at inverse temperature $\beta$
has a unique Gibbs measure if and only if the
$\hat{\phi}_{p,q}^{1}$-probability of having an infinite cluster
is $0$, where as usual $p=1-e^{-2\beta}$.

For Definitions \ref{defn:RC_type_I} and \ref{defn:RC_type_II} to be of
interest, we have to establish at least the existence of
random-cluster measures of the two types.
The following theorem tells us that at least
for $q \geq 1$, such measures do exist. (The existence problem for $q<1$
remains open in the setting of general graphs, although existence has been
established for $\Z^d$ and $\T_d$; see \cite{Gr5} and \cite{H2},
respectively.)
\begin{thm} \label{thm:RC_existence}
For $p\in [0,1]$ and $q\geq 1$, we have that (i)  $\phi_{p,q}^{0}$
is a random-cluster measure, and (ii) $\hat{\phi}_{p,q}^{1}$ is a
C-random-cluster measure.
\end{thm}
For the proof we use the following lemma which characterizes random-cluster
measures in terms of single-edge conditional probabilities.
For $e=\langle x y \rangle\in \B$ and $\xi \in \{0,1\}^{\B \setminus
\{e\}}$,
we write as usual $\{x\leftrightarrow y\}$ for the event that there exists
an open path in $\xi$ from $x$ to $y$. We also write
$\{x \stackrel{C}{\longleftrightarrow} y\}$
for the event that there either exists an open path in $\xi$ from $x$ to
$y$,
or $x$ and $y$ are both in infinite clusters of $\xi$. We think of this
C-connectivity notion $x \stackrel{C}{\longleftrightarrow} y$ as
allowing paths between $x$ and $y$ to go ``via infinity''.
\begin{lem} \label{lem:single_edge_characterizations}
Fix $p\in [0,1]$ and $q>0$, and let $\phi$ be
a probability measure on $\{0,1\}^\B$. Then $\phi$ is a
random-cluster measure for $p$ and $q$ if and only if
for each $e=\langle x y \rangle\in \B$ and
$\phi$-a.e.\ $\xi \in \{0,1\}^{\B \setminus \{e\}}$ we have
\begin{equation} \label{eq:single_edge_type_I}
\phi(e \mbox{ \em is open }|\, \xi) = \left\{
\begin{array}{ll}
p & \mbox{\em if } x\leftrightarrow y \mbox{ \em in }\xi\\
\frac{p}{p+(1-p)q} & \mbox{\em otherwise,}
\end{array} \right.
\end{equation}
Similarly, $\phi$ is a
C-random-cluster measure for $p$ and $q$ if and only if
\rf{eq:single_edge_type_I} holds with $x \stackrel{C}{\longleftrightarrow}
y$
instead of $x\leftrightarrow y$.
\end{lem}
{\bf Proof:} We consider only the first statement, as the C-case is
completely similar. For the ``only if'' part we only need to note
that the right-hand side of \rf{eq:single_edge_type_I} is equal to
$\phi_{p,q}^{\{e\},\xi}( e \mbox{ is open })$. Passing to the
``if'' part, we may restrict ourselves to the case of $p\in (0,1)$.
Assume that $\phi$ satisfies (\ref{eq:single_edge_type_I}) for each
$e=\langle x y \rangle\in \B$ and $\phi$-a.e.\ $\xi \in
\{0,1\}^{\B \setminus \{e\}}$. Let $B\subset \B$ be some finite
edge set. We need to show  that the conditional distribution
$\phi(\cdot  | \xi)$ of $\phi$ given the configuration $\xi$ on
$B^c$ equals $\phi_{p,q}^{B,\xi}$ for $\phi$-a.e.\ $\xi$. For this,
it suffices to check that for any two configurations $\eta, \eta'
\in \{0,1\}^\B$ which agree with $\xi$ on $B^c$ we have
\begin{equation} \label{eq:desired_ratio}
\frac{\phi(\eta \, | \, \xi)}{\phi(\eta' \, | \, \xi)}=
\frac{\{\prod_{e\in B}
p^{\eta(e)}(1-p)^{1-\eta(e)}\}q^{k(\eta,B)}}{\{\prod_{e\in B}
p^{\eta'(e)}(1-p)^{1-\eta'(e)}\}q^{k(\eta',B)}}
\end{equation}
with $k$ defined as in (\ref{eq:finite_volume_RC}). If $\eta$ and
$\eta'$ differ only at a single edge $e$, then $k(\eta,B)-k(\eta',B)=
k(\eta,e)-k(\eta',e)$, whence (\ref{eq:desired_ratio})
is immediate from (\ref{eq:single_edge_type_I}). In the general case, we
interpolate $\eta$ and $\eta'$ by a sequence of configurations
which successively differ in at most
one edge, and use a telescoping argument. $\Cox$

\medskip\noindent
{\bf Proof of Theorem \ref{thm:RC_existence}:}
We prove (ii) only, as (i) follows from a similar argument
and is also  better known, see e.g.\ Borgs and Chayes \cite{BC}.
Fix $e=\langle x y \rangle \in \B$. By Lemma
\ref{lem:single_edge_characterizations}, it is sufficient to establish
 \rf{eq:single_edge_type_I} with the C-connectivity relation
$\stackrel{C}{\longleftrightarrow}$ in place of the standard connectivity
relation $\leftrightarrow$.

Let $B_1, B_2, \ldots$ be an increasing sequence of finite edge sets
containing $e$ and
converging to $\B$ in the usual sense. We write, with slight
abuse of notation, $\xi({B_i})$ for the restriction of $\xi$
to $B_i\setminus \{e\}$.
We recall from the martingale convergence theorem that
\begin{equation} \label{eq:martingale_convergence}
\hat{\phi}^{1}_{p,q}(e \mbox{ is open }| \, \xi) = \lim_{j\rightarrow
\infty}
\hat{\phi}^{1}_{p,q}(e \mbox{ is open }| \, \xi(B_j))
\end{equation}
for $\hat{\phi}^{1}_{p,q}$-a.e.\ $\xi\in \{0,1\}^{\B\setminus \{e\}}$.

We suppose first that
$x \stackrel{C}{\longleftrightarrow}y$ fails in $\xi$. Then at least
one of the vertices $x$ and $y$ is in a finite cluster of $\xi$, and
consequently there is some $m$ (depending on $\xi$) such that
$\neg(x \stackrel{C}{\longleftrightarrow}y)$ can be verified by
just looking at $\xi(B_m)$. (This is a consequence of the special concept
of C-connectivity.) For any $n\ge j\geq m$ we then have
\[
\hat{\phi}^{B_n,1}_{p,q}(e \mbox{ is open }| \, \xi(B_j))=
\frac{p}{p+(1-p)q}
\]
so that by the definition of $\hat{\phi}^{1}_{p,q}$ we get
\[
\hat{\phi}^{1}_{p,q}(e \mbox{ is open }| \, \xi(B_j))=
\frac{p}{p+(1-p)q}
\]
by letting $n\rightarrow\infty$. Then we let $j\rightarrow\infty$ and use
(\ref{eq:martingale_convergence}) to deduce the C-version of
\rf{eq:single_edge_type_I} in the case $\xi\notin \{x
\stackrel{C}{\longleftrightarrow}y\}$.

We go on to the case $\xi\in \{x \stackrel{C}{\longleftrightarrow}y\}$. In
analogy to (\ref{eq:martingale_convergence}), we have
\[
\lim_{j\rightarrow\infty} \hat{\phi}^{1}_{p,q}
(x \stackrel{C}{\longleftrightarrow}y\, | \, \xi(B_j)) =1
\]
for $\hat{\phi}^{1}_{p,q}$-a.e.\ $\xi\in\{x
\stackrel{C}{\longleftrightarrow}y\}$. For such $\xi$ and any $\eps>0$,
we can thus find an $m$ (depending on $\xi$) such that
\begin{equation} \label{eq:Large_enough_set}
\hat{\phi}^{1}_{p,q}(x \stackrel{C}{\longleftrightarrow}y
\, | \, \xi(B_j))\geq 1-\eps \,
\end{equation}
for any $j \geq m$. Next we use the definition of $\hat{\phi}^{1}_{p,q}$.
For any $n=1,2,\ldots$,  let $Y$ and $Y_n$ be $\{0,1\}^\B$-valued random
edge configurations with distributions $\hat{\phi}^{1}_{p,q}$ and
$\hat{\phi}^{B_n, 1}_{p,q}$  satisfying
$Y_n \succeq Y$;
this is possible by Lemma \ref{lem:wired/free_monotone}. We write $P_n$
for the
probability measure underlying this coupling. By the same lemma and the
order relation $Y_n \succeq Y$,
we have
\begin{equation} \label{eq:setwise_limit}
\lim_{n\rightarrow\infty} P_n(Y_n(B_j) \neq Y(B_j))=0 \, .
\end{equation}
Since $Y_n\in\{x \stackrel{C}{\longleftrightarrow}y\}$ whenever
$Y\in\{x \stackrel{C}{\longleftrightarrow}y\}$, we can write
\begin{eqnarray*}
&&\hspace*{-4ex}\left|\;\hat{\phi}^{B_n, 1}_{p,q}(x
\stackrel{C}{\longleftrightarrow}y
\, , \, \xi(B_j))-\hat{\phi}^{1}_{p,q}(x \stackrel{C}{\longleftrightarrow}y
\, ,  \, \xi(B_j))\;\right|\\[1ex]
&\le& P_n\Big(\big\{x \stackrel{C}{\longleftrightarrow}y \mbox{ in }Y_n,
\; Y_n(B_j)=\xi(B_j)  \big\}\; \triangle\; \big\{x
\stackrel{C}{\longleftrightarrow}y \mbox{ in }Y,\; Y(B_j)=\xi(B_j)  \big\}\Big)
\\[1ex]
&\le&
\left(\hat{\phi}^{B_n, 1}_{p,q}(x \stackrel{C}{\longleftrightarrow}y)
-\hat{\phi}^{1}_{p,q}(x \stackrel{C}{\longleftrightarrow}y)\right) +
P_n(Y_n(B_j) \neq Y(B_j))\;,
\end{eqnarray*}
where $\triangle$ denotes symmetric difference.
Since $\{x \stackrel{C}{\longleftrightarrow}y\}$ is the decreasing limit
of the local events $\{x \leftrightarrow y \mbox{ in  }\De\}\cup
\{x \leftrightarrow \De^c,\, y \leftrightarrow \De^c\}$ as $\De\uparrow\L$,
an analogue of \rf{limit_interchange} together with
\rf{eq:setwise_limit} shows that the last expression tends to
zero as $n\ti$. It follows that
\[
\lim_{n\rightarrow\infty}
\hat{\phi}^{B_n, 1}_{p,q}(x \stackrel{C}{\longleftrightarrow}y
\, | \, \xi(B_j))=
\hat{\phi}^{1}_{p,q}(x \stackrel{C}{\longleftrightarrow}y
\, | \, \xi(B_j))
\]
which is at least $1-\eps$ by (\ref{eq:Large_enough_set}).
But since
$
\hat{\phi}^{B_n, 1}_{p,q}(e \mbox{ is open } | \,
\xi')= p
$
for each $n$ and all $\xi'\in\{x \stackrel{C}{\longleftrightarrow}y\}$,
we get
\[
p - \eps \leq \lim_{n\rightarrow\infty}
\hat{\phi}^{B_n, 1}_{p,q}(e \mbox{ is open }|\, \xi(B_j)) \leq p \, .
\]
Hence,
\[
p - \eps \leq \hat{\phi}^{1}_{p,q}(e \mbox{ is open }|\, \xi(B_j)) \leq p
\, ,
\]
and since $\eps$ was arbitrary we can use
(\ref{eq:martingale_convergence}) to deduce the C-version of
\rf{eq:single_edge_type_I} in the case $\xi\in \{x
\stackrel{C}{\longleftrightarrow}y\}$.
$\Cox$

\medskip\noindent
Let us now briefly address the issue of whether the two
types of random-cluster measures are any different. The
following result says that very often they are the same.
\begin{prop} \label{prop:types_I_II_identical}
Let $\phi$ be a probability measure on $\{0,1\}^\B$ with
\[
\phi(\exists \mbox{ \em at most one infinite open cluster})=1 \, .
\]
Then, for any $p\in [0,1]$ and $q>0$, $\phi$ is a random-cluster measure
for $p$ and $q$ if and only if it is a C-random-cluster measure for $p$
and $q$.
\end{prop}
This means that whenever ``uniqueness of the infinite cluster'' can be
verified, the two types of random-cluster measures coincide. An example
is obtained if we consider translation invariant random-cluster measures for
$\Z^d$,
since the Burton--Keane uniqueness theorem (Theorem \ref{thm:Burton_Keane})
applies in this situation. For $\Z^d$, the measures $\phi^{0}_{p,q}$
and $\hat{\phi}^{1}_{p,q}$
are translation invariant, by Lemma \ref{lem:free_and_wired_measure}.
On the other hand, uniqueness of the infinite
cluster typically fails on trees, leading to very different behavior for
the two types of random-cluster measures; see \cite{H2,H4}
for a discussion.

\medskip\noindent
{\bf Proof of Proposition \ref{prop:types_I_II_identical}:}
For $p=0$ or $1$ the result is trivial, so we may assume that $p\in (0,1)$.
The conditional probabilities in (\ref{eq:single_edge_type_I}) and
its C-counterpart differ only on the event
\[
A_{xy}=\{x \stackrel{C}{\longleftrightarrow} y\} \setminus
\{x \leftrightarrow y\} \, .
\]
Hence if $\phi$ is a random-cluster measure but not a
C-random-cluster measure (or vice versa), then $A_{xy}$ has to have
positive $\phi$-probability for some edge
$e=\langle x y \rangle\in \B$. But then the event
$A_{xy} \cap \{e \mbox{ is closed}\}$  has positive
$\phi$-probability, and since this event implies the existence of
at least two infinite clusters, we are done. $\Cox$

\medskip\noindent
Much of the study of infinite volume random-cluster measures that has
been done so far concerns the issue of uniqueness (or non-uniqueness)
of random-cluster measures. A discussion of this issue would, however,
lead us too far, so instead we advise the reader to consult Grimmett
\cite{Gr5}, H\"aggstr\"om \cite{H2} and Jonasson \cite{J} to find out
what is known and what is conjectured in this field.

\subsection{An application to percolation in the Ising model}

In Theorem \ref{thm:Ising_percolation_monotone_in_h} we have seen that
the probability of percolation of plus spins in the Ising model is an
increasing function of the external field.
A much harder question is to determine
monotonicity properties of percolation probabilities as $\beta$ (rather
than $h$) is varied. An interesting open problem
is to decide whether for $G=\Z^d$, $d\geq 2$, the probability
\[
\mu^+_{\beta}(x \stackrel{+}{\longleftrightarrow}\infty )
\]
is increasing in $\beta$. Here we write $\mu^+_{\beta}$ for the plus phase
in the Ising model at inverse temperature $\beta$ with external field
$h=0$, and $\{x \stackrel{+}{\longleftrightarrow}\infty \}$ is the
event that there exists an infinite path of plus spins starting at $x$.
At first sight, one might be seduced into
thinking that this would be a consequence of the connection between
Ising and random-cluster models, and the stochastic monotonicity
of random-cluster measures as $p$ varies; see (\ref{eq:RC_domination}).
However, such a conclusion is unwarranted. For example, in the
coupling of Theorem \ref{thm:RC_coupling} the existence of an
open path between $x$ and $y$ in the random-cluster representation is a
sufficient {\em but not necessary} condition for $x$ and $y$ to be in
the same spin cluster. In fact, H\"aggstr\"om \cite{H5}
showed, by means of a simple counterexample and in response to a question
of Cammarota \cite{Cam}, that the probability that $x$ and $y$ are in the
same
spin cluster need not be increasing in $\beta$, and similarly for
the expected size of the spin cluster containing $x$.

However,  when the underlying graph is a tree, monotonicity
in $\beta$ of the probability of plus percolation
can be established:
\begin{thm} \label{thm:tree_monotonicity}
For the Ising model on the regular tree $\T_d$, $d\geq 2$, with a
distinguished vertex $x$, the percolation probability
$\mu^+_{\beta}(x \stackrel{+}{\longleftrightarrow}\infty)$ is
increasing in $\beta$.
\end{thm}
An interesting aspect of this result is that its proof, unlike those of
the monotonicity results mentioned earlier in this section, is {\em not}
based on stochastic domination between the probability measures in
question. In fact, stochastic domination fails, i.e.\ it is not always
the case (in the setting of Theorem \ref{thm:tree_monotonicity}) that
\begin{equation} \label{eq:false_domination}
\mu^+_{\beta_1} \leqd \mu^+_{\beta_2}
\end{equation}
when $\beta_1 \leq \beta_2$. An easy way to see this is as follows.
Just as in Theorem \ref{thm:Ising}, let $\beta_c$ be the critical inverse
temperature for non-uniqueness of the Gibbs measure. (It is straightforward
to show, using either the
random-cluster approach or the methods in Section \ref{sect:unireg}, that
$\beta_c>0$ for $\L=\T_d$.)
Pick $\beta_1<\beta_2$ in $(0,\beta_c)$.
By Theorem \ref{thm:Ising_sandwich}, we then have
$\mu^+_{\beta_1}(\xi: \, \xi(y)=+1)=
\mu^+_{\beta_2}(\xi: \, \xi(y)=+1)= {1}/{2}$ for every vertex $y$. If
now (\ref{eq:false_domination}) was true we would have
$\mu^+_{\beta_1} = \mu^+_{\beta_2}$ by Proposition
\ref{prop:single_site+dom}.
This, however, is impossible because the two measures have different
conditional distributions on finite regions.

Theorem \ref{thm:tree_monotonicity} can be proved using the exact
calculations
for the Ising model on $\T_d$, which can be found in e.g.\ Spitzer
\cite{Sp}
and Georgii \cite{Geo}. Here we present a simpler proof which does not
require any exact calculation, but which exploits random-cluster
methods.

\medskip\noindent
{\bf Proof of Theorem \ref{thm:tree_monotonicity}:}
As usual we write $\L$ and $\B$ for the vertex and edge sets of
$\T_d$. Since
\[
\lim_{\La\uparrow\L}\mu^+_{\beta,\Lambda}
(x \stackrel{+}{\longleftrightarrow}\infty)=
\mu^+_{\beta} (x \stackrel{+}{\longleftrightarrow}\infty)
\]
for any $\beta$ in analogy to \rf{limit_interchange}, it suffices to show
that for $\beta_1\leq \beta_2$ and any $\La$, we have
\begin{equation} \label{eq:domination_on_finite_level}
\mu^+_{\beta_1,\Lambda}
(x \stackrel{+}{\longleftrightarrow}\infty)\le
\mu^+_{\beta_2,\Lambda}
(x\stackrel{+}{\longleftrightarrow}\infty) \, .
\end{equation}
This we will do by constructing a coupling $P$ of two $\{-1,+1\}^\L$-valued
random objects $X_1$ and $X_2$ with respective distributions
$\mu^+_{\beta_1,\Lambda}$ and $\mu^+_{\beta_2,\Lambda}$ and the property
that if $x \stackrel{+}{\longleftrightarrow}\infty$ in $X_1$, then the same
thing happens in $X_2$.

Recall the Edwards--Sokal coupling of spin and edge configurations
described ahead of Proposition \ref{prop:free_and_ordered_Potts_state}.
In the present case of the tree $\T_d$, this construction requires the
C-version of counting clusters, which corresponds to making the
complement of $\La$ connected. Therefore we will work with the
C-random-cluster distributions.
 Let $p_1=1-e^{-2\beta_1}$ and $p_2=1-e^{-2\beta_2}$,
and let $B=\B_\La^1\subset \B$ be the set of edges incident to at least
one vertex in
$\La$. We first let $Y_1$ and $Y_2$ be two $\{0,1\}^\B$-valued
random edge configurations distributed according to the random-cluster measures
$\hat{\phi}^{B,1}_{p_1, 2}$ and $\hat{\phi}^{B,1}_{p_2, 2}$,
and such that $P(Y_1 \preceq Y_2)=1$; this is possible by the
$\hat{\phi}^{B,1}_{p, 2}$-analogue of (\ref{eq:RC_domination}).
$X_1$ and $X_2$ can now be obtained by assigning spins to the connected
components of $Y_1$ and $Y_2$ in the usual way; these spin assignments are
coupled as follows. First we must assign spin $+1$ to all infinite
clusters in $Y_1$ and $Y_2$. Then we let $(Z(y))_{y\in \La}$ be
i.i.d.\ random variables taking values $+1$ and $-1$ with probability
${1}/{2}$ each, and assign to each finite cluster $ C$ of $Y_1$ and
$Y_2$ the value $Z(y)$, where $y$ is the (unique) vertex of $C$ that
minimizes the distance to $x$. This defined $X_1$ and $X_2$. A moment's
thought reveals that the set of vertices that can be reached from $x$ via
spins in $X_1$ is almost surely contained in the corresponding set for
$X_2$. Hence
(\ref{eq:domination_on_finite_level}) is established, and we are done.
$\Cox$

\medskip\noindent
Note that this proof did not use any property of $\T_d$ except for the tree
structure, so Theorem \ref{thm:tree_monotonicity} can immediately be
extended
to the setting of arbitrary trees.

\subsection{Cluster algorithms for computer simulation}
\label{sect:Swendsen_Wang}

An issue of great importance in statistical mechanics
which we have not touched upon so far is the ability to perform computer
simulations of large Gibbs systems. Many (most?) questions about
phase transition behavior etc.\ can with current knowledge
only be answered partially (or not at all)
using rigorous mathematical arguments. Computer simulations are then
important for supporting (or rejecting) heuristic arguments, or (in case
not even a good heuristic can be found) to provide ideas for what a good
conjecture might be. This topic is somewhat beside the main issue of
our survey, but since random-cluster representations have played a
key role in simulation algorithms for more than a decade we feel that it
is appropriate to describe some of these algorithms. In fact, it was the
need of efficient simulation which, in the late 1980's,
sparked the revival of the random-cluster model (Swendsen and Wang
\cite{SW})
which up to then had raised only little interest
since its introduction by Fortuin and Kasteleyn in the early 1970's.

Consider for instance the Ising model with free boundary condition on
a large cubic region $\Lambda \subset \Z^d$. Direct sampling from the Gibbs
distribution $\mu_{h, \beta, \Lambda}$ with free boundary condition is not
feasible, due to the
huge cardinality of the state space $\Omega$, and the (related)
intractability of computing the normalizing constant
for the Gibbs measure. The most
widely used way to handle this problem is the {\em Markov chain Monte
Carlo}
method, which dates back to the 1953 paper by Metropolis et al.\
\cite{MRRTT}.
The idea is to define an ergodic Markov chain having as unique stationary
distribution the target distribution
$\mu_{h, \beta, \Lambda}$. Starting the chain in an arbitrary
state and running the chain for long enough will then produce an output
with a distribution close to the target distribution. An example of
such a chain is the single-site heat bath algorithm, whose evolution is
as follows. At each integer time, a vertex $x \in \Lambda$
is chosen at random, and the spin at $x$ is replaced by a new value
according to the conditional distribution
(under $\mu_{h, \beta, \Lambda}$) of the spin at $x$ given the
spins at its neighbors. It is immediate that
$\mu_{h, \beta, \Lambda}$ is stationary for this chain, and
ergodicity of the chain follows from elementary Markov chain theory upon
checking that it is aperiodic and irreducible. The problem with this
approach is that the time taken to come close to equilibrium may be very
long. For example, let $h=0$. Then, for $\beta< \beta_c$ (with $\beta_c$
defined as in Theorem \ref{thm:Ising}), the time taken to come
within a fixed small variational distance from the target distribution
grows only like $n\log n$ in the size of the system (here $n$ is the number
of vertices in $\Lambda$) whereas in contrast the time grows
(stretched) exponentially in the size of the system
for $\beta>\beta_c$; see e.g.\ \cite{MaO1,M1}.
This means that simulation using this heat bath algorithm is
computationally feasible even for fairly large systems provided that
$\beta<\beta_c$, but not for $\beta>\beta_c$. What happens for
$\beta> \beta_c$ is that if the chain starts in a configuration
dominated by plus spins, then the plus spins continue to dominate for
an astronomical amount of time, and similarly for starting configurations
dominated by minus spins. The set of configurations where the fraction
of plus spins is around ${1}/{2}$ (rather than around the fractions
predicted by the magnetization in the infinite-volume Gibbs measures
$\mu_\beta^+$ and $\mu_\beta^-$)
has small probability and thus can be seen as a ``bottleneck''
in the state space, slowing down the convergence rate.

A way to tackle the exceedingly slow convergence rate in the phase
coexistence regime is to use the heat bath algorithm for the
corresponding random-cluster model rather than for the Ising model itself,
and only in the end go over to the Ising model by the random mapping
described in
Corollary \ref{cor:RC_to_Potts}. This has the disadvantage that the
calculation of single-site (or, rather, single-edge) conditional
probabilities
become computationally more complicated due to the possible dependence
on edges arbitrarily far away (see Lemma \ref{lem:single_edge}). This
disadvantage, however, seems to be by far outweighed by the fact that
the convergence rate of the Markov chain (for $\beta>\beta_c$) appears
to be very much faster than for the heat bath applied directly to the spin
variables. The reason for this phenomenon is that the random-cluster
representation ``doesn't see any difference'' between the plus state and
the minus state. This approach can
of course be used also for the $q\geq 3$ Potts model, and
is due to Sweeny \cite{Swe}. Later, Propp and Wilson
\cite{PW} built on this approach by coupling several such Markov chains
(i.e.\ running them in parallel) in an ingenious way, producing an
algorithm which runs for a random amount of time (determined by the
algorithm
itself) and then outputs a state which has {\em exactly} the target
distribution. The running time of this algorithm turns out
(from experiments) to be moderate except for
the case of large $q$ and $\beta$ close to the critical value.
The Propp--Wilson approach, known as exact or perfect simulation,
has received a vast amount of attention among statisticians during the
last few years (see e.g.\ the annotated bibliography \cite{Wil}) and
we believe that it has interesting potential also in physics.

There is, however, another Markov chain which appears to converge even
faster than those of Sweeny, Propp and Wilson. We are talking about
the Swendsen--Wang \cite{SW} algorithm, which runs as follows for
Ising and Potts models on a graph with vertex set $\L$ and edge set $\B$:
Starting with a spin configuration $X_0\in \{1, \ldots, q\}^\L$,
a bond configuration $Y_0\in \{0,1\}^\B$
is chosen according to the random mapping defined in
Corollary \ref{cor:Potts_to_RC}. Then another spin configuration $X_1$ is
produced from $Y_0$ by assigning random spins to the connected components,
i.e.\ by the random mapping of Corollary \ref{cor:RC_to_Potts}.
This procedure is then iterated, producing a new edge configuration
$Y_1$ and a new spin configuration $X_2$, etc.
By combining the two corollaries, we see that if $X_0$ is chosen according
to the target distribution, then the same holds for $X_1$, and consequently
for $X_2, X_3, \ldots$. In other words, the target distribution is
stationary for the chain $\{X_k\}_{k=0}^\infty$, and by the (easily
verified)
ergodicity of the chain we have a valid Markov chain Monte Carlo algorithm.
Although
it is not exact in the sense of the Propp--Wilson algorithm, it appears
to converge much faster, thus in practice allowing simulation of systems
that
are orders of magnitude larger. Heuristically, the reason for this faster
convergence is that
large chunks of spins may flip simultaneously, allowing the chain to
tunnel through any bottlenecks in the target distribution. However,
rigorous upper and lower bounds on the time taken to come close to
equilibrium
are to a large extent lacking, although Li and Sokal \cite{LiSo} have
provided
a lower bound demonstrating the phenomenon of ``critical slowing down''
as $\beta$ approaches $\beta_c$.

The Swendsen--Wang algorithm has, since its introduction in 1987,
become the standard approach to simulating Ising and Potts models.
Interesting variants and modifications of this algorithm have been
developed by Wolff \cite{Wol} and Machta et al.\ \cite{MCLSC}; the
last paper is an interesting attempt at combining the original
approach of Swendsen and Wang with ideas from so called {\em
invasion percolation} (see \cite{CCN}) to get an algorithm
specifically aimed at sampling from a Gibbs distribution at the
critical inverse temperature $\beta_c$, i.e.\ where the use of
other algorithms have proved to be most difficult. Generalizations
of the Swendsen--Wang algorithm for various models other than Ising
and Potts models have also been obtained, see e.g.\ Campbell and
Chayes \cite{CaC}, Chayes and Machta \cite{CM,ChM2}, and H\"aggstr\"om et
al.\ \cite{HvLM}.

\subsection{Random-cluster representation of the Widom--Rowlinson model}
\label{sect:WR_RC}

The random-cluster model can be seen as a perturbation of Bernoulli bond
percolation, where the probability measure is changed in favour of
configurations with many (for $q> 1$) or few (for $q<1$) connected
components. A fairly natural question is what happens if we perturb
Bernoulli site percolation in the same way. For lack of an established
name,
we call the resulting model the site-random-cluster model. Let $G$ be
a finite graph with vertex set $\L$ and edge set $\B$. For a site
configuration $\eta\in\{0,1\}^\L$, we write $k(\eta)$ for the number of
connected components in the subgraph of $G$ obtained by deleting all
vertices
$x$ with $\eta(x)=0$ and their incident edges.
\begin{defn} \label{defn:site-random-cluster}
The {\bf site-random-cluster measure} $\psi^G_{p,q}$ for $G$ with
parameters
$p\in[0,1]$ and $q>0$ is the probability measure on $\{0,1\}^\L$ which
to each $\eta\in\{0,1\}^\L$ assigns probability
\[
\psi^G_{p,q}(\eta)=\frac{1}{Z^G_{p,q}}\left\{\prod_{x \in \L}
p^{\eta(x)}(1-p)^{1-\eta(x)}\right\}q^{k(\eta)},
\]
where $Z^G_{p,q}$ is a normalizing constant.
\end{defn}
Analogously to the usual random-cluster model living on bonds, taking
$q=1$ gives the ordinary Bernoulli site percolation $\psi_p$, while other
choices of $q$ lead to dependence between vertices.

Taking $q=2$ is of particular interest because it leads to a representation
of the Widom--Rowlinson model which is similar to
(and in fact slightly simpler than)
the usual random-cluster representation of the Ising model. Let
$\mu_\lambda^G$ be the Gibbs measure for the Widom--Rowlinson model
with activity $\lambda$ on $G$, i.e.\ $\mu_\lambda^G$ is the probability
measure on $\{-1, 0, +1\}^\L$ which to each $\xi\in\{-1, 0, +1\}^\L$
assigns probability proportional to
\[
\prod_{\langle x y \rangle \in \B}
I_{\{\xi(x)\xi(y)\neq -1\}} \prod_{x\in \L} \lambda^{|\xi(x)|} \, .
\]
The following analogues of Corollaries \ref{cor:RC_to_Potts} and
\ref{cor:Potts_to_RC} are trivial to check.
\begin{prop} \label{prop:RC_to_WR}
Let $p=\frac{\lambda}{1+\lambda}$, and suppose we pick a random spin
configuration $X\in \{-1, 0, +1\}^\L$ as follows.
\begin{enumerate}
\item Pick $Y\in \{0, 1\}^\L$ according to $\psi^G_{p,2}$.
\item Set $X(x)=0$ for each $x\in \L$ such that $Y(x)=0$.
\item For each open cluster $C$
of $Y$, flip a fair coin to decide whether
to give spin $+1$ or $-1$ in $X$ to all vertices of $C$.
\end{enumerate}
Then $X$ is distributed according to the Widom--Rowlinson Gibbs measure
$\mu_\lambda^G$.
\end{prop}
\begin{prop} \label{prop:WR_to_RC}
Let $p=\frac{\lambda}{1+\lambda}$, and suppose we pick a random spin
configuration $Y\in \{0, 1\}^\L$ as follows.
\begin{enumerate}
\item Pick $X\in \{-1, 0, +1\}^\L$ according to $\mu_\lambda^G$.
\item Set $Y(x)=|X(x)|$ for each $x\in \L$.
\end{enumerate}
Then $Y$ is distributed according to the site-random-cluster measure
$\psi^G_{p,2}$.
\end{prop}
We remark that
for $q\in\{3,4,\ldots\}$, these results extend in the obvious way to
a connection between $\psi^G_{p,q}$ and the generalized Widom--Rowlinson
model whith $q$ types of particles rather than just $2$ (and strict
repulsion between all particles of different type).

Many of the arguments applied to Ising and Potts models in Section
\ref{sect:PT_Potts} can now be applied to the Widom--Rowlinson model in
a similar manner. To apply Theorem \ref{thm:holley}, we need to calculate
the conditional probability in the site-random-cluster model that
a given vertex is open given the status of all other vertices.
For $x \in \L$ and $\eta\in \{0,1\}^{\L \setminus \{x\}}$, we get
\begin{equation} \label{eq:single_site}
\psi^G_{p,q}(x \mbox{ is open }| \, \eta) =
\frac{p\,q^{1-\kappa(x,\eta)}}{p\,q^{1-\kappa(x,\eta)}+1 - p}
\end{equation}
where $\kappa(x,\eta)$ is the number of open clusters of $\eta$ which
intersect $x$'s neighborhood $\{y\in \L: \, y\sim x\}$.
If the degree of the vertices in $G$ is bounded by $N$, say, then
$0\leq\kappa(x,\eta)\leq N$ for any $x \in \L$ and
$\eta\in \{0,1\}^{\L \setminus \{x\}}$. For fixed $q$ and any
$p^*\in(0,1)$,
we can thus apply Theorem \ref{thm:holley} to show that
$\psi^G_{p,q}$ stochastically dominates $\psi_{p^*}$ for $p$ sufficiently
close to $1$, and is dominated by $\psi_{p^*}$ for $p$ small enough.
The arguments of Section \ref{sect:PT_Potts} leading to a proof of
Theorem \ref{thm:Ising}, with the random-cluster model replaced
by the site-random-cluster model, therefore go through to show
Theorem \ref{thm:WR}.

One thing that does {\em not} go through in this context, however, is the
analogue of (\ref{eq:RC_domination}). The reason for this is that,
in contrast to (\ref{eq:single_edge}), the conditional probability in
(\ref{eq:single_site}) fails to be increasing in $\eta$, so that
Theorem \ref{thm:holley} is not applicable for comparison between
site-random-cluster measures with different values of $p$. In fact, the
analogue of (\ref{eq:RC_domination}) for site-random-cluster measures
sometimes fails, and moreover the occurrence of phase transition
for the Widom--Rowlinson model on certain graphs fails to be increasing
in $\lambda$, as demonstrated by Brightwell, H\"aggstr\"om and
Winkler \cite{BHW}.

Another consequence of the failure of the conditional probability in
(\ref{eq:single_site}) to be increasing is that the FKG inequality
(Theorem \ref{thm:FKG}) cannot be applied to $\psi^G_{p,q}$. As a
consequence, the proof of Theorem \ref{thm:Potts_sandwich} cannot be
adapted to the case of the multitype ($q\geq 3$) Widom--Rowlinson model.
In fact, such a Widom--Rowlinson analogue of Theorem
\ref{thm:Potts_sandwich}
is known to be {\em false}, as shown by Runnels and Lebowitz \cite{RL};
see also \cite{CKoS} and \cite{NL}.

\section[Uniqueness and exponential mixing from
non-percolation]{Uniqueness
and exponential mixing\\ from non-percolation}
\label{sect:unireg}

In the previous section we saw examples where phase transition in one
system was equivalent to the existence of infinite
clusters in another, suitably defined, system. In this section we
shall discuss various approaches where conclusions about the phase
transition behavior can only be drawn from nonexistence (and not from
existence) of infinite clusters. On the other hand, these approaches
typically apply to a much wider range of models.
We address two problems: the uniqueness of the Gibbs measure,
and the decay of correlations for a given Gibbs measure.
In fact, the general theme of this section can be stated as follows:
To which extent
can a given spin be influenced by a configuration far away? If such an
influence disappears in the limit of infinite distance, it follows
(depending on the setting) that either there is no long-range
influence of boundary conditions at all (implying uniqueness of the
Gibbs measure), or that  a specific low temperature phase exhibits
some mixing properties. In both
cases, the decreasing influence comes from the absence of infinite
clusters of suitable type which could transport a dependence
between spins. So, both uniqueness and mixing will appear here as a
consequence of non-percolation.

In a first part, we will address the problem of uniqueness. In fact,
we will encounter conditions which not only imply the uniqueness
of the Gibbs measure, but also lead us into a regime where `all
good things' happen, i.e., where the unique Gibbs measure exhibits
nice exponential mixing properties and the free energy depends
analytically on all relevant parameters. (In general, the uniqueness
of the Gibbs measure does not imply the absence of other critical
phenomena, which might manifest themselves as singularities of the
free energy or other thermodynamic quantities.
For example, in Section \ref{sect:rain} we will see that in
the  so-called Griffiths' regime of a disordered system there is a
unique Gibbs measure, but the free energy is not analytic.)

 The `nice regime' above is usually referred to as the high
temperature, or weak coupling, low density, or also analytic
regime, and is usually studied by high
temperature cluster expansions. Dobrushin and Shlosman \cite{DS1,DS2}
developed a beautiful and general theory describing a regime of
`complete analyticity' by various equivalent properties. One of these
ranks at the top of a hierarchy of mixing properties.
While complete analyticity makes precise what actually the
`nice regime' is, and applies mainly to high temperatures or
large external fields, it is not limited to this case only
\cite{vEFSS}. The relationships between this and related notions
and also with dynamical
properties have been studied in many papers. Although some of these
have an explicit geometric flavor, we do not discuss them here
because of limitations of space. We rather refer to the sources
\cite{DS1,DS2,Ze1,M1,MaO} and also to the references following
condition (\ref{WMix}).

In Section \ref{sect:appli} we shall discuss an application of the
percolation method to
the low temperature regime, and see how percolation estimates for
the covariance between two distant observables, combined with
contour estimates, give rise to exponential mixing properties.

\subsection{Disagreement paths} \label{sect:disa}

Let $(\L,\sim)$ be an arbitrary locally finite graph, and suppose
we are given a neighbor interaction $U:S\times S\to\R$ and a
self-potential $V:S\to\R$.
Consider the associated Gibbs distributions $\mu_{\b ,\La}^{\e}$
introduced in \rf{eq:FvolGibbs}. More generally, we could consider an
arbitrary Markov specification $(G_\La)_{\La\in\E}$ in the sense of
Section \ref{sect:Gibbsmeasures}. Such specifications appear, in
particular, if we have
an  interaction of finite range $R$, say on $\Z^d$, and draw edges
between all sites of distance at most $R$.
However, for definiteness and simplicity we stick to the setting
described by the Hamiltonian \rf{eq:Hamiltonian}. We will often
consider the inverse temperature $\b$ as fixed and then simply
write $\mu_{\La}^{\e}$ instead
of $\mu_{\b,\La}^{\e}$. If $\La$ is a singleton, we use the
shorthand $x$ for $\{x\}$.

We look for a condition implying that there is only one Gibbs measure
$\mu$ for the Hamiltonian  \rf{eq:Hamiltonian}, i.e., a unique
probability measure on $\Om=S^\L$ satisfying
\[
 \mu(\,\cdot\,|\,X \equiv\e \mbox{
off } \La)=\mu_{\La}^{\e}\quad\mbox{ for $\mu$-almost all }\e\in
\Om\;.
\]
Since this property needs only to be checked for singletons
$\La=\{x\}$ (cf.\ Theorem 1.33 of \cite{Geo}), it is sufficient to
look for conditions on
the single-spin Gibbs distributions $\mu_{x}^{\e}$ with $x\in\L$.
Intuitively, we want to express that $\mu_{x}^{\e}(X(x)=a)$ depends only
weakly on $\e$ (which can be expected to hold for small $\b$).
This dependence can be measured by the maximal variation
\begin{equation} \label{dens}
p_x=\max_{\e,\e'\in \Om}
\|\mu_{x}^{\e} -\mu_{x}^{\e'} \|_x\;,
\end{equation}
where
\begin{equation}\label{TVnorm}
\|\nu\|_\De=\sup_{A\in\F_\De}|\nu(A)|
\end{equation}
is the total variation norm on the sub-$\s$-algebra $\F_\De$ of events
which depend only on the spins in $\De$. We write ${\bf p}$ as a
shorthand for the family $(p_x)_{x\in\L}$.

Given two configurations $\xi,\; \xi' \in \Om$, a path in $\L$ will be
called a {\it path of
disagreement} (for $\xi$ and $\xi'$) if $\xi(x) \neq \xi'(x)$ for all
its vertices $x$. For each finite region $\La\subset\L$ and  any two
configurations $\e,\,\e'$ on $\La^c$ we will construct a coupling
$P$ of $\mu_{\La}^{\e} $ and $\mu_{\La}^{\e'} $ describing the
difference of these measures in terms of paths of disagreement
running from the boundary $\partial\La$ into the interior of $\La$.
Intuitively, these paths of
disagreement then show how deep inside the influence of the
boundary conditions can still be felt.
We write $\{\De\stackrel{\ne}{\longleftrightarrow}\partial\La\}$ for
the event in $S^\La\times S^\La$ that there exists a path of
disagreement from some
point of a set $\De\subset\La$ to some point of $\partial\La$.

Although the coupling $P$ to be constructed is not best suited for
direct use, it has a useful special feature:
its disagreement distribution is stochastically dominated by a
Bernoulli measure. This will allow us to conclude that
absence of percolation for the latter implies uniqueness of the Gibbs
measure for the Hamiltonian \rf{eq:Hamiltonian}.

We write $\psi_{\bf p}$ for the Bernoulli measure
on $\{0,1\}^\L$ with $\psi_{\bf p}(X(x)=1)=p_x$ for all $x\in\L$, and
$\psi_{{\bf p},\La}$ for the analogous product measure on
$\{0,1\}^\La$.
As in Section \ref{sect:coupling_SD}, we use the notation $X(x)$ and
$X'(x)$ for the projections from $\Om\times\Om$ to $S$.
The following theorem is due to van den Berg and Maes \cite{vdBM}.
\begin{thm} \label{thm:dis}
For each finite $\La \subset \cal L$ and each pair $\e, \,\e' \in \Om$
there exists a coupling  $P=P_{\La,\e,\e'}$ of $\mu_{\La}^{\e} $ and
$\mu_{\La}^{\e'} $  having the following properties:
\begin{description}
\item{\rm (i)}
For each $x\in \La$,
$\{X(x) \neq X'(x)\} =
\{x\stackrel{\ne}{\longleftrightarrow}\partial\La\}$
 $P$-a.s.
\item{\rm (ii)}
For the distribution $P_\La^{\ne}$ of $(I_{\{(X(x) \neq X'(x)\}})_{
x\in\La}$ under $P$, we have
$
P_\La^{\ne}\leqd \psi_{{\bf p},\La}\;.
$
\item{\rm (iii)} For each $\De\subset\La$,
\begin{equation} \label{corr}
\|\mu_\La^{\e} - \mu_\La^{\e'}\|_\De\leq
P(\De\stackrel{\ne}{\longleftrightarrow}\partial\La)\le
\psi_{\bf p}(\De\,{\leftrightarrow}\,\partial\La)\;.
\end{equation}
\end{description}
\end{thm}
{\bf Proof:} We construct a coupling $(X,X')$ of
$\mu_{\La}^{\e}$ and $\mu_{\La}^{\e'}$ by the following algorithm.
In a preparatory step we introduce an arbitrary linear ordering on
$\La$ , set $\De=\La$, and define $X(x)=\e(x),
X'(x)=\e'(x)$ for $x\in \De^c$.

For fixing the main iteration step, suppose that $(X,X')$ is already
defined on the complement of a non-empty set
$\De\subset\La$ and is realized as a pair $(\xi,\xi')$ off $\De$,
where $(\xi,\xi')\equiv(\e,\e')$ off $\La$.
Conditional on the event that $(X,X')\equiv(\xi,\xi')$ off $\De$,
we consider the Gibbs distributions $\mu_\De^\xi$
 and $\mu_\De^{\xi'}$ obtained by conditioning $\mu_{\La}^{\e}$ and
$\mu_{\La}^{\e'}$ on $X\equiv\xi$ resp.\ $\xi'$ off $\De$, and
we pick the smallest vertex $x=x(\xi,\xi')\in \De$ for which
there exists
some vertex $y\in \De^c$ with $y\sim x$ and $\xi(y)\neq \xi'(y)$.
If such an $x$ does not exist, we have $\mu_{\De}^{\xi}
=\mu_{\De}^{\xi'}$ on $\F_\De$ by the Markov property,
so that we can take the obvious optimal coupling for which $X\equiv
X'$ on $\De$, and we are done. If such an $x$ does exist, we
consider the single vertex distributions
$\mu_{\De,x}^{\xi} = \mu_{\De}^{\xi}(X(x)=\cdot\,)$ and
$\mu_{\De,x}^{\xi'} = \mu_{\De}^{\xi'}(X(x)=\cdot\,)$ on $S$.
Conditionally on $(X,X')\equiv(\xi,\xi')$ off $\De$, we then let
$(X(x),X'(x))$ be distributed according to
an optimal coupling (as in Definition \ref{def:opt}) of
$\mu_{\De,x}^{\xi}$ and $\mu_{\De,x}^{\xi'}$. The coupling $(X,X')$
 is then defined on the set
$x\cup \De^c$, so that we can replace $\De$ by
$\De\setminus x$ and repeat the preceding iteration step.

It is clear that the algorithm above stops after finitely many
iterations and gives us a coupling of $\mu_{\La}^{\e}$ and
$\mu_{\La}^{\e'}$. Property (i) is evident from the construction,
since disagreement at a vertex is only possible if a path of disagreement
leads from this vertex to the boundary.
For (ii), we note that the measures $\mu_{\De,x}^{\xi}$ and
$\mu_{\De,x}^{\xi'}$ are mixtures of the Gibbs distributions
$\mu_x^\s$
with suitable boundary conditions $\s$, by the consistency of Gibbs
distributions. Hence
\[
\|\mu_{\De,x}^{\xi}-\mu_{\De,x}^{\xi'}\|_x\le p_x\;.
\]
By construction, this means that in each iteration of the main step
we have
\[
P(X(x)\ne X'(x)\;|\;(X,X')\equiv (\xi,\xi')\mbox{ off }\De)
\le p_x
\]
for $x= x(\xi,\xi')$, so that (ii) follows by induction. Finally,
(iii) follows
directly from (i) and (ii) because for each $\De\subset\La$
\[
\|\mu_\La^{\e} - \mu_\La^{\e'}\|_\De
\leq {P}(X(x) \neq X'(x) \mbox{ for some }x\in\De)
\]
by the coupling inequality \rf{eq:coupling_inequality}. The proof is
therefore
complete. $\Cox$

\medskip\noindent
Although the algorithm in the proof above is quite
explicit, it is  not easy to deal with directly. In particular, it is
not clear in which way the coupling depends on the
chosen ordering, because the site $x$ to be selected in each step
depends on $(\xi,\xi')$ and is therefore random. Nevertheless, if the
Gibbs distributions are monotone (in the sense of Definition
\ref{def:mon}), we get some extra properties.

\medskip\noindent{\bf Remark:}
Suppose $S$ is linearly ordered and the conditional distributions
$\mu_x^\xi$ are stochastically increasing in $\xi$. Then, if
$\e \preceq \e'$, the coupling $P$ of Theorem \ref{thm:dis}
can be chosen in such a way that, in addition to properties
(i) to (iii),  $X \preceq X'$ ${P}$--a.s.\ and, for each $x \in \La$,
\begin{equation}\label{monis}
\|\mu_\La^{\e} - \mu_\La^{\e'}\|_x\leq
P(x\stackrel{\ne}{\longleftrightarrow}\partial\La)\le
(|S| - 1) \,\|\mu_\La^{\e} - \mu_\La^{\e'}\|_x\;.
\end{equation}
This is because in each step of the algorithm proving Theorem
\ref{thm:dis} we can achieve that $X(x)\le X'(x)$, and for the second
inequality in \rf{monis} it is sufficient to note that
\[
P(x\stackrel{\ne}{\longleftrightarrow}\partial\La)=P(X(x)<X'(x))\le
\sum_{a\in S\setminus\{m\}}[P(X(x)\le a)-P(X'(x)\le a)]\;,
\]
where $m$ is the maximal element of $S$.
 For details we refer to \cite{vdBM}. In particular, for $|S| = 2$ we
have
equality in (\ref{monis}).

\medskip\noindent
Let us apply this remark to the ferromagnetic Ising model with external
field $h=0$ and any inverse temperature $\b$, with boundary
conditions $\e\equiv+1$
and $\e'\equiv-1$ outside of some finite region $\La\in \cal E$.
Then, by the spin flip symmetry and stochastic monotonicity,
\[
\mu_{\b,\La}^+(X(x))=\mu_{\b,\La}^+(X(x)=1)-\mu_{\b,\La}^-(X(x)=1)=
\|\mu_{\b,\La}^{+} - \mu_{\b,\La}^{-}\|_x
\]
and therefore, by (\ref{monis}),
\[
\mu_{\b,\La}^+(X(x)) =
P(x\stackrel{\ne}{\longleftrightarrow}\partial\La)\;.
\]
We emphasize that this relation is completely similar to what we
obtained for the random-cluster representation, viz.
\[
\mu_{\b,\La}^+(X(x)) =
\phi^1_{p,2,\La}(x\,{\leftrightarrow}\,\partial\La)
\]
for $p=1-e^{-2\b}$; cf.\ equation
\rf{eq:Potts_magnetization_connectivity}.
The coupling $P$, however, is less explicit, and
the geometric event involves site percolation rather than bond
percolation as for the  random--cluster measure, but the exact
correspondence
between the magnetization for the spin system and the percolation
probability of the geometric system is the same.

Let us now turn to the main result of this subsection, the uniqueness
theorem. Let $\mu,\,\mu'$ be any two Gibbs measures for the
Hamiltonian \rf{eq:Hamiltonian} at some inverse temperature $\b$.
Inequality \rf{corr} then shows  that
\[
\|\mu - \mu'\|_\De\leq\sup_{\e,\e'\in\Om} \|\mu_\La^{\e} -
\mu_\La^{\e'}\|_\De
\leq \psi_{\bf p}(\De\,{\leftrightarrow}\,\partial\La)
\]
whenever $\De\subset\La\in\E$. Letting $\La\uparrow\L$ we find
\[
\|\mu - \mu'\|_\De\leq \psi_{\bf p}(\De\,{\leftrightarrow}\,\infty)
\]
which gives the following uniqueness result.

\begin{thm}\label{thm:uni}
If $\, \psi_{\bf p}(\,\exists\mbox{ \em an infinite open
cluster })=0\, $ then the set ${\cal G}(\b H)$ of Gibbs measures for
the Hamiltonian \rf{eq:Hamiltonian} at inverse temperature $\b$
is a singleton.  In particular, this holds if
$\sup_x p_x < p_c$, the critical density for Bernoulli site
percolation on $(\cal L,\sim)$.
\end{thm}

\noindent
A weaker version of Theorem \ref{thm:uni} was obtained first
in \cite{vdB} using a product coupling instead of Theorem
\ref{thm:dis}; see Proposition \ref{apr} below and also \cite{vdBS}.
In some cases, the simple product coupling nevertheless gives
equivalent results; cf.\ the discussion in \cite{vdBM}.

For a large class of regular graphs such as $\Z^d$, the assumption
of Theorem \ref{thm:uni} not only implies
the uniqueness of the Gibbs measure but even yields certain exponential
mixing properties. This can
be seen almost immediately by combining inequality \rf{corr}
with Theorem \ref{thm:exponential_tail} on the exponential tail of the
distribution of the cluster diameter in sub-critical
 Bernoulli percolation. We will use similar arguments in
Section \ref{sect:G-regime} in the context of random interactions.

Let us discuss now some special cases. Clearly, the conditions of Theorem
\ref{thm:uni} hold when $\L=\Z$ with the usual graph structure, since
then $p_c=1$. This gives uniqueness of the Gibbs measure for
one-dimensional
nearest-neighbor systems. Next we consider the case $\L=\Z^d$, $d\ge2$.
Recall the bound \rf{erma} for the percolation threshold $p_c$ when
$d=2$,
and the large--dimensions asymptotics of $p_c$ in
\rf{eq:perc_high_dimensions}.
\begin{example}\label{Ising_uniq} The Ising ferromagnet.
{\rm Let $\b>0$ be any inverse temperature and $h$ an external field.
Then, for any $x$, we obtain from \rf{eq:sandwich} by a short computation
\[
p_x=
\|\mu_{h,\b,x}^{+} -\mu_{h,\b,x}^{-} \|_x=[ \tanh(\b(h+2d))
-\tanh(\b(h-2d))]/2\,.
\]
Hence, the Gibbs measure is unique when $h=0$ and $\tanh (2d\b)<p_c$,
or if $|h|>2d$ is so large that $2d<p_c\,\cosh^2(\b(|h|-2d))$, for
example.}
\end{example}
\begin{example} The hard-core lattice gas.
{\rm Setting $\b=1$, we see that $p_x=\la/(1+\la)$ for any $x$, so
that uniqueness of the Gibbs measure follows for $\lambda <
p_c/(1-p_c)$.
(This can also be obtained by using the product coupling mentioned
above, cf.\ \cite{vdBS}.)
}
\end{example}
\begin{example} The Widom--Rowlinson lattice gas.
{\rm We take again $\b=1$ and set $\la_+=\la_-=\la$. It turns out that
the maximum in equation \rf{dens} is attained for the boundary
conditions $\e\equiv0$ and $\e'$ equal to $+1$ and $-1$ on (at least)
two different neighbors
of $x$, whence $p_x=2\la/(1+2\la)$ for any $x$. It follows that the
Gibbs measure is unique when $\lambda < p_c/(2(1-p_c))$.}
\end{example}

\noindent It is interesting to compare the uniqueness condition of
Theorem \ref{thm:uni} with the celebrated Dobrushin uniqueness
condition, cf.\ \cite{Geo} and the original papers \cite{D2,D72}.
This condition reads
\begin{equation} \label{dssc}
\sup_x \sum_y \max_{\eta\equiv\eta'\mbox{ \scriptsize off }y}
\|\mu_{x}^{\e} -\mu_{x}^{\e'} \|_x<1.
\end{equation}
The constraint ``$\eta\equiv\eta'\mbox{ off }y$'' means that the
configurations
$\e,\e'$ differ only at the vertex
$y$. For systems with hard-core
exclusion or in certain antiferromagnetic models it often happens that,
for every $y\in \partial x$,  the maximum in \rf{dssc} is actually the
same as that in
\rf{dens}, see \cite{vdBM}. Dobrushin's uniqueness  condition then
takes the form  $ \sup_x |\partial x|\, p_x < 1$. For ${\cal L} =
\Z^d$ and $p_x=p_0$ independently of $x$, this means that
$p_0<1/(2d)$, while
Theorem \ref{thm:uni} only requires $p_0<p_c$, and it is known that
$p_c > 1/(2d-1)$ for $d>1$. However, if the constrained maximum in
Dobrushin's condition is much smaller than the unconstrained maximum
in \rf{dens}, then Dobrushin's condition will
be weaker than that of Theorem \ref{thm:uni}. For example, for the
Ising ferromagnet on $\Z^d$ with external field $h=0$, Dobrushin's
condition requires that $2d\tanh\b<1$ which, in view of
\rf{eq:perc_high_dimensions}, is less restrictive than the condition
obtained in Example \ref{Ising_uniq}. Thus, roughly speaking, Theorem
\ref{thm:uni} works best for ``constrained''
 systems with strong repulsive interactions and
low-dimensional lattices (or graphs with small $|\partial x|$'s) for
which reasonable lower bounds of the critical probability
$p_c$ are available.  Examples are the hard-core lattice gas and the
Widom--Rowlinson lattice gas on $\Z^2$ considered  above.

There is also another reason why Theorem \ref{thm:uni} is useful.
Namely, its condition of non-percolation is a global
condition: the absence of percolation does not depend on the value of $p_x$
at any single site $x$. In particular, $p_x$ could be large or even be
equal to 1 for all $x$'s in an infinite subset (say, a
periodic sublattice) of $\L$; once the $p_x$'s are sufficiently small
on the complementary set, there is still no infinite open cluster.
This can be applied to non-translation invariant interactions
where, in general, it is impossible to obtain uniform small
bounds on the $p_x$'s (or on the strength of the
interaction, as would be required by the Dobrushin condition or for
some standard cluster-expansion argument).  We will come back to
this point in Section \ref{sect:appli}.

\subsection{Stochastic domination by random-cluster measures}
\label{ssection:stodom}
Recently, Alexander and Chayes \cite{AC} introduced a variant of the
random-cluster technique that applies to a substantially greater
class of systems than those considered in Section
\ref{sect:random-cluster}. This approach involves a so-called
{\em graphical representation} of the original system. The graphical
representation is stochastically dominated by a random-cluster model,
and absence of infinite clusters in this random-cluster
model implies uniqueness of the Gibbs measure for the original spin
system. The price
to pay for the greater generality is that the implication goes only
one way: percolation in the random-cluster model does not, in general,
imply non-uniqueness of Gibbs measures.

We assume that the state space $S$ is a {\em finite group}
with unit element $1$; the inverse element of $a\in S$ is denoted by
$a^{-1}$, so that $a^{-1}a=aa^{-1}=1$. For simplicity we assume
that the underlying graph is $\L=\Z^d$ (although this will not really
matter). We consider the Hamiltonian
(\ref{eq:Hamiltonian}) for a  pair potential
$U$ and with no self-energy, $V=0$. By adding some constant to $U$
(which does not change the relative Hamiltonian) we can arrange that
$U\leq 0$. The basic assumption is that $U$ is {\em left-invariant}, so
that
\begin{equation}\label{pot_homogeneous}
U(a,b) = u(a^{-1}b)
\end{equation}
for all $a,b\in S$ and the even function $u=U(1,\cdot)\le0$.
Note that this setting includes the $q$-state Potts model for which
$S=\Z_{q}$ and $u=-2\,I_{\{0\}}$.
For any finite $\La \subset \Z^d$ we consider the
Gibbs distribution
\[
\mu_{\b,\La}^\eta(\s) =
\frac{I_{\{\s\equiv\eta \mbox{ \scriptsize off  } \La\}}
}{Z_\La(\b,\eta)}\;\exp\bigg[-\b \sum_{\langle xy \rangle\in\B_{\La}}
u(\s(x)^{-1}\s(y))\bigg]
\]
at inverse temperature $\b$ with boundary condition $\eta \in \Om$.
Here we write ${\cal B}_\La$ for the set of all bonds $b\in \B$ with
at least one endpoint in $\La$. The graphical representation of
$\mu_{\b,\La}^\eta$ will be based on bond configurations $\om \in
\{0,1\}^{\B_{\La}}$. Each such $\om$ will also be viewed as a subset of
${\cal B}_\La$, and the bonds in $\om$ will be called open. The key
idea of this representation is taken from the classical high
temperature expansion. For fixed $\b>0$ and any $a\in S$ we
introduce the difference
\begin{equation}\label{def_R}
R_{a} = e^{-\b\, u(a)} - 1\ge 0\,.
\end{equation}
With this notation we can write
\begin{eqnarray*}
\mu_{\b,\La}^\eta(\s) &=&
\frac{I_{\{\s\equiv\eta \mbox{ \scriptsize off  } \La\}}
}{Z_\La(\b,\eta)}\prod_{\langle xy\rangle\in\B_{\La}}
(1+R_{\s(x)^{-1}\s(y)})\\
&=&\frac{I_{\{\s\equiv\eta \mbox{ \scriptsize off  }
\La\}}}{Z_\La(\b,\eta)}\sum_{\om\in\{0,1\}^{\B_{\La}}}
\prod_{\langle xy\rangle\in\om}R_{\s(x)^{-1}\s(y)}\;.
\end{eqnarray*}
This shows that $\mu_{\b,\La}^\eta$ is the first marginal
distribution of a probability measure $P_{\b,\La}^\eta$ on $\Om\times
\{0,1\}^{\B_{\La}}$, namely
\[
P_{\b,\La}^\eta(\s,\om)=
\frac{I_{\{\s\equiv\eta \mbox{ \scriptsize off  } \La\}}}{Z_\La(\b,\eta)}
\prod_{\langle xy\rangle\in\om}R_{\s(x)^{-1}\s(y)}\;,
\]
$\s\in\Om$, $\om\in\{0,1\}^{\B_{\La}}$. The second marginal distribution
of $P_{\b,\La}^\eta$ is equal to
\[
\g_{\b,\La}^\eta(\om) =
W_{\b,\La}^\eta(\om)\bigg/ Z_\La(\b,\eta)\;,
\]
where
\begin{equation} \label{greyw}
W_{\b,\La}^\eta(\om) = \sum_{\s\equiv\e\mbox{ \scriptsize off }\La}
 \prod_{\langle
xy\rangle \in \om} R_{\s(x)^{-1}\s(y)}
\end{equation}
is the ``graphical weight'' of any $\om\in\{0,1\}^{{\cal B}_\La}$.
The probability measure $\g_{\b,\La}^\eta$ on $\{0,1\}^{\B_{\La}}$ is
called the {\em graphical distribution\/} or the {\em grey measure\/}
(since it ignores the spins which are considered as colors). The
graphical
representation of $\mu_{\b,\La}^\eta$
thus obtained is analoguous to the random-cluster
representation of the Potts model and can be summarized as follows.
\begin{lem}\label{rcdomlem1}
In the set-up described above, the Gibbs distribution $\mu_{\b,\La}^\eta$
can be derived from the graphical distribution $\g_{\b,\La}^\eta$
by means of the conditional probabilities
\[
P_{\b,\La}^\eta(\s|\om)={W_{\b,\La}^\eta(\om)}^{-1}{\prod_{b=\langle
xy\rangle\in \om} R_{\s(x)^{-1}\s(y)}}\;.
\]
That is,
\[
\mu_{\b,\La}^\eta(\s)=\sum_{\om\in\{0,1\}^{\B_{\La}}}
\g_{\b,\La}^\eta(\om)
\,P_{\b,\La}^\eta(\s|\om)\;.
\]
\end{lem}
For the Potts interaction $u=-2\,I_{\{0\}}$ with state space $S=\Z_{q}$,
the graphical representation above is easily seen to coincide
with the random-cluster representation studied in Section
\ref{sect:random-cluster}.
One important feature is that the graphical weights factorize
into cluster terms. Indeed, each bond configuration $\om$ divides $\La$
into connected components called open clusters (which may possibly
consist of isolated sites).
The set of bonds belonging to an open cluster $C$ is denoted by
$\om_{C}$. Writing ${\cal C}(\om)$ for the set of all open clusters we
then obtain that
\begin{equation}\label{factorize}
W_{\b,\La}^\eta(\om)= \prod_{C\in{\cal C}(\om)} \bar
W_{\b,\La}^\eta(C,\om_{C})
\end{equation}
with
\[
\bar W_{\b,\La}^\eta(C,\om_{C}) =\sum_{\s\in
S^{C}:\,\s\equiv\eta\mbox{ \scriptsize on }C\cap\partial\La}\ \
\prod_{\langle xy\rangle \in\om_{C}} R_{\s(x)^{-1}\s(y)}\;.
\]
(We make the usual convention that the empty product is equal to 1;
hence $\bar W_{\b,\La}^\eta(C,\om_{C}) \\=|S|$ if $C$ is an isolated
site.) Together with Lemma \ref{rcdomlem1},
equation \rf {factorize} shows that the spins belonging to disjoint open
clusters are conditionally independent.
In particular, we can simulate the spin system by first drawing a bond
configuration $\om$ with weights (\ref{greyw}) and then obtain in each
open cluster a spin configuration according to
$P_{\b,\La}^\eta(\s|\om)$.

Suppose we knew that there is no percolation in the graphical
representation,
in the sense that $\max_{\e}\g_{\b,\La}^\eta(0 \leftrightarrow
\partial\La)\to 0$ as $\La\uparrow\L$. The conditional independence of
spins in different open clusters would then suggest that there
is only one Gibbs measure for the spin system. Unfortunately,
this is not known (though weaker statements are established in
\cite{CM}). However, one can make a stochastic comparison of the
graphical distributions with wired random-cluster distributions
(Lemma \ref{rcdoml2} below), and the absence of percolation in the
dominating random-cluster distribution will then guarantee that the
original system has a unique Gibbs measure. This will be achieved
in Theorem \ref{AC} allowing to bound the dependence on
boundary conditions in terms of the connectivity probability in a
random-cluster model.

To this end we also need to consider Gibbs distributions $\mu_{\b,\La}^f$
with {\em free boundary condition}. These admit similar
graphical representations  $\g_{\b,\La}^f$ based on bond configurations
inside $\La$; that is, the bonds leading from $\La$ to $\La^c$ are
removed. In the
following, the superscript $f$ will refer to this case.

The stochastic comparison with random-cluster distributions will be
formulated using
\begin{equation}\label{R_param}
R^*=\max_{a\in S} R_a\;,\quad \bar R
=\frac 1{|S|}\sum_{a\in S} R_a\;,\quad p=R^*/(1+R^*)\;,\quad q= R^*/\bar
R\;.
\end{equation}
Note that these quantities depend on $\b$ since the $R_a$ in \rf{def_R}
do.
In the case of the $r$-state Potts model when $u=-2\,I_{\{0\}}$, we have
$R^*=1-e^{-2\b}$  and $q=r$; that is, in this case the parameters $p$ and
$q$ are nothing but the standard parameters of the random-cluster
representation. For  $p$ and $q$ as above we consider now the wired
(resp.\ free) random-cluster
distribution $\phi^1_{p,q,\La}$ (resp.\ $\phi^0_{p,q,\La}$) in $\La$ as
introduced in Section \ref{sect:Inf_vol_limits}.
\begin{lem}\label{rcdoml2} For any $\La\in\E$, $\b>0$ and $p,q$ as above,
$\g_{\b,\La}^\eta \leqd \phi^1_{p,q,\La}$ and $\g_{\b,\La}^f \leqd
\phi^0_{p,q,\La}$.
\end{lem}
{\bf Proof: } We only prove the first statement since the second is
similar
and simpler. According to Section \ref{defn:random-cluster},
the weights of the random-cluster distribution
$\phi^1_{p,q,\La}$ are proportional to
\[
\left( \frac{p}{1-p}\right)^{|\om|} q^{k(\om,\La)}
\]
with $|\om|$ the number of open bonds and
$k(\om,\La)$ the number of open clusters meeting $\La$ (where all
clusters touching $\partial\La$ are wired together into a single cluster).
Up to a constant factor, the Radon--Nikodym
density of $\g_{\b,\La}^\eta$ relative to $\phi^1_{p,q,\La}$ is thus
given by
\[
F(\om)={W_{\b,\La}^\eta(\om)}\bigg/{({R^{*})}^{|\om|}
({R^*}/{\bar R})^{k(\om,\La)}}\;.
\]
Since $\phi^1_{p,q,\La}$ has positive correlations,
the lemma will therefore be proved once we have shown that
$F$ is a decreasing function of $\om$.
To this end we let $\om\preceq\om'$ be such that
$\om' = \om \cup \{b\}$ for a bond $b\in {\cal B}_\La
\setminus \om$.

We first consider the case when $b=\langle xy\rangle$ is not connected
to $\partial\La$ and joins two open clusters $C_x,\,C_{y}\in{\cal
C}(\om)$.
For each open cluster $C$ let
$\bar W(C) =\bar W_{\b,\La}^\eta(C,\om_{C})$ be as in \rf{factorize}.
Suppose we stipulate that the spin $\s(z)$ at any site $z\in C$ is
equal to some $a\in S$. It is then easy to see that the remaining sum in
the definition of $\bar W(C)$ does not depend on $a$ and thus has the
value
$\bar W(C)/|S|$. Prescribing the values of $\s(x)$
and $\s(y)$ in this way we thus find that
\[
\bar W(C_x\cup C_y\cup b) = \bar W(C_x)\, \bar W(C_y)\;
|S|^{-2} \sum_{\s(x),\,\s(y)} R_{\s(x)^{-1}\s(y)},
\]
and therefore $W_{\b,\La}^\eta(\om') = \bar R\,W_{\b,\La}^\eta(\om)$.
Since $k(\om',\La) = k(\om,\La)-1$ and $|\om'| = |\om| + 1$, it follows
that $F(\om') = F(\om)$, proving the claim in the first case. If $b$
links some cluster to the boundary which otherwise was separated from
the boundary, then the argument above shows again that
$F(\om') = F(\om)$.

Next we consider the case when $b=\langle xy\rangle$ closes a loop in
$\om$ but is still not connected to the boundary. Since clearly
\[
R_{\s(x)^{-1}\s(y)} \leq \max_{a\in S} R_a=R^*\;,
\]
we find that $W_{\b,\La}^\eta(\om')\leq
W_{\b,\La}^\eta(\om)\, R^*$.  On the other hand, in this case we have
$|\om'|=|\om| + 1$ and $k(\om',\La) = k(\om,\La)$, so that $F(\om')
\le F(\om)$.  As we are considering the wired random-cluster measure,
this argument remains valid if $b$ joins two
clusters already attached to the boundary.
$\Cox$

\medskip\noindent
We are now in a position to state the main result of Alexander and Chayes
\cite{AC}, an estimate on the dependence of Gibbs distributions on their
boundary condition in terms of percolation in the wired random-cluster
distribution. Recall the
notation \rf{TVnorm} for the total variation norm on the
sub-$\s$-algebra $\F_\De$ of events in some $\De$.
\begin{thm}\label{AC} Consider the spin system with pair interaction
\rf{pot_homogeneous} at some inverse temperature $\b>0$, and let $p,\,q$
be given by \rf{R_param}. Then,
for any $\De \subset \La\in\E$ and any pair of boundary conditions
$\eta,\eta'\in\Om$,
\[
\|\mu_{\b,\La}^\eta -
 \mu_{\b,\La}^{\eta'}\|_\De \leq \phi^1_{p,q,\La}(\De
 \leftrightarrow \partial \La)\,.
\]
\end{thm}
{\bf Proof:} (This proof is different from the one that appeared in
\cite{AC}.)  Let $A$ be any event in $\F_\De$. From Lemma
\ref{rcdomlem1} we know that
\[
\mu_{\b,\La}^\eta(A)=
\sum_\om \g_{\b,\La}^\eta(\om)\, P_{\b,\La}^\eta(A|\om)\,.
\]
To control the $\e$-dependence of this probability we will
proceed in analogy to the argument for the implication
(ii) $\Rightarrow$ (iii) of Theorem \ref{thm:Potts_sandwich}.
If $\om\in\{\De\not\leftrightarrow\partial\La\}$ then
equation \rf{factorize} shows that the conditional
distribution $P_{\b,\La}^\eta(A|\om)$ does not depend on $\eta$.
So we need to control the $\eta$-dependence of
$\g_{\b,\La}^\eta(\De\not\leftrightarrow\partial\La)$. This, however,
does not seem possible directly. So we will replace
$\g_{\b,\La}^\eta$ by the $\e$-independent $\phi^1_{p,q,\La}$
by using a suitable coupling trick.

By Lemma \ref{rcdoml2} and Strassen's theorem (Theorem \ref{thm:strassen})
there exists a coupling $(\tilde Y,\tilde Y')$ of $\g_{\b,\La}^\eta$ and
$\phi^1_{p,q,\La}$ such that  $\tilde Y\preceq\tilde Y'$ almost surely.
If $\tilde Y'\in\{\De\not\leftrightarrow\partial\La\}$, there exists
a largest (random) set $\Gamma=\Gamma(\tilde Y')$ such that
\begin{description}
\item{(a)} $\De\subset\Gamma\subset\La$, and
\item{(b)} $\tilde Y'(b)=0$ for all bonds connecting $\Gamma$ with
$\Gamma^c$.
\end{description}
For $\tilde Y'\in\{\De\leftrightarrow\partial\La\}$ we set
$\Gamma=\emptyset$. Conditional on $\Gamma$,
Lemma \ref{rcdoml2} and Strassen's theorem provide us further with a
coupling
$(\tilde Y_\Gamma,\tilde Y_\Gamma')$ of $\g_{\b,\Gamma}^f$ and
$\phi^0_{p,q,\Gamma}$ such that  $\tilde Y_\Gamma\preceq\tilde Y_\Gamma'$.
It is then easy to see that the pair of random variables $(Y,Y')$ defined
by
\[
(Y,Y')(b)=\left\{\begin{array}{cl}
(\tilde Y_\Gamma,\tilde Y'_\Gamma)(b)&\mbox{if $b$ is contained in
$\Gamma$,}\\
(\tilde Y,\tilde Y')(b)&\mbox{otherwise}
\end{array}\right.
\]
is still a coupling of  $\g_{\b,\La}^\eta$ and
$\phi^1_{p,q,\La}$ such that  $Y\preceq Y'$ almost surely.
(Notice that also $\tilde Y(b)=0$ for all bonds from $\Gamma$ to
$\Gamma^c$.) We denote the
underlying probability measure by $Q^\e$. Now we can write
\begin{eqnarray*}
\mu_{\b,\La}^\eta(A)&=&Q^\e\bigg( P_{\b,\La}^\eta(A|Y)\bigg)\\
&=&Q^\e\bigg( P_{\b,\La}^\eta(A|Y)\,I_{\{\Gamma=\emptyset\}}\bigg)
+Q^\e\bigg( P_{\b,\La}^\eta(A|Y)\,I_{\{\Gamma\ne\emptyset\}}\bigg)\;.
\end{eqnarray*}
The first term in the last sum is at most
$Q^\e(\Gamma=\emptyset)=\phi^1_{p,q,\La}(\De\leftrightarrow\partial\La)$.
We claim that the second term
does not depend on $\e$. Indeed, if $\Gamma\ne\emptyset$ then,
by \rf{factorize}, $P_{\b,\La}^\eta(A|Y)= P_{\b,\Gamma}^f(A|Y_\Gamma)$
only depends on the restriction $Y_\Gamma$ of $Y$ to the set
of bonds inside $\Gamma$. The second term can thus be written explicitly
as
\[
\sum_{G\ne\emptyset}\phi^1_{p,q,\La}(\Gamma=G)\sum_{\om\mbox{ \scriptsize
in } G}
\g_{\b,G}^f(\om)\;
P_{\b,G}^f(A|\om)\;,
\]
which is obviously independent of $\e$. The theorem now follows
immediately.
$\Cox$

\medskip\noindent
To apply the theorem we consider the limiting random-cluster
measure $\phi^1_{p,q}$ with arbitrary parameters $p\in ]0,1[$ and
$q\geq 1$ and wired boundary condition; recall from Section
\ref{sect:Inf_vol_limits} that this limiting measure exists.
By Corollary \ref{cor:FKG_for_RC}(d), it makes sense to define
the percolation threshold
\[
p_c(q) = \inf\{p: \phi^1_{p,q}(0\leftrightarrow \infty) > 0\}.
\]
We also consider the threshold $w_c(q)$ for exponential decay of
connectivities,
which is defined as the supremum of all $p$'s for which
\[
\phi^1_{p,q}(0\leftrightarrow \partial\La) \leq C e^{-c\,d(0,\partial\La)}
\]
uniformly in $\La$ (or, at least, for $\La$ in a prescribed sequence
increasing to $\L$) with suitable constants $c>0$ and $C<\infty$. It
is evident that $w_c(q) \le p_c(q)$; for large $q$
it is known that $w_c(q) = p_c(q)$ \cite{vEFSS}.
Theorem \ref{AC} then gives us the following conditions for
high-temperature behavior; compare with Theorem \ref{thm:uni}.

\begin{cor}\label{rcdomcor}
Whenever $\b$ is so small that $p < p_c(q)$, there is a unique Gibbs
measure for the Hamiltonian $H$ with pair interaction
\rf{pot_homogeneous}. Furthermore, if in fact $p< w_c(q)$ then the
spin system is exponentially weak-mixing in the sense that there
are positive constants $C<\infty,c>0$ such that for all $\De\subset
\La\in\E$ and all boundary conditions $\eta,\eta'\in
\Omega$
\[
\|\mu_{\b,\La}^\eta - \mu_{\b,\La}^{\eta'}\|_\De
\leq C\,|\partial \De|\, e^{-c\, d(\De,\La^c)}\, .
\]
\end{cor}
In fact, Alexander and Chayes \cite{AC} go a bit further in their
exploration of `nice' high temperature behavior,
showing that for $p<p_c(1)$ the unique Gibbs measure satisfies
the condition of `complete analyticity' (investigated in \cite{DS2},
for example).

\subsection{Exponential mixing at low temperatures}
\label{sect:appli}

In the previous subsections we have seen how stochastic-geometric
methods can be used to analyze the high temperature behavior of a
spin system and, in particular, for establishing exponential decay
of correlations. Here we want to demonstrate that similar
percolation techniques can also be used in the low temperature
regime in the presence of phase transition. We will present a
method to show that, for a given phase $\mu$, the covariance
$\mu(f;g)$ of any two local observables $f$ and $g$ decays
exponentially fast with the distance between their dependence sets.
(The problem of phase transition at low temperatures will be
addressed in Section \ref{sect:agreep}.)

As a matter of fact, the problem of exponential decay of covariances
(or truncated correlation functions) arises in many physical situations.
Correlation functions are related to interesting response functions
or to fluctuations of specific order parameters.
Exponential decay of covariances also provides estimates on higher order
correlation functions,  eventually providing
infinite differentiability of the free energy with respect to
an external field \cite{DKP}. Motivated by these interests,
a variety of techniques have been developed. The most familiar
approach are cluster expansions which apply equally well
to both the high temperature (or low density) regime and the low
temperature (or high density) regime.
Although they often employ geometric concepts,
it seems useful to combine them with ideas of percolation theory
to make geometry more visible. An example of this is the method
to be described below which is taken from a paper by
Burton and Steif \cite{BS}, where it was used to show that
certain Gibbs measures exhibit a powerful mixing property called
`quite weak Bernoulli with exponential rate'.

We consider a spin system on an arbitrary graph $({\cal L},\sim)$
with Hamiltonian \rf{eq:Hamiltonian}. As before, the essential feature
of this Hamiltonian is that only adjacent spins interact, so that the
Gibbs distributions
$\mu_{\b,\La}^{\e}$ in finite regions $\La$ have the Markov property.
The inverse temperature $\b>0$ does not play any role for the moment,
so we set it equal to 1 and drop it from our notation.

Our starting point is the following estimate on the $\e$-dependence
in terms of disagreement paths for two {\em independent\/} copies
of $\mu_{\La}^{\e}$. This result is a weak version (and, in fact, a
forerunner \cite{vdB}) of Theorem  \ref{thm:dis}. It is a pleasant
surprise that although developed with high temperature situations
in mind, it also provides a useful alternative to some aspects of
the standard low temperature expansions.
\begin{prop}\label{apr} For any $\De\subset\La\in\E$
and  $\e,\e'\in\Om$,
\[
\|\mu_\La^\e-\mu_\La^{\e'} \|_{\De}   \leq \mu_\La^\e\times
\mu_\La^{\e'}(\De \stackrel{\ne}{\longleftrightarrow}\partial\La)\;.
\]
\end{prop}
{\bf Proof: } For brevity let $P=\mu_\La^\e\times\mu_\La^{\e'}$, and
write $X,X'$ for the two projections from $\Om\times\Om$ to $\Om$.
Then for any $A\in\F_{\De}$ we have
\[
\mu_\La^\e(A)-\mu_\La^{\e'}(A)=P(X\in A)-P(X'\in A)\;.
\]
We decompose the probabilities on the right-hand side
into the two contributions according to whether the event
$\{\De \stackrel{\ne}{\longleftrightarrow}\partial\La\}$
occurs or not. In the latter case,
 there exists a random set $\Gamma\subset\La$ containing
$\De$ such that $X\equiv X '$ on $\partial\Gamma$. (The union
of all  disagreement clusters in $\La$ meeting $\De$ is such a set.)
Let $\Gamma$ be the maximal random subset
of $\La$ with this property. Then for each $G$
the event $\{\Gamma=G\}$ only depends on the configuration outside $G$,
and $X\equiv X '$ on $\partial G$.
The Markov property therefore implies that, conditionally on
$\{\Gamma=G\}$
and $(X_{G^c},X'_{G^c})$,  $X_{G}$ and $X'_{G}$ are independent and
identically distributed, and this shows that
$P(X\in A,\,\De \stackrel{\mbox{\tiny$\ne$}}{\not\leftrightarrow}
\partial\La )
=P(X'\in A,\,\De \stackrel{\mbox{\tiny$\ne$}}{\not\leftrightarrow}
\partial\La )$.
The proposition now follows immediately. $\Cox$

\medskip\noindent
What is gained with the disagreement estimate above? First, let us
observe that this estimate provides bounds on covariances of local
functions in terms of disagreement percolation.
\begin{cor}\label{apr2}
Fix any $\Lambda \in \E$ and $\eta \in \Omega$.  Let $f$ and $g$ be
any two local functions depending on the spins in two disjoint
subsets $\De$ resp.\ $\De'$ of $\La$. Then
\[
|\mu_\La^\e(f;g)|
\le \de(f)\,\de(g)\;\mu_\La^\e\times\mu_\La^\e
(\De \stackrel{\ne}{\longleftrightarrow}\De'\mbox{ in }\La)\;
\]
where $\de(f)= \max_\xi f(\xi) - \min_\xi f(\xi)$ is the total oscillation
of $f$.
\end{cor}
{\bf Proof: } By rescaling and addition of suitable constants we can
assume that $0\le f\le 1$ and  $0\le g\le 1$. Proposition \ref{apr}
then shows that
\begin{eqnarray*}
&&|\mu_\La^\e(f;g)|
\le\int\mu_\La^\e(d\xi)\int\mu_\La^\e(d\xi')\; g(\xi')\,
\bigg|\mu_{\La\setminus\De'}^\xi(f)
-\mu_{\La\setminus\De'}^{\xi'}(f)\bigg|\\
&&\le\de(f)\,\de(g)\;\int\mu_\La^\e(d\xi)\int\mu_\La^\e(d\xi')\;
\mu_{\La\setminus\De'}^{\xi}\times\mu_{\La\setminus\De'}^{\xi'}
(\De \stackrel{\ne}{\longleftrightarrow}\De'\mbox{ in }\La)
\end{eqnarray*}
because $\xi\equiv\xi'\equiv\e$ on $\partial\La$. By carrying out the
last integration we obtain the result. $\Cox$

\medskip\noindent
The bounds above leave us with the task of estimating the
probability of disagreement paths in a duplicated system. In
contrast to the situation in Section \ref{sect:disa}, we are
looking now for estimates valid at low temperatures. If a cluster
expansion works, there is no need to look any further. For
instance, a low temperature analysis and estimates of
semi-invariants for the Ising model can be obtained using standard
contour representations; see \cite{rlD}. It needs to be emphasized,
however, that the main step of cluster expansions consists in
expanding the logarithm of the partition function. Only afterwards,
by taking ratios of partition functions, does one obtain
expressions for covariances and higher order correlation functions.
Therefore, a point to appreciate is that the bound of Corollary
\ref{apr2}  provides a direct geometric bound on covariances which
avoids the machinery of cluster expansions and, in particular, the
problems coming from taking logarithms. As a consequence, this
estimate also applies to some cases where standard cluster
expansions are doomed to fail.

Let us illustrate this for the case of the Ising ferromagnet. This might
not be the best example because other methods can also be applied to it;
nevertheless, it is useful to demonstrate the technique in this simple
case. Afterwards we will discuss a case where cluster expansion
techniques cannot be used equally easily.

Consider the low temperature plus phase $\mu_{\b}^+$ of the
ferromagnetic Ising model on the
square lattice $\L=\Z^2$ with zero magnetic field. By taking the
infinite volume limit in Corollary \ref{apr2} with $\eta \equiv +1$
we obtain the estimate
\begin{equation}\label{plus_cov}
|\mu_{\b}^+(f;g)|
\le \de(f)\,\de(g)\;\mu_{\b}^+\times\mu_{\b}^+
(\De \stackrel{\ne}{\longleftrightarrow}\De')
\end{equation}
for the covariance of any two local functions $f,g$ with disjoint
dependence sets $\De,\De'$. Next we observe that the event $\{\De
\stackrel{\ne}{\longleftrightarrow}\De'\}$ is clearly contained in
the event that there exists a path from $\De$ to $\De'$ along which
$(X,X') \neq (+,+)$. The latter event will be denoted by $\{\De
\stackrel{\mbox{\tiny$\ne(+,+)$ }}{\longleftrightarrow}\De'\}$.
Now, Burton and Steif \cite{BS} have shown how to estimate the
probability of this event to occur in the duplicated plus phase.
The result is the following.
\begin{thm} \label{thm:bs}
For the Ising ferromagnet on $\Z^2$ at sufficiently large
$\b$, there exist  constants $c > 0$ and
$C<\infty$ (depending on $\b$)
such that for any two disjoint sets $\De,\De'\in\E$,
\[ 
\mu_{\b}^+\times\mu_{\b}^+(\De \stackrel{\mbox{\tiny$\ne
(+,+)$}}{\longleftrightarrow}\De')
\leq C \min\{|\partial_i \De|,|\partial_i \De'|\}\; e^{-c\,d(\De,\De')}\;,
\]
where $\partial_i \De=\partial(\De^c)$ is the inner boundary of $\De$.
\end{thm}
Combining this theorem with \rf{plus_cov} we obtain
an exponential bound for the covariance of any two local observables
in the low temperature Ising plus phase. While this result is well-known,
its proof below shows how one can proceed in more general cases.

\medskip\noindent
{\bf Sketch proof of Theorem \ref{thm:bs}:} It is sufficient to prove
the statement with $\mu_{\b}^+$ replaced by $\mu_{\b,\La}^+$, where
$\La$ is a sufficiently large square box containing $\De,\De'$. For
brevity, let $P_{\La}=\mu_{\b,\La}^+\times\mu_{\b,\La}^+$.
Suppose that $|\partial \De'| \leq |\partial \De|$, fix an arbitrary $x
\in\partial_i \De'$, and suppose that the event
$\{x \stackrel{\mbox{\tiny$\ne(+,+)$}}{\longleftrightarrow}\De\}$
occurs. Let $\Gamma_{0}=\Gamma_{0}(X,X')$ be
the maximal connected set containing $x$ on which $(X, X') \neq
(+,+)$. Also, let $\Gamma$ be the union of $\Gamma_{0}$ and all
finite components of ${\Gamma_{0}}^c$. $\Gamma$ is enclosing in the
sense that both $\Gamma$ and $\Gamma^c$ are connected.
In fact, $\Gamma$ can be identified with a contour, the broken
line which separates $\Gamma$ from its complement.

Consider the set ${\cal C}$ of all enclosing sets containing $x$
and contained in $\La$. For any integer $\ell\ge2$, let ${\cal
C}_{\ell}$ be the set of all $C\in{\cal C}$ such that $|\partial_i
C| = \ell$. Since the number of contours with length $\ell$
surrounding a given site on the square lattice is bounded by $\ell
3^{\ell}$ and since each enclosing set with $|\partial_i C| =
\ell$ uniquely defines a contour with length between $\ell$ and $4\ell$,
we have that $|{\cal C}_{\ell}|
\leq
3(4\ell + 1)/2 (81)^{\ell}$ growing exponentially with $\ell$. Now
we can write
\[
P_{\La}(x \stackrel{\mbox{\tiny$\ne(+,+)$}}{\longleftrightarrow}\De)
=
P_{\La}(\Gamma\cap \De\neq \emptyset)
\le
\sum_{\ell\ge d(\De,\De')}\ \sum_{C \in {\cal C}_{\ell}}
P_{\La}(\Gamma = C)\;.
\]
Furthermore, for $C \in {\cal C}_{\ell}$ we have
\begin{eqnarray}\label{prob}
P_{\La}(\Gamma = C)
&\leq&\sum_{D \subset \partial_i C} \mu_{\b,\La}^+(X \equiv -1\mbox{ on }
D, \, X \equiv +1\mbox{ on }\partial C) \nonumber \\
&&\hspace{4em}\times\,
\mu_{\b,\La}^+(X\equiv -1\mbox{ on }\partial_i C \setminus D, \,
X\equiv +1\mbox{ on }\partial C)\;.
\end{eqnarray}
Now, a standard Peierls estimate (see e.g. \cite{Si}) shows that
\begin{equation} \label{pir}
\mu_{\b,\La}^+(X \equiv -1 \mbox{ on } D, X \equiv +1 \mbox{ on }
\partial C)\leq C(\b)\; e^{-c(\b)\, |D|}
\end{equation}
with constants $c(\b), C(\b)$ independent of $\La$ satisfying
$c(\b)\rightarrow \infty$ as $\b \to\infty$.
Substituting (\ref{pir}) into \rf{prob} we obtain the theorem by simple
combinatorics and summations. $\Cox$

\medskip\noindent
We emphasize that the specific properties of the plus phase $\mu_{\b}^+$
are used only in the last step, the Peierls estimate \rf{pir}. Before,
we needed only the Markov property. Therefore it is useful to  note that
the Peierls estimate is not limited to the Ising model; it remains
valid under the conditions of the standard Pirogov--Sinai theory \cite{Si}.
In particular, it follows that the results of Burton and Steif \cite{BS}
on the ergodic properties of the Ising model carry over to more general
Markovian models of Pirogov--Sinai type.

Let us finally discuss a case in which a Peierls estimate of the form
(\ref{pir}) is not available. Namely, we ask for covariance estimates
of local functions, still for the
ferromagnetic Ising model in a large square $\La$,  but now for
some boundary condition $\eta$ not identically equal to $+1$.
This question arises, for example, in the context of
correlations atop of a disordered surface, or in the problem
of establishing a Gibbsian description of non-Gibbsian measures. In fact,
in \cite{MaV} the method to be described below is used to
prove that the projection to a line of the low-temperature plus phase of
the two-dimensional Ising model is weakly Gibbsian.

To be specific, suppose that the boundary condition $\eta$ on $\partial
\La$ is not identically plus but contains a large proportion of
plus spins; we stipulate that  $\e \equiv +1$ on three sides of $\La$
while on the remaining side only a large fraction of the spins is plus.
In this case, (\ref{pir}) cannot be true because $D$ can be small
and close to the boundary of $\La$. For example, if $E=\partial D
\cap \partial\La \neq \emptyset$  and
$\eta$ happens to be minus on $E$ then it is not very unlikely that all
spins
in $D$ are minus,  and (\ref{pir}) will not hold.
However, one can take advantage of the
fact that for small $D$ its complement in $\partial_i C$ is
large, and vice versa. In other words, to estimate the right-hand
side of (\ref{prob}) one should not apply (\ref{pir}) separately to
each factor, but rather one can hope to estimate their product.

To make these general remarks precise we consider the Ising model on the
half-plane $\Z \times \Z_{+}$.  For any $n$ we consider the square
$\La_n =
\{(x_1,x_2) \in \Z^2: -n\leq x_1 \leq n,\, 0 < x_2
\leq n\}$ touching the boundary line
$\Z \times \{0\}$.  On this line we fix a configuration
$\xi \in \{+1,-1\}^\Z$, thereby defining a boundary condition
on one part of $\partial \La$.  On the remaining part of
$\partial \La$ we impose plus boundary condition. That is,
we choose the boundary condition
$\eta\equiv +1$ on $\partial \La_n \setminus
(\Z \times\{0\})$ and $\eta\equiv\xi$ on $\Z\times\{0\}$.
We ask for the correlation of the spins at the sites $x=(0,1)$
and $y=(k,1)$ with $0<|k|<n$.
\begin{thm} In the situation just described, suppose $\xi$ is such that
\[
\sum_{j=0}^m \xi(j,0) \geq 8m/9\quad \mbox{ and }\quad
\sum_{j=-m}^{-1} \xi(j,0) \geq 8m/9
\]
for sufficiently large $m$, and let $\e$ be defined as above. Then
there are constants $ c>0$ and $C< \infty$
(not depending on $n$) such that
\[
|\mu_{\b,\La_n}^\eta(X(0,1);X(k,1))| \leq C\, e^{-c\, |k|}
\]
whenever $\b$ and $|k|$ are sufficiently large and $n > |k|$.
\end{thm}
{\bf Sketch proof: } We proceed as in the proof of Theorem \ref{thm:bs}.
In dealing with the right-hand side of (\ref{prob}) we must take into
account
that possibly $\partial C \cap \partial\La \neq \emptyset$. We
therefore replace $\partial C$ by $\partial C \setminus \partial \La$
in the product term and also estimate the
probabilities of intersections by conditional probabilities,
yielding the upper bound
\[
\mu_{\b,C}^{+,\e}(X \equiv -1\mbox{ on }D) \;
\mu_{\b,C}^{+,\e} (X\equiv-1\mbox{ on }\partial_i C
\setminus D)
\]
for the summands on the right-hand side of (\ref{prob}).
Here, $\mu_{\b,C}^{+,\e}$ stands for the Gibbs distribution
in $C$ with boundary condition equal to $+1$ on $\partial C\cap\La$
and equal to $\eta$ on $\partial C\cap\partial \La$.
To derive the theorem
we need to replace the Peierls estimate (\ref{pir}) by a similar
bound on the last product. The exponential decay of correlations
then again follows by simple combinatorics and summations.

To make the influence of the boundary condition $\eta$ explicit  we
exploit a contour representation leading to the estimate
\begin{eqnarray} \label{form3}
&&
\mu_{\b,C}^{+,\e}(X \equiv -1\mbox{ on }D) \;
\mu_{\b,C}^{+,\e}(X\equiv -1\mbox{ on }\partial_i C\setminus D) \nonumber
\\[1ex]
&&\hspace{4em}\leq
\sum_{{\Gamma, \Gamma' \mbox{\tiny{ inside }} C
\atop \Gamma \mbox{\tiny{ compatible with }}  D} \atop
\Gamma' \mbox{\tiny{ compatible with }} \partial_i C
\setminus D}
\prod_{\g\in\Gamma}w_{\eta}(\g) \prod_{\g'\in\Gamma'}w_{\eta}(\g')\;.
\end{eqnarray}
The right-hand side is defined as follows. For any configuration
$\s \in\{+1,-1\}^\La$ we draw horizontal resp.\ vertical lines of unit
length between neighboring sites of opposite spins, doing
as if the boundary spins were all plus; we then obtain a disjoint union
of closed non-self-intersecting polygonal curves.  Each of these
curves is called a contour $\g$, and a set $\Gamma$ of contours
arising in this way is called compatible.
We thus have a one-to-one
correspondence between spin configurations $\s$ and compatible sets
$\Gamma$ of contours.
If  $\s \equiv -1$ on $D$, then each component of $D$ is surrounded by
some contour $\g$ (i.e., belongs to the interior $\mbox{Int}\,\g$ of $\g$);
the smallest contours surrounding the components of $D$ are collected into
a set $\Gamma$ of contours. Each
set $\Gamma$ arising in this way is called {\em compatible
with $D$}. The probability that a given set $\Gamma$ of contours occurs
is not larger than $\prod_{\g\in\Gamma}w_{\eta}(\g)$, where
\[
w_{\eta}(\g) = \exp \bigg[-2\b |\g| + 2\b \sum_{x \in
\mbox{\scriptsize Int}\, \g} \sum_{y\in \partial \La: y\sim x}
(1-\eta(y))\bigg]
\]
and $|\g|$ is the length of $\g$; this can be seen by comparing the
probability of a configuration containing $\Gamma$ with the probability
of the configuration obtained by flipping the spins in
$\bigcup_{\g\in\Gamma}\mbox{Int}\,\g$. These observations establish the
inequality \rf{form3}.

Note that the weight $w_{\eta}(\g)$  of a contour
$\g$ depends on the boundary configuration $\eta$; this is
because we have chosen to draw the contours for plus boundary
conditions rather than $\e$.  It follows that
$w_{\eta}(\g)$ does not necessarily tend to zero when $|\g|$
grows to infinity; this is in contrast with the case $\eta \equiv +1$.
The standard low temperature expansion would therefore become much
more complicated. However, if the density of plus spins in $\eta$ is
sufficiently large, or if $\partial\mbox{Int}\,\g \cap
\partial \La$ is rather small, the standard weight
$\exp[-2\b|\g|]$ of the Ising contours will dominate, and
the right-hand side of (\ref{form3}) can be estimated, as we will show now.

We unite the sets $\Gamma, \Gamma'$ in
(\ref{form3}) into a single set of contours
$\tilde{\Gamma} = \Gamma \cup \Gamma'$. The contours in
$\tilde{\Gamma}$ can overlap, but a site of $\partial_i C$ can only
belong to the interior of
at most two contours. On the other hand, every site of $\partial_i C$ is
in the interior of at least one contour of
$\tilde{\Gamma}$, and $|\partial_i C|\ge k$, the
distance of the two spins considered. These
ingredients allow us to control the sum on the right-hand side of
(\ref{form3}). If $k$ is so large
that the density of plus spins in $\xi$ between 0 and $(k,0)$ exceeds
$8/9$ then we find for any collection of contours $\tilde{\Gamma}=
\Gamma\cup\Gamma'$ as above
\[
\sum_{\tilde{\g}\in\tilde{\Gamma}}\
\sum_{x \in \mbox{\scriptsize Int}\,
\tilde{\g}}\ \sum_{y\in \partial \La: y\sim x} (1-\eta(y)) \leq
5/9\;\sum_{\tilde{\g}\in\tilde{\Gamma}} |\tilde{\g}|\;.
\]
This yields
\[
\prod_{\g\in\Gamma}w_{\eta}(\g) \prod_{\g'\in\Gamma'}w_{\eta}(\g') \leq
\prod_{\g\in\Gamma} \exp[-8/9 \,\b\, |\g|] \prod_{\g'\in\Gamma'}
\exp[-8/9 \,\b \,|\g'|] \;.
\]
At this point the standard arguments take over (with $\b$ replaced by $
4\b/9$), leading to an exponential estimate of (\ref{form3}).  For
example, one can conclude the proof along the
lines of Lemma 2.5 of \cite{BS}. $\Cox$

\section{Phase transition and percolation}
\label{sect:agreep}

Typically, two ends of the phase diagram are amenable to mathematical
analysis. One is the high temperature, or low density,
regime which was discussed in the previous section and in which the
system can be viewed as a small perturbation of an independent spin
system. The other end is the low temperature regime which we will
now consider.
At low temperatures, the energy dominates over the entropy which comes
from the thermal fluctuations of the spins.
One therefore expects that the spin configuration is typically similar to
some
frozen zero temperature state, which is a configuration of
minimal energy and thus called a {\em ground state}. The similarity of a
low
temperature state with a ground state is conveniently described in
geometric terms: one imagines that the spins which agree with the given
ground state form an infinitely extended sea, whereas those spins which
have chosen to deviate from the ground state are confined to interspersed
finite islands. This is, of course, a picture of percolation
theory: spins that agree with the ground state form a unique infinite
cluster. We are thus led to the concept of {\em agreement percolation},
which will be discussed in the first part of this section.

In fact, agreement
percolation is intimately related to the existence of a phase transition.
If several distinct ground states exist, we may hope to find at low
temperatures also several equilibrium phases which can be distinguished by
agreement percolation with respect to the different ground states. One may
ask further whether the geometric picture that applies to low
temperatures remains
valid throughout the whole non-uniqueness region. Physically, this
is a matter of stability of the ground states. Mathematically, it means
to look for conditions under which distinct Gibbs measures allow
distinct stochastic-geometric characterizations.

We will approach this question from two different sides. In
Sections \ref{subsect:agreement} to \ref{subsect:appl_other_models}
we investigate whether ``phase transition implies percolation''. We
study a fixed equilibrium phase $\mu$ in the non-uniqueness region
which, by its very construction, can be viewed as a random
perturbation of some ground state $\e$. We then will see that, in
many cases, spins that agree with $\e$ do percolate. After a
general discussion of agreement percolation in Section
\ref{subsect:agreement}, we investigate this concept in the
subsequent subsections for some specific models including the Ising
ferromagnet and the Potts model. In the case of the planar lattice
$\Z^2$ with its limited geometric possibilities we will also see
that conversely, the absence of phase transition sometimes implies
an absence of percolation, and that in the case of phase transition
one has restrictions on the number of phases. (Methodologically,
these results still run under the heading ``phase transition
implies percolation''.) In the last Section
\ref{ssect:ground_energy_perc} the converse will be treated more
systematically and under a different aspect: we will show that at
low temperatures one has percolation of bonds along which the
interaction energy is minimal, and we will see that such a {\em
ground-energy bond percolation} often implies a phase transition.
Taken together, these results will show that in various models
phase transition comes along with the existence of a ground-state
sea with finite islands (deviation islands) on which the spins
deviate from the ground state, and vice versa.

A theory developing this picture in much more detail is the
Pirogov--Sinai theory of phase transition which deals with the low
temperature phase diagram in the presence of several stable ground
states. One basic idea of this theory is to treat the finite
deviation islands of a low temperature system as the constituents
of a low density gas of hard-core particles. While the
Pirogov--Sinai theory is intimately related to the subject of this
section, it is much too involved to be developed here. There is,
however, a number of expositions which may serve as general
introductions and in which many additional references can be found.
We mention only \cite{BS,DOS,vEFS,Si,Sla} and \cite{Zah1} to
\cite{Zah3}. Here we will concentrate on more specific results
which are partly beyond the Pirogov--Sinai theory, in that they are
not limited to low temperatures but rather apply to the full
non-uniqueness region, and comment occasionally on some
relationships. In particular, the results of Section
\ref{ssect:ground_energy_perc} are similar in spirit to this
theory.

\subsection{Agreement percolation from phase coexistence}
\label{subsect:agreement}

We consider again the general setting of Section \ref{sect:eq_phases}.
$(\L,\sim)$ is an arbitrary locally finite graph, $S$ is a finite set,
and $\Om=S^\L$. Suppose $\mu$ is a random field and $\e\in\Om$ a
fixed configuration. We consider the event
$\{x\stackrel{\e}{\longleftrightarrow}\infty\}$ that $x\in\L$
belongs to an infinite cluster of the random
set $R(\eta) = \{y\in {\L}:X(y)=\eta(y)\}$, and we say that {\em $\mu$
exhibits agreement percolation for $\e$}
if $\mu(x\stackrel{\e}{\longleftrightarrow}\infty)>0$ for some $x\in\L$.
In short, we will then simply speak of $\e$-percolation.
To visualize such an agreement, it may be convenient to think of a
reduced description of $\mu$ in terms of its image under the map
$s_{\e}:\Omega\rightarrow\{0,1\}^{\L}$, which describes local
agreement and disagreement with $\eta$, and is defined by
\begin{equation} \label{a/dis}
(s_{\eta } (\sigma))(x)=\left\{
\begin{array}{ll}
1  &\mbox{if } \s (x) = \e(x), \\
  0   &\mbox{otherwise.} \\
\end{array} \right.
\end{equation}
With this mapping, we can write
$\{x\stackrel{\e}{\longleftrightarrow}\infty\}=s_\e^{-1}\{x
{\leftrightarrow} \infty\}$.

We are interested here in the case when $\mu$ is a Gibbs
measure for the Hamiltonian (\ref{eq:Hamiltonian}), and $\e$ is an
associated ground state. We say that a configuration $\e\in
\Omega$ is a {\it ground state} or, more explicitly, a ground state
configuration for the relative Hamiltonian
\rf{relener}, if $H(\s|\e)\geq0$  for any local modification (or
``excitation'') $\s$ of $\e$.
In other words, $\e$ is a ground state if, for any region $\Lambda
\in \cal E$, the configuration $\e$ minimizes the energy in
$\Lambda$ when $\e_{\Lambda^c}$ is fixed.
One should note in this context that, in the low temperature limit
$\beta\uparrow\infty$,
the finite-volume Gibbs distribution $\mu_{\beta,\Lambda}^\eta$
 from \rf{eq:FvolGibbs}
tends to the equidistribution on the set of all configurations $\sigma$
of minimal energy $H(\sigma|\eta)$.  This fact suggests that, at
least in some
cases, the low temperature phase diagram is only a
slight deformation of the zero-temperature phase diagram describing the
structure of ground states. This is precisely the subject of
Pirogov--Sinai theory which provides sufficient conditions for this to
hold, proposes a construction of low temperature phases as
perturbations of ground states, and also shows that the size distribution
of the deviation islands has exponential decay.

Suppose next that the Gibbs measure $\mu$ is related to the
ground state $\eta$ in some way. For example, $\mu$ might be obtained
as the infinite volume limit of the finite volume
Gibbs distributions $\mu_{\b,\La}^\e$ with boundary condition $\e$,
 possibly along some subsequence. (Under the
conditions of the Pirogov--Sinai theory such a limit always exists.)
In the case of a phase transition, when other phases than $\mu$
exist and one is interested in characteristic properties of $\mu$,
one expects that the relationship between $\mu$ and $\e$ becomes
manifest in a macroscopic pattern of the typical configurations, in
that $\mu$ shows $\e$-percolation. In short, we ask for the
validity of the hypothesis
\begin{eqnarray}\label{hypoagree}
&& |{\cal G}(\b H)|>1,\, \mu \mbox{ is extremal in } {\cal G}(\b H)
\mbox{ and related to a ground state }
 \e\in\Om\nonumber\\
 &&\Longrightarrow\
 \mu(x \stackrel{\eta}{\longleftrightarrow} \infty) >0 \quad \forall\,
 x\in\L\;.
\end{eqnarray}
In the specific cases considered below it will always be clear in what
sense $\mu$ and $\e$ are related; typically, $\mu$ will be a limiting
Gibbs measure with boundary condition $\e$.
We emphasize that (\ref{hypoagree}) does not hold in general; a
counter-example
can be constructed by combining many independent copies of the Ising
ferromagnet
to a layered system, see the discussion after Proposition
\ref{prop:reflected_plus-perc} below.
Also, even when (\ref{hypoagree}) holds,
it does not necessarily imply that the phase $\mu$ is uniquely
characterized
by the property of $\e$-percolation.

How can one establish (\ref{hypoagree})? In the context of the Ising
model, Coniglio et al.\ \cite{CNPR} and Russo \cite{Rus} developed
a convenient citerion which is based on a multidimensional analog of the
strong Markov property and thus can be used for general Markov random
fields
\cite{BLM,GLM}. One version is as follows.
\begin{thm} \label{thmrusso}
Let $(\L,\sim)$ be a locally finite graph, $\mu$ a Markov field on
$\Om=S^\L$, and $\e\in\Om$ any configuration.
Suppose there exist a constant $c\in\R$ and a local function
$f:\Om \rightarrow \R$ depending only on the configuration in
a connected set $\De$, such that $\mu(f)>c$ but
\begin{equation}\label{order}
\mu(f \,|\,  X\equiv \xi \mbox{ \em on }\partial \Gamma) \leq c
\end{equation}
for all finite connected sets $\Gamma\supset\De$ and all $\xi \in
\Om$ with $s_\e(\xi)\equiv0$ on $\partial \Gamma$. Then
\linebreak $\mu(\De\stackrel{\eta}{\longleftrightarrow}\infty)>0$,
i.e., $\mu$ exhibits agreement percolation for $\eta$.
\end{thm}
{\bf Proof:}
 Suppose by contraposition that
$\mu(\De\stackrel{\eta}{\longleftrightarrow}\infty)=0$. For any
$\eps>0$ we can then choose some finite $\La\supset\De$ such that
$\mu(\De\stackrel{\eta}{\longleftrightarrow}\La^c)<\eps$. For
$\xi\notin\{\De\stackrel{\eta}{\longleftrightarrow}\La^c\}$, there
exists a connected set $\Gamma$ such that
$\De\subset\Gamma\subset\La$ and $s_\e(\xi)\equiv 0$ on
$\partial\Gamma$; we simply let $\Gamma$ be the union of $\De$ and
all $\e$-clusters meeting $\partial\De$. As in the proof of Theorem
\ref{thm:Potts_sandwich}, we let $\Gamma(\xi)$ be the largest such
set. For $\xi\in\{\De\stackrel{\eta}{\longleftrightarrow}\La^c\}$
we put $\Gamma(\xi)=\emptyset$.
Then, for each finite connected set $\Gamma\ne\emptyset$,
the event $\{\xi:\Gamma(\xi)=\Gamma\}$ depends only on the
configuration in $\La\setminus\Gamma$, whence by the Markov property
$\mu(f\,|\,\Gamma(\cdot)=\Gamma)$ is an average of
the conditional probabilities that appear in assumption (\ref{order}).
>From this we obtain
\[
\mu(f)\le
c\,\mu(\Gamma(\cdot)\ne\emptyset) +
\mu(|f|\,I_{\{\Gamma(\cdot)=\emptyset\}})<
 c+\eps\,\|f\|\;.
\]
Letting $\eps\to0$ we find $\mu(f)\le c$, contradicting our
assumption.
$\Cox$

\medskip\noindent
In most applications we will have a natural candidate for
the function $f$. Whenever distinct phases do exist, they can be
distinguished by some order parameter, viz.\ a local function $f$ having
different expectations for the two phases. If, in addition, some
stochastic monotonicity is available then we can hope to establish
\rf{order}. In fact, the percolation phenomena stated in Examples
\ref{ex:BLM1} to \ref{ex:BLM4} can be deduced from a slight modification
of
Theorem \ref{thmrusso}. We will not go into the details of these examples
which are treated in \cite{BLM}, but rather apply Theorem \ref{thmrusso}
to our standard examples.

\subsection{Plus-clusters for the Ising ferromagnet}
\label{subs:ferro}

The idea of agreement percolation was first developed in the context of
the ferromagnetic Ising model \cite{CNPR,Rus}. Let us apply Theorem
\ref{thmrusso} to this standard case. We only consider the case of no
external field, i.e., we set $h=0$, so the only parameter is the inverse
temperature $\b>0$. We are interested in agreement percolation for the
constant
configurations $\e\equiv+1$ resp.\ $\e\equiv-1$, which are the only
periodic ground
states of the model. We write $\stackrel{+}{\longleftrightarrow}$ resp.\
$\stackrel{-}{\longleftrightarrow}$ for the corresponding connectedness
relation. Our first result shows that if there is a phase transition then
there is plus-percolation for each Gibbs measure except the minus-phase
$\mu_\b^-$; that is, assertion \rf{hypoagree} holds for $\e\equiv+1$. This
result (due to \cite{Rus})
 is valid for an arbitrary locally finite graph $(\L,\sim)$ with finite
critical inverse temperature $\b_c$.
\begin{thm}\label{thm:Ising_plus_agreep}
Let $\mu$ be an arbitrary Gibbs measure for the ferromagnetic Ising model
with parameters $\b>0$, $h=0$. If $\mu\ne\mu_\b^-$ then
$\mu(x\stackrel{+}{\longleftrightarrow}\infty)>0$ for all $x\in\L$. In
particular, if $\b>\b_c$ then
$\mu_\b^+(x\stackrel{+}{\longleftrightarrow}\infty)>0$ for all $x\in\L$.
\end{thm}
{\bf Proof: } By the sandwiching inequality \rf{eq:Ising_sandwich} and
Proposition
\ref{prop:single_site+dom}, there exists a site $x\in\L$ such that
$\mu(X(x)=1)>c\equiv\mu_\b^-(X(x)=1)$. On the other hand, the
analogue of inequality \rf{eq:decreasing} for the minus boundary
condition shows that
\[
\mu(X(x)=1\,|\,X\equiv-1\mbox{ on
}\partial\Gamma)=\mu_{\b,\Gamma}^-(X(x)=1)\le c
\]
for every finite $\Gamma\ni x$. Theorem \ref{thmrusso} thus gives the
result for the $x$ at hand. In view of the finite energy property of
$\mu$, this extends easily to all other $x\in\L$. $\Cox$

\medskip\noindent
Let us rephrase the last statement of Theorem \ref{thm:Ising_plus_agreep}
as
follows: below the critical temperature
the plus spins percolate in the plus phase and, by symmetry,
the minus spins percolate in the minus phase $\mu_\b^-$. In the case of
graphs with symmetry axes, this statement allows an interesting
refinement.
\begin{prop}\label{prop:reflected_plus-perc}
Suppose $(\L,\sim)$ admits an involutive graph automorphism $r$  which
maps a subset ${\cal H}\subset\L$ onto its complement ${\cal H}^c$, and
that for $x\in{\cal H}$, $y\in{\cal H}^c$ either $x\not\sim y$ or $y=rx$.
For $x\in{\cal H}$ let
$\{x\stackrel{\mbox{\tiny+$r$+}}{\longleftrightarrow}\infty\mbox{ \em in
}{\cal H}\}$ be the event that there exists an infinite path $\gamma$ in
${\cal H}$ starting from $x$ such that all spins along both $\gamma$ and
its
reflection image $r\gamma$ are positive. If $\mu_\b^+\ne\mu_\b^-$ then
$\mu_\b^+(x\stackrel{\mbox{\tiny+$r$+}}{\longleftrightarrow}\infty \mbox{
\em in }{\cal H})>0$ for all $x\in{\cal H}$.
\end{prop}
One natural case to think of is when $\L=\Z^d$ for $d\ge2$, ${\cal H}$ a
halfspace with boundary orthogonal to an axis, and $r$ the associated
reflection. The proposition then asserts that $\mu_\b^+$-almost surely
there exists an infinite connected mirror-symmetric set of plus spins.
Another interesting case is when $\L$ consists of two disjoint copies of a
graph
${\cal H}$ which are not connected to each other by any bond. In this
case, $\mu_{\b}^+=\mu_{\b,{\cal H}}^+\times\mu_{\b,{\cal H}}^+$, and the
statement is that two independent realizations of $\mu_{\b,{\cal H}}^+$
exhibit simultaneous plus-percolation; in this case, the preceding
proposition was observed by Giacomin et al.\ \cite{GLM}.
It is, however, not possible to take an arbitrarily large number $k$ of
independent realizations $X_1,\ldots,X_k$ of $\mu_{\b,{\cal H}}^+$,
at least when ${\cal H}$ has bounded degree $N$. For, if $p_c$ is
the Bernoulli site percolation threshold of ${\cal H}$ and $k$ is
so large that
$$\sup_{x\in{\cal H}}\mu_{\b,\{x\}}^+(X(x)=1)^k<p_c$$ then the set
$\{x\in{\cal H}:X_1(x)=\cdots=X_k(x)=1\}$ does not percolate.
This follows from a standard domination argument. Since the layered
system consisting of $k$ independent copies of the Ising model with
$\b>\b_c$ certainly exhibits a phase transition, we see that
hypothesis \rf{hypoagree} does not hold in general.

\medskip\noindent
{\bf Proof of Proposition \ref{prop:reflected_plus-perc}: }
We identify each $\xi\in\Om$ with $(\xi(x),\xi(rx))_{x\in{\cal H}}\in
S^{\cal H}$, where $S=\{-1,1\}^2$. The event under consideration then
corresponds to $\e$-percolation for the configuration $\e\in S^{\cal H}$
with $\e(x)=(1,1)$ for all $x\in{\cal H}$. Let $f=X(x)+X(rx)$. Then
$\mu_\b^+(f)=2\,
\mu_\b^+(X(x))>0$ by the $r$-invariance of $\mu_\b^+$. On the other hand,
let $\Gamma\subset{\cal H}$ be a finite set containing $x$, and
$\tilde\Gamma=\Gamma\cup r\Gamma$. If
$(\xi,\xi')\in S^{\cal H}=\Om$ with $s_\e(\xi,\xi')\equiv 0$ on
$\partial_{\cal H}\Gamma=\partial\Gamma\cap{\cal H}$,
then $\xi'\preceq-\xi$ on $\partial_{\cal H}\Gamma$, and therefore
$(\xi,\xi') \preceq(\xi,-\xi)$ (as elements of $\Om$) on
$\partial\tilde\Gamma=\partial_{\cal H}\Gamma\cup r\,\partial_{\cal
H}\Gamma$. We can thus write
\begin{eqnarray*}
\mu_\b^+(f\,|\, (X,X')\equiv (\xi,\xi')\mbox{ on }\partial\Gamma)
&=&\mu_{\b,\tilde\Gamma}^{(\xi,\xi')}(X(x))+
\mu_{\b,\tilde\Gamma}^{(\xi,\xi')}(X(rx))\\
&\le&\mu_{\b,\tilde\Gamma}^{(\xi,-\xi)}(X(x))+
\mu_{\b,\tilde\Gamma}^{(\xi,-\xi)}(X(rx))\ =\ 0
\end{eqnarray*}
by Lemma \ref{lem:boundary_domination} and the symmetry under $r$ and
simultaneous spin flip. The proposition thus follows from Theorem
\ref{thmrusso}. $\Cox$

\medskip\noindent
In the remaining part of this subsection we consider the case of the
square lattice $\L=\Z^2$, in which we can obtain much stronger
conclusions. The following result gives a complete
characterization of the non-uniqueness regime of the parameter space
in terms of percolation of plus spins in the Gibbs measure $\mu_\beta^+$.
It is due to Coniglio et al.\ \cite{CNPR}; see also \cite{Hig}.
\begin{cor} \label{cor:CNPR}
For the Ising ferromagnet on the square lattice $\Z^2$ with no external
field and inverse temperature $\beta$, the $\mu_\b^+$-probability of
having an infinite plus-cluster is $0$ in the
uniqueness regime $\beta\leq \beta_c$, and $1$ in the non-uniqueness
regime
$\beta>\beta_c$.
\end{cor}
{\bf Proof: } The existence of an infinite cluster of plus spins is a tail
event and thus, by the extremality of $\mu_\beta^+$, has probability $0$
or $1$.
The case $\b>\b_c$ is thus covered by Theorem \ref{thm:Ising_plus_agreep}.
For $\b\le\b_c$, $\mu_\beta^+$ coincides with $\mu_\beta^-$. Thus, if an
infinite plus-cluster existed with probability $1$ then, by symmetry, an
infinite minus-cluster would also exist, in contradiction to Theorem
\ref{thm:GKR};
 the assumptions of this theorem are satisfied by Proposition
\ref{prop:plus_measure}.
$\Cox$

\medskip\noindent
Combining the corollary above with Proposition
\ref{prop:Ising_external_field},  we can also obtain some bounds for the
percolative region of the Ising model for $h\ne0$; see \cite{ABL} for a
detailed discussion.

The equivalence of non-uniqueness and percolation just observed for the
Ising model on $\Z^2$ cannot be expected to hold for non-planar graphs.
Consider, for example, the Ising model on the cubic lattice $\Z^3$. For
$\beta=0$ uniqueness certainly holds, and plus-percolation is equivalent
to Bernoulli site percolation on $\Z^3$ with parameter ${1}/{2}$. But
a result of \cite{CR} states that $p_c(\Z^3)<{1}/{2}$. The
plus spins thus percolate at $\b=0$.
In view of Proposition \ref{prop:Ising_external_field}, this is still the
case for sufficiently small $\beta$, so that plus-percolation does occur
in a non-trivial part of the uniqueness region.

For the planar graph $\Z^2$, however, Theorem \ref{thm:GKR} does not only
imply the equivalence of phase transition and percolation, but also gives
some information on the number of phases in the non-uniqueness region. As
a warm-up let us show that, for the Ising ferromagnet on $\Z^2$ at inverse
temperature $\b>\b_c$, there are no other translation and rotation
invariant
extremal Gibbs measures than $\mu_\b^+$ and $\mu_\b^-$. For,
suppose another such phase $\mu$ existed. By Theorem
\ref{thm:Ising_plus_agreep} and the Burton-Keane uniqueness theorem
\ref{thm:Burton_Keane}, there exist unique infinite plus- and
minus-clusters
with $\mu$-probability $1$. As an extremal Gibbs measure, $\mu$ has
positive
correlations; recall the paragraph below Proposition
\ref{prop:plus_measure}. Proposition \ref{prop:GKR2} thus shows that $\mu$
cannot exist.

The statement just shown is a weak version of the following result which
characterizes all translation invariant Gibbs measures. In fact, it is
sufficient to assume periodicity, which means invariance under the
translation
subgroup $(\theta_x)_{x\in p\Z^2}$ for some $p>1$.
\begin{prop}\label{prop:MM} Any periodic Gibbs measure $\mu$ for
the Ising ferromagnet on $\L=\Z^2$ with no external field and inverse
temperature $\beta>\b_c$ is a mixture of the two phases $\mu_\b^+$ and
$\mu_\b^-$.
\end{prop}
Under the condition of translation invariance, this proposition was first
derived for large $\b$ by Gallavotti and Miracle-Sole \cite{GalMS}, and
later
for all $\b>\b_c$ by Messager and Miracle-Sole \cite{MesMS} using some
specific
correlation inequalities; it follows also from the Onsager-formula
for the free energy density and a result of Lebowitz \cite{Leb}. We will
give a
geometric proof below.

Remarkably enough, one can go one step further: each (not necessarily
periodic) Gibbs measure for the Ising model on the square
lattice is a mixture of the plus-phase and the minus-phase, and thus
automatically automorphism invariant. This beautiful
result was obtained independently by Aizenman \cite{A} and Higuchi
\cite{Hi} based on the work of Russo \cite{Rus}. For more general
two-dimensional systems the absence of non-translation-invariant Gibbs
measures at sufficiently low temperatures was proved in \cite{DOS}. In
three or more dimensions, however, non-translation invariant phases of the
Ising model do exist; this is a famous result of Dobrushin \cite{D72}, see
also \cite{vBe} for a short proof.
\begin{thm}\label{thm:AH}
{\bf (Aizenman--Higuchi)}
For the Ising ferromagnet on $\L=\Z^2$ with no external field and inverse
temperature $\beta>\b_c$,  $\mu_\b^+$ and $\mu_\b^-$ are the only phases,
and any other Gibbs measure is a mixture of these two.
\end{thm}
The proof is a masterpiece of random-geometric analysis of equilibrium
phases and contains various ingenious ideas, but unfortunately it is too
long to be sketched here. For the full result we thus need to refer to
the original papers cited above, as well as to the survey \cite{A1}.
However,
to provide an idea of some of the geometric ideas involved we will now
give
a (new) geometric proof of Proposition \ref{prop:MM}. This proof resulted
from discussions of the first author with Y.\ Higuchi.

In this proof we need to consider infinite clusters in halfplanes. Here,
we
say that a set ${\cal H}\subset\Z^2$ is a halfplane if ${\cal H}$ is a
translate of either the upper halfplane $\{x=(x_1,x_2)\in\Z^2:x_2\ge0\}$
or
its complement, the lower halfplane, or a translate of the right and left
halfplanes which are similarly defined. The next lemma provides a first
step
in the proof of Proposition \ref{prop:MM}.
\begin{lem}\label{lem:halfplanes}
Consider the Ising ferromagnet on $\Z^2$, and let $D$ be the event that
for
at least one halfplane ${\cal H}$ in $\Z^2$, both ${\cal H}$ and ${\cal
H}^c$
contain an infinite cluster of the same sign. Then $\mu(D)=1$ for all
$\mu\in{\cal G}(\b H)$ and $\b>\b_c$.
\end{lem}
{\bf Proof: } Since each Gibbs measure is a mixture of extremal Gibbs
measures,
we only need to show that $\mu(D)=1$ for any extremal $\mu$. Suppose the
contrary. Since $D$ is tail measurable, it then follows that $\mu(D)=0$
for
some extremal $\mu$. We will show that this is impossible.

{\em Step 1: }Let ${\cal H}$ be any halfplane, $r$ the reflection of
$\Z^2$
which maps ${\cal H}$ onto ${\cal H}^c$, and $\tau:\s\to-\s$ the spin flip
on $\Om$. We show that $\mu=\mu\circ r\circ\tau$. Since $\mu(D)=0$, at
least
one of the halfplanes ${\cal H}$ and ${\cal H}^c$ contains no infinite
minus-cluster, and this or the other halfplane contains no infinite
plus-cluster. In view of the tail triviality of $\mu$, we can assume that
${\cal H}$ contains no infinite minus-cluster $\mu$-almost surely. Hence,
for any given $\De\in\E$ and $\mu$-almost every $\xi\in\Om$, there exists
an $r$-symmetric region $\Gamma(\xi)\in\E$ such that
$\Gamma(\xi)\supset\De$ and $\xi\equiv1$ on $\partial\Gamma(\xi)\cap{\cal
H}$.
The last property implies that
$\xi\succeq r\circ\tau(\xi)$ on $\partial\Gamma(\xi)$, and using Lemma
\ref{lem:boundary_domination} and the flip-reflection symmetry of $H$ we
find that
\[
\mu_{\b,\Gamma(\xi)}^\xi \geqd \mu_{\b,\Gamma(\xi)}^{r\circ\tau(\xi)} =
\mu_{\b,\Gamma(\xi)}^\xi\circ r\circ \tau \mbox{ on } \F_\De\;.
\]
Assuming that $\Gamma(\xi)$ is maximal in a large box $\La\supset\De$,
we can apply the Markov property of $\mu$ in the same way as in the
proof of Theorem \ref{thmrusso}. This yields that $\mu\geqd\mu\circ
r\circ\tau$
on $\F_\De$ for any $\De$, and thus $\mu\geqd\mu\circ r\circ\tau$.
(The preceding argument is a variant of an idea of Russo \cite{Rus}.)
Using the absence of infinite plus-clusters in ${\cal H}$ or ${\cal H}^c$
we find analogously that $\mu\leqd\mu\circ r\circ\tau$. Hence
$\mu=\mu\circ r\circ\tau$ as claimed.

{\em Step 2: } Here we use a variant of Zhang's argument which was
explained in the proof of Theorem \ref{thm:GKR}. To begin, we observe
that the composition of two reflections in parallel axes is a translation.
Step 1 therefore implies that $\mu$ is periodic with period $2$. The
flip-reflection symmetry of $\mu$ implies further that $\mu$ is different
from $\mu_\b^+$ and $\mu_\b^-$, so that (by Theorem
\ref{thm:Ising_plus_agreep})
there exist both an infinite plus- and an infinite minus-cluster
$\mu$-almost surely. By the Burton--Keane uniqueness
theorem \ref{thm:Burton_Keane}, these infinite clusters are almost
surely unique. We now choose a square $\La=[-n,n-1]^2\cap\Z^2$ so large
that
$\mu(\La\stackrel{+}{\longleftrightarrow}\infty)>1-10^{-3}$. Since $\mu$
is extremal, $\mu$ has positive correlations. By the argument leading to
\rf{eq:A-events} we thus obtain that $\mu(\partial_k\La
\stackrel{+}{\longleftrightarrow}\infty)>1-10^{-3/4}$ for some
$k\in\{1,\ldots,4\}$, where $\partial_k\La$ is the intersection of
$\partial\La$ with the $k$'th quadrant (relative to the axes
$\{x_2=-1/2\}$
and $\{x_1=-1/2\}$). For definiteness, we assume that $k=1$. By the
flip-reflection symmetry, it follows that the intersection
\[
\{\partial_1\La \stackrel{+}{\longleftrightarrow}\infty,\
\partial_2\La \stackrel{-}{\longleftrightarrow}\infty,\
\partial_3\La \stackrel{+}{\longleftrightarrow}\infty,\
\partial_4\La \stackrel{-}{\longleftrightarrow}\infty\}
\]
has probability at least $1 - 4\cdot 10^{-3/4}>0$, which is impossible
because of the uniqueness of the infinite clusters. This contradiction
concludes the proof of the lemma.
$\Cox$

\medskip\noindent
{\bf Proof of Proposition \ref{prop:MM}: } Let $\mu$ be any Gibbs
measure invariant under $(\theta_x)_{x\in p\Z^2}$ for some $p>1$.
Using the ergodic decomposition, we can assume that $\mu$ is in fact
ergodic
with respect to this group of translations. By Lemma \ref{lem:halfplanes},
there exists a pair $({\cal H},{\cal H}^c)$ of halfplanes such that, with
positive probability, both ${\cal H}$ and ${\cal H}^c$ contain infinite
clusters of spins of the same constant sign. For definiteness, suppose
${\cal H}$ is the upper halfplane, and the sign is plus. In view of the
finite energy property, it then follows that also $\mu(A_0)>0$, where for
$k\in\Z$
\[
A_k=\{(k,0)\stackrel{+}{\longleftrightarrow}\infty\mbox{ both in ${\cal
H}$
and ${\cal H}^c$}\}\;.
\]
Let $A$ be the event that $A_k$ occurs for infinitely many $k<0$ and
infinitely many $k>0$. The horizontal periodicity and Poincar\'e's
recurrence theorem (or the ergodic theorem) then show that
$\mu(A_0\setminus A)=0$,
and therefore $\mu(A)>0$.

Next, let $B$ be the event that there exists an infinite minus-cluster.
We claim that $\mu(A\cap B)=0$. Indeed, suppose $\mu(A\cap B)>0$. Since
$A$
is tail measurable and horizontally periodic, we can use the finite
energy property
and horizontal periodicity of $\mu$ as above to show that the event
\[
C=A\cap \{(k,0)\stackrel{-}{\longleftrightarrow}\infty\mbox{ for
infinitely
many $k<0$ and infinitely many $k>0$}\}
\]
has positive probability. But on $C$ there exist infinitely many
minus-clusters,
which is impossible by the Burton--Keane theorem.

To complete the proof, we note that $\mu(B)\le \mu(A^c)<1$, and thus
$\mu(B)=0$ by ergodicity. In view of Theorem \ref{thm:Ising_plus_agreep},
this means that $\mu=\mu_\b^+$. In the case considered, the proposition is
thus proved. The other cases are similar; in particular, in the case of
negative sign we find that $\mu=\mu_\b^-$.
$\Cox$

\subsection{Constant-spin clusters in the Potts model}

Consider the $q$-state Potts model on the lattice $\L=\Z^d$ introduced
in Section \ref{sect:Potts}, $q,d\geq 2$, and recall the results of
Section
\ref{sect:PT_Potts} on the phase transition in this model. The periodic
ground states are the constant configurations $\e_i \equiv i$, $1\le i\le
q$.  We write
$\stackrel{i}{\longleftrightarrow}$ for the agreement connectivity
relation relative to $\e_i$, and we consider
the limiting Gibbs measure $\mu^i_{\beta,q}$ at inverse
temperature $\b$ associated to $\e_i$, which exists by Proposition
\ref{prop:free_and_ordered_Potts_state}. As a further illustration of
assertion
\rf{hypoagree}, we show that  $\mu^i_{\beta,q}$ exhibits $i$-percolation
whenever there is a phase transition. This is a Potts-counterpart of
Theorem  \ref{thm:Ising_plus_agreep}. For its proof, we use the
random-cluster
representation rather than Theorem \ref{thmrusso} because for $q>2$ there
is no stochastic monotonicity available in the spin configuration.
\begin{thm}\label{thm:Potts_i_agreep}
For the Potts model on $\Z^d$ at any inverse temperature $\b$ with
$|{\cal G}(\b H)|>1$, we have
$
\mu^i_{\b,q}(x\stackrel{i}{\longleftrightarrow}\infty)>0
$
for all $x\in\Z^d$ and $i\in \{1,\ldots,q\}$.
\end{thm}
{\bf Proof: } By translation invariance we can choose $x=0$. In
Theorem \ref{thm:Potts_sandwich} we have seen that
$\phi^1_{p,q}(0\leftrightarrow \infty)=c>0$ for $\b>\b_c$, where
$p=1-e^{-2\b}$ as usual. In view of \rf{limit_interchange}, this
means that $\phi^1_{p,q,\La}(0\leftrightarrow \La^c)\ge c$ for all
$\La\ni0$. But for the Edwards-Sokal coupling $P^{i}_{p,q,\La}$ of
$\mu_{\b,q,\La}^i$ and  $\phi^1_{p,q,\La}$ (defined before
Proposition \ref{prop:free_and_ordered_Potts_state}) we have
$\{0\leftrightarrow \La^c\}\subset
\{0\stackrel{i}{\longleftrightarrow}\La^c\}$ almost surely, so that
$\mu_{\b,q,\La}^i(0\stackrel{i}{\longleftrightarrow}\La^c)\ge c$. In
particular,
$\mu_{\b,q,\La}^i(0\stackrel{i}{\longleftrightarrow}\Delta^c)\ge c$
whenever $0\in\Delta\subset\La$. Letting first $\La\uparrow\Z^d$ and then
$\Delta\uparrow\Z^d$ we find that
$\mu^i_{\b,q}(0\stackrel{i}{\longleftrightarrow}\infty)\ge c$, and the
theorem follows. $\Cox$

\medskip\noindent
Next we ask for a converse stating that ``agreement percolation implies
phase transition''. As
we already noticed in the case of the Ising model, this can
be expected to hold only in the case of a planar lattice. But then a
counterpart of Corollary \ref{cor:CNPR} does indeed hold, as was shown by
L.\ Chayes \cite{Lincoln}.
\begin{thm} \label{thm:CNPR_Chayes}
For the unique Gibbs measure $\mu_{\b,q}$ of the $q$-state Potts model on
the square lattice $\Z^2$ at inverse temperature $\b<\b_c$, we have
$
\mu_{\b,q}(\exists\mbox{ \rm an infinite $i$-cluster})=0
$ for all $i\in \{1,\ldots,q\}$.
\end{thm}
The strategy of proving this theorem is the same as that in the proof of
Corollary \ref{cor:CNPR}.
Suppose the $i$-spins percolate in $\mu_{\b,q}$ for some $i$. Then, by
symmetry, this holds for all $i$, so that in particular the $1$-spins
{\em and} the other spins percolate.
Hence, Theorem \ref{thm:GKR} leads to a contradiction, provided we can
show
that the set of $1$'s has positive
correlations. Theorem \ref{thm:CNPR_Chayes} thus follows
from the following lemma.
\begin{lem}
Consider the phase $\mu^i_{\b,q}$ of the $q$-state Potts model at any
inverse temperature $\b>0$, and let the mapping $s_i$ be defined by
\rf{a/dis} with $\e=\e_i$, $i\in \{1,\ldots,q\}$. Then the measure
$\nu_{\b,q}=\mu^i_{\b,q}\circ s_i^{-1}$ has positive correlations.
\end{lem}
{\bf Sketch of proof: } By symmetry, $\nu_{\b,q}$ does not depend on $i$.
For
definiteness we set $i=1$ in the following.
Since the property of positive correlations is preserved under weak
limits,
it is sufficient to consider the finite volume Gibbs distribution
$\mu^1_{\b,q, \La}$ and its image $\nu_{\b,q, \La}=
\mu^1_{\b,q, \La}\circ s_1^{-1}$.
By the FKG inequality, Theorem \ref{thm:FKG}, it is further sufficient
to show that $\nu_{\b,q,\La}$ is monotone. In terms of $\mu^1_{\b,q,
\La}$ and the random field $Y=s_1(X)$, this means that the conditional
probability
\[
q_x(\xi)=\mu^1_{\b,q,\La}(Y(x)=1 \, | \, Y\equiv\xi \mbox{ off }x)
\]
is increasing in $\xi\in\{0,1\}^{\Z^d}$ for any $x\in \La$. Since the
boundary
condition is fixed to be equal to 1 off $\La$, we can assume that
$\xi$ is equal to $1$ off $\La$, and it is sufficient to prove the
inequality $q_x(\xi)\le q_x(\xi')$ for any two such $\xi,\,\xi'$ that
differ
only at a single site $y\in \La$ and are such that
$\xi(y)=0$ and $\xi'(y)=1$. For such $\xi,\xi'$, the inequality
$q_x(\xi)\le q_x(\xi')$ simply means that $Y(x)$ and $Y(y)$
are positively correlated under the conditional distribution
$\mu_{x,y|\xi}$ of
$\mu^1_{\b,q,\La}$ given that $Y\equiv\xi$ off $\{x,y\}$.

To show this we fix $x,\,y,\,\xi$. For $\mu_{x,y|\xi}$, we have
$X\equiv1$ on the complement of $\Delta=\{x,y\}\cup\{v\in
\La\setminus\{x,y\}
:\, \xi(v)=0\}$. We thus consider the graph $G$ with vertex set $\Delta$
and
edge set $\B(\Delta)$ consisting of all edges of $\B$ with both endpoints
lying in
$\Delta$. If we knew that $Y(x)=Y(y)=0$, then
$\mu_{x,y|\xi}$ would be the distribution of a $(q-1)$-state Potts model
on $G$ with state space $\{2,3, \ldots, q\}$.
Now that we don't know $Y(x)$ and $Y(y)$, $\mu_{x,y|\xi}$ is still a
modification of this $(q-1)$-state Potts model,
in which $x$ and $y$ are allowed to have the $q$'th spin value $1$.

To describe this modification we suppose first that $x$ and $y$ are not
adjacent. Let $n_x$ be the number of neighbors
$v$ of $x$ with $\xi(v)=1$, and define $n_y$ accordingly.
The probability weight of $\mu_{x,y|\xi}$ then contains
the additional biasing factor
$$\exp[2\b (n_x I_{\{X(x)=1\}} + n_y I_{\{X(y)=1\}})]$$ which
acts like an external field at $x$ and $y$. For this modified Potts model,
we can still define a modified
random-cluster representation which gives any edge configuration
$\zeta\in\{0,1\}^{\B(\Delta)}$ a probability proportional to
\[
(q-1)^{k(\zeta)}
(q-1+e^{2\beta n_x})^{k_x(\zeta)} (q-1+e^{2\beta n_y})^{k_y(\zeta)}
\prod_{e\in \B(\Delta)} p^{\zeta(e)} (1-p)^{1-\zeta(e)} \,.
\]
Here $p= 1-e^{-2\beta}$, $k(\zeta)$ is the number of connected components
{\em excluding singletons at $x$ or $y$}, and $k_x(\zeta)$ and
$k_y(\zeta)$ are the
indicator functions of having a singleton connected component at $x$
resp.\ $y$. A spin configuration with distribution $\mu_{x,y|\xi}$ is then
obtained from the edge configuration by assigning spins at random
uniformly
from $\{2, \ldots, q\}$ to connected components, except for a singleton at
$x$,
where the spin is taken from $\{1, \ldots, q\}$ with probabilities
proportional to $(e^{2\beta n_x}, 1,\ldots, 1)$, and similarly for
a singleton at $y$. Just as in Corollary \ref{cor:Potts_pos_cor}, this
representation gives the desired
positive correlation of $Y(x)$ and $Y(y)$ under $\mu_{x,y|\xi}$,
provided we can show that $k_x$ and $k_y$ are positively
correlated in the modified random-cluster model. Since these indicator
variables are decreasing, it suffices to check that
the modified random-cluster model has positive correlations, which
follows from Theorem \ref{thm:FKG} by verifying that it is
monotone; this, however, is similar to Lemma \ref{lem:single_edge}.

The case when $x$ and $y$ are neighbors is handled analogously;
in fact, the positive correlation can only become stronger when $x$ and
$y$
have an edge in common.
$\Cox$

\subsection{Further examples of agreement percolation}
\label{subsect:appl_other_models}

Here we treat the Ising antiferromagnet, the hard-core lattice gas, and
the Widom--Rowlinson lattice model, and shortly mention the
Ashkin--Teller model.

\medskip\noindent
{\em The Ising antiferromagnet}. Consider the setting of
Section \ref{sect:anti}. We need to assume that the underlying lattice
$\L$
is bipartite, and thus splits off into two parts, $\L_{even}$ and
$\L_{odd}$.
If $|h|<2d$, there exist two periodic ground
states, $\e_{even}$ and $\e_{odd}=-\e_{even}$, where
$\e_{even}\equiv1$ on $\L_{even}$ and $\e_{even}\equiv-1$ on $\L_{odd}$.
(There are no other periodic ground states,
see for example \cite{DKS}.) The phase transition in
this model has been studied in \cite{D2} and \cite{Hei}. Because of
the bipartite structure, we can flip all spins on a sublattice as in
(\ref{flipping_the_odd_lattice}), which turns the model into an Ising
ferromagnet in a staggered magnetic field of alternating sign on
$\L_{even}$
and $\L_{odd}$. The latter model still satisfies the FKG
inequality. As pointed out in
Section \ref{sect:anti}, for $h=0$ there is a one-to-one correspondence
between all Gibbs measures for the Ising ferromagnet and the Ising
antiferromagnet. In particular, both models then have the same
critical inverse temperature $\b_c$. For general $|h|<2d$, we still have
two limiting Gibbs measures $\mu_\b^{\e_{even}}$ and $\mu_\b^{\e_{odd}}$,
and  these measures have positive correlations relative to the
``staggered'' ordering
$\s\preceq \s'$ iff $\s(x)\e_{even}(x) \le \s'(x) \e_{even}(x)$ for all
$x \in \L$. Relative to this ordering, an analogue of the sandwiching
inequality (\ref{eq:Ising_sandwich}) holds; for more details see Section 9
of \cite{Pr1}. Here is a version of statement \rf{hypoagree} for this
model.
\begin{thm} \label{antiper}
Consider the Ising antiferromagnet on a bipartite graph $(\L,\sim)$ in
an external field $h$ at any inverse
temperature $\b>0$. If $|{\cal G}(\b H)|>1$, we have
\[
\mu_\b^{\e_{even}}(x\stackrel{\e_{even}}{\longleftrightarrow}\infty)>0
\]
for all $x\in\L$.
\end{thm}
This follows from Theorem \ref{thmrusso} in the same way as
Theorem \ref{thm:Ising_plus_agreep}. For $\L=\Z^2$, the obvious
counterparts of Corollary \ref{cor:CNPR} and Proposition \ref{prop:MM}
are also
valid since the proofs of these results carry over to the
case of a staggered external field.

\medskip\noindent
{\em The hard-core lattice gas. } As we have seen in Section
\ref{sect:hard-core}, this model has state space $S=\{0,1\}$ and
corresponds
to setting $U(a,b)=\infty I_{\{a=b=1\}}$ and $V(a)=-a\log\la$ in
\rf{eq:Hamiltonian}, $a,b\in S$; $\la>0$ is an activity parameter.
The hard-core model is the limit of the Ising antiferromagnet for
$\b\rightarrow \infty$ and $h\rightarrow 2d$ along
$\beta (2d-h)=\frac12\log\la$, provided a configuration $\s\in \{-1, +1
\}^{\L}$ is mapped to $(1-\s)/2 \in \{ 0,1 \}^{\L}$; see \cite{DKS} for
details. (The phase diagram point $h=2d,\beta =+\infty$ of the
Ising antiferromagnet is highly degenerate since there are
infinitely many, in general nonperiodic, ground states.) For
$\L=\Z^d$, the hard-core lattice gas can be seen as a gas of hard
(i.e., non-overlapping) diamonds. In general, we still assume that
$\L$ is bipartite. For $\la>1$, the hard-core model then has two
periodic ground states of chessboard type, namely $\e_{even}$ which
is equal to 1 on $\L_{even}$ and 0 otherwise, and
$\e_{odd}=1-\e_{even}$. As noticed in Section
\ref{sect:other_appl}, the associated limiting Gibbs states
$\mu_\lambda^{even}$ and $\mu_\lambda^{odd}$ exist. So, following
the program stated in \rf{hypoagree}, we may ask whether these
Gibbs measures exibit agreement percolation in the case of phase
transition. The answer is again positive:
\begin{thm} \label{perha}
For the hard-core model on a bipartite graph $\L$ we have for any
activity $\la>0$: If $\mu_\lambda^{even}\ne\mu_\lambda^{odd}$ then
$
\mu_\lambda^{even}(\,x\stackrel{\e_{even}}{\longleftrightarrow}\infty)>0
$
for all $x\in\L$, and similarly with `odd' in place of `even'.
\end{thm}
This result is completely analogous to Theorem \ref{antiper}, and was
conjectured by Hu and Mak \cite{Hu,Hu1} from
computer simulations. In these papers, the authors also
discuss the case of hard-core particles on a triangular lattice,
the hard hexagon model. While Theorem \ref{perha} does apply to
the hard triangle model on the hexagonal lattice (which is bipartite),
the non-bipartite triangular lattice with nearest-neighbor bonds is
excluded. The results of \cite{Hu,Hu1} suggest that Theorem \ref{perha}
still holds for the triangular lattice. A geometric proof of this
conjecture would be of particular interest.

The hard-core model on the square lattice $\Z^2$ admits an analogue
to Corollary \ref{cor:CNPR}, in that nonuniqueness of the Gibbs
measure is equivalent to $\eta_{even}$-percolation for the Gibbs
measure $\mu_\lambda^{even}$; see \cite{GLM} or \cite{H3} for more
details.

\medskip\noindent
{\em The Widom--Rowlinson lattice model. } Consider the set-up of
Section \ref{sect:WR}, with equal activities $\la_+=\la_-=\la>0$
for the plus and minus particles. For $\la>1$ we have two distinct
periodic ground states $\e_+\equiv +1$ and $\e_-\equiv -1$. From
Section \ref{sect:other_appl} we know that the associated limiting
Gibbs measures $\mu_\la^+=\lim_{\La\uparrow\L}\mu_{\la,\La}^{\e_+}$
and $\mu_\la^-$ exist. Moreover, Theorem \ref{WR_sandwich} asserts
that a phase transition occurs for some activity $\la$ if and only
if $\mu_\la^+(X(x)=1)>\mu_\la^+(X(x)=-1)$ for some $x\in\L$. Now,
it turns out that in this model not only hypothesis
\rf{hypoagree} holds, but that the nonuniqueness
of the Gibbs measure is in fact equivalent to agreement
percolation, not only for the square lattice but {\em for any
graph}. This comes from the nature of the random-cluster
representation of Section \ref{sect:WR_RC}, which is related to the
sites rather than the bonds of the lattice, and is a curious
exception from the fact that, on the whole, the Widom--Rowlinson
model is less amenable to sharp results than the Ising model.
However, by the reasons discussed in Section \ref{sect:WR_RC}, this
result does {\em not} carry over to the multitype Widom--Rowlinson
lattice model with $q\ge3$ types of particles.
\begin{thm} \label{WR_agree}
Consider the Widom--Rowlinson lattice model on an arbitrary graph
$(\L,\sim)$ for any activity $\la>0$. Then the following statements
are equivalent.
\begin{description}
\item{\rm (i)} The Gibbs measure for the parameter $\la$ is non-unique.
\item{\rm (ii)} $\mu_\la^+(x\stackrel{\e_+}{\longleftrightarrow}\infty)>0$
for some, and thus all $x\in\L$.
\end{description}
\end{thm}
{\bf Sketch of Proof: } Consider $\mu_{\la,\La}^{\e_+}$ for some
finite $\La$. In the same way as the random-cluster representation
of Section \ref{sect:FK} was modified in Section
\ref{sect:Inf_vol_limits} to deal with boundary conditions, we can
modify the site-random-cluster representation of Section
\ref{sect:WR_RC} to obtain a coupling of $\mu_{\la,\La}^{\e_+}$ and
a wired site-random cluster distribution $\psi^1_{p,2,\La}$, so
that analogues of Propositions \ref{prop:RC_to_WR} and
\ref{prop:WR_to_RC} hold. As a counterpart to  equation
\rf{eq:Potts_magnetization_connectivity} and by the
specific nature of the site-random-cluster representation, we then
find that
\[
\mu_{\la,\La}^{\e_+}(X(x)=1)-\mu_{\la,\La}^{\e_+}(X(x)=-1)=
\psi^1_{p,2,\La}(x\leftrightarrow\partial\La)
= \mu_{\la,\La}^{\e_+}(x\stackrel{\e_+}{\longleftrightarrow}\partial\La)
\]
for all $x\in\L$. In the limit $\La\uparrow\L$ we obtain by an
analogue to \rf{limit_interchange}
\[
\mu_\la^+(X(x)=1)-\mu_\la^+(X(x)=-1)=
\mu_\la^+(x\stackrel{\e_+}{\longleftrightarrow}\infty)\;,
\]
and the theorem follows immediately. $\Cox$

\medskip\noindent
To conclude this subsection, we note that hypothesis \rf{hypoagree} also
holds in other models. We mention here only the
{\em Ashkin--Teller model} \cite{AT}, a
$4$-state model which interpolates in an interesting way between the
$4$-state Potts and the so called $4$-state clock model, which is also
accessible to random-cluster methods; we refer to \cite{CM,CMW,PVe,SS}.

\subsection{Percolation of ground-energy bonds}
\label{ssect:ground_energy_perc}

So far in this section we considered a number of models which are
known to show a phase transition, and asked whether this phase
transition goes hand in hand with agreement percolation. These results
run under the heading ``phase transition implies percolation'', even
though for the square lattice we established results of converse type
coming from the impossibility of simultaneous occupied and vacant
percolation on $\Z^2$.

We now take the opposite point of view and ask whether ``percolation
implies phase transition''. More precisely, we want to deduce the
existence of a phase transition (at low temperatures or high densities)
from a percolation result.
In fact, such an idea is already implicit in Peierls' \cite{Pe} and
Dobrushin's \cite{D1} proof of phase transition in the Ising model,
and is an integral part of the Pirogov--Sinai theory. For models
with neighbor interaction as in the Hamiltonian \rf{eq:Hamiltonian},
the underlying principle can be sketched as follows. At low temperatures
(or high densities), each pair of adjacent spins (or particles) tries
to minimize its pair interaction energy. Note that this minimization
involves the bonds rather than the sites of the lattice. So, one
expects that bonds of minimal energy -- the {\em ground-energy bonds}
 -- prevail, forming regions separated by boundaries that consist of bonds
of higher energy. Such boundaries, which are known as {\em contours},
cost an energy proportional to their size, and are therefore typically
small when $\b$ is large. This implies that the ground-energy bonds should
percolate. Now, the point is that if the spins along a bond can choose
between different states of minimal energy then this ambiguity can be
transmitted to the macroscopic level by an infinite ground-energy cluster,
and this gives rise to phase transition.
In other words, the classical contour argument for the existence of phase
transition can be summarized in the phrase: ground-energy bond percolation
together with a (clear-cut) non-uniqueness of the local ground state
implies
non-uniqueness of Gibbs measures. We will now describe this picture
in detail.

We consider the cubic lattice $\L=\Z^d$ of dimension
$d\ge2$ with its usual graph structure. For definiteness we consider the
Hamiltonian
\rf{eq:Hamiltonian} for some pair potential $U:S\times S\to\R$. We can and
will assume that the self-potential $V$ vanishes; this is because
otherwise we can replace $U$ by
\begin{equation}\label{eq:U-transform}
U'(a,b)=U(a,b)+\frac1{2d}\,[V(a)+V(b)]\;,\quad a,b\in S,
\end{equation}
which, together with the self-potential $V'\equiv0$, leads to the same
Hamiltonian. Let
\begin{equation}\label{eq:minimum}
m=\min_{a,b\in S} U(a,b)
\end{equation}
 be the minimal value of $U$.

Given an arbitrary configuration $\s\in\Om$, we will say that an edge
$e=\{x,y\}\in\B$ is a {\em ground-energy bond} for $\s$ if
$U(\s(x),\s(y))=m$. The subgraph of $\Z^d$ consisting of all vertices of
$\Z^d$ and only the ground-energy bonds for $\s$  splits then
off into connected components which will be called {\em ground-energy
clusters} for $\s$. We are interested in the existence of infinite
ground-energy clusters, and we also need to identify specific such
clusters. Unfortunately, the Burton--Keane uniqueness theorem
\ref{thm:Burton_Keane} does not apply here because, for any Gibbs
measure,
the distribution of the set of ground-energy bonds fails to have the
finite-energy property. We therefore resort to considering ground-energy
clusters in any fixed two-dimensional layer of $\Z^d$; the uniqueness of
planar infinite clusters can be shown in our case. (An alternative
argument avoiding the use of planar layers but requiring stronger
conditions on the temperature has been suggested in
\cite{Fuku}.) In fact, we have the following result.
\begin{thm}\label{thm:ground_energy_perc}
Consider the Hamiltonian \rf{eq:Hamiltonian} on the lattice $\L=\Z^d$,
$d\ge2$, with neighbor interaction $U$ and no self-potential, and let
${\cal P}$ be any planar layer in $\L$. (So ${\cal P}=\L$ for $d=2$.)
If $\b$ is large enough, there exists a Gibbs measure $\mu\in{\cal G}(\b H
)$ which is invariant under all automorphisms of $\L$ and all symmetries
of $U$ such that
\[
\mu(\exists\mbox{ \em a unique infinite ground-energy cluster in ${\cal
P}$})=1\;.
\]
In the above, a symmetry of $U$ is a transformation $\tau$ of $S$ such
that $U(\tau a,\tau b)=U(a,b)$ for all $a,b\in S$; such a $\tau$ acts
coordinatewise on configurations.
\end{thm}
Theorem  \ref{thm:ground_energy_perc} is a particular case of a result
first derived in \cite{Ge1} and presented in detail in Chapter 18 of
\cite{Geo}. We will sketch its proof below. The remarkable fact is that
this type of percolation often implies that $\mu$ has a non-trivial
extremal decomposition, so that there must be a phase transition. In fact,
this happens whenever the set
\begin{equation}\label{eq:bond_ground_states}
G_U=\{(a,b)\in S\times S:U(a,b)=m\}
\end{equation}
of bond ground states splits into sufficiently disjoint parts. To explain
the underlying mechanism (which may be viewed as the core of the classical
Peierls argument, and a rudimentary version of Pirogov--Sinai theory)
we consider first the standard Ising model.
\begin{example} The Ising ferromagnet at zero external field.
{\rm In this model, we have as usual $S=\{-1,1\}$, $U(a,b)=-ab$ for
$a,b\in S$, $m=-1$, and
$G_U=\{(-1,-1),(1,1)\}$. Hence, either all spins of a ground-energy
cluster are negative, or else all these spins are positive. In other
words,
each ground-energy cluster is either a minus-cluster or a plus-cluster.
This implies that
\[
\{\exists\mbox{ a unique infinite ground-energy cluster in ${\cal P}$}\}
\subset A_-\cup A_+\;,
\]
where $A_-$ and $A_+$ are the events that there exists an infinite cluster
of negative, resp.\ positive, spins in ${\cal P}$. For the Gibbs measure
$\mu$ of Theorem \ref{thm:ground_energy_perc} we thus have $\mu(A_-\cup
A_+)=1$ and,
by the spin-flip symmetry of $U$ and thus $\mu$, $\mu(A_-)=\mu(A_+)$.
Hence
$\mu(A_-)>0$ and $\mu(A_+)>0$, so that the measures
$\mu^-=\mu(\cdot\,|A_-)$ and
$\mu^+=\mu(\cdot\,|A_+)$ are well-defined. Since $A_-$, $A_+$ are tail
events, it follows immediately that  $\mu^-$, $\mu^+$ are Gibbs measures
for $\b H$. Also,
$A_-\cap A_+$ is contained in the event that there are two distinct
ground-energy clusters in ${\cal P}$, and therefore has $\mu$-measure 0.
Hence $\mu^-$ and $\mu^+$ are mutually singular, whence $|{\cal G}(\b
H)|>1$.}
\end{example}
The same argument as in the preceding example yields the following theorem
on phase transition by symmetry breaking. A detailed proof (in a slightly
different setting) can be found in Section 18.2 of \cite{Geo}.
\begin{thm}\label{thm:symmetry_breaking}
Under the conditions of Theorem \ref{thm:ground_energy_perc} suppose that
the set $G_U$ defined by \rf{eq:bond_ground_states} admits a decomposition
$G_U=G_1\cup\ldots\cup G_N$ such that
\begin{enumerate}
\item the sets $G_n$, $1\le n\le N$, have disjoint projections, i.e., if
$(a,b)\in G_n$, $(a',b')\in G_{n'}$, and $n\ne n'$, then $a\ne a', b\ne
b'$, and
\item for any two indices $n,n'\in\{1,\ldots,N\}$ we have
$\bar\tau(G_n)=G_{n'}$  for some transformation $\bar\tau$ of $S\times S$
which is either the reflection, or the coordinatewise application of some
symmetry of $U$, or a composition of both.
\end{enumerate}
Then, if $\b$ is sufficiently large, there exist $N$ mutually singular
Gibbs measures $\mu^1,\ldots\mu^N$ $\in{\cal G}(\b H)$, invariant
under all even automorphisms of $\Z^d$ and such that
\[
\mu^n(\exists\mbox{ \em an infinite $n$-cluster in ${\cal P}$})=1
\]
for all $1\le n\le N$. In particular, there exist $N$ distinct phases for
$\b H$.
\end{thm}
In the statement above, an infinite $n$-cluster for a configuration $\s$
is an infinite cluster of the subgraph of $\Z^d$ obtained by keeping only
those edges
 $e\in\B$ with $(\s(x),\s(y))\in G_n$, where $x$ is the endpoint of $e$ in
the even sublattice $\L_{even}$ and $y\in\L_{odd}$ is the other endpoint
of $e$. Also, an even automorphism of $\Z^d$ is an automorphism leaving
$\L_{even}$ invariant.

We illustrate this theorem by applying it to our other standard examples.
\begin{example}\label{ex:IsingAF}
The Ising antiferromagnet in an external field. {\rm We have again
$S=\{-1,1\}$, but the interaction is now $U(a,b)=ab-\frac h{2d}(a+b)$ for
some constant $h\in\R$. (Here we applied the transformation
\rf{eq:U-transform}.) If $|h|<2d$ then $m=-1$ and
$G_U=\{(-1,1),(1,-1)\}$. $G_U$ splits up into the singletons
$G_1=\{(1,-1)\}$
and $G_2=\{(-1,1)\}$. This decomposition meets the conditions of the
theorem;
in particular, $G_1$ and $G_2$ are related to each other by the reflection
of $S\times S$. Consequently, there exist two mutually singular Gibbs
measures
$\mu^1$ and $\mu^2$ which are invariant under even automorphisms and have
an infinite cluster of chessboard type, either with plus spins on the even
cluster sites and minus spins at the odd cluster sites, or vice versa.
}\end{example}
\begin{example} \label{ex:Potts}
The Potts model. {\rm In this case, $S=\{1,\ldots,q\}$ for some integer
$q\ge 2$ and $U(a,b)=1-2I_{\{a=b\}}$. Again $m=-1$, and $G_U=\{(n,n):1\le
n\le q\}$. Theorem \ref{thm:symmetry_breaking} is obviously applicable,
and we recover
the result that for sufficiently large $\b$ there exist $q$ mutually
singular, automorphism invariant Gibbs measures, the $n$th of which has an
infinite cluster of spins with value $n$.
}\end{example}
\begin{example}
The hard-core lattice gas. {\rm This model has state space $S=\{0,1\}$ and
neighbor interaction $U$ of the form $U(a,b)=\infty$ if $a=b=1$, and
$U(a,b)=
-\frac{\log \la}{2d}(a+b)$ for all other $(a,b)\in S^2$. Here $\la>0$ is
an
activity parameter, and we have again used the transformation
\rf{eq:U-transform}. For $\la>1$ we have $G_U=\{(0,1),(1,0)\}$, so that
Theorem \ref{thm:symmetry_breaking} applies. Since multiplying $U$ with a
factor $\b$ amounts to changing $\la$,
 we obtain that for sufficiently large $\la$ there exist two distinct
Gibbs measures with infinite clusters of chessboard type, just as for the
Ising antiferromagnet at low temperatures.
}\end{example}
\begin{example}
The Widom--Rowlinson lattice model. {\rm Here we have $S=\{-1,0,1\}$ and
$U(a,b)=\infty$ if $ab=-1$, $U(a,b)=-\frac{\log\la}{2d}(|a|+|b|)$
otherwise, $a,b\in S$. If $\la>1$ then $G_U=\{(-1,-1),(1,1)\}$.
Theorem \ref{thm:symmetry_breaking} thus shows that for sufficiently large
$\la$ there exist two translation invariant Gibbs measures having infinite
clusters of plus- resp.\ minus-particles.
}\end{example}
Although the results in the examples above are weaker than those obtained
by the random-cluster methods of Section \ref{sect:random-cluster} (when
these apply), the ideas presented here have the advantage of providing a
general picture of
 the geometric mechanisms that imply a phase transition, and Theorem
\ref{thm:symmetry_breaking} can quite easily be applied. Moreover, the
ideas
can be extended immediately to systems with arbitrary state space and
suitable
 interactions. In this way one obtains phase transitions in anisotropic
plane rotor models, classical Heisenberg ferromagnets or antiferromagnets,
and related $N$-vector models; see Chapter 16
of \cite{Geo}. One can also consider next-nearest neighbor interactions,
and thus obtain various other interesting examples; for this one has to
consider
percolation of ground-energy plaquettes rather than ground-energy bonds,
which is in fact the set-up of \cite{Geo}. Last but not least, the
symmetry assumption of Theorem \ref{thm:symmetry_breaking} can often be
replaced by either some direct argument, or a Peierls condition in the
spirit of the Pirogov--Sinai theory; see Chapter 19 of \cite{Geo}. One
such extension will be used in our next example.
\begin{example}\label{ex:Potts_first_order}
First-order phase transition in the Potts model. {\rm Consider again the
Potts model of Example \ref{ex:Potts}, and suppose for simplicity that
$d=2$. Any translate in $\Z^2$ of the quadratic cell $\{0,1\}^2$ is
called a
plaquette. For a given configuration
$\s\in\Om$, a plaquette $P$ is called {\em ordered\/} if all spins in $P$
agree,
{\em disordered\/} if no two adjacent spins in $P$ agree, and {\em pure\/}
if one of these two cases occurs. If $q$ (the number of distinct spin
 values) is large enough then, for arbitrary $\b$, there exists an
 automorphism invariant Gibbs measure $\mu$
supported on configurations with a unique infinite cluster of
pure plaquettes. This variant of Theorem \ref{thm:ground_energy_perc} is
due to Koteck\'y and Shlosman \cite{KS}, see also Section 19.3.2 of
\cite{Geo}. Clearly, each cluster of pure plaquettes only contains
plaquettes of the same type, either ordered or disordered. For some
specific critical value $\b_c(q)$ both possibilities must occur with
positive probability;
this follows from thermodynamic considerations, namely by convexity of
the free energy as a function of $\b$ \cite{KS,Geo}. Conditioning on each
of these two cases yields two mutually singular Gibbs measures with an
infinite cluster of ordered
resp.\ disordered plaquettes. Furthermore, all spins of a cluster of
ordered
plaquettes must have the same value, so that by symmetry the ``ordered''
Gibbs
measure can be decomposed further into $q$ Gibbs measures with infinite
clusters of constant spin value. As a result,
 for large $q$ and $\b=\b_c(q)$ there exist $q+1$ mutually singular Gibbs
measures which behave qualitatively similar to the disordered phase for
$\b<\b_c(q)$ resp.\ the $q$ ordered phases for $\b>\b_c(q)$. This is the
first-order phase transition in the
Potts model for large $q$. For further discussions we refer to
\cite{Wu,KS,LMR,LMMRS} and the references therein.
}\end{example}
We now give an outline of the proof of Theorem
\ref{thm:ground_energy_perc}.

\medskip\noindent
{\bf Sketch proof of Theorem \ref{thm:ground_energy_perc}:}
For simplicity we stick to the case $d=2$. For any inverse temperature
$\b>0$ and any square box $\La_n=[-n,n-1]^2\cap\Z^2$ we write
$\mu^{per}_{\b,n}$ for the Gibbs distribution relative to $\b H$ in the
box $\La_n$ with periodic boundary condition. The
latter means that $\La_n$ is viewed as a torus, so that
$(i,n-1)\sim(i,-n)$ and $(n-1,i)\sim(-n,i)$ for $i\in[-n,n-1]\cap\Z$; the
Hamiltonian $H_n^{per}$ in $\La_n$ with periodic boundary condition is
then defined in the natural way. Let $\mu^{per}$ be an arbitrary limit
point of the sequence
$(\mu^{per}_{\b,n})_{n\ge1}$. Evidently, $\mu^{per}$ has the symmetry
properties required
of $\mu$ in Theorem \ref{thm:ground_energy_perc}, and $\mu^{per}\in{\cal
G}(\b H)$.

To establish percolation of ground-energy bonds we fix some $\a<1$ and
consider the wedge ${\cal W}=\{x=(x_1,x_2)\in\L:x_1\ge0,\,|x_2|\le\a
x_1\}$.
Let $A_{\cal W}$ be the event that there exists an infinite path of
ground-energy bonds in ${\cal W}$ starting from
the origin. We want to show that $\mu^{per}(A_{\cal W})>3/4$ when $\b$ is
large enough. Suppose $\xi\notin A_{\cal W}$. Then there exists a contour
crossing ${\cal W}$, i.e., a path $\gamma$ in the dual lattice
$\L^*=\Z^2+(\frac12,\frac12)$ which crosses
no ground-energy bond for $\xi$ and connects the two half-lines bordering
${\cal W}$. For each such path $\gamma$ we will establish the contour
estimate
\begin{equation}\label{eq:contour_est}
\mu^{per}(\gamma \mbox{ is a contour})\le (|S|e^{-\b\de})^{|\gamma|}\,,
\end{equation}
where $|\gamma|$ is the length (the number of vertices) of $\gamma$, and
$\de>0$ is such that $m+2\de$ is the second lowest value of $U$.

Assuming \rf{eq:contour_est} we obtain the theorem as follows. The number
of paths of length $k$ crossing ${\cal W}$ is at most $ck\,3^{k}$ for some
$c<\infty$ depending on $\a$. Hence
\[
\mu^{per}(A_{\cal W}^c)\le c \sum_{k\ge1} k\,(3|S|e^{-\b\de})^k<\frac14
\]
for sufficiently large $\b$. By the rotation invariance of $\mu^{per}$, it
follows that $\mu^{per}(A_0)>0$, where $A_0$ is the intersection of
$A_{\cal W}$ with its three counterparts obtained by lattice rotations.
Roughly speaking, $A_0$ is the event that the origin belongs to two
doubly infinite ground-energy paths, one being
quasi-horizontal and the other quasi-vertical. Since $\mu^{per}$ is
invariant under horizontal and vertical translations, the Poincar\'e
recurrence theorem (or the ergodic theorem)
 implies that the event
\[
A_\infty=\{\xi\in\Om:\theta_x\xi\in A_0\mbox{ for infinitely many $x$ in
each of the four half-axes}\}
\]
has also positive $\mu^{per}$-probability. Each configuration in
$A_\infty$ has infinitely many  quasi-horizontal and quasi-vertical
ground-energy paths in each of the four directions of the compass, and by
planarity all paths of different orientation must intersect. Therefore
all these paths
belong to the same infinite ground-energy cluster which has
only finite holes, and is therefore unique. Hence $A_\infty$ is contained
in the event $B$ that there exists a ground-energy cluster surrounding
each finite set of $\L$, and $\mu^{per}(B)>0$.
As $B$ is a tail-event and invariant under all automorphisms
of $\L$ and all symmetries of $U$, the theorem follows by setting
$\mu=\mu^{per}(\,\cdot\,|B)$.

It remains to establish the contour estimate \rf{eq:contour_est}. For this
it is sufficient to show that
\begin{equation}\label{eq:contour_est2}
\mu^{per}_{\b,n}(\gamma \mbox{ is a contour})\le
(|S|e^{-\b\de})^{|\gamma|}
\end{equation}
when $n$ is so large that $\gamma$ is contained in $\La_n$. This bound is
based on reflection positivity and the chessboard estimate, which are
treated at length in Chapter 17 of \cite{Geo}. Here we give only the
principal ideas.
The basic observation is the following consequence of the toroidal
symmetry of $\mu^{per}_{\b,n}$: for any $i\in\{0,\ldots,n-1\}$, the
configurations on the two parts $\La^+_{n,i}=\{x\in\La_n:x_1\ge i\mbox{ or
}x_1\le i-n\}$ and $\La^-_{n,i}=\{x\in\La_n:i
-n\le x_1\le i\}$ of $\La_n$ are conditionally independent and, up to
reflection, identically distributed given the spin values on the two
separating lines $\{x_1=i\}$ and $\{x_1=i-n\}$. Hence, if $f,g$ are real
functions on $S^{\La_n}$ depending only on
the configuration in $\La^+_{n,i}$, and $g^{(i)}$ is the function obtained
from $g$ by reflection in these two separating lines (and thus depending
on the
configuration in $\La^-_{n,i}$), then the bilinear form
$(f,g)\to\mu^{per}_{\b,n}(fg^{(i)})$ is positive definite and thus
satisfies
the Cauchy--Schwarz inequality. Similar
Cauchy--Schwarz inequalities hold for vertical reflections. The chessboard
inequality is obtained by suitable combinations of all these, as we will
illustrate next.

Let us mark the plaquettes around the vertices of $\gamma$ with a
$\bullet$;
this gives $|\gamma|$ marked pla\-quet\-tes. Marking a plaquette
indicates that at least one of its bonds has non-minimal energy; leaving
it unmarked does not say that it consists of ground-energy bonds, but
that we don't need any information on its spins
. In the case $n=2$, this might lead
to the picture
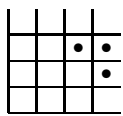
\begin{center}\scriptsize\begin{tabular}[b]{|c|c|c|c}
$\!$&$\!$&$\!$&$\!$\\ \hline
$\!$&$\!$&$\!\!\bullet\!\!$&$\!\!\bullet\!\!$\\ \hline
$\!$&$\!$&$\!\!\!\!$&$\!\!\bullet\!\!$\\ \hline
$\!$&$\!$&$\!\!\!\!$&$\!\!\!\!$\\ \hline
\end{tabular}\quad .\end{center}
To estimate its probability we use repeatedly the Cauchy--Schwarz
inequality
relative to suitable pairs of reflection lines.
Indicating each time only the pair of lines used next, and omitting the
event that no plaquette is marked (which has probability 1), we obtain
\begin{eqnarray*}
&&\mu^{per}_{\b,n} \left(\;\mbox{\scriptsize\begin{tabular}{|cc|cc}
$\!$&$\!$&$\!$&$\!$\\
$\!$&$\!$&$\!\!\bullet\!\!$&$\!\!\bullet\!\!$\\
$\!$&$\!$&$\!$&$\!\!\bullet\!\!$\\
$\!$&$\!$&$\!$&$\!$\\
\end{tabular}}\,\right)
\ \le\ \
\mu^{per}_{\b,n} \left(\;\mbox{\scriptsize\begin{tabular}{cccc}
$\!$&$\!$&$\!$&$\!$\\ \hline
$\!\!\bullet\!\!$&$\!\!\bullet\!\!$&$\!\!\bullet\!\!$&$\!\!\bullet\!\!$\\
$\!\!\bullet\!\!$&$\!$&$\!$&$\!\!\bullet\!\!$\\ \hline
$\!$&$\!$&$\!$&$\!$\\
\end{tabular}}\,\right)^{1/2}
\ \le\ \
\mu^{per}_{\b,n} \left(\;\mbox{\scriptsize\begin{tabular}{cccc}
$\!\!\bullet\!\!$&$\!\!\bullet\!\!$&$\!\!\bullet\!\!$&$\!\!\bullet\!\!$\\
$\!\!\bullet\!\!$&$\!\!\bullet\!\!$&$\!\!\bullet\!\!$&$\!\!\bullet\!\!$\\
\hline
$\!\!\bullet\!\!$&$\!$&$\!$&$\!\!\bullet\!\!$\\
$\!\!\bullet\!\!$&$\!$&$\!$&$\!\!\bullet\!\!$\\ \hline
\end{tabular}}\,\right)^{1/4}
\\[2ex]
&&\ \le\ \
\mu^{per}_{\b,n} \left(\,\mbox{\scriptsize\begin{tabular}{cccc}
$\!\!\bullet\!\!$&$\!\!\bullet\!\!$ & $\!\!\bullet\!\!$&$\!\!\bullet\!\!$
\\
$\!\!\bullet\!\!$&$\!\!\bullet\!\!$ &
$\!\!\bullet\!\!$&$\!\!\bullet\!\!$\\
$\!\!\bullet\!\!$&$\!\!\bullet\!\!$ &
$\!\!\bullet\!\!$&$\!\!\bullet\!\!$\\
$\!\!\bullet\!\!$&$\!\!\bullet\!\!$ & $\!\!\bullet\!\!$&$\!\!\bullet\!\!$
\\
\end{tabular}}\,\right)^{1/8}
\mu^{per}_{\b,n} \left(\,\mbox{\scriptsize\begin{tabular}{c|cc|c}
$\!\!\bullet\!\!$&$\!$&$\!$&$\!\!\bullet\!\!$\\
$\!\!\bullet\!\!$&$\!$&$\!$&$\!\!\bullet\!\!$\\
$\!\!\bullet\!\!$&$\!$&$\!$&$\!\!\bullet\!\!$\\
$\!\!\bullet\!\!$&$\!$&$\!$&$\!\!\bullet\!\!$\\
\end{tabular}}\,\right)^{1/8}
\ \le\ \
\mu^{per}_{\b,n} \left(\,\mbox{\scriptsize\begin{tabular}{cccc}
$\!\!\bullet\!\!$&$\!\!\bullet\!\!$ &
$\!\!\bullet\!\!$&$\!\!\bullet\!\!$\\
$\!\!\bullet\!\!$&$\!\!\bullet\!\!$ &
$\!\!\bullet\!\!$&$\!\!\bullet\!\!$\\
$\!\!\bullet\!\!$&$\!\!\bullet\!\!$ &
$\!\!\bullet\!\!$&$\!\!\bullet\!\!$\\
$\!\!\bullet\!\!$&$\!\!\bullet\!\!$ &
$\!\!\bullet\!\!$&$\!\!\bullet\!\!$\\
\end{tabular}}\,\right)^{3/16}\;.
\end{eqnarray*}
In general, we obtain in this way
\[
\mu^{per}_{\b,n}(\gamma \mbox{ is a contour})\le
\mu^{per}_{\b,n}(C_n)^{|\gamma|/|\La_n|}\;,
\]
where $C_n$ is the event that all plaquettes in $\La_n$ contain at least
one bond of non-minimal energy. But if $C_n$ occurs then at least
$|\La_n|/2$ of
the $2|\La_n|$ edges in $\La_n$ are no ground-energy bonds. The
Hamiltonian $H_n^{per}$ with periodic boundary condition is therefore
at least $(2m+\de)|\La_n|$ on $C_n$. Since there is at
least one $\s\in S^{\La_n}$ with
$H_n^{per}(\s)=2|\La_n|m$, it follows that
\[
\mu^{per}_{\b,n}(C_n)\le \sum_{\xi\in C_n}e^{-\b(2m+\de)|\La_n|}\Big/
e^{-\b 2m|\La_n|}\le (|S|\,e^{-\b\de})^{|\La_n|}\;.
\]
This gives  estimate \rf{eq:contour_est2} and completes the proof of
Theorem \ref{thm:ground_energy_perc}. $\Cox$

\section{Random interactions} \label{sect:rain}

So far in this review, the spin systems considered had an interaction
which was invariant under all automorphisms of the underlying graph
$\L$.
Here we will assume for convenience that $\L$ is the cubic lattice
$(\Z^d,\B)$, $d\geq 2$, but the interaction between adjacent spins
will no longer be translation invariant. That is, instead of the
Hamiltonian \rf{eq:Hamiltonian} we now consider
a modified Hamiltonian of the form
\begin{equation} \label{eq:ranHam}
H(\s)= \sum_{b=\langle xy\rangle\in\B } J_b \;U (\s(x),\s(y)) +
\sum_{x\in\Z^d} h_x\, V(\s(x))\;
\end{equation}
where the $J_b$ and the $h_x$ may vary from bond to bond, resp.\ from
site to site.
In fact, we are interested in the case where these coupling coefficients
show no regular structure, and thus assume that they are {\em
random}. Such systems of spins interacting differently depending on their
position and in a way governed by chance are known as
{\em disordered systems}.  We will not elaborate on the physical
origins of such random interactions. We merely mention that they can
be related to the presence of impurities or defects in an originally
homogeneous system, and are used to model quenched alloys of magnetic
and nonmagnetic materials like {\sl FeAu}.
For details we refer to \cite{Fr,BY,FHe}.

We assume that the family $\J=(J_b)_{b\in\B}$ of
coupling coefficients and the external fields $\h=(h_x)_{x\in \Z^d}$
are independent,  and each collection constitutes a family of mutually
independent and identically distributed real random variables.
Hence, while no realization of the coupling coefficients is translation
invariant, we still have translation invariance in a statistical sense.
We will not specify the underlying probability space,
except that the letter $P$ will be used to denote the probability
measure and the associated expectation.  The random families $\J$ and $\h$
are often referred to as the
{\em disorder}. The disorder is called bounded if $P(|J_b| >
c) = 0$ for some finite $c$. Physically, this is the
most relevant case.

(In some physical applications it is natural to assume that the
$J_b$  are not independent but rather have some finite-range dependence
structure, but we will not include this case here. We also assumed
for simplicity that the disorder is real valued, although some of
the following also applies to the case when $J_b$ or $h_x$ are
allowed to take the value $+\infty$ with positive probability.)

In Section \ref{sect:dilute} we will discuss {\em diluted
ferromagnets}. A bond-diluted Ising or Potts ferromagnet on $\Z^d$
can be viewed as an Ising or Potts model on the open clusters for
Bernoulli bond percolation on $\Z^d$. As observed in \cite{ACCN1},
these models can quite easily be understood using their
random-cluster representation. They form just about the only class of
disordered systems where the phase transition can be investigated
in such detail.

Then, in Section \ref{sect:G-regime}, we study the so-called {\em
Griffiths regime} which is the only non-trivial regime for
disordered systems where by now quite general results are
available. It occurs at intermediate temperatures if the disorder
is bounded, or at arbitrary high temperatures if the disorder is
unbounded, and is characterized by the fact that the Gibbs measure
is still unique but fails to have a nice high temperature behavior
uniformly in the disorder. The study of random Gibbs measures in
the Griffiths regime has started in the early 1980's and has
reached a satisfactory stage only recently. The simplest and also
most powerful methods use stochastic-geometric representations and
will be presented here.

As representative for the  large literature on the subject we refer
to \cite{AKN1,BD,Be,CK,DKP,Fr,FI,FZ,GM1,GrL,Per,Pet}. Dynamical
problems (which are not touched upon here) are treated e.g.\ in
\cite{CM1,CM2,CM3,AC1,GMV,GZ2,Gie}.

\subsection{Diluted and random Ising and Potts ferromagnets}
\label{sect:dilute}

The random Potts model is defined as follows. Spins take values in
the state space $S=\{1,\ldots,q\}$, and the interaction is given by
the Hamiltonian \rf{eq:ranHam} with
\[
U(\s(x),\s(y))= 2\,I_{\{\s(x)\ne\s(y)\}}
\]
and $V\equiv0$. Note that this choice of $U$ and $V$ coincides with
that used in Section \ref{sect:FK} for the standard Potts model and
differs from that in Section \ref{sect:Potts} only by constants which
cancel in the definition of Gibbs distributions.
The Ising model corresponds to the choice $q=2$. As for the disorder,
we make the essential assumption that the random coupling coefficients
$J_b$
are nonnegative, so that the interaction is still ferromagnetic. Of
course, we also make the general assumption of this section that the
$J_b$ are independent with the same distribution, say $\pi$.
A particular case of special interest is that of {\em dilution},
in which the $J_b$ take the values $1$ and $0$ with probabilities $p$
and $1-p$, respectively, which means that $\pi=p\de_{1}+(1-p)\de_{0}$.
For $p=1$ we then recover the homogeneous Potts
model of Section \ref{sect:Potts}.

In the following, the distribution $\pi$ of the $J_b$ will enter only
through the quantities
\[
\bar{p}(\b,\pi)=P\bigg(1-e^{-2\b J_b}\bigg) \;,\quad
\underline{p}(\b,\pi)=P\bigg(\frac{1-e^{-2\b J_b}}{1+ (q-1)e^{-2\b
J_b}}\bigg)
\]
for $\b>0$. Note that they do not depend on $b\in\B$, and
that $0\leq \underline{p}(\b,\pi) \leq \bar{p}(\b,\pi) \leq p$ with
$p= P(J_b>0)$.

For a given realization $\J=(J_b)_{b\in\B}$ of the disorder and
inverse temperature $\b$, we can introduce the Gibbs measure
$\mu^i_{\b \J,q}$
obtained from the Gibbs distributions with constant boundary
condition $i\in\{1,\ldots,q\}$ in the infinite volume limit. This
limit exists by the arguments of  Proposition
\ref{prop:free_and_ordered_Potts_state}, since these use only the
stochastic monotonicity coming from Corollary
\ref{cor:FKG_for_RC} (a) as well as the random-cluster representation,
which
both remain valid in the non-homogeneous case.

The key quantity for phase transition, the order parameter,
is the ``quench\-ed magnetization''
\[
m(\b,\pi) = \frac q{q-1}\;P\bigg(\mu^i_{\b\J,q}(X(0)=i) - \frac
1{q}\bigg)\;.
\]
Indeed, an inhomogeneous version of equation
\rf{eq:Potts_magnetization_connectivity} shows that
$m(\b,\pi)=P(\theta_q(\b\J))$,
where $\theta_q(\b\J)= \phi^1_{\p,q}(0\leftrightarrow\infty)$ is
the percolation probability
for the wired infinite-volume random-cluster measure with
bond probabilities $p_b=1-e^{-2\b J_b}$. Hence, we have $m(\b,\pi)>0$
if and only if $\theta_q(\b\J)>0$ with positive $P$-probability. But
whether or not $\theta_q(\b\J)>0$ does not depend on the value of
$J_b$ for a single bond $b$.
So, Kolmogorov's zero--one law implies that $m(\b,\pi)>0$ if and only if
$\theta_q(\b\J)>0$ $P$-almost surely, and by an inhomogeneous version
of Theorem \ref{thm:Potts_sandwich} the latter means that multiple Gibbs
measures for $\b\J$ exist with $P$-probability $1$. Moreover, an
inhomogeneous version of relation \rf{eq:RC_domination} shows that
$\theta_q(\b\J)$ is an increasing function of $\b\J$.
It follows that $m(\b,\pi)$ is increasing in $\b$ and also in $\pi$
(relative to $\leqd$). In particular, for each $\pi$ there  exists
a critical inverse temperature $\b_{c}(\pi)$, possibly $=+\infty$,
such that $m(\b,\pi)>0$ for $\b>\b_{c}(\pi)$ and $m(\b,\pi)=0$ for
$\b<\b_{c}(\pi)$.

It remains to investigate the quenched magnetization $m(\b,\pi)$.
The following lemma shows how $m(\b,\pi)$ can be estimated
in terms of Bernoulli bond percolation; recall the end of Section
\ref{sect:bond_percolation}.
\begin{lem}\label{ranPlem}
Let $\theta(p) = \phi_p(0\leftrightarrow\infty)$ be the Bernoulli
bond percolation probability on $\Z^d$ with parameter $p$. Then
\[
\theta(\underline{p}(\b,\pi)) \leq m(\b,\pi) \leq
\theta(\bar{p}(\b,\pi))\,.
\]
\end{lem}
{\bf Proof:} We will use an inhomogeneous limiting version of the
domination bounds (b) and (c) of Corollary \ref{cor:FKG_for_RC}.
Although they were stated only in the case of a homogeneous interaction,
they do extend also to the inhomogeneous case. Define two families
$\p=(p_b)_{b\in\B}$ and $\p'=(p'_b)_{b\in\B}$ in terms of a realization
$\J$ by $p_b=1-e^{-2\b J_b}$,
$p_b'=(1-e^{-2\b J_b})/(1+(q-1)e^{-2\b J_b})= p_b/(p_b+q(1-p_b))$.
Let $\phi^1_{\p,q}$ be the associated wired random-cluster measure,
and $\phi_\p$, $\phi_{\p'}$ the corresponding product measures on
$\{0,1\}^\B$. An inhomogeneous version of Corollary \ref{cor:FKG_for_RC}
then shows that
\[
\phi_{\p'}(0\leftrightarrow \infty) \leq
\phi^1_{\p,q}(0\leftrightarrow
\infty) \leq \phi_{\p}(0\leftrightarrow \infty)\;.
\]
We now take the expectation with respect to $P$. In view of the
preceding remarks, the middle term
has $P$-expectation $m(\b,\pi)$, while the $P$-integration of the
Bernoulli measures $\phi_{\p'}$ and $\phi_{\p}$ again leads to
Bernoulli measures,
namely the homogeneous Bernoulli measures
$\phi_{\underline{p}(\b,\pi)}$ and $\phi_{\bar{p}(\b,\pi)}$. The lemma
follows immediately.
$\Cox$

\medskip\noindent
Combining the lemma with the discussion before we arrive at the
following result on phase transition in the random Potts model.
\begin{thm}\label{thm:randomPotts}
Consider the random Potts model on $\Z^d$ at inverse temperature
$\b>0$ with coupling distribution $\pi$. Set $p(\pi)=\pi(]0,\infty[)
=P(J_b>0)$, and
let $p_c$ be the  Bernoulli bond percolation threshold of $\Z^d$; so
$p_{c}=1/2$ when $d=2$.
\begin{itemize}
\item[{\rm(i)}]
If $\bar{p}(\b,\pi) < p_c$ then with $P$-probability $1$ there exists
 only one Gibbs measure with interaction $\b\J$. In particular, this
holds when $p(\pi)<p_c$ or $\b$ is small enough.
\item[{\rm(ii)}]
If $\underline{p}(\b,\pi) > p_c$ then $m(\b,\pi) > 0$, and with
$P$-probability $1$ there exist $q$ distinct phases for the
interaction $\b\J$. In particular, this holds when $p(\pi)>p_c$ and
 $\b$ is large enough.
\end{itemize}
\end{thm}
Another way of stating this result is the following. Suppose
$\pi=(1-p)\de_{0}+p\,\pi_{+}$ with $\pi_{+}=\pi(\cdot\,|]0,\infty[)$,
and let $L_{+}(\b)=\int_{0}^{\infty}e^{-\b t}\pi_{+}(dt)$ be the
Laplace transform of $\pi_{+}$. (Note that then
$\bar{p}(\b,\pi)=p(1-L_{+}(2\b))$ and $\underline{p}(\b,\pi)\ge
p(1-q\,L_{+}(2\b))$.)
Then there is no phase transition for $p<p_{c}$, whereas for $p>p_{c}$
the critical inverse temperature $\b_{c}(p,\pi_{+})\equiv\b_{c}(\pi)$
is finite (and decreasing in $p$) and satifies the bounds
\[
\frac{p-p_{c}}{pq}\le L_{+}(2\,\b_{c}(p,\pi_{+}))\le
\frac{p-p_{c}}{p}\;.
\]
If $\theta(p_{c})=0$ (which
is known to hold for $d=2$, and is expected to hold for all
dimensions) then uniqueness holds when $p=p_{c}$ or
$\b=\b_{c}(p,\pi_{+})$. In physical terminology,
the preceding bounds on  $\b_{c}(p,\pi_{+})$ imply that the
so-called crossover exponent is 1.
\begin{example}\label{ex:dil_Ising}
The case of dilution. {\rm If $J_b$ is 1 or 0 with probability $p$
resp.\ $1-p$ then $L_{+}(\b)=e^{-\b}$. Hence, for $p>p_c$ the critical
inverse temperature satisfies the logarithmic bounds
\[
-\ln\frac{p-p_{c}}{p}\le 2\,\b_{c}(p,\de_{1})\le -\ln\frac{p-p_{c}}{pq}\;.
\]
For $q=2$, the diluted Ising model, assertion (ii) of Theorem
\ref{thm:randomPotts} gives the slightly sharper  upper bound
$\b_{c}(p,\de_{1}) \le \tanh^{-1} (p_c/p)$.
}\end{example}
\begin{example}
The case of power law singularities. {\rm Suppose
$\pi_{+}(dt)=\Gamma(a)^{-1}t^{a-1}e^{-t}dt$ is the Gamma distribution
with parameter $a>0$. Then $L_{+}(\b)=(1+\b)^{-a}$, so that for $p>p_c$
the critical inverse temperature satisfies a power law with
exponent $-1/a$:
\[
\bigg(\frac{p_{c}}{p-p_{c}}\bigg)^{1/a}-1\ \le\ 2\,\b_{c}(p,\pi_{+})\ \le\
\bigg(\frac{p_{c}q}{p-p_{c}}\bigg)^{1/a}\;.
\]
}\end{example}
Examples with other kinds of singularities can easily be produced
\cite{Ge84}.

Theorem \ref{thm:randomPotts} is due to \cite{ACCN1}. Earlier,
a generalized Peierls argument was used in  \cite{Ge81} to show
that for the diluted Ising model ($q=2$) in $d=2$
dimensions a phase transition occurs almost surely when
$p>p_c =1/2$ and $\b$ is large enough. In fact, this paper dealt
mainly with the case of site dilution, in which sites rather than bonds
are randomly removed from the lattice, and which in the present
framework can be described by setting $J_{\langle xy\rangle}=
\xi(x)\xi(y)$ for a family $(\xi(x))_{x\in\Z^d}$ of Bernoulli variables;
the $J_{b}$ are thus 1-dependent. This was continued in \cite{Ge84,Ge85}
for a class of random interaction models including the random-bond Ising
model as considered here, obtaining improved bounds on $\beta_c(p,
\pi_{+})$ for $d=2$ as $p\downarrow 1/2$.
Extensions, in particular to $d\geq 3$, were obtained in \cite{CCF}.
The diluted Ising model with a non-random external field $h\ne0$ does not
exhibit a phase transition; this was shown in \cite{Ge81} for
$\L=\Z^d$ and recently extended to quite general graphs in \cite{HSSt}.

For the diluted Ising model there is also a dynamical phase
transition at the point $p=p_c$. For $p > p_c$ and
$\b_{c}(1,\de_{1})<\b<\b_{c}(p,\de_{1})$ the
relaxation to equilibrium is no longer exponentially fast
\cite{AC1}. This illustrates that uniqueness of the Gibbs measure
does not in itself imply the absence of a critical phenomenon.
Beside such dynamical phenomena, there are also some static effects
of the disorder in the uniqueness regime,
albeit these are perhaps less remarkable. These are the subject of the
next
subsection.

\subsection{Mixing properties in the Griffiths
regime}\label{sect:G-regime}

As we have seen above, the diluted Ising ferromagnet shows spontaneous
magnetization when $p>p_c$ and $\b >\b_c(p)\equiv \b_c(p,\de_{1})$,
and multiple Gibbs measures for $\b\J$ exist almost surely.
In the uniqueness region when still $p>p_c$ but
$\b <\b_c(p)$ we need to distinguish between two different regimes.
At high temperatures when actually $\b < \b_c\equiv\b_c(1)$,
the critical inverse temperature of the undiluted system,
we are in the so-called paramagnetic case. This is
comparable to the usual uniqueness regime for translation
invariant Ising models.  At intermediate temperatures, namely when
$\b_c< \b < \b_c(p)$, we encounter different behavior arising from the
fact that the system starts to feel the disorder. This regime is called
the {\em Griffiths regime}, since it was he \cite{Grif} who discovered in
this parameter region the phenomenon now called Griffiths' singularities.
He studied site-diluted ferromagnets , but the arguments remain valid
also in the bond-diluted case. The basic fact is the following:
adding a complex magnetic field $h$ to the Hamiltonian of the
diluted Ising model we find that the partition function
in a box with plus boundary conditions, as a function of $h$,
can take values arbitrarily close to zero.
The reason is that typically a large part of the box is filled by a
huge cluster of interacting bonds, giving a contribution
corresponding to an Ising partition function in the phase transition
region. The radius of
analyticity of the free energy around $h=0$ is thus zero. In other
words, the magnetization $m(\b,p,h)$ cannot be continued analytically
from $h>0$ to $h<0$ through $h=0$ when $p>p_{c}$ and $\b >\b_c$.
So, the presence of macrosopic clusters of strongly interacting spins
(on which the spins show the low temperature
behavior of the corresponding translation invariant system)
gives rise to singular behavior.
Related phenomena show up in a large variety of other random
models, though not necessarily in the form of non-analyticity in the
uniqueness regime; in general it may be difficult to pinpoint their
precise nature. Nevertheless, we will speak of Griffiths' phase or
Griffiths' regime whenever such singularities are expected to occur,
even when a proof is still lacking. These terms then simply
indicate that the usual high temperature techniques
cannot be applied as such.

As another illustration we consider a random Ising model with
unbounded, say Gaussian coupling variables $J_{b}$. Then
$\b J_b$ is also unbounded, even for arbitrarily small $\b$, and
with high probability a large box contains a positive fraction of
strongly interacting spins. In particular, there is no paramagnetic
regime, and the whole uniqueness region belongs to the Griffiths
phase. For this reason, it is a non-trivial problem to show the
uniqueness of the Gibbs measure. For example, the
standard Dobrushin uniqueness condition encountered in \rf{dssc} (cf.\
\cite{D2,DS1}) is useless in this case; similarly, a naive use of
standard cluster expansion techniques fails. In fact, these methods
are bound to fail since they also imply analyticity which is probably
too much to hope for (even though we cannot disprove it).

In the following we will not deal with the singular behavior in the
Griffiths phase. Instead, we address the problem of showing nice
behavior, which we specify here as good mixing properties of the system.
We shall present two techniques: the use of random-cluster
representations, and the use of disagreement percolation.

\medskip\noindent
{\em Application of random-cluster representations. }
Consider a random Ising model. Spins take values $\s(x)=\pm 1$,
and the formal Hamiltonian is
\begin{equation}\label{rais}
H(\s) = -\sum_{b=\langle xy\rangle\in\B} J_b\; \s(x) \s(y)\;.
\end{equation}
We set $\b=1$.  Let $\mu_{\J,\La}^{\e}$ be the
associated Gibbs distribution in $\La\in {\cal E}$ with boundary condition
$\e\in\Om$. For many applications it is important to have good estimates
on
the variational distance $\|\cdot\|_{\De}$ on $\De
\subset \La$ (see \rf{TVnorm}) of these
measures with different boundary conditions.

\begin{defn}
The random spin system above is called {\bf exponentially
weak-mixing} with rate  $m>0$ if for some $C<\infty$ and all $\La
\in \cal E$ and $\De\subset \La$
\begin{equation}\label{WMix}
P\bigg(
\max_{\e,\e'\in \Om}
\|\mu_{\J,\La}^{\e} -\mu_{\J,\La}^{\e'} \|_\De\bigg)
 \leq C|\De|\,e^{-m\,d(\De,\La^c)}\;,
 \end{equation}
where $d(\De,\La^c)$ is the Euclidean distance of $\De$ and $\La^c$.
\end{defn}
Various other mixing conditions can also be considered. A stronger
condition requires that the variational distance in (\ref{WMix}) is
exponentially small in the distance between the set $\De$ and
the region where the boundary conditions $\e$ and $\e'$ really
differ. One could also restrict $\La$ and/or $\De$ to
regular boxes. See \cite{MaO,Ze1,DS1,DS2,vB,CM1,CM2,AC1}.

Let us comment on the significance of the exponential weak-mixing
condition above.

\medskip\noindent
{\bf Remarks:} (1) Suppose condition \rf{WMix} holds. A straightforward
application of the Borel-Cantelli lemma then shows that for any $m'<m$
\[
\max_{\e,\e'\in \Om}
\|\mu_{\J,\La}^{\e} -\mu_{\J,\La}^{\e'} \|_\De \leq
C_\J\,|\De|\,e^{-m'\,d(\De,\La^c)}
\]
with some realization-dependent $C_\J<\infty$ $P$-almost surely.
Integrating over $\e'$ for any
Gibbs measure $\mu_{\J}$ we find in particular that
$\mu_\J=\lim_{\La\uparrow\Z^d}
\mu_{\J,\La}^{\e}$ for all $\e$, implying that $\mu_{\J}$ is the only
Gibbs measure
(and depends measurably on $\J$). Moreover, noting that
$\mu_\J(A|B)=\int \mu_{\J,\La}^{\e}(A)\, \mu_\J(d\e|B)$ for $A\in
\F_\De$ and $B\in\F_{\La^c}$, we see that this
realization-dependent Gibbs measure $\mu_{\J}$ satisfies the
exponential weak-mixing condition
\begin{equation}\label{eq:wmix2}
\sup_{A\in \F_\De,\, B\in\F_{\La^c},\, \mu_\J(B)>0}
|\mu_\J(A|B) -\mu_\J(A)|
 \leq C_\J\,|\De|\,e^{-m'\,d(\De,\La^c)}\;.
\end{equation}

(2) Condition \rf{eq:wmix2} above also implies an exponential decay of
covariances.
Let $f$ be any local observable with dependence set $\De\in\E$ and
$g$ be any bounded observable
depending only on the spins off $\La$, where $\De\subset\La$.
Also, let $\de(f)=\max_{\s,\s'} | f(\s) - f(\s')|$ be the total
oscillation of $f$ and $\de(g)$ that of $g$. The covariance
$\mu_{\J}(f;g)$ of $f$ and $g$ then satisfies $P$-almost surely
the inequality
\[
|\mu_{\J}(f;g)|
       \leq C_{\J}\,|\De|\, \de(f)\, \de(g)\, e^{-m' \,d(\De,\La^c)}
       /2\;.
\]
Indeed, a short computation shows that
$|\mu_{\J,\h}(f;g)|$ is not larger than the  left-hand side of
inequality \rf{eq:wmix2} times $\de(f)\,\de(g)/2$,
cf.\ inequality (8.33) in \cite{Geo}. If $g$ is local, a similar
inequality holds for covariances relative to
finite volume Gibbs distributions in sufficiently large
regions with arbitrary boundary conditions.

\medskip\noindent
We will now investigate under which conditions the random Ising
system with Hamiltonian \rf{rais} is exponentially weak-mixing. We
start from the estimate
\[
\|\mu_{\J,\La}^{\e} -\mu_{\J,\La}^{\e'} \|_\De
\leq  \phi_{\p,2,\La}^1 (\De \leftrightarrow \partial\La)
\]
obtained in Theorem \ref{AC}. As before, this
bound is also valid in the inhomogeneous case considered here, and
$\phi_{\p,2,\La}^1$ stands for the wired random-cluster distribution
in $\La$ with bond-probabilities $\p\equiv\p(|\J|)=(p_{b})_{b\in\B}$
given by $p_b=1 -\exp[-2|J_b|]$. Next we can use a recent concavity
result of \cite{AC1}:
\begin{lem}\label{concave}
Let\/ $\K=(K_b)_{b\in\B}$ be a collection of positive real
numbers and
$\p=\p(\K)$ denote the collection of densities $p_b
= 1-e^{-2K_b}$. For any increasing function $f$, the expectation
$\phi_{\p,2,\La}^1(f)$ is then a concave function of each $K_b$.
\end{lem}
{\bf Proof:} For brevity we set $\phi_{\p,2,\La}^1 = \phi$. For any
fixed bond $b$ we consider the functions
$F(p_b) = \phi(f)$ and $G(K_b)=F(1-e^{-2K_b})$.
Using equation (\ref{eq:single_edge}) of Lemma \ref{lem:single_edge}
we then find that
\[
F'(p_b) = \phi(f\,g_{p_b}) - \phi(f)\,\phi(g_{p_b})
\]
and
\[
F''(p_b)= -2F'(p_b)\phi(g_{p_b})\;,
\]
where for $\eta_b\in\{0,1\}$
\[
g_{p_b}(\eta_b) = \frac{2\eta_b-1}{p_b^{\eta_b}(1-p_b)^{1-\eta_b}}\;.
\]
This implies
\[
G''(K_b)=-4e^{-2K_b}F'(p_b)\,[1+2e^{-2K_b}\phi(g_{p_b})].
\]
Now, $F'$ is nonnegative because $f$ and $g_{p_b}$ are
increasing and $\phi$ has positive correlations by
Corollary \ref{cor:FKG_for_RC}(a).
Another explicit computation shows that
\[
1+2e^{-2K_b}\phi(g_{p_b}) = 2\phi(\eta_b=1)/p_b - 1 \geq 2/(2-p_b) - 1
\geq 0,
\]
where the first inequality uses the fact that the
random-cluster distribution for $q=2$ dominates an independent percolation
model with densities $p_b'= p_b/(2 - p_b)$, see Corollary
\ref{cor:FKG_for_RC}(c). It follows that $G'' \leq 0$,
proving the claimed concavity. $\Cox$

\medskip\noindent
The preceding lemma implies that
\[
P\bigg( \phi^1_{\p(|\J|),2,\La}(\De \leftrightarrow \partial\La)\bigg)
\leq
\phi^1_{p,2,\La}(\De \leftrightarrow \partial\La)
\le \sum_{x\in\partial\De}\phi^1_{p,2,\La}(x\leftrightarrow \partial\La),
\]
where $p=1-\exp[-2\,P(|J_{b}|]$. The weak
mixing property \rf{WMix} thus follows provided we can show an
exponential decay of connectivity in the wired
random-cluster distribution $\phi^1_{p,2,\La}$.
The latter certainly holds when $p < p_c$ because
$\phi^1_{p,2,\La}$ is dominated (on $\La$) by the Bernoulli measure
$\phi_{p}$ (cf.\ Corollary \ref{cor:FKG_for_RC}(b)),
and the connectivity of a subcritical Bernoulli model decays
exponentially fast (recall Theorem \ref{thm:exponential_tail}).
So we arrive in particular at the following result.
\begin{thm}\label{thm:wmix}
Consider a random Ising model on $\Z^d$ with Hamiltonian \rf{rais}.
If \linebreak $2\,P(|J_{b}|)<-\ln(1-p_{c})$ for the Bernoulli bond
percolation threshold $p_{c}$ of $\Z^d$ then the system is
exponentially weak-mixing.
\end{thm}
As the proof above shows, the factor $|\De|$ on the right-hand side of
\rf{WMix} can
actually be replaced by $|\partial\De|$.

In the diluted Ising model at inverse temperature $\b$ (see Example
\ref{ex:dil_Ising}), the condition of the preceding theorem reads $2\b
p<-\ln(1-p_c)$.
Therefore, if $d$ is so large that $2\b_c\, p_c<-\ln(1-p_c)$ then the
theorem
covers part of the Griffiths regime. This fact is evident in the case of
an unbounded, say Gaussian, disorder because then (as explained above)
there is no paramagnetic phase.

Exponential weak-mixing for random Ising models can also be shown
by other applications of random-cluster domination. Let us sketch
such an alternative route. Using the random-cluster representation,
Newman \cite{N} showed that (pointwise in the disorder $\J$)
\begin{equation}\label{new}
\max_{\e, \e'}
\|\mu_{\J,\La}^{\e} -\mu_{\J,\La}^{\e'} \|_\De
 \leq
2 \sum_{x\in \De} \mu_{|\J|,\La}^{+}(X(x)),
\end{equation}
where $\mu_{|\J|,\La}^{+}$ is the Gibbs distribution
in $\La$ with plus boundary conditions for the Hamiltonian
(\ref{rais}) with $J_b$ replaced by
$|J_b|$. On the other hand,
Higuchi \cite{Hi1} obtained the estimate
\[
\mu_{|\J|,\La}^{+}(X(x)) \leq \sum_{y\notin \La}
  \sum_{z\in \La :\, z\sim y}\mu_{|\J|,\La}^{f}(X(x)X(z))
\]
(the superscript `$f$' referring to the free boundary condition), while
Olivieri, Perez and Goulart Rosa \cite{OPR} proved that
\begin{equation} \label{oli}
P\bigg(\mu_{|\J|,\La}^{f}(X(x)X(z))\bigg) \leq
\mu_{J,\La}^{f}(X(x)X(z))
\end{equation}
with $J=P(|J_{b}|)$. We are thus back to the standard Ising Gibbs
distribution in $\La$ with zero external field, free boundary
condition and coupling constant $J$. Now we can take advantage of
the second Griffiths inequality (which we did not discuss so far in
this text) stating that correlation functions such as on the
right-hand side above are monotone in the coupling coefficients,
see e.g.\ \cite{L}. This implies that the right-hand side of
(\ref{oli}) is an increasing function of $\La$ and thus bounded
above by its infinite volume limit. But for $J$ less than $J_c$,
the critical coupling, the Gibbs measure is unique and has an
exponential decay of correlations, see \cite{ABF}. We thus find
that for $J<J_c$ the right-hand side of (\ref{oli}) has an
exponential upper bound $C\exp[-m|x-z|]$ for suitable constants $C
< \infty$ and $m > 0$.  Together with the previous estimates, we
conclude that {\em under the condition $P(|J_{b}|)<J_{c}$, the
random Ising model on $\Z^d$ with Hamiltonian \rf{rais} is
exponentially weak-mixing}. Again, this condition covers part of
the Griffiths regime for the diluted Ising ferromagnet.

As an alternative to the use of Griffiths' inequalities above
we can also apply Theorem \ref{thm:RC_coupling} and
Corollary \ref{cor:FKG_for_RC}(b), giving
\[ 
P\bigg(\mu_{|\J|,\La}^{f}(X(x)X(z))\bigg) \leq
\phi_p(x\leftrightarrow z)
\]
where $\phi_p$ is the bond Bernoulli measure with density $p= P(1 -
\exp[-2|J_b|])$. As mentioned above, its connectivity function
decays exponentially fast when $p < p_c$.  We therefore conclude that
{\em the random Ising model on $\Z^d$ is also exponentially
weak-mixing when $P(1 -\exp[-2|J_b|]) < p_c$}.

The above estimate (\ref{new}) does not hold when we add a random magnetic
field to the Hamiltonian (\ref{rais}),
\begin{equation} \label{raish}
H(\s) = -\sum_{b=\langle x y\rangle\in\B} J_b\, \s(x) \s(y) -
\sum_{x\in\Z^d} h_x\,\s(x),
\end{equation}
with i.i.d.\ real random variables $h_x$ independent from the
$J_b$. However, if  $J_b, h_x \geq 0$  then we can take advantage of
Section 2 of \cite{Hi1} to replace (\ref{new}) with
\[
\max_{\e, \e'}
\|\mu_{\J,\h,\La}^{\e} -\mu_{\J,\h,\La}^{\e'} \|_\De
  \leq
2 \sum_{x\in \De} \mu_{\J,0,\La}^{+}(X(x)),
\]
and we can continue as above.

\medskip\noindent
{\em Application of disagreement percolation. }
As we have indicated at the end of Section \ref{sect:disa},
the idea of disagreement percolation can be applied
to study the Griffiths regime for rather general
random-interaction systems. To be specific we consider Ising spins with
the Hamiltonian (\ref{raish}).  We consider the
finite volume Gibbs distribution $\mu_{\J,\h,\La}^\eta$ in a box
$\La$ with  boundary condition $\e\,$; as before, the
subscripts $\J,\h$ describe the random interaction.
We are going to apply Theorems
\ref{thm:dis} and \ref{thm:uni} pointwise in the
disorder. As in (\ref{dens}) we thus have to consider the variational
single-spin oscillations
\begin{equation}\label{dens3}
p_x^{\J,\h}=\max_{\e,\e'\in \Om}
\|\mu_{\J,\h,x}^{\e} -\mu_{\J,\h,x}^{\e'} \|_x
\end{equation}
which depend on the disorder $\J,\h$. Under the present
assumptions, $(p_x^{\J,\h})_{x\in\Z^d}$ is a 1-dependent random
field; this is the only property of the disorder needed below. Now,
from Theorems \ref{thm:dis} and \ref{thm:uni} we can conclude that
if with $P$-probability 1 there is no Bernoulli site-percolation
with densities (\ref{dens3}) then, again with $P$-probability 1,
there is a unique Gibbs measure for $\J,\h$.  (In fact, similarly
to the results of Section \ref{ssection:stodom} the random system
can be dominated by a random percolation system which more or less
coincides with a stochastic-geometric representation of the diluted
ferromagnet \cite{Gie}.) The following theorem is an immediate
consequence of the results of Section \ref{sect:disa}.
\begin{thm}\label{twop} An Ising system with random Hamiltonian
\rf{raish} satisfying
\[
P(p_x^{\J,\h}) < \frac{1}{(2d-1)^2}
\]
is exponentially weak-mixing.
\end{thm}
Of course, we now have to add in \rf{WMix} the subscript $\h$
referring to the random external field.

\medskip\noindent
{\bf Proof:}  Consider Theorem \ref{thm:dis}. Evidently, the right
hand side of (\ref{corr}) can only increase if the density for an
open site is set to be $1$ on the odd sublattice $\L_{odd}$ of
$\Z^d$. Due to the 1-dependence of the random field
$\p=(p_x^{\J,\h})_{x\in\Z^d}$ noticed above, the remaining
$p_x^{\J,\h}$ with $x\in\L_{even}$ are mutually independent. Taking
the $P$-expectation in (\ref{corr}) we thus obtain on the
right-hand side the Bernoulli percolation probability
$\psi_{p,\,even}(\De\,{\leftrightarrow}\,\partial\La)$, where
$\psi_{p,\,even}$ is the Bernoulli measure with density
$p=P(p_x^{\J,\h})$ on the even sublattice $\L_{even}$ and density 1
on $\L_{odd}$. The exponential weak-mixing property now follows by
an argument similar to that following
(\ref{eq:percolation_needs_open_paths}); note that, relative to
$\psi_{p,\,even}$, a path of length $k$ is open with probability at
most $p^{\lfloor k/2\rfloor}$. $\Cox$

\medskip\noindent
{\bf Remarks:} (1) Suppose we add a uniform magnetic field $h$
to the random Hamiltonian \rf{raish}. Under the conditions of
Theorem \ref{twop}
it is then not too difficult to show that the disorder-averaged
expectation
$P(\mu_{\J,\h+h}(f))$ of any local function $f$ is an infinitely
differentiable function of $h$, see e.g.\ \cite{DKP}.

(2) Following Theorem \ref{thm:uni}, Theorem \ref{twop} requires
that the single-point densities (\ref{dens3}) are globally small enough
to prevent percolation in the associated Bernoulli model.
This condition can be extended to so-called constructive
conditions involving finite boxes rather than single sites \cite{vB}.

(3) Another very powerful and transparent treatment of the
Griffiths regime (based on very similar percolation ideas) has been
developed in \cite{DKP}.  This paper proposes a technique similar
to \cite{BD} and \cite{FI} of a resummation in the high temperature
cluster expansion. The bounds then allow a probabilistic
interpretation of the expansion linking it with a bond percolation
process.

\section{Continuum models}

All models considered so far lived on a lattice.
Physical systems like real gases, however, are more realistically
modelled by particles living in continuous space. This section is an
outline of how
some of the stochastic-geometric ideas developed in previous sections can
be
applied to a continuum setting. First, in Section \ref{sect:cont_perc},
we consider the natural continuum analogues, based on Poisson processes,
of the Bernoulli percolation models introduced in Section
\ref{sect:percolation}.
Then, in Section \ref{sect:cont_WR}, we consider a continuum variant
of the Widom--Rowlinson model introduced in Section \ref{sect:WR}
and discuss its phase transition behavior. (As mentioned in
Section \ref{sect:cont_WR}, this continuum variant
is the one originally
considered by Widom and Rowlinson \cite{WR}, so it predates the
lattice model.)

\subsection{Continuum percolation} \label{sect:cont_perc}

Here we consider the basic models of continuum percolation. For a thorough
treatment of the mathematical theory of continuum percolation we
refer to Meester and Roy \cite{MR}.

We first need to introduce the Poisson process on $\R^d$ and its subsets.
Heuristically, a Poisson process with intensity $\la>0$
on $\R^d$
is a random set $X$ of points in $\R^d$ with the properties that
\begin{itemize}
\item[(i) ] for any bounded Borel set $\La$ of $\R^d$ with volume
$|\La|$, the number of points of $X$ in $\La$
is Poisson distributed with mean $\la |\La|$, i.e., for $n=0,1,2,\ldots$
the probability of seeing exactly $n$ points in $\La$ equals
${ \exp(-\la |\La|)(\la |\La|)^n}/{n!}\,$;
\item[(ii) ] for any two disjoint such subsets $\La_1$ and $\La_2$,
the numbers of points observed in $\La_1$ and in $\La_2$ are independent.
\end{itemize}
For a construction of such a process, we first consider
a bounded Borel set $\La$ of $\R^d$. Let $\Om_\La$
be the set of all finite subsets of $\La$.
A Poisson process
on $\La$ with intensity $\la>0$ is then given by a random element of
$\Om_\La$ having distribution $\pi_{\la,\La}$, where
\[
\pi_{\la,\La}(F)= e^{-\la |\La|}\;\sum_{n=0}^\infty
\,\frac{\la^n }{n!} \int \cdots \int
I_{\{\{x_1, \ldots, x_n\}\in F\}}\, dx_1\cdots dx_n
\]
for all $F\in {\cal F}_\La$, the smallest $\sigma$-field
which allows us to count the number of points in each Borel subset of
$\La$.

Next, let $\Om$ be the set of all locally finite point
configurations on $ \R^d$; locally finite means that any bounded set
contains only finitely many points.
The Poisson process on $ \R^d$ with intensity $\la$
is a random point configuration $X$ distributed
according to the unique probability measure $\pi_\la $ on $\Om$
which, when projected on any bounded Borel set $\La \subset  \R^d$, yields
$\pi_{\la,\La}$. Properties (i) and (ii) above are easily verified,
and make Poisson processes the natural analogues of Bernoulli measures
for lattice models.

To study percolation properties of the Poisson process $X$, we need to
introduce some notion of connectivity. A natural way is to imagine a
closed
Euclidean ball $B(x,R)$ of fixed radius $R$ around each point $x$ of
the Poisson process, which leads us to considering the random subset
$\bar X=\bigcup_{x\in X}B(x,R)$ of $\R^d$.
Such random subsets are widely known as {\em Boolean models}.
(More generally, $B(x,R)$ could be replaced  by a
closed compact random set centered at $x$.)
This particular Boolean model is often referred to as the
{\em Poisson blob model} or {\em lily pond model}.
Two points $x,\,y\in X$ are then considered as connected to each other
if they are connected in $\bar X$, meaning that there
exists a continuous path from $x$ to $y$ in $\bar X$. By scaling, there
is no
loss of generality in setting $R={1}/{2}$, so that two balls centered at
$x$ and $y$ intersect if and only if $|x-y| \leq 1$, where $| \cdot |$
denotes
Euclidean distance.
 The basic result on Boolean continuum percolation,
analogous to Theorem \ref{thm:perc_nontrivial} for ordinary site
percolation,
is the following.
\begin{thm} \label{thm:cont_perc_nontrivial}
Pick a  Poisson process  $X$ on $ \R^d$, $d\geq 2$, with intensity $\la$.
Let $\bar{X}= \bigcup_{x\in X}B(x, {1}/{2})$ be the associated Boolean
model,
and let $\theta(\la)$ denote the probability that the origin
belongs to an unbounded connected component of $\bar{X}$. Then there
exists a critical value $\la_c=\la_c(d)\in (0, \infty)$
such that
$\theta(\la)=0$ if $\la<\la_c$ and
$\theta(\la)> 0$ if $\la>\la_c $.
\end{thm}
The standard proof of this result is based on a partitioning of $\R^d$
into
small cubes, reducing the problem to its lattice analogue,
Theorem \ref{thm:perc_nontrivial}; see \cite{MR} or \cite{Gr1}.

Another continuum percolation model is the
so called {\em random connection model}, or {\em Poisson random edge
model},
which was introduced by M.\ Penrose \cite{MPenrose}. Let
$g: \, [0, \infty) \rightarrow [0,1]$ be a decreasing function with
bounded support (that is $g(x)=0$ when $x$ exceeds some $R<\infty$).
The random connection model with intensity $\la$ and
connectivity function $g$ arises by taking a Poisson process $X$ in $
\R^d$
with intensity $\la$ and independently drawing an edge between
each pair of points $x$ and $y$ of $X$ with probability
$g(|x-y|)$. This setting includes the Boolean model,
which corresponds to the choice $g=I_{[0,1]}$. Theorem
\ref{thm:cont_perc_nontrivial} extends to this model: For $g$
as above with $\int_0^\infty g(x)dx>0$ and dimensions $d\geq 2$,
there is a critical value $\la_c=\la_c(d,g)$ such that
infinite connected components a.s.\ occur (resp.\ do not occur)
whenever $\la>\la_c$ (resp.\ $\la<\la_c$) in the
random connection model with intensity $\la$ and
connectivity function $g$.

Much of the theory of standard (lattice) percolation has analogues for
these continuum models. An example is the uniqueness of the infinite
cluster
(Theorem \ref{thm:Burton_Keane}), which goes through for the Boolean and
random connection models. See \cite{MR} for this and much more.

\subsection{The continuum Widom--Rowlinson model} \label{sect:cont_WR}

The continuum Widom--Rowlinson model  is a marked point process
where the points are of two types: we call them plus-points and
minus-points.
(From now on and throughout this section,
we drop the term ``continuum''
when refering to this model, and instead add the word ``lattice'' when
talking about the model of Section \ref{sect:WR}.)
For the model defined on a region $\La \subseteq  \R^d$
realizations take values in $\Om_\La \times \Om_\La$;
the first coordinate describes the locations of the plus-points, and the
second coordinate the minus-points. There is a hard sphere interaction
preventing two points from coming within Euclidean distance $R$ from each
other; again, we set $R=1$ without loss of generality. This interaction
corresponds to the Hamiltonian
\[
H({\bf x},{\bf y})=\sum_{x\in{\bf x},y\in{\bf y}} \infty\, I_{\{|x-y|\le
1\}}\;,
\]
${\bf x},{\bf y}\in\Om_\La\,$.

When $\La$ is bounded, the Widom--Rowlinson model on $\La$
with intensity $\la$ is obtained by conditioning the Poisson product
measure $\pi_{\la,\La}\times\pi_{\la,\La}$
on the event that there is no pair of points of opposite type within
unit distance from each other. The extension to $\R^d$ is done in the
usual DLR fashion: a probability measure $\mu$ on
$\Om\times\Om$ is a Gibbs measure for the
Widom--Rowlinson model at intensity $\la$ if,
for any bounded Borel set $\La$, the conditional distribution of the point
configuration on $\La$ given the point
configuration on $ \R^d \setminus \La$ is that of two independent
Poisson processes conditioned on the event that no point in $\La$ is
placed within
unit distance from a point of the opposite type, either inside or outside
$\La$.
The resemblence with the lattice Widom--Rowlinson model of
Section \ref{sect:WR} is evident. We have the following
analogue of Theorem \ref{thm:WR}.
\begin{thm} \label{thm:cont_WR}
For the Widom--Rowlinson model on $ \R^d$, $d\geq 2$,
with activity $\la$,
there exist constants $0<\la'_c\leq \la''_c<\infty$
(depending on $d$) such that for $\la<\la'_c$
the model has a unique Gibbs measure, while for
$\la>\la''_c$ there are multiple Gibbs measures.
\end{thm}
The proof of this result splits naturally into two parts: first, we need
to demonstrate uniqueness of Gibbs measures for $\la$ sufficiently small,
and secondly we need to show non-uniqueness for $\la$ sufficiently large.
The first half can be done by a variety of techniques. For instance, one
can partition $ \R^d$ into cubes of unit sidelength and apply
disagreement percolation (Theorem \ref{thm:dis}).
Two observations are crucial in
order to make this work: that the conditional distribution of the
configuration in such a cube given everything else only depends on
the configurations in its neighboring cubes, and that the conditional
probability of seeing no point at all
in a cube tends to $1$ as $\la \rightarrow 0$, uniformly in the
neighbors' configurations.

The more difficult part, the nonuniqueness for large $\la$,
was first obtained by Ruelle \cite{Ru} using a Peierls-type argument.
Here we shall sketch a modern stochastic-geometric approach using
a random-cluster representation. This approach is due mainly to
Chayes, Chayes and Koteck\'y \cite{CCK} (but see also \cite{GLM}),
and works in showing both parts of Theorem \ref{thm:WR}.
The so called {\em continuum random-cluster model} is defined as follows.

\begin{defn}
The {\bf continuum random-cluster distribution} $\phi_{\la,\La}$ with
intensity $\la$ for the compact region $\La \subset  \R^d$
is the probability measure on $\Om_\La$ with density
\begin{equation} \label{eq:CRC}
f({\bf x})=
\frac{1}{Z_{\la,\La}} \; 2^{k({\bf x})}\;,\quad {\bf x}\in\Om_\La
\end{equation}
with respect to the  Poisson process
$\pi_{\la,\La}$ of intensity $\la$;
here $Z_{\la,\La}$ is a normalizing constant and $k({\bf x})$
is the number of connected components of the set
$\bar{\bf x}=\bigcup_{x\in {\bf x}}B(x,{1}/{2})$.
\end{defn}
In analogy to the correspondence between the lattice
Widom--Rowlinson model and its random-cluster representation in
Propositions \ref{prop:RC_to_WR} and \ref{prop:WR_to_RC}, we obtain
the continuum random-cluster model by simply disregarding the types
of the points in the Widom--Rowlinson model, with the same choice
of the parameter $\la$. Conversely, the Widom--Rowlinson model is
obtained when the connected components in the continuum
random-cluster model are assigned independent types, plus or minus
with probability ${1}/{2}$ each. To see why this is true, note that
for any ${\bf x}\in \Om_\La$ there are exactly $2^{k({\bf x})}$
elements of $\Om_\La \times \Om_\La$ which do not contradict the
hard sphere condition of the Widom--Rowlinson model and which map
into ${\bf x}$ when we disregard the types of the points.

Besides the random-cluster representation, we can also take advantage
of stochastic monotonicity properties.  Let us define a
partial order $\preceq$ on $\Om \times \Om$
by setting
\begin{equation}\label{eq:CWR-order}
({\bf x}, {\bf y}) \preceq ({\bf x}', {\bf y}') \hspace{4 mm}
\mbox{if} \hspace{4 mm} {\bf x} \subseteq {\bf x}' \mbox{ and }
{\bf y} \supseteq {\bf y}'\;,
\end{equation}
so that in other words a configuration increases with respect to
$\preceq$ if plus-points are added and minus-points are deleted.
A straightforward extension of Theorem \ref{thm:stochdom_pointproc}
below then implies that the Gibbs distributions for the Widom--Rowlinson
model
have positive correlations relative to this order.
The methods of Sections \ref{sect:ineq} and \ref{sect:Ising_monotone}
can therefore be adapted to show that the Widom--Rowlinson model on
$\R^d$ at
intensity $\la$ admits two particular phases $\mu^+_\la$ and $\mu^-_\la$,
where $\mu^+_\la$ is obtained as a weak limit of Gibbs measures on
compact sets
(tending to $\R^d$) with boundary condition consisting of a dense crowd
of plus-points, and $\mu^-_\la$ is obtained similarly. We also
have the sandwiching relation
\begin{equation}\label{eq:CWR_sandwich}
\mu^-_\la \leqd \mu \leqd \mu^+_\la
\end{equation}
for any Gibbs measure $\mu$ for the intensity $\la$
Widom--Rowlinson model on $\R^d$, so that
uniqueness of Gibbs measures is equivalent to having $\mu^-_\la=
\mu^+_\la$.

The Gibbs measure for the Widom--Rowlinson model on a box $\La$
with ``plus'' (or ``minus'') boundary condition corresponds to the wired
continuum random-cluster model $\phi^1_{\la,\La}$ on $\La$ where all
connected components
within distance ${1}/{2}$ from the boundary count as a single
component. Arguing as in Sections \ref{sect:PT_Potts} and
\ref{sect:WR_RC} we find that uniqueness of Gibbs
measures for the Widom--Rowlinson model is equivalent to not having
any infinite connected components in the continuum random-cluster model.
Let $\la_c$ be as in Theorem \ref{thm:cont_perc_nontrivial}.
Theorem \ref{thm:cont_WR} follows if we can show
that the continuum random-cluster model $\phi^1_{\la,\La}$ with
sufficiently large intensity $\la$
stochastically dominates $\pi_{\la_1,\La}$ for some $\la_1>\la_c$,
whereas $\phi^1_{\la,\La}\leqd\pi_{\la_2,\La}$ for some $\la_2<\la_c$
when $\la$ is sufficiently small.

To this end we need a point process analogue of
Theorem \ref{thm:holley},
which is based on the concept of {\em Papangelou
(conditional) intensities} for point processes. Suppose $\mu$ is a
probability measure on $\Om_\La$ which is absolutely continuous with
density $f({\bf x})$ relative to
the unit intensity Poisson process $\pi_{1,\La}$. For $x\in \La$
and a point
configuration ${\bf x}\in \Om_\La$ not containing $x$, the
Papangelou intensity of $\mu$ at $x$ given ${\bf x}$
is, if it exists,
\begin{equation} \label{eq:Pap}
\la(x|{\bf x})=\frac{f({\bf x}\cup \{x\})}{f({\bf x})}\;.
\end{equation}
Heuristically, $\la(x|{\bf x})dx$ can be interpreted as the
probability of finding a point inside an infinitesimal region $dx$
around $x$,
given that the point configuration outside this region is ${\bf x}$.
Alternatively, $\la(\cdot\,|\,\cdot)$ can be characterized as
the Radon--Nikodym
density of the measure
$\int\mu(d{\bf x})\sum_{x\in{\bf x}}\de_{(x,{\bf x}\setminus\{x\})}$ on
$\La\times\Om_\La$, the so-called reduced Campbell measure of $\mu$,
relative to the Lebesgue measure times $\mu$ \cite{GK}. It is easily
checked that the Poisson process $\pi_{\la,\La}$ has Papangelou
intensity $\la(x|{\bf x})=\la$.

The following point process analogue of Theorem \ref{thm:holley} was
proved by Preston \cite{Pr20} under an additional technical assumption,
using a coupling of so called spatial
birth-and-death processes similar to the coupling used in the proof of
Theorem \ref{thm:holley}. Later, the full result was proved in Georgii and
K\"uneth \cite{GK} by a discretization argument.
\begin{thm} \label{thm:stochdom_pointproc}
Suppose $\mu$ and $\tilde{\mu}$ are probability measures on
$\Om_\La$ with Papangelou intensities
$\la(\cdot\,|\,\cdot)$ and $\tilde{\la}(\cdot\,|\,\cdot)$ satisfying
\[
\la(x|{\bf x})\leq \tilde{\la}(x|\tilde{\bf x})
\]
whenever $x\in \La$ and ${\bf x}, \tilde{\bf x}\in\Om_\La$
are such that
${\bf x}\subseteq\tilde{\bf x}$. Then $\mu\leqd \tilde{\mu}$, in that
there exists a coupling $(X,\tilde{X})$ of $\mu$ and $\tilde{\mu}$
such that $X\subseteq \tilde{X}$ a.s.
\end{thm}
Plugging (\ref{eq:CRC}) into (\ref{eq:Pap}) we find that
the continuum random-cluster measure $\phi_{\la,\La}$
has Papangelou intensity
\begin{equation} \label{eq:Pap-CRC}
\la(x|{\bf x})=\la\,2^{1-\kappa(x,{\bf x})}\;,
\end{equation}
where $\kappa(x,{\bf x})$ is the number of connected components of
$\bigcup_{y\in {\bf x}}B(y,{1}/{2})$ intersecting
$B(x, {1}/{2})$. It is a simple geometric fact that there
exists a constant $\kappa_{max}=\kappa_{max}(d)<\infty$ such that
$\kappa(x,{\bf x})\leq\kappa_{max}$ for all $x$ and ${\bf x}$; for $d=2$
we may take $\kappa_{max}=5$. It follows that
\begin{equation} \label{eq:stochdom-CRC}
\la\,2^{1-\kappa_{max}}\leq\la(x|{\bf x})\leq 2\la
\end{equation}
for all $x$ and ${\bf x}$. Hence, applying Theorems
\ref{thm:stochdom_pointproc} and \ref{thm:cont_perc_nontrivial}
we find that $\phi^1_{\la,\La}\leqd\pi_{2\la,\La}$, so that for
$\la<{\la_c}/{2}$ we obtain
the absence of unbounded connected components of
$\bigcup_{x\in {\bf x}}B(x,{1}/{2})$ in the limit $\La\uparrow\R^d$
of continuum random-cluster measures. On the other hand, taking
$\la>\la_c\,2^{\kappa_{max}-1}$ yields the presence of unbounded
connected components in the same limit. Theorem \ref{thm:cont_WR} follows
immediately.

It is important to
note that this approach does not allow us to show that nonuniqueness of
Gibbs
measures depends monotonically on $\la$. The reason is similar to what
we saw for the lattice Widom--Rowlinson model in Section \ref{sect:WR_RC}:
the right hand side
of (\ref{eq:Pap-CRC}) fails to be increasing in ${\bf x}$. It thus remains
an open problem whether
one can actually take $\la'_c=\la''_c$ in Theorem \ref{thm:cont_WR}.

There are several interesting generalizations of the Widom--Rowlinson
model.
Let us mention one of them, in which neighboring pairs of particles of the
opposite type are not forbidden, but merely discouraged. Let
$h: [0, \infty) \rightarrow [0, \infty]$ be an
``interspecies repulsion function''  which is
decreasing and has bounded support. For $\La \subset \R^d$
compact and $\la>0$, the associated Gibbs distribution
$\mu_{h, \la,\La}$ on $\Om_\La \times \Om_\La$ is given by its density
\[
f({\bf x}, {\bf y}) = \frac{1}{Z_{h, \la,\La}}\;
 \exp\bigg( - \sum_{x \in {\bf x}, y \in {\bf y}} h(|x-y|) \bigg) \, .
\]
relative to $\pi_{\la,\La}\times \pi_{\la,\La}$.
Infinite volume Gibbs measures on $\Om\times\Om$ are
then defined in the usual way. Lebowitz and
Lieb \cite{LL} proved nonuniqueness of Gibbs measures for large
$\la$ when $h(x)$ is large enough in a neighborhood of the origin.
Georgii and H\"aggstr\"om \cite{GH} later established the same behavior
without
this condition, and for
a larger class of systems, using the random-cluster approach. This
involves
a generalization of the continuum random-cluster model, which arises by
taking
the random connection model of Section \ref{sect:cont_perc}  with
connectivity function $g(x)= 1-e^{-h(x)}$ and biasing it with a factor
$2^{k({\bf z})}$,
where $k({\bf z})$ is the number of connected components of a
configuration
${\bf z}$ of points and edges. To establish the phase transition behavior
of this ``soft-core Widom--Rowlinson model''
(i.e., uniqueness of Gibbs measures for small $\la$
and nonuniqueness for large $\la$) one can basically use the same
arguments as the ones sketched above for the standard Widom--Rowlinson
model.
However, due to the extra randomness of the edges some parts of the
argument become
more involved. In particular, there is no longer a deterministic
bound (corresponding to $\kappa_{max}$) on how much the number of
connected components can decrease when a point is added to the
random-cluster configuration; thus more work is needed to obtain an
analogue of the
first inequality in (\ref{eq:stochdom-CRC}).

To conclude, we note that the Widom--Rowlinson model on $\R^d$ has an
obvious
multitype analogue with $q\ge3$ different types of particles. This
multitype model
still admits a random-cluster representation from which the existence of
a phase
transition can be derived \cite{GH}. There is, however, no partial
ordering
like \rf{eq:CWR-order} giving rise to
stochastic monotonicity or an analogue of \rf{eq:CWR_sandwich}.

\end{document}